\documentclass[11pt,a4paper]{article}
\usepackage[utf8]{inputenc}

\usepackage{mathrsfs,amsthm,xcolor,verbatim,bbm,amsmath,amsfonts,amssymb,nicefrac,hyperref,bm,mathtools,xparse,etoolbox,enumitem,derivative,cite,tzplot,pgfplots,float,mdframed}
\usepgfplotslibrary{fillbetween}
\usepackage[margin=1in]{geometry}
\usepackage[capitalise,sort]{cleveref} %

\crefname{enumi}{item}{items}
\crefname{equation}{}{}
\crefname{subsection}{Subsection}{Subsections}
\crefname{figure}{Figure}{Figures}

\hypersetup{
	colorlinks,
	linkcolor={blue!80!black},
	citecolor={green},
	urlcolor={blue!80!black},
	linktoc=section
}

\theoremstyle{plain}
\newtheorem{theorem}{Theorem}[section]
\newtheorem{lemma}[theorem]{Lemma}
\newtheorem{prop}[theorem]{Proposition}

\newtheorem{definition}[theorem]{Definition}

\theoremstyle{definition}

\newmdenv[
skipabove=\topsep,
skipbelow=\topsep,
backgroundcolor=gray!10,
linecolor=gray!80,
linewidth=1pt,
innertopmargin=4pt,
innerbottommargin=0.5\baselineskip,
innerleftmargin=0.5\baselineskip,
innerrightmargin=0.5\baselineskip
]{boxedthm}

\newmdenv[
skipabove=\topsep,
skipbelow=\topsep,
backgroundcolor=blue!10,
linecolor=gray!80,
linewidth=1pt,
innertopmargin=4pt,
innerbottommargin=0.5\baselineskip,
innerleftmargin=0.5\baselineskip,
innerrightmargin=0.5\baselineskip
]{boxeddef}

\DeclareMathAlphabet{\mathpzc}{OT1}{pzc}{m}{it}

\DeclareFontEncoding{LS1}{}{}
\DeclareFontSubstitution{LS1}{stix}{m}{n}
\DeclareMathAlphabet{\mathscr}{LS1}{stixscr}{m}{n}

\newcommand{\E}{\mathbb{E}}
\renewcommand{\P}{\mathbb{P}}

\newcommand{\C}{\mathbb{C}}
\newcommand{\R}{\mathbb{R}}
\newcommand{\N}{\mathbb{N}}
\newcommand{\Z}{\mathbb{Z}}

\renewcommand{\c}[1]{\mathfrak{c}^{#1}}

\newcommand{\cF}{\mathcal{F}}

\newcommand{\cK}{\mathcal{K}}

\newcommand{\fC}{\mathfrak{C}}

\newcommand{\fc}{\mathfrak{c}}

\newcommand{\fl}{\mathfrak{l}}

\newcommand{\sr}{\rho\cfadd{def:sr}}
\newcommand{\fd}{d}

\renewcommand{\emptyset}{\varnothing}
\newcommand{\eps}{\varepsilon}

\newcommand{\sprod}{\textstyle\prod}

\newcommand{\m}{\mathbf{m}}

\newcommand{\lp}{\vartheta}

\DeclarePairedDelimiter{\norm}{\lVert}{\rVert}
\DeclarePairedDelimiter{\abs}{\lvert}{\rvert}
\DeclarePairedDelimiter{\rbr}{(}{)}
\DeclarePairedDelimiter{\br}{[}{]}
\DeclarePairedDelimiter{\cu}{\{}{\}}

\makeatletter
\def\@tvsp{\mathchoice{{}\mkern-4.5mu}{{}\mkern-4.5mu}{{}\mkern-2.5mu}{}}
\def\ltrivert{\left|\@tvsp\left|\@tvsp\left|}
\def\rtrivert{\right|\@tvsp\right|\@tvsp\right|}
\makeatother
\newcommand{\mnorm}[1]{\vass{\kern-0.25ex\vass{\kern-0.25ex\vass{#1}\kern-0.25ex}\kern-0.25ex}}
\newcommand{\bmnorm}[1]{\vass[\big]{\kern-0.25ex\vass[\big]{\kern-0.25ex\vass[\big]{#1}\kern-0.25ex}\kern-0.25ex}}
\newcommand{\bbmnorm}[1]{\vass[\bigg]{\kern-0.25ex\vass[\bigg]{\kern-0.25ex\vass[\bigg]{#1}\kern-0.25ex}\kern-0.25ex}}

\makeatletter
\ExplSyntaxOn
\seq_new:N \g_abbrs
\prop_new:N \g_abbr_counts
\tl_new:N \l_abbr_count_tl

\ExplSyntaxOn

\bool_new:N \g_forexample

\NewDocumentCommand{\eg}{ o }{
	\IfValueT{#1}{
		\str_if_eq:noTF {fe} {#1} {
			\bool_gset_true:N \g_forexample
		} {\bool_gset_false:N \g_forexample}
	}
	\bool_if:nTF { \g_forexample } {
		\bool_gset_false:N \g_forexample
		for~example
	}{
		\bool_gset_true:N \g_forexample
		for~instance
	}
}

\NewDocumentCommand{\abbr}{m m O{#1} m m O{#4} m}{
	\expandafter\newcommand\csname#3\endcsname[1][]{
		\seq_if_in:NnTF \g_abbrs {#1} {
			\prop_get:NnN \g_abbr_counts {#1} \l_abbr_count_tl
			\prop_gput:Nnx \g_abbr_counts {#1} {\int_eval:n {\l_abbr_count_tl + 1}}
			\hyperref[#1]{#7}
		} {
			\seq_gput_left:Nn \g_abbrs {#1}
			\prop_gput:Nnn \g_abbr_counts {#1} {1}
			\expandafter\gdef\csname#1@def\endcsname{#2}
			\phantomsection\label{#1}
			\str_if_eq:nnTF{##1}{}{\emph{#2}}{##1}~(\hyperref[#1]{#7})
		}
	}
	\expandafter\newcommand\csname#6\endcsname[1][]{
		\seq_if_in:NnTF \g_abbrs {#1} {
			\prop_get:NnN \g_abbr_counts {#1} \l_abbr_count_tl
			\prop_gput:Nnx \g_abbr_counts {#1} {\int_eval:n {\l_abbr_count_tl + 1}}
			\hyperref[#1]{#4}
		} {
			\expandafter\gdef\csname#1@def\endcsname{#5}
			\seq_gput_left:Nn \g_abbrs {#1}
			\prop_gput:Nnn \g_abbr_counts {#1} {1}
			\phantomsection\label{#1}
			\str_if_eq:nnTF{##1}{}{\emph{#5}}{##1}~(\hyperref[#1]{#4})
		}
	}
}

\ExplSyntaxOff
\makeatother

\abbr{ANN}{artificial neural network}{ANNs}{artificial neural networks}{ANN}
\abbr{DNN}{deep artificial neural network}{DNNs}{deep artificial neural networks}{DNN}
\abbr{SGD}{Stochastic gradient descent}{SGDs}{Stochastic gradient descent}{SGD}
\abbr{GD}{gradient descent}{GDs}{gradient descent}{GD}
\abbr{RMSprop}{root mean square propagation}{RMSprops}{root mean square propagation}{RMSprop}
\abbr{Adam}{adaptive moment estimation}{Adams}{Adaptive moment estimation}{Adam}
\abbr{iid}{independent and identically distributed}{i.i.d}{independent and identically distributed}{i.i.d.}
\abbr{DL}{Deep learning}{DLs}{Deep learning}{DL}
\abbr{AI}{artificial intelligence}{AIs}{artificial intelligence}{AI}
\abbr{LLM}{large language model}{LLMs}{large language models}{LLM}
\abbr{PDE}{partial differantial equation}{PDEs}{partial differantial equations}{PDE}
\abbr{as}{almost surely}{ass}{almost surely}{a.s.}
\abbr{pas}{$\P$-almost surely}{pass}{almost surely}{$\P$-a.s.}
\abbr{OP}{optimization problem}{OPs}{optimization problems}{OP}
\abbr{Adagrad}{adaptive gradient}{Adagrads}{Adaptive gradients}{Adagrad}

\newcommand{\qandq}{\quad \text{and} \quad }
\newcommand{\qqandqq}{\qquad\text{and}\qquad}

\DeclarePairedDelimiter{\vass}{\lvert}{\rvert}
\DeclarePairedDelimiter{\pr}{(}{)}
\DeclarePairedDelimiter{\PR}{[}{]}
\DeclarePairedDelimiter{\pR}{\{}{\}}

\newcommand{\vertiii}[1]{\vass{\kern-0.25ex\vass{\kern-0.25ex\vass{#1}\kern-0.25ex}\kern-0.25ex}}

\newcommand{\prb}[1]{\pr[\big]{ #1 }}
\newcommand{\PRb}[1]{\PR[\big]{ #1 }}
\newcommand{\pRb}[1]{\pR[\big]{ #1 }}
\newcommand{\prbbb}[1]{\pr[\bigg]{ #1 }}
\newcommand{\PRbbb}[1]{\PR[\bigg]{ #1 }}
\newcommand{\pRbbb}[1]{\pR[\bigg]{ #1 }}
\newcommand{\prbb}[1]{\pr[\Big]{ #1 }}
\newcommand{\PRbb}[1]{\PR[\Big]{ #1 }}
\newcommand{\pRbb}[1]{\pR[\Big]{ #1 }}

\newcommand{\indicator}[1]{\mathbbm{1}_{\smash{#1}}}

\newcommand{\diag}{\operatorname{diag}}

\ExplSyntaxOn

\bool_new:N \g_noteobserve

\NewDocumentCommand{\setnote}{}{
	\bool_gset_true:N \g_noteobserve
}

\NewDocumentCommand{\setobserve}{}{
	\bool_gset_false:N \g_noteobserve
}

\NewDocumentCommand{\nobs}{ o }{
	\IfValueT{#1}{
		\str_if_eq:noTF {note} {#1} {
			\bool_gset_true:N \g_noteobserve
		} {
			\str_if_eq:noTF {Note} {#1} {
				\bool_gset_true:N \g_noteobserve
			} {
				\bool_gset_false:N \g_noteobserve
			}
		}
	}
	\bool_if:nTF { \g_noteobserve } {
		\bool_gset_false:N \g_noteobserve
		note
	} {
		\bool_gset_true:N \g_noteobserve
		observe
	}
	\IfValueF{#1}{~}
}

\NewDocumentCommand{\Nobs}{ o }{
	\IfValueT{#1}{
		\str_if_eq:noTF {note} {#1} {
			\bool_gset_true:N \g_noteobserve
		} {
			\str_if_eq:noTF {Note} {#1} {
				\bool_gset_true:N \g_noteobserve
			} {
				\bool_gset_false:N \g_noteobserve
			}
		}
	}
	\bool_if:nTF { \g_noteobserve } {
		\bool_gset_false:N \g_noteobserve
		Note
	} {
		\bool_gset_true:N \g_noteobserve
		Observe
	}
	\IfValueF{#1}{~}
}

\int_new:N \g_furthermore

\NewDocumentCommand{\Moreover}{ o o }{
	\IfValueT{#1}{
		\str_case:nn {#1} {
			{Furthermore} {\int_set:Nn {\g_furthermore} {0}}
			{Moreover} {\int_set:Nn {\g_furthermore} {1}}
			{In~addition} {\int_set:Nn {\g_furthermore} {2}}
			{note} {\bool_gset_true:N \g_noteobserve}
			{observe} {\bool_gset_false:N \g_noteobserve}
		}
		\IfValueT{#2}{
			\str_case:nn {#2} {
				{Furthermore} {\int_set:Nn {\g_furthermore} {0}}
				{Moreover} {\int_set:Nn {\g_furthermore} {1}}
				{In~addition} {\int_set:Nn {\g_furthermore} {2}}
				{note} {\bool_gset_true:N \g_noteobserve}
				{observe} {\bool_gset_false:N \g_noteobserve}
			}
		}
	}
	\int_case:nn { \int_mod:nn {\g_furthermore} {3} } {
		{ 0 } { Furthermore,~\nobs that}
		{ 1 } { Moreover,~\nobs that}
		{ 2 } { In~addition,~\nobs that}
	}
	\int_incr:N \g_furthermore
	\IfValueF{#1}{~}
}

\bool_new:N \g_hencetherefore

\NewDocumentCommand{\hence}{}{
	\bool_if:nTF { \g_hencetherefore } {
		\bool_gset_false:N \g_hencetherefore
		hence~
	} {
		\bool_gset_true:N \g_hencetherefore
		therefore~
	}
}

\NewDocumentCommand{\Hence}{}{
	\bool_if:nTF { \g_hencetherefore } {
		\bool_gset_false:N \g_hencetherefore
		Hence,~we~obtain~
	} {
		\bool_gset_true:N \g_hencetherefore
		Therefore,~we~obtain~
	}
}

\usepackage{cleveref}

\ExplSyntaxOn

\seq_new:N \g_cflist_loaded
\seq_new:N \g_cflist_pending

\NewDocumentCommand{\cfadd} { m } {
	\seq_if_in:NnF \g_cflist_loaded { #1 } {
		\seq_if_in:NnF \g_cflist_pending { #1 } {
			\seq_gput_right:Nn \g_cflist_pending { #1 }
		}
	}
}

\NewDocumentCommand{\cfconsiderloaded} { m } {
	\seq_gput_right:Nn \g_cflist_loaded {#1}
}

\NewDocumentCommand{\cfremove} { m } {
	\seq_gremove_all:Nn \g_cflist_pending { #1 }
}

\NewDocumentCommand{\cfload} { o } {
	\seq_if_empty:NTF \g_cflist_pending {
		\IfValueTF{#1}{\ignorespaces}{\unskip}
	} {
		(cf.\ \cref{\seq_use:Nn \g_cflist_pending {,}})\IfValueTF{#1}{#1~}{\unskip}
		\seq_gconcat:NNN \g_cflist_loaded \g_cflist_loaded \g_cflist_pending
		\seq_gclear:N \g_cflist_pending
		\IfValueT{#1}{\ignorespaces}
	}
}

\NewDocumentCommand{\cfclear} {} {
	\seq_gclear:N \g_cflist_loaded
	\seq_gclear:N \g_cflist_pending
}

\NewDocumentCommand{\cfout} { o } {
	\seq_if_empty:NTF \g_cflist_pending {\unskip\IfValueT{#1}{\ignorespaces}} {
		(cf.\ \cref{\seq_use:Nn \g_cflist_pending {,}})\IfValueTF{#1}{#1~}{\unskip}
		\seq_gclear:N \g_cflist_pending
		\IfValueT{#1}{\ignorespaces}
	}
}

\NewDocumentCommand{\ifnocf} { m } {
	\seq_if_empty:NT \g_cflist_pending { #1 }
}

\ExplSyntaxOff

\ExplSyntaxOff

\ExplSyntaxOn

\seq_const_from_clist:Nn \g_prove_mru {
	establish,
	demonstrate,
	prove,
	show,
	imply,
	ensure
}

\prop_new:N \l__verbs
\prop_put:Nnn \l__verbs {show} {shows}
\prop_put:Nnn \l__verbs {imply} {implies}
\prop_put:Nnn \l__verbs {demonstrate} {demonstrates}
\prop_put:Nnn \l__verbs {prove} {proves}
\prop_put:Nnn \l__verbs {establish} {establishes}
\prop_put:Nnn \l__verbs {ensure} {ensures}
\prop_put:Nnn \l__verbs {assure} {assures}

\tl_new:N \g_wordtmp
\seq_new:N \l_mytmps

\cs_generate_variant:Nn \str_if_in:nnTF { nVTF }
\cs_generate_variant:Nn \str_if_in:nnTF { xVTF }

\NewDocumentCommand{\prove}{ o }{
	\IfValueTF{#1}{
		\seq_clear:N \l_mytmps
		\seq_map_inline:Nn \g_prove_mru {
			\str_if_eq:nnTF {##1} {ensure} {
				\str_set:Nn \l_temps {n}
			} {
				\str_set:Nx \l_temps {\str_head_ignore_spaces:n {##1}}
			}
			\str_if_in:xVTF {#1} \l_temps {
				\seq_put_right:Nn \l_mytmps {##1}
			} { }
		}
		\seq_get_right:NN \l_mytmps \g_wordtmp
	} {
		\seq_get_right:NN \g_prove_mru \g_wordtmp
	}
	\tl_use:N \g_wordtmp
	\IfValueTF{#1}{}{~}
	\seq_gput_left:NV \g_prove_mru \g_wordtmp
	\seq_gremove_duplicates:N \g_prove_mru
}

\NewDocumentCommand{\proves}{ o }{
	\IfValueTF{#1}{
		\seq_clear:N \l_mytmps
		\seq_map_inline:Nn \g_prove_mru {
			\str_if_eq:nnTF {##1} {ensure} {
				\str_set:Nn \l_temps {n}
			} {
				\str_set:Nx \l_temps {\str_head_ignore_spaces:n {##1}}
			}
			\str_if_in:xVTF {#1} \l_temps {
				\seq_put_right:Nn \l_mytmps {##1}
			} { }
		}
		\seq_get_right:NN \l_mytmps \g_wordtmp
	} {
		\seq_get_right:NN \g_prove_mru \g_wordtmp
	}
	\str_set:NV \l_tmpa_str \g_wordtmp
	\prop_get:NVN \l__verbs \l_tmpa_str \l_tmpa_tl
	\tl_use:N \l_tmpa_tl
	\IfValueTF{#1}{}{~}
	\seq_gput_left:NV \g_prove_mru \g_wordtmp
	\seq_gremove_duplicates:N \g_prove_mru
}

\newcommand{\llabel}[1]{\savelabel{#1}\label{\loc.#1}\ignorespaces}

\clist_new:N \l_localreflist
\clist_new:N \l_reflist

\NewDocumentCommand{\lref} { m } {
	\clist_set:No \l_localreflist {#1}
	\clist_clear:N \l_reflist
	\clist_map_inline:Nn \l_localreflist { \clist_put_right:Nn \l_reflist {\loc.##1} }
	\cref{\l_reflist}
}

\NewDocumentCommand{\Lref} { m } {
	\clist_set:No \l_localreflist {#1}
	\clist_clear:N \l_reflist
	\clist_map_inline:Nn \l_localreflist { \clist_put_right:Nn \l_reflist {\loc.##1} }
	\Cref{\l_reflist}
}

\NewDocumentCommand{\itref}{ m m }{
	\clist_set:No \l_localreflist {#2}
	\clist_clear:N \l_reflist
	\clist_map_inline:Nn \l_localreflist { \clist_put_right:Nn \l_reflist {#1.##1} }
	\cref{\l_reflist}~in~\cref{#1}
}

\seq_new:N \l_enum_seq
\int_new:N \l_num_items

\bool_new:N \g_commaused_bool

\providecommand{\comma}{}

\cs_new:Nn \enum_it:nn {
	\int_case:nnF {\l_num_items - #1} {
		{0} {
			\renewcommand{\comma}{}
			#2\space
		}
		{1} {
			\bool_gset_false:N \g_commaused_bool
			\renewcommand{\comma}{,~\bool_gset_true:N \g_commaused_bool}
			#2
			\bool_if:NTF \g_commaused_bool {} {,~}
			and~
		}
	} {
		\bool_gset_false:N \g_commaused_bool
		\renewcommand{\comma}{,~\bool_gset_true:N \g_commaused_bool}
		#2
		\bool_if:NTF \g_commaused_bool {} {,~}
	}
}

\cs_new:Nn \enum_it_U:nn {
	\int_case:nnF {\l_num_items - #1} {
		{0} {
			\renewcommand{\comma}{}
			#2
			\space
		}
		{1} {
			\bool_gset_false:N \g_commaused_bool
			\renewcommand{\comma}{,~\bool_gset_true:N \g_commaused_bool}
			#2
			\bool_if:NTF \g_commaused_bool {} {,~}
			and~
		}
	} {
		\bool_gset_false:N \g_commaused_bool
		\renewcommand{\comma}{,~\bool_gset_true:N \g_commaused_bool}
		\int_compare:nTF {#1=1} {
			\text_titlecase_first:n {#2}
		} {
			#2
		}
		\bool_if:NTF \g_commaused_bool {} {,~}
	}
}

\cs_generate_variant:Nn \tl_if_eq:nnTF {onTF}

\cs_new:Nn \enum:nnnn {
	\seq_set_split:Nnn \l_enum_seq ; {#1}
	\seq_remove_all:Nn \l_enum_seq { }
	\seq_remove_all:Nn \l_enum_seq {#2}
	\seq_log:N \l_enum_seq
	\int_set:Nn \l_num_items {\seq_count:N \l_enum_seq}
	\int_log:N \l_num_items
	\int_case:nnF {\l_num_items} {
		{ 0 } { 0 }
		{ 1 } {
			\IfBooleanTF{#4} {
				\tl_set:Nn \l_text_case_exclude_arg_tl {\cref}
				\bool_set_false:N \l_text_titlecase_check_letter_bool
				\text_titlecase_first:n {\seq_use:Nn \l_enum_seq {}}
			} {
				\seq_use:Nn \l_enum_seq {}
			}
			\space
			\tl_if_eq:onTF{#3}{-}{}{
				\bool_if:NTF \l_plural_bool {
					\prove[#3]~
				} {
					\proves[#3]~
				}
			}
		}
		{ 2 } {
			\IfBooleanTF{#4} {
				\tl_set:Nn \l_text_case_exclude_arg_tl {\cref}
				\bool_set_false:N \l_text_titlecase_check_letter_bool
				\text_titlecase_first:n {\seq_item:Nn \l_enum_seq {1}}
				{} ~and~
				\seq_item:Nn \l_enum_seq {2}
			} {
				\seq_use:Nn \l_enum_seq {~and~}
			}
			\space
			\tl_if_eq:onTF{#3}{-}{}{
				\prove[#3]~
			}
		}
	} {
		\IfBooleanTF{#4} {
			\tl_set:Nn \l_text_case_exclude_arg_tl {\cref}
			\seq_indexed_map_function:NN \l_enum_seq \enum_it_U:nn
		} {
			\seq_indexed_map_function:NN \l_enum_seq \enum_it:nn
		}
		\tl_if_eq:onTF{#3}{-}{}{
			\prove[#3]~
		}
	}
}

\cs_generate_variant:Nn \enum:nnnn {nxnn}
\cs_generate_variant:Nn \enum:nnnn {nxxn}

\NewDocumentCommand{\enum}{O{} m O{-} s}{
	\IfBooleanTF{#4}{
		\enum:nxnn {#2} {#1} {sindep} \BooleanFalse
	} {
		\enum:nxxn {#2} {#1} {#3} \BooleanFalse
	}
}

\bool_new:N \g_arg_start_bool
\bool_gset_true:N \g_arg_start_bool

\NewDocumentCommand{\startnewargseq}{}{\bool_gset_true:N \g_arg_start_bool \tl_set:Nn \g_label_tl {}}

\cs_generate_variant:Nn \seq_if_in:NnTF {NxTF}
\cs_generate_variant:Nn \seq_remove_all:Nn {Nx}

\int_new:N \l_random_int

\cs_generate_variant:Nn \tl_if_head_eq_catcode:nNTF {oNTF}
\cs_generate_variant:Nn \tl_if_head_eq_catcode:nNTF {VNTF}

\cs_generate_variant:Nn \tl_if_head_eq_catcode:nNTF {eNTF}
\cs_generate_variant:Nn \tl_log:n {o}
\cs_generate_variant:Nn \tl_log:n {f}
\cs_generate_variant:Nn \tl_log:n {x}
\cs_generate_variant:Nn \tl_log:n {e}

\cs_generate_variant:Nn \tl_if_in:nnTF {onTF}
\cs_generate_variant:Nn \tl_if_in:NnTF {NeTF}

\cs_generate_variant:Nn \tl_if_head_eq_meaning:nNTF {VNTF}

\seq_const_from_clist:Nn \g_arg_mru_this {
	Ahpr,
	Tapr,
	Ctapr,
	H
}

\seq_const_from_clist:Nn \g_arg_mru_nothis {
	Ia,
	Nwc,
	N,
	Itns,
	Fm,
	Itnswc,
	Mo
}

\seq_new:N \l_arg_seq
\tl_new:N \l_cons_tl
\tl_new:N \l_dummy_tl

\bool_new:N \g_debug_bool
\bool_gset_false:N \g_debug_bool

\bool_new:N \l_insidearg_bool

\bool_new:N \g_firstargletter_bool

\sys_gset_rand_seed:n {0903}

\bool_new:N \l_plural_bool
\tl_new:N \l_arg_verbs_tl

\NewDocumentCommand{\argument}{mom}{
	\bool_set_false:N \l_plural_bool
	\tl_set:Nn \l_arg_verbs_tl {sindep}
	\keys_define:nn { benno/argument } {
		plural .value_forbidden:n = true,
		plural .code:n = {\bool_set_true:N \l_plural_bool},
		verbs .value_required:n = false,
		verbs .tl_set:N = \l_arg_verbs_tl,
	}
	\IfValueT{#2}{
		\keys_set:nn { benno/argument } {#2}
	}
	\bool_log:N \l_plural_bool
	\bool_gset_true:N \l_insidearg_bool
	\seq_set_split:Nnn \l_arg_seq ; {#1}
	\seq_remove_all:Nn \l_arg_seq { }
	\seq_log:N \l_arg_seq
	\tl_set:Nn \l_cons_tl {#3}
	\tl_trim_spaces:N \l_cons_tl
	\seq_if_in:NxTF \l_arg_seq {\lref{\g_label_tl}} {
		\seq_remove_all:Nx \l_arg_seq {\lref{\g_label_tl}}
		\seq_get_left:NNTF \l_arg_seq \l_dummy_tl {
			\tl_trim_spaces:N \l_dummy_tl
			\bool_gset_false:N \g_firstargletter_bool
			\tl_if_head_eq_catcode:nNTF \l_dummy_tl a {
				\bool_gset_true:N \g_firstargletter_bool
			} {
				\tl_if_head_eq_meaning:VNTF \l_dummy_tl {\cref} {
					\tl_set:Nx \l_tmpa_tl {\tl_tail:N \l_dummy_tl}
					\tl_set:Nx \l_tmpb_tl {\tl_head:N \l_tmpa_tl}
					\bool_gset_true:N \g_firstargletter_bool
					\tl_if_in:NeTF \l_tmpb_tl {lem\c_colon_str} {} {
						\tl_if_in:NeTF \l_tmpb_tl {thm\c_colon_str} {} {
							\tl_if_in:NeTF \l_tmpb_tl {prop\c_colon_str} {} {
								\tl_if_in:NeTF \l_tmpb_tl {cor\c_colon_str} {} {
									\bool_gset_false:N \g_firstargletter_bool
								}
							}
						}
					}
				} {
				}
			}
			\bool_if:NTF \g_firstargletter_bool {
				\seq_set_eq:NN \l_tmpa_seq \g_arg_mru_this
				\seq_remove_all:Nn \l_tmpa_seq {H}
				\seq_get_right:NN \l_tmpa_seq \l_tmpa_tl
				\int_case:nnF {\seq_count:N \l_arg_seq} {
					{1} {
						\str_case:VnF {\l_tmpa_tl} {
							{Ahpr} {
								\bool_if:NT \g_debug_bool {C1.1}
								\seq_gput_left:Nn \g_arg_mru_this {Ahpr}
								\seq_gremove_duplicates:N \g_arg_mru_this
								\enum:nxnn {#1} {\lref{\g_label_tl}} {-} {\BooleanTrue}
								\hence~
								\bool_if:NTF \l_plural_bool {
									\prove[\l_arg_verbs_tl]~\ignorespaces #3
								} {
									\proves[\l_arg_verbs_tl]~\ignorespaces #3
								}
							}
							{Tapr} {
								\bool_if:NT \g_debug_bool {C1.2}
								\seq_gput_left:Nn \g_arg_mru_this {Tapr}
								\seq_gremove_duplicates:N \g_arg_mru_this
								\enum[\lref{\g_label_tl}]{
									This;
									#1
								}[\l_arg_verbs_tl]\ignorespaces #3
							}
							{Ctapr} {
								\bool_if:NT \g_debug_bool {C1.3}
								\seq_gput_left:Nn \g_arg_mru_this {Ctapr}
								\seq_gremove_duplicates:N \g_arg_mru_this
								Combining~
								\enum[\lref{\g_label_tl}]{
									this;
									#1
								} \proves[\l_arg_verbs_tl]~\ignorespaces #3
							}
						} {}
					}
				} {
					\str_case:VnF {\l_tmpa_tl} {
						{Ahpr} {
							\bool_if:NT \g_debug_bool {C2.1}
							\seq_gput_left:Nn \g_arg_mru_this {Ahpr}
							\seq_gremove_duplicates:N \g_arg_mru_this
							\enum:nxnn {#1} {\lref{\g_label_tl}} {-} {\BooleanTrue}
							\hence~
							\prove[\l_arg_verbs_tl]~\ignorespaces #3
						}
						{Tapr} {
							\bool_if:NT \g_debug_bool {C2.2}
							\seq_gput_left:Nn \g_arg_mru_this {Tapr}
							\seq_gremove_duplicates:N \g_arg_mru_this
							\enum[\lref{\g_label_tl}]{
								This;
								#1
							}[\l_arg_verbs_tl]\ignorespaces #3
						}
						{Ctapr} {
							\int_case:nn {\int_rand:nn {0} {1}} {
								{0} {
									\bool_if:NT \g_debug_bool {C2.3}
									\seq_gput_left:Nn \g_arg_mru_this {Ctapr}
									\seq_gremove_duplicates:N \g_arg_mru_this
									Combining~
									\enum[\lref{\g_label_tl}]{
										this;
										#1
									} \proves[\l_arg_verbs_tl]~\ignorespaces #3
								}
								{1} {
									\bool_if:NT \g_debug_bool {C2.4}
									\seq_gput_left:Nn \g_arg_mru_this {Ctapr}
									\seq_gremove_duplicates:N \g_arg_mru_this
									Combining~
									\enum:nxnn {#1} {\lref{\g_label_tl}} {-} {\BooleanFalse}
									\hence~
									\proves[\l_arg_verbs_tl]~\ignorespaces #3
								}
							}
						}
					} {}
				}
			} {
				\seq_set_eq:NN \l_tmpa_seq \g_arg_mru_this
				\seq_remove_all:Nn \l_tmpa_seq {H}
				\seq_remove_all:Nn \l_tmpa_seq {Ahpr}
				\seq_get_right:NN \l_tmpa_seq \l_tmpa_tl
				\int_case:nnF {\seq_count:N \l_arg_seq} {
					{1} {
						\str_case:VnF {\l_tmpa_tl} {
							{Tapr} {
								\bool_if:NT \g_debug_bool {C3.1}
								\seq_gput_left:Nn \g_arg_mru_this {Tapr}
								\seq_gremove_duplicates:N \g_arg_mru_this
								\enum[\lref{\g_label_tl}]{
									This;
									#1
								}[\l_arg_verbs_tl]\ignorespaces #3
							}
							{Ctapr} {
								\bool_if:NT \g_debug_bool {C3.2}
								\seq_gput_left:Nn \g_arg_mru_this {Ctapr}
								\seq_gremove_duplicates:N \g_arg_mru_this
								Combining~
								\enum[\lref{\g_label_tl}]{
									this;
									#1
								} \proves[\l_arg_verbs_tl]~\ignorespaces #3
							}
						} {}
					}
				} {
					\str_case:VnF {\l_tmpa_tl} {
						{Tapr} {
							\bool_if:NT \g_debug_bool {C4.1}
							\seq_gput_left:Nn \g_arg_mru_this {Tapr}
							\seq_gremove_duplicates:N \g_arg_mru_this
							\enum[\lref{\g_label_tl}]{
								This;
								#1
							}[\l_arg_verbs_tl]\ignorespaces #3		
						}
						{Ctapr} {
							\int_case:nn {\int_rand:nn {0} {1}} {
								{0} {
									\bool_if:NT \g_debug_bool {C4.2}
									\seq_gput_left:Nn \g_arg_mru_this {Ctapr}
									\seq_gremove_duplicates:N \g_arg_mru_this
									Combining~
									\enum[\lref{\g_label_tl}]{
										this;
										#1
									} \proves[\l_arg_verbs_tl]~\ignorespaces #3		
								}
								{1} {
									\bool_if:NT \g_debug_bool {C4.3}
									\seq_gput_left:Nn \g_arg_mru_this {Ctapr}
									\seq_gremove_duplicates:N \g_arg_mru_this
									Combining~
									\enum:nxnn {#1} {\lref{\g_label_tl}} {-} {\BooleanFalse}
									\hence~
									\proves[\l_arg_verbs_tl]~\ignorespaces #3    
								}
							}
						}
					} {}
				}
			}
		} {
			\tl_if_head_eq_catcode:oNTF \l_cons_tl a {
				\seq_set_eq:NN \l_tmpa_seq \g_arg_mru_this
				\seq_remove_all:Nn \l_tmpa_seq {Ctapr}
				\seq_remove_all:Nn \l_tmpa_seq {Ahpr}
				\seq_get_right:NN \l_tmpa_seq \l_tmpa_tl
				\str_case:VnF {\l_tmpa_tl} {
					{H} {
						\bool_if:NT \g_debug_bool {C5.1}
						\seq_gput_left:Nn \g_arg_mru_this {H}
						\seq_gremove_duplicates:N \g_arg_mru_this
						Hence,~we~obtain~\ignorespaces #3
					}
					{Tapr} {
						\bool_if:NT \g_debug_bool {C5.2}
						\seq_gput_left:Nn \g_arg_mru_this {Tapr}
						\seq_gremove_duplicates:N \g_arg_mru_this
						This~\proves[\l_arg_verbs_tl]~\ignorespaces #3
					}
				} {}
			} {
				\bool_if:NT \g_debug_bool {C6.1}
				\seq_gput_left:Nn \g_arg_mru_this {Tapr}
				\seq_gremove_duplicates:N \g_arg_mru_this
				This~\proves[\l_arg_verbs_tl]~\ignorespaces #3
			}
		} 
	} {
		\int_compare:nNnTF {\seq_count:N \l_arg_seq} = {0} {
			\bool_if:NTF \g_arg_start_bool {
				\bool_if:NT \g_debug_bool {C7.1}
				\Nobs\unskip
				#3
			} {
				\bool_if:NT \g_debug_bool {C7.2}
				\Moreover~
				#3
			}
		} {
			\bool_if:NTF \g_arg_start_bool {
				\bool_if:NT \g_debug_bool {C8.1}
				\tl_log:N \l_arg_verbs_tl
				\Nobs~that~
				\enum{
					#1
				}[\l_arg_verbs_tl]\ignorespaces #3
			} {
				\int_compare:nNnTF {\seq_count:N \l_arg_seq} = {1} {
					\seq_set_eq:NN \l_tmpa_seq \g_arg_mru_nothis
					\seq_remove_all:Nn \l_tmpa_seq {Nwc}
					\seq_remove_all:Nn \l_tmpa_seq {Itnswc}
					\seq_get_right:NN \l_tmpa_seq \l_tmpa_tl
				} {
					\seq_get_right:NN \g_arg_mru_nothis \l_tmpa_tl
				}
				\str_case:VnF {\l_tmpa_tl} {
					{Mo} {
						\bool_if:NT \g_debug_bool {C9.1}
						\seq_gput_left:Nn \g_arg_mru_nothis {Mo}
						\seq_gremove_duplicates:N \g_arg_mru_nothis
						Moreover,~\nobs~that~
						\enum{
							#1
						}[\l_arg_verbs_tl]\ignorespaces #3		
					}
					{Fm} {
						\bool_if:NT \g_debug_bool {C9.2}
						\seq_gput_left:Nn \g_arg_mru_nothis {Fm}
						\seq_gremove_duplicates:N \g_arg_mru_nothis
						Furthermore,~\nobs~that~
						\enum{
							#1
						}[\l_arg_verbs_tl]\ignorespaces #3		
					}
					{Ia} {
						\bool_if:NT \g_debug_bool {C9.3}
						\seq_gput_left:Nn \g_arg_mru_nothis {Ia}
						\seq_gremove_duplicates:N \g_arg_mru_nothis
						In~addition,~\nobs~that~
						\enum{
							#1
						}[\l_arg_verbs_tl]\ignorespaces #3		
					}
					{N} {
						\bool_if:NT \g_debug_bool {C9.4}
						\seq_gput_left:Nn \g_arg_mru_nothis {N}
						\seq_gremove_duplicates:N \g_arg_mru_nothis
						Next,~\nobs~that~
						\enum{
							#1
						}[\l_arg_verbs_tl]\ignorespaces #3		
					}
					{Itns} {
						\bool_if:NT \g_debug_bool {C9.5}
						\seq_gput_left:Nn \g_arg_mru_nothis {Itnswc}
						\seq_gput_left:Nn \g_arg_mru_nothis {Itns}
						\seq_gremove_duplicates:N \g_arg_mru_nothis
						In~the~next~step~we~\nobs~that~
						\enum{
							#1
						}[\l_arg_verbs_tl]\ignorespaces #3		
					}
					{Nwc} {
						\bool_if:NT \g_debug_bool {C9.6}
						\seq_gput_left:Nn \g_arg_mru_nothis {Nwc}
						\seq_gremove_duplicates:N \g_arg_mru_nothis
						Next~we~combine~
						\enum{
							#1
						}to~obtain~\ignorespaces #3
					}
					{Itnswc} {
						\bool_if:NT \g_debug_bool {C9.7}
						\seq_gput_left:Nn \g_arg_mru_nothis {Itns}
						\seq_gput_left:Nn \g_arg_mru_nothis {Itnswc}
						\seq_gremove_duplicates:N \g_arg_mru_nothis
						In~the~next~step~we~combine~
						\enum{
							#1
						}to~obtain~\ignorespaces #3
					}
				} {}
			}
		}
	}
	\bool_gset_false:N \g_arg_start_bool
	\bool_gset_false:N \l_insidearg_bool
	\cfload[.]
	\color{black}
}

\tl_new:N \g_label_tl
\tl_gset:Nn \g_label_tl { }

\NewDocumentCommand{\savelabel}{m}{
	\bool_if:NTF \l_insidearg_bool {
		\tl_gset:Nn \g_label_tl {#1}
	} {
		\tl_gset:Nn \g_label_tl { }
	}
}

\ExplSyntaxOff

\ExplSyntaxOn

\NewDocumentEnvironment {athm} {m m o} {
	\str_if_eq:noTF {example} {#1} {
		\bool_gset_true:N \g_example_bool
	} {
		\bool_gset_false:N \g_example_bool
	}
	\cfclear
	\IfNoValueTF{#3}{
		\begin{#1}\label{#2}\global\def\loc{#2}
		}{
			\begin{#1}[#3]\label{#2}\global\def\loc{#2}
			}
		}{
		\end{#1}
	}
	
	\NewDocumentEnvironment {adef} {m} {
		\begin{definition}\label{#1}\global\def\loc{#1}
		}{
		\end{definition}
	}
	
	\NewDocumentEnvironment{aproof} {} {
		\bool_if:NTF \g_example_bool {
			\bool_gset_true:N \g_arg_start_bool
			\begin{proof}[Proof~for~\cref{\loc}]
			} {
				\bool_gset_true:N \g_arg_start_bool
				\begin{proof}[Proof~of~\cref{\loc}]
				}
				\bool_gset_false:N \g_finishproof_bool
			}{
				\bool_if:NTF \g_finishproof_bool {}
				{\finishproofthus}
			\end{proof}
		}
		
		\NewDocumentCommand{\finishproofthus} {} {
			\bool_gset_true:N \g_finishproof_bool 
			\bool_if:NTF \g_example_bool {
				The~proof~for~\cref{\loc}~is~thus~complete.
			} {
				The~proof~of~\cref{\loc}~is~thus~complete.
			}
		}
		\NewDocumentCommand{\finishproofthis} {} {
			\bool_gset_true:N \g_finishproof_bool 
			\bool_if:NTF \g_example_bool {
				This~completes~the~proof~for~\cref{\loc}.
			} {
				This~completes~the~proof~of~\cref{\loc}.
			}
		}
		
		\ExplSyntaxOff
		
		\ExplSyntaxOff
		
		\NewDocumentEnvironment{cproof}{m}
		{\begin{proof}[Proof of \cref{#1}]}%
			{\noindent The proof of \cref{#1} is thus complete.
		\end{proof}}
		
		\NewDocumentEnvironment{cproof2}{m}
		{\begin{proof}[Proof of \cref{#1}]}%
			{\noindent This completes the proof of \cref{#1}.
		\end{proof}}

		\ExplSyntaxOff

		\ExplSyntaxOn

		\ExplSyntaxOff

		\newcommand{\stable}[1]{is $#1$-asymptotically stable \cfadd{def:asymptotically:stable}}
		\newcommand{\sstable}[1]{is strongly $#1$-asymptotically stable \cfadd{def:strongly:asymptotically:stable}}
		\newcommand{\ssstable}[1]{is super strongly $#1$-asymptotically stable \cfadd{def:super:strongly:asymptotically:stable}}
		\newcommand{\ustable}{is uniformly stable \cfadd{def:uniformly:stable}}
		\newcommand{\grad}[2]{g_{#1}^{(#2)}}
		
		\newcommand{\mom}{\mathbb{M}}
		\newcommand{\MOM}{\mathbb{V}}
		\newcommand{\MOMs}{\mathfrak{M}}
		\newcommand{\batch}{J}
		\newcommand{\Cst}{\mu}
		\newcommand{\cst}{\nu}
		
		\newcommand{\pars}{\mathfrak{d}}
		
		\newcommand{\bscl}{\mathfrak{b}}
		\newcommand{\loss}{\ell}
		\newcommand{\LR}{\Gamma}
		\newcommand{\GRAD}{G}
		\newcommand{\SPR}{\varTheta}

		\title{Asymptotic stability properties and a priori bounds for Adam and other gradient descent optimization methods}

		\author{Steffen Dereich$^{1}$, Robin Graeber$^{2}$, Arnulf Jentzen$^{3,4}$, and Adrian Riekert$^5$
				\bigskip
			\\
			\small{$^1$ Applied Mathematics: Institute for Mathematical Stochastics, Faculty of Mathematics and}
			\vspace{-0.1cm}\\
			\small{Computer Science, University of M{\"u}nster, Germany, e-mail: \texttt{steffen.dereich@uni-muenster.de}}
			\smallskip
			\\
			\small{$^2$ School of Data Science,	The Chinese University of Hong Kong, Shenzhen}
			\vspace{-0.1cm}\\
			\small{(CUHK-Shenzhen), China, e-mail:
				\texttt{223040041@link.cuhk.edu.cn}}
			\smallskip
			\\
			\small{$^3$ School of Data Science and School of Artificial Intelligence, The Chinese University}
			\vspace{-0.1cm}\\
			\small{of Hong Kong, Shenzhen (CUHK-Shenzhen), China, e-mail: \texttt{ajentzen@cuhk.edu.cn}}
			\smallskip
			\\
			\small{$^4$ Applied Mathematics: Institute for Analysis and Numerics, Faculty of Mathematics and}
			\vspace{-0.1cm}\\
			\small{Computer Science, University of M{\"u}nster, Germany, e-mail: \texttt{ajentzen@uni-muenster.de}}
			\smallskip
			\\
			\small{$^5$ Applied Mathematics: Institute for Analysis and Numerics, Faculty of Mathematics and}
			\vspace{-0.1cm}\\
			\small{Computer Science, University of M{\"u}nster, Germany, e-mail: \texttt{ariekert@uni-muenster.de}}
			\smallskip
			\\
		}
		
\begin{document}
		
		\maketitle
		
		\begin{abstract}
		Gradient descent (GD) based optimization methods are these days the standard tools to train deep neural networks in artificial intelligence systems.
		In optimization procedures in deep learning the employed optimizer is often not the standard GD method but instead suitable adaptive and accelerated variants of standard GD (including the momentum and the root mean square propagation (RMSprop) optimizers) are considered.
		The adaptive moment estimation (Adam) optimizer proposed in 2014 by Kingma \& Ba is presumably the most popular variant of such adaptive and accelerated GD based optimization methods.
		Despite the popularity of such sophisticated optimization methods, it remains a fundamental open problem of research to provide a rigorous mathematical analysis for such accelerated and adaptive optimization methods.
		In particular, it remains an open problem of research to establish boundedness of the Adam optimizer. In this work we solve this problem in the case of a simple class of quadratic strongly convex stochastic optimization problems.
		Specifically, for the considered class of stochastic optimization problems we reveal a priori bounds for momentum, RMSprop, and Adam.
		In particular, we prove for the considered class of strongly convex stochastic optimization problems, for the first time, that Adam does \emph{not explode} but \emph{stays bounded}
		for any choices of the learning rates.
		In this work we also introduce certain stability concepts -- such as the notion of the \emph{stability region} -- for deep learning optimizers and we discover that among standard GD, momentum, RMSprop, and Adam we have that Adam is the \emph{only} optimizer that achieves the optimal higher order convergence speed and also has the \emph{maximal stability region}. 
		Furthermore, we prove that the stability region of Nesterov momentum is strictly smaller than the stability region of standard GD, that the stability region of standard GD is strictly smaller than the stability region of momentum, and that the stability region of momentum is strictly smaller than the stability region of RMSprop and Adam, which both have the maximal stability region.
		\end{abstract}
		
		\newpage
		\tableofcontents

		\section{Introduction}
		\label{section:1}
		\SGD\ optimization schemes are nowdays the standard instuments to train \ANNs\ in deep learning.
		Often not the plain vanilla standard \SGD\ method is the employed optimization scheme but instead suitable adaptive and/or accelerated \SGD\ methods such as 
		the momentum optimizer \cite{Poljak_momentum_SGD}, 
		the Nesterov optimizer \cite{Nesterov1983}, 
		and
		the \RMSprop\ optimizer \cite{Hinton24_RMSprop} are considered
		(cf., \eg, 
		\cite{ruder2017overviewgradientdescentoptimization}, 
		\cite[Chapters 5, 6, and 7]{jentzen2023mathematical}, and
		\cite[Chapters 4, 5, 6, and 8]{DeepLearningBook}).
		The most popular variant of such adaptive and/or accelerated \SGD\ methods is the \Adam\ optimizer introduced in 2014 by Kingma and Ba \cite{KingmaBa2024_Adam}. 
		%
		%
		Despite the popularity and practical relevance of these methods in the training of \AI\ systems, an open research question is to provide a convergence analysis for such sophisticated adaptive and accelerated \SGD\ methods, even for the simplest class of quadratic strongly convex stochastic \OPs.
		In particular, even in the situation of such a simple class of quadratic \OPs, it remains an open research question to show that such adaptive and accelerated \SGD\ methods do not escape to \emph{infinity} but stay bounded and satisfy a priori moment bounds when applied to such \OPs.
		
		It is precisely the topic of this project to answer this question and to establish uniform a priori bounds for such adaptive and accelerated \SGD\ optimization methods such as \Adam\ and, in general, to develop 
		a theory of stability properties for such methods.
		More specifically, in \cref{thm:stoch:ADAM:bounded} below we establish explicit uniform a priori bounds for the \Adam\ and the \RMSprop\ optimizers applied to such simple quadratic stochastic \OPs.
		We illustrate this contribution within this introductory section by means of \cref{thm:stoch:ADAM:bounded} below, which is a direct consequence of \cref{cor:stoch:ADAM:bounded} in \cref{section:5} below.
		
		In addition, we propose different \emph{stability concepts} for adaptive and accelerated optimization methods and, in particular, in \cref{def:stab:region} below we propose the concept, which we refer to as \emph{stability region}, characterizing the set of all possible values of learning rates and eigenvalues of the Hessian of the objective function such that the considered optimization method \emph{does not diverge} to infinity but \emph{stays bounded}.
		We illustrate this theory by explicitly specifying the stabilitiy regions for 
		the standard \GD\ optimizer, 
		the Nesterov optimizer, 
		the momentum optimizer, 
		the \RMSprop\ optimizer, and
		the \Adam\ optimizer in \cref{intro:thm:1.1} below.
		After this brief informal description of the contributions of this work, we now introduce in \cref{def:stab:region}
		the notion of the stability region with all mathematical details and, thereafter, we explain the proposed concept in words.

	\subsection{Introduction of the notion of the stability region}
	\label{subsection:1:1}
	The natural number $d\in\N=\{1,2,3,\dots\}$ in \cref{def:stab:region} specifies the dimension of the \OP\ under consideration (the number of parameters/degrees of freedom that need to be optimized),
	the functions $\Phi_n\colon(\R^d)^{n}\to\R^d$, $n \in \N$, in \cref{def:stab:region} specify the optimization method under consideration,
	and
	the set $\mathcal{A}\subseteq[0,\infty)^{d+1}$ in \cref{def:stab:region} serves as the stability region.
	
			\renewcommand{\thempfootnote}{\arabic{mpfootnote}}
		\begin{definition}[Stability region]
			\begin{boxeddef}
				\label{def:stab:region}
				Let $d \in \N$,
				let $\Phi_n\colon(\R^d)^{n}\to\R^d$, $n \in \N$, be functions,
				and let $\mathcal{A}\subseteq [0,\infty)^{d+1}$ be a set.
				Then we say that the stability region of $(\Phi_n)_{n\in \N}$ is $\mathcal{A}$ if and only if it holds\footnotemark\ for all 
				$\gamma\in[0,\infty)$, 
				$\lambda=(\lambda_1,\dots,\lambda_d)\in[0,\infty)^d$,
				$\vartheta\in\R^d$ and all 
				$\Theta\colon\N_0  \to \R^d$ with
				\begin{equation}
					\label{eq:def:stab:reg:1}
					\forall\, n \in \N\colon \Theta_{n}=\Theta_{n-1}
					-
					\gamma\Phi_n\prb{
						\operatorname{diag}(\lambda)(\Theta_0-\vartheta),
						\operatorname{diag}(\lambda)(\Theta_1-\vartheta),
						\dots,
						\operatorname{diag}(\lambda)(\Theta_{n-1}-\vartheta)}
				\end{equation}
				that
				\begin{equation}
					\label{eq:def:stab:reg:2}
					\limsup_{n\to \infty} \norm{\Theta_n}
					\in
					\begin{cases}
						\R & \colon (\gamma,\lambda_1,\lambda_2,\dots,\lambda_d)\in \mathcal{A}\\
						\{\infty\} & \colon (\gamma,\lambda_1,\lambda_2,\dots,\lambda_d)\notin \mathcal{A}.
					\end{cases}
				\end{equation}
			\end{boxeddef}
		\end{definition}
		\footnotetext{Note that for all $d \in \N$, $x=(x_1, \dots, x_d)\in \C^d$ it holds that 
			$\norm{x} = ( \sum_{i=1}^d \abs{x_i}^2 ) ^{1/2}$ (standard norm)
			and note that for all
			$d \in \N$, $x=(x_1, \dots, x_d)$, $y=(y_1, \dots, y_d)\in \R^d$ it holds that
			$\operatorname{diag}(x) y = (x_1 y_1, x_2 y_2, \dots, x_d y_d )$ (diagonal matrix associated to a vector).}

		A large class of deep learning optimizers can be formulated using the functions $\Phi_n\colon(\R^d)^{n}\to\R^d$, $n \in \N$, in \cref{def:stab:region} (cf., \eg, \cite[Definitions~1.1, 2.1, 2.2, 3.1, and 4.1]{dereich2025sharphigherorderconvergence} and \cite[Sections~6.4 and 6.5]{MR4688424}).
		For example, in the case of the standard \GD\ optimization we have for all $n \in \N$, $g_1,g_2,\dots,g_n\in\R^d$ that $\Phi_n(g_1,g_2,\dots,g_n)=g_n$ (cf.\ \cref{def:SGD}).
		Similarly,
		the momentum optimizer (cf.\ \cref{def:MOM}), 
		the Nesterov optimizer (cf.\ \cref{def:Nest}), 
		the \mbox{\RMSprop} optimizer (cf.\ \cref{def:RMSprop}),
		and the \Adam\ optimizer (cf.\ \cref{def:ADAM}) can be described in a full history recursive manner (cf.\ \cref{eq:def:stab:reg:1} above) by employing such functions $\Phi_n\colon(\R^d)^{n}\to\R^d$, $n\in\N$, from \cref{def:stab:region} above.
		Note that in \cref{def:stab:region} we consider the \OP\ to approximately compute the global minimizer $\vartheta=(\vartheta_1,\dots,\vartheta_d) \in \R^d$ of the minimization problem 
		\begin{equation}
			\label{eq:min:prob}
			\textstyle
			\min_{\theta=(\theta_1,\dots,\theta_d)\in\R^d}\prb{\sum_{i=1}^d \frac{\lambda_i}{2} (\theta_i-\vartheta_i)^2}
		\end{equation}
		(cf.\ \eqref{eq:def:stab:reg:1} above).
		We note that \cref{def:stab:region} assures that the stability region of a deep learning optimizer is a subset of $[0,\infty)^{d+1}$ that contains exactly those tuples of learning rates and eigenvalues of the Hessian of the objective function in \eqref{eq:min:prob} 
		such that the optimization process does not have a divergent subsequence but stays bounded.

		\subsection{Stability regions for deep learning optimizers}
		After having presented the notion of stability region in \cref{def:stab:region} above we are now in the position to state \cref{intro:thm:1.1} that explicitly specifies the stability region for 
		the standard \GD\ optimizer, 
		the Nesterov optimizer, 
		the momentum optimizer, 
		the \RMSprop\ optimizer, and
		the \Adam\ optimizer.
		
		\begin{athm}{theorem}{intro:thm:1.1}[Stability for optimizers]
			\begin{boxedthm}
				Let $d\in\N$, $\alpha\in[0,1)$, $\beta\in(\alpha^2,1)$, $\eps\in(0,\infty)$. Then
				\begin{enumerate}[label=(\roman*)]
					\item \llabel{it:NEST}\label{it:Nest:perm}
					it holds for every $\alpha$-Nesterov optimizer $\Phi_n\colon(\R^d)^{n}\to\R^d$, $n \in \N$,\cfadd{def:Nest} 
					that the stability region\cfadd{def:stab:region} of $(\Phi_n)_{n\in \N}$ is
					\begin{equation}
						\label{eq:stab:Nest:intro}
						\color{magenta}
						\pRb{\lambda=(\lambda_0,\lambda_1,\dots,\lambda_{d})\in[0,\infty)^{d+1}\colon \textstyle\max_{i\in\{1,2,\dots,d\}}\pr{\lambda_0\lambda_{i}}\leq2\PRb{\frac{1-\alpha^2}{1+2\alpha}}}\color{black},
					\end{equation}
					\vspace{-0.7cm}
					\item \llabel{it:SGD}\label{it:SGD:perm}
					it holds
					for every \GD\ optimizer $\Phi_n\colon(\R^d)^{n}\to\R^d$, $n \in \N$,\cfadd{def:SGD} 
					that the stability region\cfadd{def:stab:region} of $(\Phi_n)_{n\in \N}$ is
					\begin{equation}
						\color{magenta}
						\pRb{\lambda=(\lambda_0,\lambda_1,\dots,\lambda_{d})\in[0,\infty)^{d+1}\colon \textstyle\max_{i\in\{1,2,\dots,d\}}\pr{\lambda_0\lambda_{i}}\leq2}\color{black},
					\end{equation}
					\vspace{-0.7cm}
					\item \llabel{it:MOM}\label{it:MOM:perm}
					it holds for every $\alpha$-momentum optimizer $\Phi_n\colon(\R^d)^{n}\to\R^d$, $n \in \N$,\cfadd{def:MOM} 
					that the stability region\cfadd{def:stab:region} of $(\Phi_n)_{n\in \N}$ is
					\begin{equation}
						\label{eq:stab:MOM:intro}
						\color{magenta}
						\pRb{\lambda=(\lambda_0,\lambda_1,\dots,\lambda_{d})\in[0,\infty)^{d+1}\colon \textstyle\max_{i\in\{1,2,\dots,d\}}\pr{\lambda_0\lambda_{i}}\leq2\PRb{\frac{1+\alpha}{1-\alpha}}}\color{black},
					\end{equation}
					\vspace{-0.7cm}
					\item \llabel{it:RMS}\label{it:RMS:perm}
					it holds for every $\beta$-$\eps$-\RMSprop\ optimizer $\Phi_n\colon(\R^d)^{n}\to\R^d$, $n \in \N$,\cfadd{def:RMSprop}
					that the stability region\cfadd{def:stab:region} of $(\Phi_n)_{n\in \N}$ is
					\color{magenta}$[0,\infty)^{d+1}$\color{black}, and
					\item \llabel{it:ADAM}\label{it:ADAM:perm}
					it holds for every $\alpha$-$\beta$-$\eps$-\Adam\ optimizer  $\Phi_n\colon(\R^d)^{n}\to\R^d$, $n \in \N$,\cfadd{def:ADAM} 
					that the stability region\cfadd{def:stab:region} of $(\Phi_n)_{n\in \N}$ is
					\color{magenta}$[0,\infty)^{d+1}$\color{black}
				\end{enumerate}
				\cfload.
			\end{boxedthm}
		\end{athm}
		
		\Cref{intro:thm:1.1} is a direct consequence of \cref{final:thm:1.1} in \cref{section:5} below.
		The natural number $d\in\N$ in \cref{intro:thm:1.1} specifies again the dimension of the \OP\ under consideration (the number of parameters/degrees of freedom that need to be optimized).
		The parameter $\alpha\in[0,1)$ in \cref{intro:thm:1.1} describes the momentum decay parameter in 
		the Nesterov optimizer, 
		the momentum optimizer, and the \Adam\ optimizer,
		the parameter $\beta\in(\alpha^2,1)$ in \cref{intro:thm:1.1} specifies the second moment decay parameter in the adaptive optimization methods \RMSprop\ and \Adam,
		and the real number $\eps\in(0,\infty)$ in \cref{intro:thm:1.1} specifies the regularizing parameter in the adaptive \GD\ methods \RMSprop\ and \Adam\ that avoids dividing by zero.
		\Cref{intro:thm:1.1} explicitly specifies the stability region of different optimizers.
		In particular, we note that 
		\begin{itemize}
			\item the stability region of the Nesterov optimizer is a proper subset of the stability region of the standard \GD\ optimizer, 
			\item the stability region of the standard \GD\ optimizer is a proper subset of the stability region of the momentum optimizer, 
		and 
		\item the stability region of the momentum optimizer is a proper subset of the stability region of the \RMSprop\ and the \Adam\ optimizers, which both have maximal stability region.
				\end{itemize}
		We also note that, without the employment of the notion of the stability region, parts of the conclusion of \cref{it:Nest:perm,it:SGD:perm,it:MOM:perm} are already well-known in the literature.
		In particular, without employing the notion of the stability region, the elementary conclusion of \cref{it:SGD:perm} can be found, \eg, in \cite[Theorem~6.1.12]{jentzen2023mathematical}.
		Furthermore, without employing the concept of the stability region, lower bounds for the stability regions in \cref{it:Nest:perm,it:MOM:perm}, which are slightly smaller than \eqref{eq:stab:Nest:intro} and \eqref{eq:stab:MOM:intro}, respectively, have been established, \eg, in \cite{TrungMomLR}.
				
		\subsection{A priori bounds for the Adam optimizer}
		The concept of the stability region and \cref{intro:thm:1.1}, respectively, only offers a conclusion for adaptive and/or accelerated gradient based optimization methods applied to deterministic \OPs.
		Many of the findings in this work are, however, also applicable to the gradient based optimization methods with possibly non-constant learning rates applied to simple stochastic \OPs.
		This is precisely the subject of the next result, \cref{thm:stoch:ADAM:bounded}, in which we establish a priori bounds for the \Adam\ and the \RMSprop\ optimizers applied to a simple class of quadratic stochastic \OPs.

			\newcommand{\Grad}{\mathcal{G}}
		\begin{athm}{theorem}{thm:stoch:ADAM:bounded}[A priori bounds for \Adam]
			\begin{boxedthm}
				Let 
				$d\in\N$, $\alpha\in[0,1)$, $\beta\in(\alpha^2,1)$, $\eps\in(0,\infty)$, 
				let $\Phi_n\colon(\R^d)^{n}\to\R^d$, $n \in \N$,\cfadd{def:ADAM} be the $\alpha$-$\beta$-$\eps$-\Adam\ optimizer,
				let $\gamma\colon\N\to[0,\infty)$ be bounded,
				let $\batch\colon\N\to\N$ be a function,
				let $\lambda\in\R^d$, $\fc\in[0,\infty)$,
				let
				$\loss\colon\R^d\times\R^d\to\R$ satisfy for all $x,\theta\in\R^d$ that
				$
				\loss(\theta,x)
				=\norm{\diag(\lambda)(\theta-x	)}^2$,
				let $(\Omega,\cF,\P)$ be a probability space,
				for every $n,j\in\N$
				let $X_{n,j}\colon\Omega\to[-\fc,\fc]^d$ be a random variable, and
				let
				$ \Grad\colon\N\times\Omega \to \R^d$ 
				and
				$ \Theta \colon\N_0\times\Omega \to \R^d$ 
				satisfy for all
				$n\in\N$ that
				\begin{equation}
					\llabel{eq:ADAM:3}
					\begin{split}
						\textstyle
						\Grad_n=
						\frac{1}{\batch_n}\sum_{j=1}^{\batch_n}\pr{\nabla_{\theta}\loss}\pr{\Theta_{n-1},X_{n,j}}
						\qqandqq
						\Theta_n
						= 
						\Theta_{ n - 1 }
						-
						\gamma_n\Phi_n\prb{\Grad_1,\Grad_2,\dots,\Grad_n}
					\end{split}
				\end{equation}
				\cfload.
				Then there exists $c\in\R$ such that $\sup_{n\in\N_0}\norm{\Theta_n}\leq c\norm{\Theta_0}+ c$.
			\end{boxedthm}
		\end{athm}

		\Cref{thm:stoch:ADAM:bounded} follows immediately from
		\cref{cor:stoch:ADAM:bounded} in \cref{section:5} below. \Cref{cor:stoch:ADAM:bounded}, in turn,
		is a direct consequence of \cref{lem:ADAM:bounded} in \cref{section:Aprioribounds} below.
		
		As before, 
		the natural number $d\in\N$ in \cref{thm:stoch:ADAM:bounded} specifies the dimension of the \OP\ under consideration,
		the parameter $\alpha\in[0,1)$ in \cref{thm:stoch:ADAM:bounded} describes the momentum decay parameter of \Adam,
		the parameter $\beta\in(\alpha^2,1)$ in \cref{thm:stoch:ADAM:bounded} specifies the second moment decay parameter of \Adam,
		the real number $\eps\in(0,\infty)$ in \cref{thm:stoch:ADAM:bounded} specifies the regularizing parameter of \Adam\ that avoids dividing by zero,
		and the functions $\Phi_n\colon(\R^d)^{n}\to\R^d$, $n \in \N$, in \cref{thm:stoch:ADAM:bounded} specify the full history recursion dynamics of \Adam\ (cf.\ \lref{eq:ADAM:3} and \cref{def:ADAM}).
		Furthermore, we note that in \cref{thm:stoch:ADAM:bounded} for every $n\in\N$ we have that
		$\gamma_n\in[0,\infty)$ specifies the learning rate of \Adam\ and
		$\batch_n\in\N$ specifies the size of the mini-batches of \Adam.
		Moreover, the function $\loss\colon\R^d\times\R^d\to\R$ in \cref{thm:stoch:ADAM:bounded} specifies the loss function of the \OP\ under consideration (cf.\ \eqref{eq:min:prob}),
		the random variables $X_{n,j}\colon\Omega\to[-\fc,\fc]^d$, $(n,j)\in\N^2$, represent the data of the considered stochastic \OP,
		and the process $\Theta=(\Theta_n)_{n\in\N_0}\colon\N_0\to\R^d$ represents exactly the \Adam\ optimization process.
		We note that \cref{thm:stoch:ADAM:bounded} establishes boundedness of the \Adam\ optimization process and, in particular, we observe that \cref{thm:stoch:ADAM:bounded} ensures that if the \Adam\ optimization process at initial time is an integrable random variable also the whole process is integrable in the sense that $\E\PRb{\sup_{n\in\N_0}\norm{\Theta_n}}<\infty$.

		\subsection{Literature review}
		In this subsection we provide a concise overview of selected works in the literature that address the convergence and/or boundedness properties of the gradient based optimization methods analyzed in \cref{intro:thm:1.1} and \cref{thm:stoch:ADAM:bounded} above.
		
		Sufficient conditions for the convergence, which can be translated into the description of a subset of the stability region, of the heavy-ball method when applied to a simple class of quadratic optimization problems are established, \eg, in \cite[Section 2]{TrungMomLR} and \cite[Theorem 1]{GhadimiFeyzJohansson2015}. 
		The heavy-ball method (cf., \eg, \cite{TrungMomLR} and \cite[Theorem 1]{GhadimiFeyzJohansson2015}), in turn, can in a straightforward way be reformulated as the classical momentum method (cf., \eg, \cite[Lemma 6.3.12]{jentzen2023mathematical}).
		Error estimates for the momentum method applied to certain stochastic \OPs\ can, \eg, be found in \cite[Theorem 1]{YangLinLi2016} and
		further convergence analyses for the momentum method applied to certain deterministic \OPs\ can be found, \eg, in \cite[Theorem 1]{Poljak_momentum_SGD} and \cite[Theorem 1]{ZavrievKostyuk1993}.

		Sufficient conditions for convergence of the Nesterov optimizer when applied to a simple class of quadratic \OPs\ are presented in \cite[Section 3]{TrungMomLR}.
		In \cite[Theorem 3]{GhadimiFeyzJohansson2015} a priori bounds are derived for the Nesterov method when applied to a class of convex continuously differentiable objective functions in the case of a fixed learning rate determined by the underlying objective function.

		%
		Error estimates for \Adam\ when applied to a class of \OPs\ with strongly convex objective functions with uniformly bounded second order moments of the gradients can, \eg, be found in \cite[Theorem 1 and Theorem 2]{Mazumder2023AdamExactStepSize}.
		For a certain class of learning rates \cite[Theorem 4]{Defossez2022} proves that the second moments of the gradients can be found to be arbitrarily small when \Adam\ is applied to stochastic \OPs\ where the gradient of the objective function is globally bounded.
		An upper bound for the expected norm of a randomly chosen gradient during the application of \Adam\ in a non-convex setting is provided in \cite[Theorem 4]{Zou2019SufficientConditionAdamRMSProp} under the assumption of the boundedness of the second moments of the gradients.
		Moreover, \cite[Theorem 1 and Theorem 2]{Zhang2024RMSPropAdamGeneralizedSmooth} establishes an upper bound for the mean of the first-order moments of the gradients along the sequence of iterates generated by \Adam\ and \RMSprop, respectively, in a non-convex stochastic setting.
		In \cite[Theorem 4.3]{barakat2021convergence} it is shown that \Adam\ can approximate the solution of an ordinary differential equation in the sense that the probability to exceed a certain error over a compact set is arbitrarily low if the learning rate is sufficiently small.
		Furthermore,
		 \cite[Theorem 1.2]{dereich2025sharphigherorderconvergence} proves
		that the order of convergence of 
		the \Adam\ optimizer and
		the momentum optimizer
		exceeds the order of convergence of
		the \RMSprop\ optimizer and
		the \GD\ optimizer, respectively, 
		when applied to a certain class of deterministic \OPs.

		More general theoretical frameworks for the analysis of optimization algorithms are proposed, \eg, in \cite{LelucPortier2020,GodichonBaggioniTarrago2025}. These works introduce a unified representation for a broad class of optimizers via so-called conditioning matrices and, thereby, provide both upper bounds and convergence guarantees for the considered optimizers. While both \Adam\ and \RMSprop\ can, in principle, be represented within this framework, the specific assumptions required for the derived results are in general not satisfied in the case of simple quadratic stochastic \OPs.
		
		For further reviews on \Adam\ and other \GD\ optimization methods we refer, \eg, to the monograph \cite{jentzen2023mathematical} and the survey article \cite{ruder2017overviewgradientdescentoptimization}.

		\subsection{Stucture of this article}
		\label{subsection:1.5}
		The remainder of this article is structured as follows.
		In \cref{section:Aprioribounds} we establish a priori bounds for the \Adam\ and 
		the \RMSprop\ optimizers
		when applied to a certain class of simple quadratic stochastic \OPs.
		In \cref{section:apriori:momentum} we explicitly calculate the set of tupels of learning rates and eigenvalues of the Hessian for which the momentum method does not explode but stays bounded when applied to the class of simple quadratic \OPs\ in \eqref{eq:min:prob} in \cref{subsection:1:1} above.
		In \cref{section:apriori:Nest} we explicitly calculate the set of tupels of learning rates and eigenvalues of the Hessian for which the Nesterov method does not explode but stays bounded when applied to the class of simple quadratic \OPs\ in \eqref{eq:min:prob}.
		In \cref{section:5} we combine the findings from 
		\cref{section:Aprioribounds,section:apriori:momentum,section:apriori:Nest}  to explicitly specify the stability region (cf.\ \cref{subsection:1:1}) for 
		the Nesterov optimizer,
		the \GD\ optimizer,
		the momentum optimizer,
		the \Adam\ optimizer, and the
		\RMSprop\ optimizer.\\

		\newpage


		\definecolor{brightpink}{rgb}{1, 0.4, 0.09}

		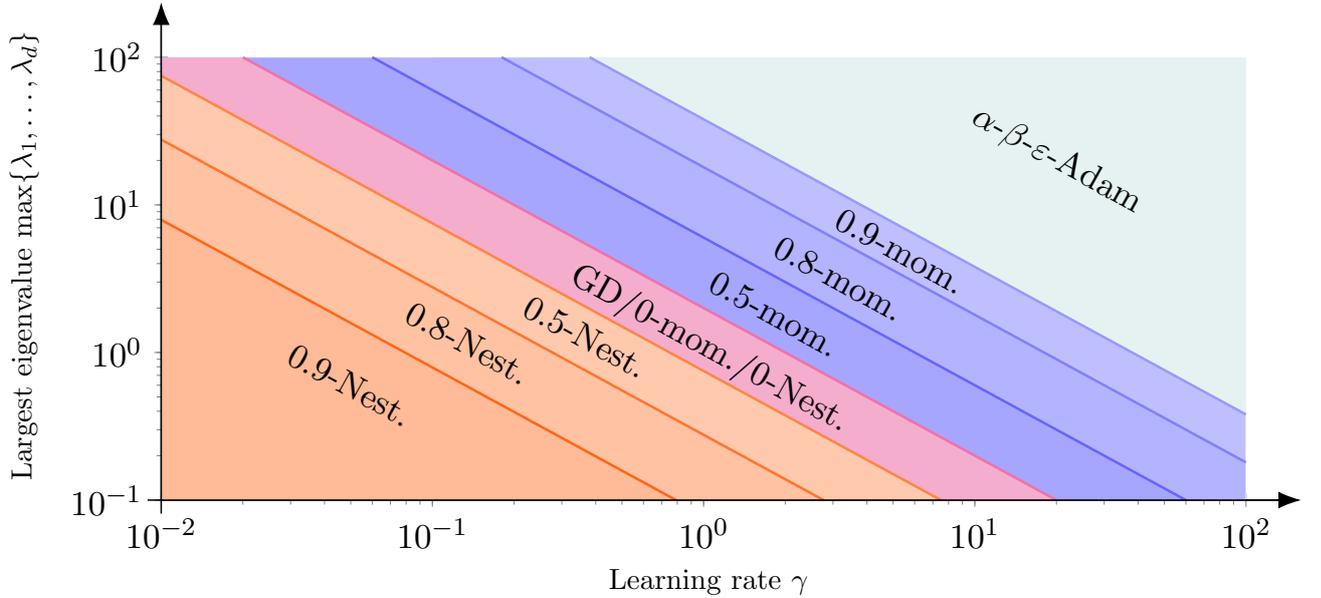
\begin{figure}[H]
			\centering
			\begin{tikzpicture}[scale=1.2]
				\begin{axis}[
					axis lines = middle,
					ymode=log,
					xmode=log,
					ymax=210,xmax=150,
					ymin=0.1,xmin=0.01,
					width=14cm,
					height=7cm,
					axis x line* = bottom,
					axis y line* = left
					]
					\fill [teal!10, domain=0.1:100, variable=\x]
					(-4.6,4.6)
					-- (4.6,4.6)
					-- (4.6,-2.3)
					-- (-4.6,-2.3)
					-- cycle;
					\fill [blue!25, domain=0.1:100, variable=\x]
					(-2.3,4.6)
					-- (-0.95,4.6)
					-- (4.6,-0.95)
					-- (4.6,-2.3)
					-- (-2.3,-2.3)
					-- cycle;
					\fill [blue!30, domain=0.1:100, variable=\x]
					(-2.8,4.6)
					-- (-1.7,4.6)
					-- (4.6,-1.7)
					-- (4.6,-2.3)
					-- (-2.3,-2.3)
					-- cycle;
					\fill [blue!35, domain=0.1:100, variable=\x]
					(-4.6,4.6)
					-- (-2.8,4.6)
					-- (4.1,-2.3)
					-- (-4.6,-2.3)
					-- cycle;
					\fill [magenta!40, domain=0.1:100, variable=\x]
					(-4.6,4.6)
					-- (-3.9,4.6)
					-- (3,-2.3)
					-- (-4.6,-2.3)
					-- cycle;
					\fill [brightpink!35, domain=0.1:100, variable=\x]
					(-4.6,4.3)
					-- (2,-2.3)
					-- (-4.6,-2.3)
					-- cycle;
					\fill [brightpink!40, domain=0.1:100, variable=\x]
					(-4.6,3.34)
					-- (1,-2.3)
					-- (-4.6,-2.3)
					-- cycle;
					\fill [brightpink!45, domain=0.1:100, variable=\x]
					(-4.6,2.08)
					-- (-0.25,-2.3)
					-- (-4.6,-2.3)
					-- cycle;
					\addplot+[domain=0.38:100,mark=none,smooth,samples=200,blue!40,thick] (\x, {38/\x}); 
					\addplot+[domain=0.18:100,mark=none,smooth,samples=200,blue!50,thick] (\x, {18/\x}); 
					\addplot+[domain=0.06:100,mark=none,smooth,samples=200,blue!60,thick] (\x, {6/\x}); 
					\addplot+[domain=0.02:100,mark=none,smooth,samples=200,magenta!70,thick] (\x, {2/\x}); 
					\addplot+[domain=0.01:100,mark=none,smooth,samples=200,brightpink!80,thick] (\x, {0.75/\x}); 
					\addplot[domain=0.01:100,mark=none,smooth,samples=200,brightpink!90,thick] (\x, {0.2769231/\x}); 
					\addplot[domain=0.00791666:100,mark=none,smooth,samples=200,brightpink,thick] (\x, {0.0791666/\x}); 
					\tznode(3,3){$\alpha$-$\beta$-$\eps$-Adam}[rotate=331]
					\tznode(1.67,1.58){$0.9$-mom.}[rotate=331]
					\tznode(1.15,1.15){$0.8$-mom.}[rotate=331]
					\tznode(0.6,0.6){$0.5$-mom.}[rotate=331]
					\tznode(0.05,0.05){GD/$0$-mom./$0$-Nest.}[rotate=331]
					\tznode(-1,0.25){$0.5$-Nest.}[rotate=331]
					\tznode(-2,0.15){$0.8$-Nest.}[rotate=331]
					\tznode(-3,-0.5){$0.9$-Nest.}[rotate=331]
				\end{axis}
				\tznode(6,-0.9){Learning rate $\gamma$}
				\tznode(-1.5,2.7){Largest eigenvalue $\max\{\lambda_1,\dots,\lambda_d\}$}[rotate=90]
				\draw [line width=0.2mm, -{Latex[length=3mm]}]
				(-0.15,0) -- (12.5,0);
				\draw [line width=0.2mm, -{Latex[length=3mm]}]
				(0,-0.15) -- (0,5.5);
			\end{tikzpicture}
			\caption{
				In this figure we graphically represent for every 
				$\alpha\in[0,1)$, $\beta\in(\alpha^2,1)$, $\eps\in(0,\infty)$ the stability region of
				the $0.9$-Nesterov optimizer, 
				the $0.8$-Nesterov optimizer, 
				the $0.5$-Nesterov optimizer,
				the \GD\ optimizer (the $0$-momentum optimizer or the $0$-Nesterov optimizer), 
				the $0.5$-momentum optimizer, 
				the $0.8$-momentum optimizer, 
				the $0.9$-momentum optimizer, and
				the $\alpha$-$\beta$-$\eps$-\Adam\ optimizer.
			}
			\label{fig:1}
		\end{figure}

\begin{figure}[H]
	\centering
	\begin{tikzpicture}[shorten >=1pt,-latex,draw=black!100, node distance=\layersep,auto]
		\tznode(-2.5,0.7){Convergence speed}[rotate=90]
		\tznode(3,-2){Stability}
		\tzplot[fill,blue!20,smooth cycle,thick,opacity=.5]{1} 
		(4,-0.5)(2.7,0.5)(4,1.5)(5.2,0.5);
		\tzplot[fill,red!20,smooth cycle,thick,opacity=.5]{1} 
		(2.3,0.5)(-1,1.1)(3,1.5)(5,1);
		\tzplot[draw,red,smooth cycle,thick]{1} 
		(2.3,0.5)(-1,1.1)(3,1.5)(5,1);
		\tzplot[draw,blue,smooth cycle,thick]{1} 
		(4,-0.5)(2.7,0.5)(4,1.5)(5.2,0.5);
		
		\node [blue,text width=12em, text centered,font=\small](uniformly stable) at (6,-0.5) {uniformly stable};
		\node [red,text width=12em, text centered,font=\small](optimal convergence rate) at (2,1.7) {optimal convergence rate};
		
		\node (SGD) at (1,0) {GD};
		\node (momentum-SGD) at (1,1) {momentum};
		\node (Adam) at (4,0) {RMSprop};
		\node (RMSprop) at (4,1) {Adam};
		\draw [line width=0.2mm, -{Latex[length=3mm]}]
		(-2.5,-1.5) -- (8,-1.5);
		\draw [line width=0.2mm, -{Latex[length=3mm]}]
		(-2,-2) -- (-2,3);
	\end{tikzpicture}
	\caption{Graphical illustration of stability and convergence speed properties of the momentum optimizer,
	the \Adam\ optimizer,
	the \GD\ optimizer, and
	the \RMSprop\ optimizer
	}
\end{figure}
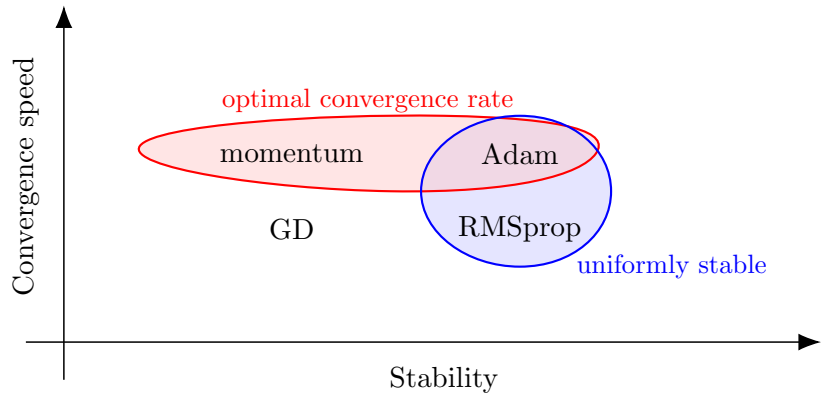

		\section{A priori bounds for Adam and other gradient based methods}
		\label{section:Aprioribounds}

		In this section we establish in \cref{lem:ADAM:bounded} a priori bounds for the \Adam\ optimizer
		(cf.\ \cite{KingmaBa2024_Adam} and, \eg,  \cite[Definitions~6.8.1 and 7.9.1]{jentzen2023mathematical}) and 
		the \RMSprop\ optimizer
		(cf.\ \cite{Hinton24_RMSprop} and, \eg,  \cite[Definitions~6.6.5 and 7.7.3]{jentzen2023mathematical}) 
		when applied to a certain class of simple quadratic stochastic \OPs.
		Our proof of \cref{lem:ADAM:bounded} employs the more general a priori estimates in \cref{lem:ADAM:form} in which we establish general a priori bounds
		that also apply to 
		the AMSGrad optimizer (cf.\ \cite{AdamandBeyond2019} and, \eg,  \cite[Definitions~6.13.1 and 7.14.1]{jentzen2023mathematical}) and 
		the \Adagrad\ optimizer (cf.\ \cite{duchi2011adaptive} and, \eg,  \cite[Definitions~6.5.1 and 7.6.1]{jentzen2023mathematical}).
		In our proof of \cref{lem:ADAM:form} we employ, among other things, the elementary and well-known bounds for the increments of \Adam\ in \cref{lem:hoelder:grad:prep} below (cf., \eg, \cite[Lemma~3.1]{dereich2024convergenceratesadamoptimizer}) as well as the elementary and well-known explicit representation for affine one-step recursions in \cref{momentum:representation} below.

		\subsection{A priori bounds for standard gradient descent (GD) optimization}
		\begin{athm}{prop}{prop:a_priori_bound_one_dimensional}
			Let 
			$ \gamma \colon \N \to \R $, 
			$
			X \colon \N \to \R
			$, 
			and 
				$ \Theta \colon \N_0 \to \R $ satisfy 
			for all $ n \in \N $ that
			\begin{equation}
				\label{eq:recursion}
				\Theta_n
				= 
				\Theta_{ n - 1 }
				-
				\gamma_n 
				( \Theta_{ n - 1 } - X_n )
			\end{equation}
			and let $ \delta \in \N $, $ \fc \in (0,\infty) $
			satisfy 
			for all $ n \in \N \cap [ \delta, \infty) $ with 
			$ \min_{ m \in \N \cap [1,\delta] } | \Theta_{ n - m } | \geq \fc $ 
			that 
			\begin{equation}
				\llabel{eq:learning_rate_smallness}
				\textstyle 
				0 \leq \gamma_n \leq 1
				\qquad 
				\text{and}
				\qquad 
				| X_n |
				\leq 
				\fc 
				.
			\end{equation}
			Then 
			\begin{equation}
				\label{eq:a_priori_to_prove}
				\textstyle 
				\sup_{ n \in \N_0 }
				| \Theta_n |
				\leq
				\bigl[ 
				1 
				+ 
				\sup_{ n \in \N } | \gamma_n | 
				\bigr]^{ \delta }
				\bigl( 
				\max\{ \fc, | \Theta_0 | \}
				+
				\sup_{ n \in \N }
				| X_n |
				\bigr)
				.
			\end{equation}
		\end{athm}

		\begin{aproof}
			Throughout this proof let 
			$ \Gamma, \fC \in [0,\infty] $
			satisfy 
			\begin{equation}
				\label{eq:definition_of_sup_gamma_and_sup_X}
				\textstyle 
				\Gamma = 
				\sup_{ n \in \N } | \gamma_n | 
				\qqandqq
				\fC = \sup_{ n \in \N } | X_n |
			\end{equation}
			and assume without loss of generality that $ \Gamma + \fC < \infty $. 
			\argument{
			\cref{eq:recursion};}
			{that for all $ m \in \N $ 
			it holds that 
			\begin{equation}
				\llabel{eq:1}
				\begin{split}
					| \Theta_m | 
					& \leq 
					| \Theta_{ m - 1 } |
					+
					| \gamma_m |
					| \Theta_{ m - 1 } - X_m |
					\\&\leq 
					| \Theta_{ m - 1 } |
					+
					| \gamma_m |
					\bigl[ 
					| \Theta_{ m - 1 } | + | X_m |
					\bigr]
					\\ &
					\leq 
					| \Theta_{ m - 1 } |
					+
					\Gamma
					\bigl[ 
					| \Theta_{ m - 1 } | + \fC
					\bigr]
					=
					( 1 + \Gamma ) | \Theta_{ m - 1 } | 
					+ \Gamma \fC
					.
				\end{split}
			\end{equation}
		}
		\argument{\lref{eq:1};}{
		for all $ n, m \in \N $ 
			with $ n - m \geq 0 $
			that
			\begin{equation}
				\llabel{eq:2}
				\begin{split}
					| \Theta_n | 
					& \leq 
					( 1 + \Gamma ) | \Theta_{ n - 1 } | 
					+ 
					\Gamma \fC
					\\&\leq 
					( 1 + \Gamma )^2 | \Theta_{ n - 2 } | 
					+
					( 1 + \Gamma ) \Gamma \fC 
					+ \Gamma \fC
					\\ &
					\leq 
					( 1 + \Gamma )^3 | \Theta_{ n - 3 } | 
					+
					( 1 + \Gamma )^2 \Gamma \fC
					+
					( 1 + \Gamma ) \Gamma \fC 
					+ \Gamma \fC
					\\ & \leq 
					\dots 
					\\&\leq \textstyle
					( 1 + \Gamma )^m | \Theta_{ n - m } | 
					+
					\left[ 
					\sum_{ k = 0 }^{ m - 1 }
					( 1 + \Gamma )^k \Gamma \fC 
					\right]
					\\&=\textstyle
					( 1 + \Gamma )^m | \Theta_{ n - m } | 
					+
					\left[ 
					\sum_{ k = 0 }^{ m - 1 }
					( 1 + \Gamma )^k 
					\right]
					\Gamma \fC 
					\\ & =
					( 1 + \Gamma )^m | \Theta_{ n - m } | 
					+
					\bigl( 
					( 1 + \Gamma )^m - 1 
					\bigr)
					\fC 
					\\&\leq
					( 1 + \Gamma )^m 
					\bigl( 
					| \Theta_{ n - m } | 
					+ \fC
					\bigr) 
					.
				\end{split}
			\end{equation}
		}
			\argument{\lref{eq:2};}{for all $ n, m \in \N_0 $ with $ n - m \geq 0 $ that
			\begin{equation}
				\llabel{eq:recursive_growth_bound_with_large_learning_rates}
				| \Theta_n | 
				\leq 
				( 1 + \Gamma )^m 
				\bigl( 
				| \Theta_{ n - m } | 
				+ \fC
				\bigr)
				.
			\end{equation}
		}
			\argument{\lref{eq:recursive_growth_bound_with_large_learning_rates};}{
			for all $ n \in \N_0 $ that 
			\begin{equation}
				\llabel{eq:3}
				| \Theta_n | 
				\leq 
				( 1 + \Gamma )^n 
				\bigl( 
				| \Theta_0 | 
				+ \fC
				\bigr)
				.
			\end{equation}
		}
		\argument{\lref{eq:3};}{
			for all $ n \in \N_0 \cap [0,\delta] $ that 
			\begin{equation}
				\label{eq:estimate_for_small_values_of_n}
				| \Theta_n | 
				\leq 
				( 1 + \Gamma )^n 
				\bigl( 
				| \Theta_0 | 
				+ \fC
				\bigr)
				\leq 
				( 1 + \Gamma )^{ \delta }
				\bigl( 
				| \Theta_0 | 
				+ \fC
				\bigr)
				\leq 
				( 1 + \Gamma )^{ \delta } 
				\bigl( 
				\max\{ \fc, | \Theta_0 | \}
				+ \fC
				\bigr)
				.
			\end{equation}
		}
			\argument{\lref{eq:recursive_growth_bound_with_large_learning_rates}}{
			that for all $ n \in \N_0 \cap [ \delta, \infty ) $, 
			$ m \in \N_0 \cap [0,\delta] $
			it holds that 
			\begin{equation}
				\llabel{eq:4}
				| \Theta_n | 
				\leq 
				( 1 + \Gamma )^m 
				\bigl( 
				| \Theta_{ n - m } | 
				+ \fC
				\bigr)
				\leq 
				( 1 + \Gamma )^{ \delta }
				\bigl( 
				| \Theta_{ n - m } | 
				+ \fC
				\bigr)
				.
			\end{equation}
		}
			\argument{\lref{eq:4}}{
			for all $ n \in \N_0 \cap [ \delta, \infty ) $, 
			$ m \in \N_0 \cap [0, \delta] $
			with 
			$
			| \Theta_{ n - m } | \leq \fc 
			$
			that
			\begin{equation}
				\llabel{eq:5}
				| \Theta_n | 
				\leq 
				( 1 + \Gamma )^{ \delta }
				\bigl( 
				| \Theta_{ n - m } | 
				+ \fC
				\bigr)
				\leq 
				( 1 + \Gamma )^{ \delta } 
				\bigl( 
				\fc 
				+ \fC
				\bigr)
				.
			\end{equation}
		}
			\argument{\lref{eq:5}}{ 
			for all $ n \in \N_0 \cap [ \delta, \infty ) $
			with 
			$
			\min_{ m \in \N_0 \cap [0, \delta] }
			| \Theta_{ n - m } | \leq \fc 
			$
			that
			\begin{equation}
				\llabel{eq:bound_if_min_is_small}
				| \Theta_n | 
				\leq 
				( 1 + \Gamma )^{ \delta } 
				\bigl( 
				\fc 
				+ \fC
				\bigr)
				\leq 
				( 1 + \Gamma )^{ \delta } 
				\bigl( 
				\max\{ \fc, | \Theta_0 | \}
				+ \fC
				\bigr)
				.
			\end{equation}
		}
			\argument{
			\lref{eq:learning_rate_smallness} 
		}{for all
			$ n \in \N \cap [ \delta, \infty) $ with 
			$ \min_{ m \in \N \cap [ 1, \delta ] } | \Theta_{n-m} | \geq \fc $ 
			it holds that 
			\begin{equation}
				\llabel{eq:6}
				\begin{split}
					| \Theta_n |
					=
					\left|
					\Theta_{ n - 1 }
					- 
					\gamma_n
					(
					\Theta_{ n - 1 } - X_n
					)
					\right|
					&=
					\left|
					( 1 - \gamma_n )
					\Theta_{ n - 1 }
					+
					\gamma_n
					X_n
					\right|
					\\&\leq 
					\left|
					1 - \gamma_n 
					\right|
					\left| 
					\Theta_{ n - 1 }
					\right|
					+
					\left|
					\gamma_n
					\right|
					\left|
					X_n
					\right|
					\\ &
					=
					\left(
					1 - \gamma_n 
					\right)
					\left| 
					\Theta_{ n - 1 }
					\right|
					+
					\gamma_n
					\left|
					X_n
					\right|
					\\&\leq 
					\left(
					1 - \gamma_n 
					\right)
					\left| 
					\Theta_{ n - 1 }
					\right|
					+
					\gamma_n
					\fc 
					\\ & 
					\textstyle 
					\leq 
					\left(
					1 - \gamma_n 
					\right)
					\left| 
					\Theta_{ n - 1 }
					\right|
					+
					\gamma_n 
					\left[ 
					\min_{ m \in \N \cap [1,\delta] } 
					\left| 
					\Theta_{ n - m }
					\right|
					\right]
					\\&\leq 
					\left(
					1 - \gamma_n 
					\right)
					\left| 
					\Theta_{ n - 1 }
					\right|
					+
					\gamma_n 
					\left| 
					\Theta_{ n - 1 }
					\right|
					.
				\end{split}
			\end{equation}
		}
			\argument{\lref{eq:6}}{
			for all 
			$ n \in \N \cap [ \delta, \infty) $ with 
			$ \min_{ m \in \N_0 \cap [ 0, \delta ] } | \Theta_{n-m} | \geq \fc $ 
			that 
			\begin{equation}
				\llabel{eq:7}
				\begin{split}
					| \Theta_n |
					\leq 
					\left(
					1 - \gamma_n 
					\right)
					\left| 
					\Theta_{ n - 1 }
					\right|
					+
					\gamma_n 
					\left| 
					\Theta_{ n - 1 }
					\right|
					=
					\left| 
					\Theta_{ n - 1 }
					\right|
					.
				\end{split}
			\end{equation}
		}
			\argument{\lref{eq:7};\lref{eq:bound_if_min_is_small}}{
			for all 
			$ n \in \N \cap [ \delta, \infty) $ 
			it holds that 
			\begin{equation}
				\llabel{eq:8}
				| \Theta_n |
				\leq 
				\max\bigl\{ 
				\left| 
				\Theta_{ n - 1 }
				\right|
				,
				( 1 + \Gamma )^{ \delta } 
				\bigl( 
				\max\{ \fc, | \Theta_0 | \}
				+ \fC
				\bigr)
				\bigr\} 
				.
			\end{equation}
		}
			\argument{\lref{eq:8};\cref{eq:estimate_for_small_values_of_n};induction}{that for all $ n \in \N_0 $ it holds that
			\begin{equation}
				\llabel{eq:9}
				| \Theta_n |
				\leq 
				( 1 + \Gamma )^{ \delta } 
				\bigl( 
				\max\{ \fc, | \Theta_0 | \}
				+ \fC
				\bigr)
				.
			\end{equation}
		}
			\argument{\lref{eq:9};\cref{eq:definition_of_sup_gamma_and_sup_X};}{\cref{eq:a_priori_to_prove}}.
		\end{aproof}

		\subsection{A priori bounds for momentum optimization}

		\begin{athm}{prop}{prop:a_priori_bound_one_dim:Adam:1}
			Let $\alpha\in[0,1)$, $\fc \in [0,\infty) $, $\cst\in(0,\infty)$, $ \Cst\in[\cst,\infty)$, $N\in\N$, $M\in\N_0$,
			let
			$ \Theta\colon\N_0\to \R $,
			$ \gamma \colon \N \to [0,\infty) $, and
			$g\colon\N_0\to\R$ satisfy for all $n\in\N\cap[N,N+M]$ that
			\begin{equation}
				\llabel{eq:learning_rate_limited}
					\Theta_n
				= 
				\Theta_{ n - 1 }
				-
				\gamma_n \PRb{
					\textstyle\sum_{k=0}^n(1-\alpha)\alpha^{n-k}g_k} 
					,
				\quad
				\gamma_n\leq \frac{1-\alpha}{(1+2\alpha)\max\{1,\Cst\}}
				,
				\qandq
				\vass{g_0}\leq \Cst(\vass{\Theta_0}+\fc)
				,
			\end{equation}
			and assume for all $n\in\N_0$ that
			\begin{equation}
				\label{eq:setup:gen_grad}
				\pr{\Theta_{n}-\fc}
				\prb{\cst+(\Cst-\cst)\indicator{(-\infty,\fc]}(\Theta_{n})}
				\leq
				g_{n+1}
				\leq
				\pr{\Theta_{n}+\fc}
				\prb{\cst+(\Cst-\cst)\indicator{[-\fc,\infty)}(\Theta_{n})}.
			\end{equation}
			Then
			\begin{equation}
				\label{it:upper:bound:theta:base}
				\begin{split}
					\max_{n\in\N\cap[N,N+M]}\vass{\Theta_{n}}
					&\leq
					\max\pRbbb{4\fc+\frac{3\fc\alpha\Cst}{(1-\alpha)\cst},\fc+3\vass{\Theta_{N-1}},\max_{n\in\N_0\cap[0,N)}\vass{\Theta_{n}}}
					\\&\leq
					4\fc+\frac{3\fc\alpha\Cst}{(1-\alpha)\cst}
					+3\PRbbb{\max_{n\in\N_0\cap[0,N)}\vass{\Theta_{n}}}
					.
				\end{split}
			\end{equation}
		\end{athm}
		
		\begin{aproof}
			Throughout this proof let $\SPR\colon\Z\to \R$ and $\GRAD\colon\Z\to\R$ satisfy for all $n\in\Z$ that
			\begin{equation}
				\begin{split}
					\llabel{eq:recursion:Adam:prep:1}
									\textstyle
					\SPR_{n}=\Theta_{\max\{n,0\}}
					\qqandqq
					\GRAD_{n}=\sum_{k=0}^{n}(1-\alpha)\alpha^{n-k}g_k
				\end{split}
			\end{equation}
			and let $\fC,\LR\in\R$ satisfy
			\begin{equation}
				\llabel{eq:setup:1st:part:Adam:aprior0}
				\fC=\max\pRbbb{\fc+\frac{\fc\alpha\Cst}{(1-\alpha)\cst},\vass{\SPR_{N-1}},\frac{1}{3}\prbb{\max_{n\in\N_0\cap[0,N)}\vass{\SPR_{n}}}-\frac{\fc}{3}}
				\qandq\LR=\frac{1-\alpha}{(1+2\alpha)\max\{1,\Cst\}}
				.
			\end{equation}
			\argument{
				\lref{eq:recursion:Adam:prep:1};}
			{that for all $n\in\N$ it holds that
				\begin{equation}
					\begin{split}
						\llabel{eq:lambda:recursion}
						\GRAD_n
						&=\textstyle(1-\alpha)g_n
						+
						\alpha\PRb{\sum_{k=0}^{n-1}(1-\alpha)\alpha^{n-1-k}g_k}
						=(1-\alpha)g_n+\alpha\GRAD_{n-1}
						.
					\end{split}
				\end{equation}
			}
						\argument{
				\lref{eq:learning_rate_limited};
				\lref{eq:recursion:Adam:prep:1};
				\lref{eq:lambda:recursion};}{
				that for all $n\in\N\cap[N,N+M]$ it holds that
				\begin{equation}
					\begin{split}
						\llabel{eq:new:exp:recursion}
						\SPR_{n}
						=\Theta_n
						=\Theta_{n-1}-\gamma_n\PRb{
							\textstyle\sum_{k=0}^n(1-\alpha)\alpha^{n-k}g_k} 
						&=\SPR_{n-1}-\gamma_{n}\GRAD_{n}
						\\&=\SPR_{n-1}-\gamma_{n}(1-\alpha)g_{n}-\gamma_{n}\alpha\GRAD_{n-1}
						.
					\end{split}
				\end{equation}
			}
			\argument{
			\lref{eq:setup:1st:part:Adam:aprior0};
			the assumption that $0\leq\alpha< 1$;
		}{\llabel{L:001}that
			\begin{equation}
					\label{eq:global:G:estimate:0}
				\begin{split}
				\max\{ 3 \LR\alpha, 3 \LR \alpha \Cst ( 1 - \alpha )^{ - 1 } \}
				&=
				3 \LR \alpha \max\{ \Cst( 1 - \alpha )^{ - 1 }, 1 \} 
				\leq
				3 \LR \alpha \max\{ \Cst, 1 \} ( 1 - \alpha )^{ - 1 }
				\\&=
				3 \alpha ( 1 + 2 \alpha )^{ - 1 } 
				\leq ( 1 + 2 \alpha ) ( 1 + 2 \alpha )^{ - 1 } = 1.
				\end{split}
			\end{equation}
		}
		\argument{\eqref{eq:setup:gen_grad};
			\lref{eq:recursion:Adam:prep:1};
			the fact that for all $u,v,w\in\R$ with $u\leq v\leq w$ it holds that $\vass{v}\leq\max\{\vass{u},\vass{w}\}$}{that for every $n\in\N$, $z\in\{-1,1\}$ and every $\psi\colon\R\to\R$ with $\forall\,x\in\R,y\in[x,\infty)\colon z\psi(x)\leq z\psi(y)$ it holds that
			\begin{equation}
				\begin{split}
					\llabel{eq:pre:basis:first:incr}
					\vass{\psi(g_n)}
					=
					\vass{z\psi(g_n)}
					&\leq
					\textstyle\max_{s\in\{-1,1\}}\max_{t\in\{\cst,\Cst\}}\vass{z\psi(t(\Theta_{n-1}+s\fc))}
					\\&=
					\textstyle\max_{s\in\{-1,1\}}\max_{t\in\{\cst,\Cst\}}\vass{\psi(t(\Theta_{n-1}+s\fc))}
					\\&=
					\textstyle\max_{s\in\{-1,1\}}\max_{t\in\{\cst,\Cst\}}\vass{\psi(t(\SPR_{n-1}+s\fc))}
					.
				\end{split}
		\end{equation}}
					\argument{
			\lref{eq:pre:basis:first:incr};
		}{\llabel{L:002}that for all $n\in\N$ it holds that
			\begin{equation}
				\label{eq:global:G:estimate:1}		
				\vass{g_n}
				\leq \textstyle\max_{s\in\{-1,1\}}\max_{t\in\{\cst,\Cst\}}\vass{t(\SPR_{n-1}+s\fc)}
				\leq\Cst(\vass{\SPR_{n-1}}+\fc).
			\end{equation}
		}
			\argument{\lref{eq:learning_rate_limited};
				\lref{eq:recursion:Adam:prep:1};}{that
				\begin{equation}		
										\llabel{eq:grad_0}
				\vass{g_0}
				\leq \Cst(\vass{\Theta_{\max\{-1,0\}}}+\fc)
				= \Cst(\vass{\SPR_{-1}}+\fc)
				.
			\end{equation}
		}
		\argument{\lref{eq:grad_0};
			\eqref{eq:global:G:estimate:0};
			\eqref{eq:global:G:estimate:1};}{that for all $n\in\N_0$ it holds that
			\begin{equation}
				\llabel{eq:global:G:estimate}
				3\LR\alpha\leq 1
				,
				\quad
				3\LR\alpha\Cst(1-\alpha)^{-1}\leq 1,
				\quad
				\LR\alpha\Cst(3\fC+\fc)\leq \fC+\fc,
				\qandq
				\vass{g_n}
				\leq \Cst(\vass{\SPR_{n-1}}+\fc).
			\end{equation}
		}
			\argument{
			\lref{eq:recursion:Adam:prep:1};
			\lref{eq:global:G:estimate};}
			{that for all $n\in\N_0$ it holds that
			\begin{equation}
				\llabel{eq:apriori:bound:increments}
				\begin{split}
					\vass{\GRAD_n}
					\leq\textstyle\sum_{k=0}^n(1-\alpha)\alpha^{n-k}\vass{g_k}
					&\leq\textstyle\sum_{k=0}^n(1-\alpha)\alpha^{n-k}\Cst(\vass{\SPR_{k-1}}+\fc)
					\\&\textstyle\leq 
					\PRb{\sum_{k=0}^n(1-\alpha)\alpha^{n-k}}
					\PRb{\Cst\fc+ \Cst\max_{k\in\{0,1,\dots,n\}}\vass{\SPR_{k-1}}}
					\\&\textstyle\leq 
					\PRb{(1-\alpha)\sum_{k=0}^\infty\alpha^{k}}
					\PRb{\Cst\fc+ \Cst\max_{k\in\{0,1,\dots,n\}}\vass{\SPR_{k-1}}}
					\\&\textstyle= \Cst\prb{\fc+ \max_{k\in\{0,1,\dots,n\}}\vass{\SPR_{k-1}}}
					\\&\textstyle= \Cst\prb{\fc+ \max_{k\in\{1,2,\dots,n+1\}}\vass{\SPR_{k-2}}}
					.
				\end{split}
			\end{equation}
		}
	\argument{
		\lref{eq:new:exp:recursion};
		\lref{eq:apriori:bound:increments};
		}{ 
		that for all $n\in\N\cap[N,N+M]$ it holds that
		\begin{equation}
			\llabel{eq:one:step:for:small:theta:pre}
			\begin{split}
				\vass{\SPR_{n}}
				&=\vass{\SPR_{n-1}-\gamma_{n}(1-\alpha)g_{n}-\gamma_{n}\alpha\GRAD_{n-1}}
				\\&\leq\vass{\SPR_{n-1}-\gamma_{n}(1-\alpha)g_{n}}
				+\gamma_{n}\alpha\vass{\GRAD_{n-1}}
				\\&\leq \vass{\SPR_{n-1}-\gamma_{n}(1-\alpha)g_{n}}
				+\gamma_{n}\alpha\Cst\prb{\fc+\textstyle\max_{k\in\{1,2,\dots,n\}}\vass{\SPR_{k-2}}}
				\\&= \vass{\SPR_{n-1}-\gamma_{n}(1-\alpha)g_{n}}
				+\gamma_{n}\alpha\Cst\fc
				+\textstyle\gamma_{n}\alpha\Cst\max_{k\in\{1,2,\dots,n\}}\vass{\SPR_{k-2}}
				.
			\end{split}
		\end{equation}
	}
	\argument{\lref{eq:learning_rate_limited};the fact that $\cst\leq\Cst$;}{that for all $n\in\N\cap[N,N+M]$, $t\in\{\cst,\Cst\}$ it holds that
			\begin{equation}
				\llabel{eq:mod:LR:bound}
				\begin{split}
				\gamma_{n}(1-\alpha)t
				\leq
				\gamma_{n}(1-\alpha)\Cst
				\leq
				\frac{(1-\alpha)^2\Cst}{(1+2\alpha)\max\{1,\Cst\}}
					\leq 1
					.
				\end{split}
			\end{equation}
		}
		\argument{\lref{eq:mod:LR:bound};\lref{eq:pre:basis:first:incr};}{that for all $n\in\N\cap[N,N+M]$ it holds that
		\begin{equation}
			\llabel{eq:basis:first:incr}
			\begin{split}
				\vass{\SPR_{n-1}-\gamma_{n}(1-\alpha)g_{n}
				}
				&\leq
				\textstyle\max_{s\in\{-1,1\}}\max_{t\in\{\cst,\Cst\}}\vass{\SPR_{n-1}-\gamma_{n}(1-\alpha)t(\SPR_{n-1}+s\fc)
				}
			\\&=
			\textstyle\max_{s\in\{-1,1\}}\max_{t\in\{\cst,\Cst\}}\vass{(1-\gamma_{n}(1-\alpha)t)\SPR_{n-1}-\gamma_{n}(1-\alpha)ts\fc
			}
			\\&\leq
			\textstyle\max_{s\in\{-1,1\}}\max_{t\in\{\cst,\Cst\}}\vass{(1-\gamma_{n}(1-\alpha)t)\SPR_{n-1}}+\vass{\gamma_{n}(1-\alpha)ts\fc}
			\\&=
		\textstyle\max_{t\in\{\cst,\Cst\}}\PRb{(1-\gamma_{n}(1-\alpha)t)\vass{\SPR_{n-1}}+\gamma_{n}(1-\alpha)t\fc}
		.
			\end{split}
		\end{equation}
	}
	\argument{\lref{eq:setup:1st:part:Adam:aprior0}}{that
	\begin{equation}
		\llabel{eq:pre:setup}
		\fC\geq \fc+\frac{\fc\alpha\Cst}{(1-\alpha)\cst}\geq \fc
		\qqandqq
		(1-\alpha)\cst\fC\geq(1-\alpha)\cst\prbbb{\fc+\frac{\fc\alpha\Cst}{(1-\alpha)\cst}}=(\alpha\Cst+(1-\alpha)\cst) \fc.
		\end{equation}}
			\argument{
			\lref{eq:mod:LR:bound};
			\lref{eq:basis:first:incr};
			\lref{eq:pre:setup};
						}{
			that for all $n\in\N\cap[N,N+M]$ with $\vass{\SPR_{n-1}}\leq\fC$ it holds that
			\begin{equation}
				\llabel{eq:gradient:step2}
				\begin{split}
					&\vass{\SPR_{n-1}-\gamma_{n}(1-\alpha)g_{n}
					}+\gamma_{n}\alpha\Cst\fc		
					\\		&\leq\textstyle\max_{t\in\{\cst,\Cst\}}\PRb{(1-\gamma_{n}(1-\alpha)t)\vass{\SPR_{n-1}}+\gamma_{n}(1-\alpha)t\fc}
					+\gamma_{n}\alpha\Cst\fc	
					\\		&\leq\textstyle\max_{t\in\{\cst,\Cst\}}\PRb{(1-\gamma_{n}(1-\alpha)t)\fC+\gamma_{n}(1-\alpha)t\fc}
					+\gamma_{n}\alpha\Cst\fc	
					\\		&=\textstyle\max_{t\in\{\cst,\Cst\}}\PRb{\fC-\gamma_{n}(1-\alpha)t(\fC-\fc)}
					+\gamma_{n}\alpha\Cst\fc	
					\\&=\fC-\gamma_{n}(1-\alpha)\cst(\fC-\fc)
					+\gamma_{n}\alpha\Cst\fc	
					\\&=\fC-\gamma_{n}\PR{(1-\alpha)\cst\fC-(1-\alpha)\cst\fc
						-\alpha\Cst\fc	}
					\\&\leq\fC-\gamma_{n}\PR{(\alpha\Cst+(1-\alpha)\cst) \fc-(1-\alpha)\cst\fc
							-\alpha\Cst\fc	}
					\leq \fC
					.
				\end{split}
			\end{equation}
		}
			\argument{\lref{eq:gradient:step2};
				\lref{eq:learning_rate_limited};
				\lref{eq:setup:1st:part:Adam:aprior0};
				\lref{eq:one:step:for:small:theta:pre};}{ 
			that for all $n\in\N\cap[N,N+M]$ with $\vass{\SPR_{n-1}}\leq\fC$ it holds that
			\begin{equation}
				\label{eq:one:step:for:small:theta}
				\begin{split}
					\vass{\SPR_{n}}
					&\leq \vass{\SPR_{n-1}-\gamma_{n}(1-\alpha)g_{n}}
					+\gamma_{n}\alpha\Cst\fc
					+\textstyle\gamma_{n}\alpha\Cst\max_{k\in\{1,2,\dots,n\}}\vass{\SPR_{k-2}}
					\\&\leq \fC
					+\gamma_{n}\alpha\Cst\textstyle\max_{k\in\{1,2,\dots,n\}}\vass{\SPR_{k-2}}
					\\&\leq \fC
					+\LR\alpha\Cst\textstyle\max_{k\in\{1,2,\dots,n\}}\vass{\SPR_{k-2}}
					.
				\end{split}
			\end{equation}
		}
			In our proof of \cref{it:upper:bound:theta:base} we distinguish between the case $\alpha=0$ and the case $\alpha>0$. We first \prove[ps] \cref{it:upper:bound:theta:base} in the case $\alpha=0$.
			\startnewargseq
			\argument{
			\lref{eq:recursion:Adam:prep:1};
			\lref{eq:setup:1st:part:Adam:aprior0};			
			\eqref{eq:one:step:for:small:theta};
			induction}{that for all $n\in\N_0\cap[N-1,N+M]$ it holds that
			\begin{equation}
				\llabel{eq:alpha:zero}
				\vass{\Theta_n}=\vass{\SPR_n}\leq\fC\leq3\fC+\fc.
			\end{equation}
		}
		\argument{\lref{eq:alpha:zero};\lref{eq:setup:1st:part:Adam:aprior0};}{\cref{it:upper:bound:theta:base} in the case $\alpha=0$}.
			In the next step we \prove[ps] \cref{it:upper:bound:theta:base} in the case $\alpha>0$.
			\startnewargseq
			\argument{
			\eqref{eq:setup:gen_grad}; 
			\lref{eq:recursion:Adam:prep:1};
			\lref{eq:new:exp:recursion};
			\lref{eq:apriori:bound:increments};
			\lref{eq:mod:LR:bound};
			}
			{that for all $n\in\N\cap[N,N+M]$ with 
			$\SPR_{n-1}\geq\fC$, 			 
			$\max_{k\in\{1,2,\dots,n\}}\vass{\SPR_{k-2}}\leq 3\fC+\fc$ it holds that
			\begin{equation}
				\llabel{eq:sign:change:value:limit_pos:pre}
				\begin{split}
					\SPR_{n}
					&=\SPR_{n-1}-\gamma_{n}(1-\alpha)g_{n}-\gamma_{n}\alpha\GRAD_{n-1}
					\\&\geq
					\SPR_{n-1}-\gamma_{n}(1-\alpha)\Cst(\Theta_{n-1}+\fc)
					-\gamma_{n}\alpha\Cst(\fc+\textstyle\max_{k\in\{1,2,\dots,n\}}\vass{\SPR_{k-2}})
					\\&=
					\SPR_{n-1}-\gamma_{n}(1-\alpha)\Cst(\SPR_{n-1}+\fc)
					-\gamma_{n}\alpha\Cst(\fc+\textstyle\max_{k\in\{1,2,\dots,n\}}\vass{\SPR_{k-2}})
					\\&\geq
					(1-\gamma_{n}(1-\alpha)\Cst)\SPR_{n-1}-\gamma_{n}(1-\alpha)\Cst \fc
					-\gamma_{n}\alpha\Cst(3\fC+2\fc)
					\\&\geq
					(1-\gamma_{n}(1-\alpha)\Cst)\fC
					-\gamma_{n}(1-\alpha)\Cst \fc
					-\gamma_{n}\alpha\Cst(3\fC+2\fc)
					\\&=
					\fC
					-\gamma_n\Cst\pr{(1-\alpha)\fC+(1-\alpha) \fc
						+\alpha(3\fC+2\fc)}
					.
				\end{split}
			\end{equation}
		}
					\argument{
			\lref{eq:sign:change:value:limit_pos:pre};
			\lref{eq:learning_rate_limited};
			\lref{eq:setup:1st:part:Adam:aprior0}; 
			\lref{eq:pre:setup};
			the assumption that $\alpha>0$;}
		{that for all $n\in\N\cap[N,N+M]$ with 
			$\SPR_{n-1}\geq\fC$, 			 
			$\max_{k\in\{1,2,\dots,n\}}\vass{\SPR_{k-2}}\leq 3\fC+\fc$ it holds that
			\begin{equation}
				\llabel{eq:sign:change:value:limit_pos}
				\begin{split}
					\SPR_{n}
					&\geq
					\fC
					-\gamma_n\Cst\pr{(1-\alpha)\fC+(1-\alpha) \fc
						+\alpha(3\fC+2\fc)}
					\\&\geq
					\fC
					-\LR\Cst\pr{(1-\alpha)\fC+(1-\alpha) \fc
						+\alpha(3\fC+2\fc)}
					\\&=
					\fC
					-\LR\Cst\pr{\fC+2\alpha\fC+\fc+\alpha\fc}
					\\&\geq\textstyle
					\fC-\LR\max\{1,\Cst\}(\fC+\fc)(1+2\alpha)
					=
					\fC-(1-\alpha)(\fC+\fc)
					\geq \fC-2(1-\alpha)\fC
					> -\fC
					.
				\end{split}
			\end{equation}
		}
			\argument{
			\eqref{eq:setup:gen_grad};
			\lref{eq:recursion:Adam:prep:1};
			\lref{eq:new:exp:recursion};
			\lref{eq:pre:setup};
			}{that for all $n\in\N\cap[N,N+M]$ with 
			$\SPR_{n-1}\geq\fC$, 
			$\GRAD_{n-1}\geq0$ it holds that
			\begin{equation}
				\label{eq:upper:half:no:increase:case}
				\begin{split}
					\SPR_{n}
					=\SPR_{n-1}-\gamma_{n}(1-\alpha)g_{n}-\gamma_{n}\alpha\GRAD_{n-1}
					&\leq
					\SPR_{n-1}-\gamma_{n}(1-\alpha)g_{n}
					\\&\leq 
					\SPR_{n-1}-\gamma_{n}(1-\alpha)\cst(\Theta_{n-1}-\fc)
					\\&=
					\SPR_{n-1}-\gamma_{n}(1-\alpha)\cst(\SPR_{n-1}-\fc)
					\\&\leq 
					\SPR_{n-1}-\gamma_{n}(1-\alpha)\cst(\fC-\fc)
					\leq \SPR_{n-1}
					.
				\end{split}
			\end{equation}
		}
			\argument{
			\lref{eq:lambda:recursion};induction;}{that for all $n\in\N_0$, $k\in\N$ it holds that
			\begin{equation}
				\llabel{eq:rephrase:gradient}
				\begin{split}
					\GRAD_{n+k}
					&=\alpha\GRAD_{n+k-1}+(1-\alpha)g_{n+k}
					=\textstyle\alpha^k\GRAD_n+\sum_{j=0}^{k-1}\alpha^j(1-\alpha)g_{n+k-j}
					.
				\end{split}
			\end{equation}
		}
		\argument{
			\lref{eq:rephrase:gradient};
			\lref{eq:new:exp:recursion};
			induction}
			{that for all $n\in\N_0\cap[N-1,N+M)$, $m\in\N\cap(0,N+M-n]$ it holds that
			\begin{equation}
				\llabel{eq:large:step:recursion:for:Adam:gen}
				\begin{split}
					\SPR_{n+m}
					=\SPR_{n+m-1}-\gamma_{n+m}\GRAD_{n+m}
					&=\textstyle\SPR_n-\sum_{k=1}^{m}\gamma_{n+k}\GRAD_{n+k}
					\\&=\textstyle
					\SPR_n-
					\sum_{k=1}^{m}\gamma_{n+k}\PRb{\alpha^k\GRAD_n+\sum_{j=0}^{k-1}\alpha^j(1-\alpha)g_{n+k-j}}
					.
				\end{split}
			\end{equation}
		}
			\argument{
			\lref{eq:large:step:recursion:for:Adam:gen};
			\lref{eq:learning_rate_limited};
			\lref{eq:recursion:Adam:prep:1};
			\eqref{eq:setup:gen_grad};
			\lref{eq:setup:1st:part:Adam:aprior0};
			\lref{eq:pre:setup};}{
			that for all $n\in\N_0\cap[N-1,N+M)$, $m\in\N\cap(0,N+M-n]$ with $\min\{\SPR_n,\SPR_{n+1},\dots,\SPR_{n+m}\}\geq\fC$, $\GRAD_n<0$ it holds that
			\begin{equation}
				\llabel{eq:large:step:recursion:for:Adam:upper}
				\begin{split}
					\SPR_{n+m}
					&=\textstyle\SPR_n-
					\sum_{k=1}^{m}\gamma_{n+k}\PRb{\alpha^k\GRAD_n+\sum_{j=0}^{k-1}\alpha^j(1-\alpha)g_{n+k-j}}
					\\&\leq\textstyle\SPR_n-
					\sum_{k=1}^{m}\gamma_{n+k}\PRb{\alpha^k\GRAD_n+\sum_{j=0}^{k-1}\alpha^j(1-\alpha)\cst(\Theta_{n+k-j-1}-\fc)}
					\\&=\textstyle\SPR_n-
					\sum_{k=1}^{m}\gamma_{n+k}\PRb{\alpha^k\GRAD_n+\sum_{j=0}^{k-1}\alpha^j(1-\alpha)\cst(\SPR_{n+k-j-1}-\fc)}
					\\&\leq\textstyle
					\SPR_n-
					\sum_{k=1}^{m}\gamma_{n+k}\alpha^k\GRAD_n
					\\&=\textstyle
					\SPR_n+
					\vass{\GRAD_n}\sum_{k=1}^{m}\gamma_{n+k}\alpha^k
					\\&\leq\textstyle
					\SPR_n+\LR
					\vass{\GRAD_n}\sum_{k=1}^{m}\alpha^k
					\leq\textstyle
					\SPR_n+\LR\alpha
					\vass{\GRAD_n}\sum_{k=0}^{\infty}\alpha^k
					=\textstyle
					\SPR_n+\LR\alpha
					\vass{\GRAD_n}(1-\alpha)^{-1}
					.
				\end{split}
			\end{equation}
		}
			\argument{\lref{eq:large:step:recursion:for:Adam:upper};\lref{eq:global:G:estimate};}{that for all 
			$n\in\N_0\cap[N-1,N+M)$, $m\in\N\cap(0,N+M-n]$ with 
			$\SPR_{n}\leq\fC+\LR\alpha\Cst(3\fC+\fc)$, 
			$\min\{\SPR_n,\SPR_{n+1},\dots,\SPR_{n+m}\}\geq\fC$, 
			and $-\Cst(3\fC+2\fc)\leq\GRAD_n<0$ it holds that
			\begin{equation}
				\label{eq:upper:half:part:1}
				\begin{split}
					\SPR_{n+m}
					\leq \SPR_n+\LR\alpha
					\vass{\GRAD_n}(1-\alpha)^{-1}
					&\leq\fC+\LR\alpha\Cst(3\fC+\fc)+{\LR\alpha\Cst
						(3\fC+2\fc)}\pr{1-\alpha}^{-1}
					\\&\leq \fC +\LR\alpha\Cst(1-\alpha)^{-1}(6\fC+3\fc)
					\leq
					3\fC+\fc
					.
				\end{split}
			\end{equation}
		}
		\argument{\eqref{eq:upper:half:no:increase:case};induction;}{that for all $n\in\N_0\cap[N-1,N+M)$, $m\in\N\cap(0,N+M-n]$ with 
				$\SPR_{n}\leq\fC+\LR\alpha\Cst(3\fC+\fc)$,
				$\min\{\GRAD_n,\GRAD_{n+1},\dots,\GRAD_{n+m-1}\}\geq0$, and
				$\min\{\SPR_n,\SPR_{n+1},\dots,\SPR_{n+m}\}\geq\fC$ it holds that
				\begin{equation}
					\llabel{eq:upper:half:part:2}
					\begin{split}
						\SPR_{n+m}
						&\leq \SPR_n
						\leq \fC+\LR\alpha\Cst(3\fC+\fc)
						.
					\end{split}
				\end{equation}
			}
			\argument{\lref{eq:upper:half:part:2};
			\lref{eq:global:G:estimate};
			\eqref{eq:upper:half:part:1};}{that 
			for all $n\in\N_0\cap[N-1,N+M)$, $m\in\N\cap(0,N+M-n]$ with 
			$\SPR_{n}\leq\fC+\LR\alpha\Cst(3\fC+\fc)$,
			$\vass{\GRAD_n}\leq\Cst(3\fC+2\fc)$, and
			$\min\{\SPR_n,\SPR_{n+1},\dots,\SPR_{n+m}\}\geq\fC$
			it holds that
			\begin{equation}
				\label{eq:upper:half}
				\begin{split}
					\vass{\SPR_{n+m}}
					=
					\SPR_{n+m}
					&\leq 
					\max\{\fC+\LR\alpha\Cst(3\fC+\fc),3\fC+\fc\}
					\leq 
					3\fC+\fc
					.
				\end{split}
			\end{equation}
		}
					\argument{
			\eqref{eq:setup:gen_grad}; 
			\lref{eq:recursion:Adam:prep:1};
			\lref{eq:new:exp:recursion};
			\lref{eq:apriori:bound:increments};
			\lref{eq:mod:LR:bound};
		}
		{that for all $n\in\N\cap[N,N+M]$ with 
			$\SPR_{n-1}\leq-\fC$, 			 
			$\max_{k\in\{1,2,\dots,n\}}\vass{\SPR_{k-2}}\leq 3\fC+\fc$ it holds that
			\begin{equation}
				\llabel{eq:sign:change:value:limit_neg:pre}
				\begin{split}
					\SPR_{n}
					&=\SPR_{n-1}-\gamma_{n}(1-\alpha)g_{n}-\gamma_{n}\alpha\GRAD_{n-1}
					\\&\leq
					\SPR_{n-1}-\gamma_{n}(1-\alpha)\Cst(\Theta_{n-1}-\fc)
					+\gamma_{n}\alpha\Cst(\fc+\textstyle\max_{k\in\{1,2,\dots,n\}}\vass{\SPR_{k-2}})
					\\&=
					\SPR_{n-1}-\gamma_{n}(1-\alpha)\Cst(\SPR_{n-1}-\fc)
					+\gamma_{n}\alpha\Cst(\fc+\textstyle\max_{k\in\{1,2,\dots,n\}}\vass{\SPR_{k-2}})
					\\&\leq
					(1-\gamma_{n}(1-\alpha)\Cst)\SPR_{n-1}
					+\gamma_{n}(1-\alpha)\Cst \fc
					+\gamma_{n}\alpha\Cst(3\fC+2\fc)
					\\&\leq
					-(1-\gamma_{n}(1-\alpha)\Cst)\fC
					+\gamma_{n}(1-\alpha)\Cst \fc
					+\gamma_{n}\alpha\Cst(3\fC+2\fc)
					\\&=
					-\fC
					+\gamma_n\Cst\pr{(1-\alpha)\fC+(1-\alpha) \fc
						+\alpha(3\fC+2\fc)}
					.
				\end{split}
			\end{equation}
		}
			\argument{
			\lref{eq:sign:change:value:limit_neg:pre};
			\lref{eq:learning_rate_limited};
			\lref{eq:setup:1st:part:Adam:aprior0}; 
			\lref{eq:pre:setup};
			the assumption that $\alpha>0$;}
		{that for all $n\in\N\cap[N,N+M]$ with 
			$\SPR_{n-1}\leq-\fC$, 			 
			$\max_{k\in\{1,2,\dots,n\}}\vass{\SPR_{k-2}}\leq 3\fC+\fc$ it holds that
			\begin{equation}
				\llabel{eq:sign:change:value:limit_neg}
				\begin{split}
					\SPR_{n}
					&\leq
					-\fC
					+\gamma_n\Cst\pr{(1-\alpha)\fC+(1-\alpha) \fc
						+\alpha(3\fC+2\fc)}
					\\&\leq
					-\fC
					+\LR\Cst\pr{(1-\alpha)\fC+(1-\alpha) \fc
						+\alpha(3\fC+2\fc)}
					\\&=
					-\fC
					+\LR\Cst\pr{\fC+2\alpha\fC+\fc+\alpha\fc}
					\\&\leq\textstyle
					-\fC+\LR\max\{1,\Cst\}(\fC+\fc)(1+2\alpha)
					=
					-\fC+(1-\alpha)(\fC+\fc)
					\leq -\fC+2(1-\alpha)\fC
					< \fC
					.
				\end{split}
			\end{equation}
		}
			\argument{ 
			\eqref{eq:setup:gen_grad};
			\lref{eq:recursion:Adam:prep:1};
			\lref{eq:new:exp:recursion};
			\lref{eq:pre:setup};}{that for all $n\in\N\cap[N,N+M]$ with $\SPR_{n-1}\leq-\fC$, $\GRAD_{n-1}\leq0$ it holds that
			\begin{equation}
				\llabel{eq:lower:half:no:decrease:case}
				\begin{split}
					\SPR_{n}
					=\SPR_{n-1}-\gamma_{n}(1-\alpha)g_{n}-\gamma_{n}\alpha\GRAD_{n-1}
					&\geq\SPR_{n-1}-\gamma_{n}(1-\alpha)g_{n}
					\\&\geq 
					\SPR_{n-1}-\gamma_{n}(1-\alpha)\cst(\Theta_{n-1}+\fc)
					\\&=
					\SPR_{n-1}-\gamma_{n}(1-\alpha)\cst(\SPR_{n-1}+\fc)
					\\&\geq 
					\SPR_{n-1}+\gamma_{n}(1-\alpha)\cst(\fC-\fc)
					\geq \SPR_{n-1}
					.
				\end{split}
			\end{equation}
		}
			\argument{
			\lref{eq:learning_rate_limited};
			\eqref{eq:setup:gen_grad};
			\lref{eq:recursion:Adam:prep:1};
			\lref{eq:setup:1st:part:Adam:aprior0};
			\lref{eq:large:step:recursion:for:Adam:gen};
			the fact that $-\fC\leq -\fc$;}{that for all $n\in\N_0\cap[N-1,N+M)$, $m\in\N\cap(0,N+M-n]$ with $\max\{\SPR_n,\SPR_{n+1},\dots,\SPR_{n+m}\}\leq-\fC$, $\GRAD_n>0$ it holds that
			\begin{equation}
				\llabel{eq:large:step:recursion:for:Adam:lower}
				\begin{split}
					\SPR_{n+m}
					&=\textstyle\SPR_n-
					\sum_{k=1}^{m}\gamma_{n+k}\PRb{\alpha^k\GRAD_n+\sum_{j=0}^{k-1}\alpha^j(1-\alpha)g_{n+k-j}}
					\\&\geq\textstyle\SPR_n-
					\sum_{k=1}^{m}\gamma_{n+k}\PRb{\alpha^k\GRAD_n+\sum_{j=0}^{k-1}\alpha^j(1-\alpha)\cst(\Theta_{n+k-j-1}+\fc)}
					\\&=\textstyle\SPR_n-
					\sum_{k=1}^{m}\gamma_{n+k}\PRb{\alpha^k\GRAD_n+\sum_{j=0}^{k-1}\alpha^j(1-\alpha)\cst(\SPR_{n+k-j-1}+\fc)}
					\\&\geq\textstyle
					\SPR_n-
					\sum_{k=1}^{m}\gamma_{n+k}\alpha^k\GRAD_n
					\\&=\textstyle
					\SPR_n-
					\vass{\GRAD_n}\sum_{k=1}^{m}\gamma_{n+k}\alpha^k
					\\&\geq\textstyle
					\SPR_n-\LR
					\vass{\GRAD_n}\sum_{k=1}^{m}\alpha^k
					\geq\textstyle
					\SPR_n-\LR\alpha
					\vass{\GRAD_n}\sum_{k=0}^{\infty}\alpha^k
					=\textstyle
					\SPR_n-\LR\alpha
					\vass{\GRAD_n}(1-\alpha)^{-1}
					.
				\end{split}
			\end{equation}
		}
		\argument{\lref{eq:large:step:recursion:for:Adam:lower};
		\lref{eq:global:G:estimate};}{that for all 
			$n\in\N_0\cap[N-1,N+M)$, $m\in\N\cap(0,N+M-n]$ with 
			$\SPR_{n}\geq-\fC-\LR\alpha\Cst(3\fC+\fc)$, 
			$\max\{\SPR_n,\SPR_{n+1},\dots,\SPR_{n+m}\}\leq-\fC$, 
			and $\Cst(3\fC+2\fc)\geq\GRAD_n>0$ it holds that
			\begin{equation}
				\label{eq:lower:half:part:1}
				\begin{split}
					\SPR_{n+m}
					\geq \SPR_n-\LR\alpha
					\vass{\GRAD_n}(1-\alpha)^{-1}
					&\geq-\fC-\LR\alpha\Cst (3\fC+\fc)-{\LR\alpha\Cst
						(3\fC+2\fc)}\pr{1-\alpha}^{-1}
					\\&\geq -\fC -\LR\alpha\Cst(1-\alpha)^{-1}(6\fC+3\fc)
					\geq
					-3\fC-\fc
					.
				\end{split}
			\end{equation}
		}
			\argument{\lref{eq:lower:half:no:decrease:case};induction;}{that for all $n\in\N_0\cap[N-1,N+M)$, $m\in\N\cap(0,N+M-n]$ with 
			$\SPR_{n}\geq-\fC-\LR\alpha\Cst(3\fC+\fc)$,
			$\max\{\GRAD_n,\GRAD_{n+1},\dots,\GRAD_{n+m-1}\}\leq0$, and
			$\max\{\SPR_n,\SPR_{n+1},\dots,\SPR_{n+m}\}\leq-\fC$ it holds that
			\begin{equation}
				\llabel{eq:lower:half:part:2}
				\begin{split}
					\SPR_{n+m}
					&\geq \SPR_n
					\geq -\fC-\LR\alpha\Cst(3\fC+\fc)
					.
				\end{split}
			\end{equation}
		}
		\argument{
			\lref{eq:lower:half:part:2};
			\lref{eq:global:G:estimate};
			\eqref{eq:lower:half:part:1};}{that 
			for all $n\in\N_0\cap[N-1,N+M)$, $m\in\N\cap(0,N+M-n]$ with 
			$\SPR_{n}\geq-\fC-\LR\alpha\Cst(3\fC+\fc)$,
			$\vass{\GRAD_n}\leq\Cst(3\fC+2\fc)$, and 
			$\max\{\SPR_n,\SPR_{n+1},\dots,\SPR_{n+m}\}\leq-\fC$
			it holds that
			\begin{equation}
				\llabel{eq:lower:half}
				\begin{split}
					\vass{\SPR_{n+m}}
					=
					-\SPR_{n+m}
					&\leq 
					\max\{\fC+\LR\alpha\Cst(3\fC+\fc),3\fC+\fc\}
					\leq 
					3\fC+\fc
					.
				\end{split}
			\end{equation}
		}
			\argument{\lref{eq:lower:half}; 
			\eqref{eq:upper:half}}{
			that
			for all $n\in\N_0\cap[N-1,N+M)$, 
			$m\in\N\cap(0,N+M-n]$,
			$s\in\{-1,1\}$ with 
			$\vass{\SPR_{n}}\leq\fC+\LR\alpha\Cst(3\fC+\fc)$, 
			$\vass{\GRAD_n}\leq\Cst(3\fC+2\fc)$,
			and 
			$\min\{s\SPR_n,s\SPR_{n+1},\dots,s\SPR_{n+m}\}\geq\fC$
			it holds that
			\begin{equation}
				\llabel{eq:both:halfs}
				\begin{split}
					\textstyle
					\vass{\SPR_{n+m}}
					&\leq 
					3\fC+\fc
					.
				\end{split}
			\end{equation}
		}
			\argument{
			\lref{eq:sign:change:value:limit_pos};
			\lref{eq:sign:change:value:limit_neg};
			}{
			that 
			for all $n\in\N\cap[N,N+M]$ with  
			$\min\{\vass{\SPR_{n-1}},\vass{\SPR_{n}}\}\geq\fC$
			and
			$\max_{k\in\{1,2,\dots,n\}}\vass{\SPR_{k-2}}\leq 3\fC+\fc$
			there exists $s\in\{-1,1\}$ such that
			\begin{equation}
				\llabel{pre1:eq:no:sign:change:maintaining:high-values}
				\begin{split}
					\textstyle
					\min\{s\SPR_{n-1},s\SPR_{n}\}\geq \fC
					.
				\end{split}
			\end{equation}
		}
		\argument{
			\lref{pre1:eq:no:sign:change:maintaining:high-values};
			\lref{eq:recursion:Adam:prep:1};
		}{
			that 
			for all $n\in\N_0\cap[N-1,N+M)$ with  
			$\min\{\vass{\SPR_{n}},\vass{\SPR_{n+1}}\}\geq\fC$
			and
			$\max\{\vass{\SPR_{0}},\vass{\SPR_{1}},\dots,\vass{\SPR_{n}}\}\leq 3\fC+\fc$
			there exists $s\in\{-1,1\}$ such that
			\begin{equation}
				\llabel{pre:eq:no:sign:change:maintaining:high-values}
				\begin{split}
					\textstyle
					\min\{s\SPR_{n},s\SPR_{n+1}\}\geq \fC
					.
				\end{split}
			\end{equation}
		}
			\argument{
				\lref{eq:both:halfs};
				\lref{pre:eq:no:sign:change:maintaining:high-values};
			}
			{that 
			for all $n\in\N_0\cap[N-1,N+M)$ with 
			$\vass{\SPR_{n}}\leq\fC+\LR\alpha\Cst(3\fC+\fc)$, 
			$\vass{\GRAD_n}\leq\Cst(3\fC+2\fc)$,
			$\max\{\vass{\SPR_{0}},\vass{\SPR_{1}},\dots,\vass{\SPR_{n}}\}\leq 3\fC+\fc$,
			and 
			$\min\{\vass{\SPR_n},\vass{\SPR_{n+1}}\}\geq\fC$
			there exists $s\in\{-1,1\}$
			such that
			\begin{equation}
				\llabel{eq:both:halfs:pre}
				\begin{split}
					\textstyle
					\min\{s\SPR_{n},s\SPR_{n+1}\}\geq \fC
					\qqandqq
					\vass{\SPR_{n+1}}
					&\leq 
					3\fC+\fc
					.
				\end{split}
			\end{equation}
		}
				\argument{
					\lref{eq:both:halfs:pre};
					\lref{eq:both:halfs};
					\lref{pre:eq:no:sign:change:maintaining:high-values};
					induction;}
		{that 
			for all $n\in\N_0\cap[N-1,N+M)$, $m\in\N\cap(0,N+M-n]$ with 
			$\vass{\SPR_{n}}\leq\fC+\LR\alpha\Cst(3\fC+\fc)$, 
			$\vass{\GRAD_n}\leq\Cst(3\fC+2\fc)$,
			$\max\{\vass{\SPR_{0}},\vass{\SPR_{1}},\dots,\vass{\SPR_{n}}\}\leq 3\fC+\fc$,
			and 
			$\min\{\vass{\SPR_n},\vass{\SPR_{n+1}},\dots,\vass{\SPR_{n+m}}\}\geq\fC$
			there exists $s\in\{-1,1\}$
			such that
			\begin{equation}
				\llabel{eq:both:halfs:pre2}
				\begin{split}
					\textstyle
					\min\{s\SPR_{n},s\SPR_{n+1},\dots,s\SPR_{n+m}\}\geq \fC
					\qqandqq
					\max\{\vass{\SPR_{0}},\vass{\SPR_{1}},\dots,\vass{\SPR_{n+m}}\}
					&\leq 
					3\fC+\fc
					.
				\end{split}
			\end{equation}
		}
			\argument{
				\lref{eq:both:halfs:pre2};
			\lref{eq:apriori:bound:increments};}{
			that 
			for all $n\in\N_0\cap[N-1,N+M)$, $m\in\N\cap(0,N+M-n]$ with  
			$\vass{\SPR_{n}}\leq\fC+\LR\alpha\Cst(3\fC+\fc)$, 
			$\vass{\GRAD_n}\leq\Cst(3\fC+2\fc)$,
			$\max\{\vass{\SPR_{0}},\vass{\SPR_{1}},\dots,\vass{\SPR_{n}}\}\leq 3\fC+\fc$, and
			$\min\{\vass{\SPR_n},\vass{\SPR_{n+1}},\allowbreak\dots,\vass{\SPR_{n+m}}\}\geq\fC$
			it holds that
			\begin{equation}
				\llabel{eq:both:halfs:pre4}
				\begin{split}
					\textstyle
					\vass{\SPR_{n+m}}
					&\leq 
					3\fC+\fc
					\qqandqq
					\textstyle
					\vass{\GRAD_{n+m}}\leq\Cst\prb{\fc+\max_{k\in\{1,2,\dots,n+m+1\}}\vass{\SPR_{k-2}}}\leq\Cst(3\fC+2\fc)
					.
				\end{split}
			\end{equation}
		}
			\argument{
			\eqref{eq:one:step:for:small:theta};
			}{that for all $n\in\N_0\cap[N-1,N+M)$ with 
			$\vass{\SPR_n}\leq\fC$ and
			$\max\{\vass{\SPR_0},\vass{\SPR_1},\dots,\vass{\SPR_n}\}\leq 3\fC+\fc$
			it holds that
			\begin{equation}
				\llabel{eq:step:out:of:tiny:box}
				\begin{split}
					\textstyle
					\vass{\SPR_{n+1}}
					\leq\fC+\LR\alpha\Cst\max_{k\in\{1,2,\dots,n+1\}}\vass{\SPR_{k-2}}
					\leq\fC+\LR\alpha\Cst(3\fC+ \fc)
				.
				\end{split}
			\end{equation}
		}
			\argument{
				\lref{eq:step:out:of:tiny:box};
			\lref{eq:apriori:bound:increments};}{
			that for all $n\in\N_0\cap[N-1,N+M)$ with 
			$\vass{\SPR_n}\leq\fC$
			 and
			$\max\{\vass{\SPR_0},\vass{\SPR_1},\dots,\vass{\SPR_n}\}\leq 3\fC+\fc$
			it holds that
			\begin{equation}
				\label{eq:leaving:small:area}
				\begin{split}
					\textstyle
					\vass{\SPR_{n+1}}
					\leq \fC+\LR\alpha\Cst(3\fC+ \fc)
					\qandq
					\vass{\GRAD_{n+1}}
				 \leq\Cst\prb{\fc+\max_{k\in\{1,2,\dots,n+2\}}\vass{\SPR_{k-2}}}
					\leq\Cst(3\fC+2\fc).
				\end{split}
			\end{equation}
		}
			\argument{\lref{eq:setup:1st:part:Adam:aprior0}}{that 
			\begin{equation}
				\llabel{eq:induction:base}
				\vass{\SPR_{N-1}}\leq \fC
				\qqandqq
				\max\{\vass{\SPR_0},\vass{\SPR_1},\dots,\vass{\SPR_{N-1}}\}\leq 3\fC+\fc
				.
			\end{equation}
		}
			\argument{
				\lref{eq:global:G:estimate};
				\lref{eq:induction:base};
				\lref{eq:both:halfs:pre4};
				\eqref{eq:leaving:small:area};
				induction;}
			{for all $n\in\N\cap[N,N+M]$ that
			\begin{equation}
				\llabel{123784}
				\max\{\vass{\SPR_0},\vass{\SPR_1},\dots,\vass{\SPR_n}\}\leq 3\fC+\fc
				\qqandqq
				\vass{\GRAD_{n}}
				\leq\Cst(3\fC+2\fc)
				.
			\end{equation}
		}
		\argument{
			\lref{123784};
			\lref{eq:recursion:Adam:prep:1};
			\lref{eq:setup:1st:part:Adam:aprior0};
			}{\cref{it:upper:bound:theta:base} in the case $\alpha>0$}.
		\end{aproof}

		\subsection{A priori bounds for Adam and other adaptive GD optimization methods}

		\begin{athm}{lemma}{lem:hoelder:grad:prep}
			Let $\alpha\in[0,1)$, $\beta\in(\alpha^2,1)$, $\bscl,\eps\in(0,\infty)$, $n\in\N$, $g_1,g_2,\dots,g_n\in\R$.
			Then
			\begin{equation}
				\label{eq:result:lem:hoelder:grad:prep}
				\frac{\vass{ \sum_{k=1}^n \alpha^{n-k}						g_k}}{\eps+\PR{\sum_{k=1}^n\bscl\beta^{n-k}\pr{g_k}^2 }^{\nicefrac{1}{2}}} 
				\leq 
				\bscl^{-\nicefrac{1}{2}}(1-\alpha^2\beta^{-1})^{-\nicefrac{1}{2}}.
			\end{equation}
		\end{athm}

		\begin{aproof}
		\startnewargseq
		\argument{
				the fact that $0\leq\alpha^2<\beta$;
				the Hölder inequality;}
			{that 
				\begin{equation}
					\llabel{eq:hoelder}
					\begin{split}
						\vass{\textstyle\sum_{k=1}^n\alpha^{n-k}g_k}
						&\leq
						\textstyle\sum_{k=1}^n\alpha^{n-k}\vass{g_k}
						\\&=\textstyle
						\bscl^{-\nicefrac{1}{2}}
						\textstyle\PR{\sum_{k=1}^n\alpha^{n-k}\beta^{\frac{k-n}{2}}\bscl^{\nicefrac{1}{2}}\beta^{\frac{n-k}{2}}\vass{g_k}}
						\\&\textstyle\leq
						\bscl^{-\nicefrac{1}{2}}
						\PR{\textstyle\sum_{k=1}^n\alpha^{2n-2k}\beta^{k-n}}^{\nicefrac{1}{2}}\PR{\sum_{k=1}^n\bscl\beta^{n-k}\pr{g_k}^2 }^{\nicefrac{1}{2}}
						\\&\textstyle=
						\bscl^{-\nicefrac{1}{2}}\PR{\sum_{k=1}^n\bscl\beta^{n-k}\pr{g_k}^2 }^{\nicefrac{1}{2}}\PR{\textstyle\sum_{k=0}^{n-1}(\alpha^2\beta^{-1})^{k}}^{\nicefrac{1}{2}}
						\\&\textstyle\leq
						\bscl^{-\nicefrac{1}{2}}\PR{\sum_{k=1}^n\bscl\beta^{n-k}\pr{g_k}^2 }^{\nicefrac{1}{2}}\PR{\textstyle\sum_{k=0}^{\infty}(\alpha^2\beta^{-1})^{k}}^{\nicefrac{1}{2}}
						\\&=\textstyle
						\bscl^{-\nicefrac{1}{2}}
						\PR{\sum_{k=1}^n\bscl\beta^{n-k}\pr{g_k}^2 }^{\nicefrac{1}{2}}\pr{1-\alpha^2\beta^{-1}}^{\nicefrac{-1}{2}}
						.
					\end{split}
				\end{equation}
			}
			\argument{\lref{eq:hoelder};the fact that $\eps>0$;}{that
			\begin{equation}
				\llabel{eq:incremet:upper:bound:class}
				\begin{split}
					\frac{ \vass{\sum_{k=1}^n \alpha^{n-k}						g_k}}{\eps+\PR{\sum_{k=1}^n\bscl\beta^{n-k}\pr{g_k}^2 }^{\nicefrac{1}{2}}} 
					&\leq 
					\frac{\bscl^{-\nicefrac{1}{2}}\PR{\sum_{k=1}^n\bscl\beta^{n-k}\pr{g_k}^2 }^{\nicefrac{1}{2}}\pr{1-\alpha^2\beta^{-1}}^{\nicefrac{-1}{2}}}{\eps+\PR{\sum_{k=1}^n\bscl\beta^{n-k}\pr{g_k}^2 }^{\nicefrac{1}{2}}}
					\\&\leq 
					\bscl^{-\nicefrac{1}{2}}(1-\alpha^2\beta^{-1})^{-\nicefrac{1}{2}}.
				\end{split}
			\end{equation}
		}
			\argument{\lref{eq:incremet:upper:bound:class};}{\eqref{eq:result:lem:hoelder:grad:prep}}.
		\end{aproof}
		
		\begin{athm}{prop}{prop:one_dim:Adam}
			Let $\alpha\in[0,1)$, $\beta\in(\alpha^2,1)$, $\bscl,\eps\in(0,\infty)$,  $\gamma,\MOM_0,\MOM_1,\fc,\cst\in[0,\infty)$, 
			$n\in\N$,
			$g_0,g_1,\dots,g_n,\Theta_0,\Theta_1\in\R$  satisfy 
			\begin{equation}
				\llabel{eq:recursion:Adam:prep:class}
				\begin{split}
					\Theta_1
					&=  
					\Theta_{0 }
					-
					\frac{\gamma\PR{\sum_{k=0}^n (1-\alpha)\alpha^{n-k}g_k}}{\eps+\PR{\MOM_1}^{\nicefrac{1}{2}}}
					,
			\qquad\qquad
			\MOM_1\geq\beta^n\MOM_0+\sum_{k=1}^n\bscl\beta^{n-k}\pr{g_k}^2,
				\end{split}
			\end{equation}
%
%
			and $\vass{\Theta_{0}}	
			\leq	
			\fc +
			\cst\vass{g_n}$. Then 
			\begin{equation}
				\label{eq:result:prop:a_priori_bound_one_dimensional:Adam:prep:class}
				\begin{split}
					&\vass{\Theta_1}
					\leq
					\fc
					+\frac{ \cst\PR{\MOM_1}^{\nicefrac{1}{2}}}{\bscl	^{\nicefrac{1}{2}}}
					+\frac{\gamma(1-\alpha)\alpha^n\vass{g_0}}{\eps+\PR{\MOM_1}^{\nicefrac{1}{2}}}
					+\frac{\gamma (1-\alpha)}{\bscl^{\nicefrac{1}{2}}(1-\alpha^2\beta^{-1})^{\nicefrac{1}{2}}}
					.
				\end{split}
			\end{equation}
		\end{athm}

		\begin{aproof}
			\argument{\lref{eq:recursion:Adam:prep:class};\cref{lem:hoelder:grad:prep};}{that
			\begin{equation}
				\label{eq:incremet:upper:bound:class}
				\begin{split}
					\frac{\vass{\gamma \sum_{k=1}^n (1-\alpha)\alpha^{n-k}						g_k}}{\eps+\PR{\MOM_1}^{\nicefrac{1}{2}}} 
					&\leq
					\frac{\gamma \vass{\sum_{k=1}^n (1-\alpha)\alpha^{n-k}						g_k}}{\eps+\PR{\beta^n\MOM_0+\sum_{k=1}^n\bscl\beta^{n-k}\pr{g_k}^2}^{\nicefrac{1}{2}}} 
					\\&\leq
					\frac{ \gamma (1-\alpha)\vass{\sum_{k=1}^n\alpha^{n-k}						g_k}}{\eps+\PR{\sum_{k=1}^n\bscl\beta^{n-k}\pr{g_k}^2}^{\nicefrac{1}{2}}} 
					\leq
					\frac{\gamma (1-\alpha)}{\bscl^{\nicefrac{1}{2}}(1-\alpha^2\beta^{-1})^{\nicefrac{1}{2}}}.
				\end{split}
			\end{equation}
		}
			\argument{\lref{eq:recursion:Adam:prep:class}}{that 
			\begin{equation}
				\llabel{eq:TEST}
				\begin{split}
													\vass{\Theta_{0}}
											\leq\fc+\cst \vass{g_n}
											&=
											\fc+\cst \PRb{\bscl\pr{g_n}^2}^{\nicefrac{1}{2}} \bscl^{\nicefrac{-1}{2}}
											\leq
											\fc+\cst\PR{\MOM_1}^{\nicefrac{1}{2}}\bscl^{\nicefrac{-1}{2}}.
				\end{split}
			\end{equation}
		}
		\argument{
			\lref{eq:TEST};
			\lref{eq:recursion:Adam:prep:class};
			\eqref{eq:incremet:upper:bound:class};}{
			that 
			\begin{equation}
				\llabel{eq:rough:UB}
				\begin{split}
					\vass{\Theta_1}
					&=\vass[\bigg]{\Theta_{ 0 }
						-
						\frac{\gamma\PR{\sum_{k=0}^n (1-\alpha)\alpha^{n-k}g_k}}{\eps+\PR{\MOM_1}^{\nicefrac{1}{2}}}}
					\\&\leq		\vass{\Theta_{0}}
					+
					\frac{\vass{\gamma(1-\alpha)\alpha^n g_0}+\vass{\gamma \sum_{k=1}^n (1-\alpha)\alpha^{n-k}g_k}}{\eps+\PR{\MOM_1}^{\nicefrac{1}{2}}} 
					\\&\leq		
					\fc
					+\frac{ \cst\PR{\MOM_1}^{\nicefrac{1}{2}}}{\bscl	^{\nicefrac{1}{2}}}
					+\frac{\gamma(1-\alpha)\alpha^n\vass{g_0}}{\eps+\PR{\MOM_1}^{\nicefrac{1}{2}}}
					+\frac{\gamma (1-\alpha)}{\bscl^{\nicefrac{1}{2}}(1-\alpha^2\beta^{-1})^{\nicefrac{1}{2}}}
					.
				\end{split}
			\end{equation}
		}
		\argument{\lref{eq:rough:UB}}{\eqref{eq:result:prop:a_priori_bound_one_dimensional:Adam:prep:class}}.
		\end{aproof}
		
		\begin{athm}{cor}{cor:a_priori_bound_gen:momentum}
			Let 
			$\alpha\in[0,1)$, $\beta\in(\alpha^2,1)$, 
			$\eps,\cst\in(0,\infty)$,
			$ \Cst\in[\cst,\infty)$,
			$\fc\in\R$,
			let
			$g_n\colon\R\to\R$, $n\in \N_0$,
			satisfy for all $n\in\N$, $\theta\in\R$ that
			\begin{equation}
				\label{eq:setup:gen_grad:2.1}
				\pr{\theta-\fc}
				\prb{\cst+(\Cst-\cst)\indicator{(-\infty,\fc]}(\theta)}
				\leq
				g_n(\theta)
				\leq
				\pr{\theta+\fc}
				\prb{\cst+(\Cst-\cst)\indicator{[-\fc,\infty)}(\theta)},
			\end{equation}
			let 
			$\bscl\colon\N\to(0,\infty)$ satisfy $\inf_{n\in\N}\bscl_n>0$, and let
			$
			\gamma \colon \N \to [0,\infty)
			$, 
			$\MOM\colon\N_0\to[0,\infty)$,
			and 
			$ \Theta \colon \Z \to \R $ satisfy 
			for all $ n \in \N $ that
			\begin{equation}
				\label{eq:recursion_cor:momentum2} 
				\Theta_n
				= 
				\Theta_{ n - 1 }
				-
				\frac{\gamma_n\PR{\sum_{k=0}^n (1-\alpha)\alpha^{n-k}g_k\pr
						{\Theta_{k-1}}}}{\eps+\PR{\MOM_n}^{\nicefrac{1}{2}}}
				,
				\qquad
				\MOM_n\geq\beta^n\MOM_0+\sum_{k=1}^n \bscl_n\beta^{n-k}\pr
				{g_k(\Theta_{k-1})}^2
,
			\end{equation}
			and
			$\vass{g_0(\Theta_{-1})}\leq \Cst(\vass{\Theta_0}+\fc)$.
			Then 
			\begin{multline}
				\label{eq:a_priori_to_prove_cor:momentum2} 
				\sup_{ n \in \N_0 }
				\vass{\Theta_n}\leq
				\fc
				+3\max\pRbbb{
					\vass{\Theta_0},
					\fc+\frac{\fc\alpha\Cst}{(1-\alpha)\cst},
					\fc
					+
					\PRbbb{\sup_{n\in\N}\gamma_n}
					\prbbb{\frac{\pr{1-\alpha}\vass{g_0(\Theta_{-1})}}{\eps+\PR{\MOM_0}^{\nicefrac{1}{2}}}
						\\+\frac{\max\{1,\Cst\}(2+\alpha)\beta^{\nicefrac{1}{2}}}{\PR{\textstyle\inf_{n\in\N}\bscl_n}^{\nicefrac{1}{2}}\cst\pr{\beta^{\nicefrac{1}{2}}-\alpha}}}
				}
				.
			\end{multline}
		\end{athm}

		\begin{aproof}
			Throughout this proof 
			assume without loss of generality that
			$\eps\pr{1-\alpha}\leq\PR{\sup_{n\in\N}\gamma_n}(1+2\alpha)\max\{1,\Cst\}$ (cf.\ \cref{prop:a_priori_bound_one_dim:Adam:1}),
			let $D\in\R$ satisfy
			\begin{multline}
				\label{eq:setup:divine:bound2}
				D=3\max\pRbbb{
					\vass{\Theta_0},
					\fc+\frac{\fc\alpha\Cst}{(1-\alpha)\cst}
					,
					\fc
					+
					\PRbbb{\sup_{n\in\N}\gamma_n}
					\prbbb{\frac{\pr{1-\alpha}\vass{g_0(\Theta_{0})}}{\eps+\PR{\MOM_0}^{\nicefrac{1}{2}}}
					\\+\frac{\max\{1,\Cst\}(2+\alpha)\beta^{\nicefrac{1}{2}}}{\PR{\textstyle\inf_{n\in\N}\bscl_n}^{\nicefrac{1}{2}}\cst\pr{\beta^{\nicefrac{1}{2}}-\alpha}}}
				}
				,
			\end{multline}
			and 
			let $S\in\R$ satisfy
			\begin{equation}
				\label{eq:setup:arbitrary:gammas:2}
				S
				= \frac{\PR{\sup_{m\in\N}\gamma_m}(1+2\alpha)\max\{1,\Cst\}}{1-\alpha}-\eps
				.
			\end{equation}
			\startnewargseq
			\argument{
				\eqref{eq:setup:gen_grad:2.1};
				\eqref{eq:recursion_cor:momentum2};}{that for all $n\in\N$, $\theta\in\R$ it holds that
			\begin{equation}
				\llabel{eq:basic:prop:setup}
				\fc\geq 0
				\qqandqq
				\vass{\theta}\leq \fc+\cst^{-1}\vass{g_n(\theta)}.
			\end{equation}
		}
			\argument{\lref{eq:basic:prop:setup};
			\eqref{eq:recursion_cor:momentum2};
			\eqref{eq:setup:arbitrary:gammas:2};
			\cref{prop:one_dim:Adam} (applied for every $n\in\N$ with
			$\eps\curvearrowleft\eps$,
			$\cst\curvearrowleft\cst^{-1}$,
			$\bscl\curvearrowleft \bscl_n$,
			$\alpha\curvearrowleft\alpha$,
			$\beta\curvearrowleft\beta$,
			$\gamma\curvearrowleft\gamma_n$,
			$\fc\curvearrowleft\fc$,
			$\MOM_0\curvearrowleft\MOM_0$,
			$\MOM_1\curvearrowleft\MOM_n$,
			$n\curvearrowleft n$,
			$(g_0,g_1,\dots,g_n) \curvearrowleft (g_0(\Theta_{-1}),g_1(\Theta_{0}),\dots,g_n(\Theta_{n-1}))$, 
			$\Theta_0\curvearrowleft \Theta_{n-1}$,
			$\Theta_1\curvearrowleft \Theta_n$
			in the notation of \cref{prop:one_dim:Adam});}{that for all $n\in\N$ it holds that
			\begin{equation}
				\llabel{eq:bound:for:large:lamda}
				\begin{split}
				\vass{
					\Theta_{n}}
				&\leq
					\fc
					+\frac{ \PR{\MOM_n}^{\nicefrac{1}{2}}}{\cst(\bscl_n)	^{\nicefrac{1}{2}}}
					+\frac{\gamma_n(1-\alpha)\alpha^n\vass{g_0(\Theta_{-1})}}{\eps+\PR{\MOM_n}^{\nicefrac{1}{2}}}
					+\frac{\gamma_n (1-\alpha)}{(\bscl_n)^{\nicefrac{1}{2}}(1-\alpha^2\beta^{-1})^{\nicefrac{1}{2}}}
				.
				\end{split}
			\end{equation}
		}
		\argument{\lref{eq:bound:for:large:lamda};
			\eqref{eq:recursion_cor:momentum2};
			\eqref{eq:setup:arbitrary:gammas:2};
			the fact that $\alpha^2<\beta$;}{that for all $n\in\N$ with $\MOM_n\leq S^2$ it holds that
			\begin{equation}
				\llabel{eq:lambda:well:bounded}
				\begin{split}
					\vass{
						\Theta_{n}}
					&\leq
					\fc
					+\frac{ \PR{\MOM_n}^{\nicefrac{1}{2}}}{\cst(\bscl_n)	^{\nicefrac{1}{2}}}
					+\frac{\gamma_n(1-\alpha)\alpha^n\vass{g_0(\Theta_{0})}}{\eps+\PR{\MOM_n}^{\nicefrac{1}{2}}}
					+\frac{\gamma_n (1-\alpha)\beta^{\nicefrac{1}{2}}}{(\bscl_n)^{\nicefrac{1}{2}}(\beta-\alpha^2)^{\nicefrac{1}{2}}}
					\\&\leq
					\fc
					+\frac{\eps+S}{\cst\PR{\textstyle\inf_{m\in\N}\bscl_m}^{\nicefrac{1}{2}}}
					+\frac{\gamma_n(1-\alpha)\alpha^n\vass{g_0(\Theta_{0})}}{\eps+\PR{\beta^n\MOM_0}^{\nicefrac{1}{2}}}
					+\frac{\PR{\sup_{m\in\N}\gamma_m} (1-\alpha)\beta^{\nicefrac{1}{2}}}{\PR{\textstyle\inf_{m\in\N}\bscl_m}^{\nicefrac{1}{2}}(\beta-\alpha^2)^{\nicefrac{1}{2}}}
					\\&=
					\fc
					+\frac{\gamma_n(1-\alpha)(\alpha^2\beta^{-1})^{\nicefrac{n}{2}}\vass{g_0(\Theta_{0})}}{\beta^{-\nicefrac{n}{2}}\eps+\PR{\MOM_0}^{\nicefrac{1}{2}}}
					+\frac{\PR{\sup_{m\in\N}\gamma_m}(1+2\alpha)\max\{1,\Cst\}}{\cst(1-\alpha)\PR{\textstyle\inf_{m\in\N}\bscl_m}^{\nicefrac{1}{2}}}
					\\&\quad+\frac{\PR{\sup_{m\in\N}\gamma_m} (1-\alpha)\beta^{\nicefrac{1}{2}}}{\PR{\textstyle\inf_{m\in\N}\bscl_m}^{\nicefrac{1}{2}}(\beta-\alpha^2)^{\nicefrac{1}{2}}}
					\\&\leq
					\fc
					+
					\PRbbb{\sup_{m\in\N}\gamma_m}\prbbb{\frac{(1-\alpha)\vass{g_0(\Theta_{0})}}{\eps+\PR{\MOM_0}^{\nicefrac{1}{2}}}
					+\frac{\max\{1,\Cst\}}{\cst\PR{\textstyle\inf_{m\in\N}\bscl_m}^{\nicefrac{1}{2}}}
					\prbbb{
						\frac{1+2\alpha}{1-\alpha}
						+\frac{ (1-\alpha)\beta^{\nicefrac{1}{2}}}{(\beta-\alpha^2)^{\nicefrac{1}{2}}}}}
					\\&\leq
					\fc
					+
					\PRbbb{\sup_{m\in\N}\gamma_m}\prbbb{
					\frac{(1-\alpha)\vass{g_0(\Theta_{0})}}{\eps+\PR{\MOM_0}^{\nicefrac{1}{2}}}
					+\frac{\max\{1,\Cst\}}{\cst\PR{\textstyle\inf_{m\in\N}\bscl_m}^{\nicefrac{1}{2}}}
					\prbbb{
						\frac{(1+2\alpha)\beta^{\nicefrac{1}{2}}}{\beta^{\nicefrac{1}{2}}-\alpha}
						+\frac{ (1-\alpha)\beta^{\nicefrac{1}{2}}}{\beta^{\nicefrac{1}{2}}-\alpha}}
					}
					\\&=
					\fc
					+
					\PRbbb{\sup_{m\in\N}\gamma_m}\prbbb{
					\frac{(1-\alpha)\vass{g_0(\Theta_{0})}}{\eps+\PR{\MOM_0}^{\nicefrac{1}{2}}}
					+
					\frac{\max\{1,\Cst\}(2+\alpha)\beta^{\nicefrac{1}{2}}}{\PR{\textstyle\inf_{m\in\N}\bscl_m}^{\nicefrac{1}{2}}\cst\pr{\beta^{\nicefrac{1}{2}}-\alpha}}
				}
					\leq \frac{D}{3}
					.
				\end{split}
			\end{equation}
		}
		\argument{\lref{eq:lambda:well:bounded};
			\eqref{eq:setup:divine:bound2};}{that 
			\begin{equation}
				\label{eq:divine:small:mus}
				3\vass{\Theta_0}\leq D
				\qqandqq
				\forall\, n\in\{m\in\N\colon \MOM_m\leq S^2\}\colon
				3\vass{\Theta_n}\leq D.
			\end{equation}
		}
			\argument{\eqref{eq:setup:arbitrary:gammas:2};}
			{for all $n\in\N$ with $\MOM_n> S^2$ it holds that
			\begin{equation}
				\llabel{eq:LR:prop}
				\begin{split}
					\frac{\gamma_n }{\eps+\PR{\MOM_n}^{\nicefrac{1}{2}}} 
					\leq\frac{\gamma_n }{\eps+S} 
					&	=
					\frac{\gamma_n(1-\alpha)}{\PR{\sup_{m\in\N}\gamma_m}(1+2\alpha)\max\{1,\Cst\}}
					\leq 
					\frac{1-\alpha}{(1+2\alpha)\max\{1,\Cst\}}
					.
				\end{split}
			\end{equation}
		}
			\argument{\lref{eq:LR:prop};
			\eqref{eq:setup:gen_grad:2.1};
			\eqref{eq:recursion_cor:momentum2};
			\eqref{eq:setup:arbitrary:gammas:2};
			\lref{eq:basic:prop:setup};
			\cref{prop:a_priori_bound_one_dim:Adam:1} (applied for every $N\in\N$, $M\in\N_0$ with
			$\alpha\curvearrowleft\alpha$,
			$\fc\curvearrowleft\fc$,
			$\cst\curvearrowleft\cst$,
			$\Cst\curvearrowleft\Cst$,
			$\pars\curvearrowleft\pars$,
			$N\curvearrowleft N$,
			$M\curvearrowleft M$,
			$\gamma\curvearrowleft\prb{\N\ni n\mapsto \gamma_n\pr{\eps+\PR{\MOM_n}^{\nicefrac{1}{2}}}^{-1}\in [0,\infty)}$,	
			$\Theta\curvearrowleft\pr{\N_0\ni n\mapsto \Theta_n\in\R}$,
			$g\curvearrowleft \pr{\N\ni n\mapsto g_n(\Theta_{n-1})\in\R}$
			in the notation of \cref{prop:a_priori_bound_one_dim:Adam:1});}{
			that for all $N\in\N$, $M\in\{m\in\N_0\colon\forall\, n\in\N\cap[N,N+m]\colon \MOM_n>S^2\}$ it holds that
			\begin{equation}
				\llabel{eq:bound:for:large:M}
				\textstyle\max_{n\in\N\cap[N,N+M]}\vass{\Theta_n}
				\leq
				\max\pRbb{4\fc+\frac{3\fc\alpha\Cst}{(1-\alpha)\cst},
					\fc+3\vass{\Theta_{N-1}},\textstyle\max_{k\in\{1,2,\dots,N\}}\vass{\Theta_{k-1}}}
				.
			\end{equation}
		}
			\argument{
				\lref{eq:bound:for:large:M};
				\eqref{eq:setup:divine:bound2};}
				{for all $N\in\N$, 
			$M\in\{m\in\N_0\colon\forall\, n\in\N\cap[N,N+m]\colon \pr{\MOM_n>S^2}\wedge\pr{3\vass{\Theta_{N-1}}\leq D}\wedge\pr{\max_{k\in\{1,2,\dots,N\}}\vass{\Theta_{k-1}}\leq \fc+D}\}$ that
			\begin{equation}
				\label{eq:final:for:a:priori:upper}
				\begin{split}
					\textstyle\max_{n\in\N\cap[N,N+M]}\vass{\Theta_n}
					&\textstyle\leq \max\pRbb{4\fc+\frac{3\fc\alpha\Cst}{(1-\alpha)\cst},\fc+3\vass{\Theta_{N-1}},\textstyle\max_{k\in\{1,2,\dots,N\}}\vass{\Theta_{k-1}}}
					\leq \fc + D.
				\end{split}
			\end{equation}
		}
			\argument{the fact that for all $N\in\{n\in\N\colon\MOM_n>S^2\}$ it holds that
			\begin{equation}
				\max\{M\in\N_0\cap[0,N)\colon(\forall\, m\in\N\cap[N-M,N]\colon \MOM_m>S^2)\}\in\N_0
			\end{equation}
			;
			\eqref{eq:divine:small:mus} }
			{
			that for all $N\in\{n\in\N\colon\MOM_n>S^2\}$ there exists $M\in\N_0$
			such that for all $n\in\N\cap[N-M,(N-M)+M]$ it holds that
			\begin{equation}
				\llabel{eq:decomposition:possible}
				\MOM_n>S^2
				\qqandqq
				3\vass{\Theta_{\max\{N-M-1,0\}}}\leq D.
			\end{equation}
		}
		\argument{
			\lref{eq:decomposition:possible};
			\eqref{eq:final:for:a:priori:upper};
			induction;}
			{that for all $N\in\{n\in\N\colon\pr{\max_{k\in\{1,2,\dots,n\}}\vass{\Theta_{k-1}}\leq \fc+D}\wedge\pr{\MOM_n>S^2}\}$ it holds that
			\begin{equation}
				\llabel{eq:base:case}
				\vass{\Theta_N}\leq\fc+D.
			\end{equation}
		}
		\argument{\lref{eq:base:case};
			\eqref{eq:divine:small:mus};
			induction;}
			{that for all $n\in\N_0$ it holds that
			\begin{equation}
				\llabel{eq:fin:UB}
				\textstyle
				\max_{k\in\{0,1,\dots,n\}}\vass{\Theta_{k}}\leq \fc+D.
			\end{equation}
		}
			 \argument{\lref{eq:fin:UB};}{
			\eqref{eq:a_priori_to_prove_cor:momentum2}}.
		\end{aproof}
		
		\newcommand{\ind}{m}
		
		\begin{athm}{cor}{cor:a_priori_bound_gen:momentum:tilde}
			Let 
			$\pars\in\N$,
			$i\in\{1,2,\dots,\pars\}$, 
			$\alpha\in[0,1)$, $\beta\in(\alpha^2,1)$,
			$\eps,\cst\in(0,\infty)$,
			$ \Cst\in[\cst,\infty)$,
			$\fc\in\R$,
			let
			$g_n\colon\R^{\pars}\to\R$, $n\in\N_0$,
			satisfy for all $n\in\N$, $\theta=(\theta_1,\dots,\theta_\pars)\in\R^{\pars}$ that
			\begin{equation}
				\label{eq:setup:gen_grad:2:tilde}
				\pr{\theta_i-\fc}
				\prb{\cst+(\Cst-\cst)\indicator{(-\infty,\fc]}(\theta_i)}
				\leq
				g_n(\theta)
				\leq
				\pr{\theta_i+\fc}
				\prb{\cst+(\Cst-\cst)\indicator{[-\fc,\infty)}(\theta_i)},
			\end{equation}
			let 
			$\bscl\colon\N\to(0,\infty)$ satisfy $\inf_{n\in\N}\bscl_n>0$, and let
			$
			\gamma \colon \N \to [0,\infty)
			$, 
			$\MOM\colon\N_0\to[0,\infty)$,
			and 
			$ \Theta=(\Theta^{(1)},\dots,\Theta^{(\pars)}) \colon \Z \to \R^{\pars} $ satisfy 
			for all $ n \in \N $ that
						\begin{equation}
				\label{eq:recursion_cor:momentum2:tilde} 
				\Theta_n^{(i)}
				= 
				\Theta_{ n - 1 }^{(i)}
				-
				\frac{\gamma_n\PR{\sum_{k=0}^n (1-\alpha)\alpha^{n-k}g_k\pr
						{\Theta_{k-1}}}}{\eps+\PR{\MOM_n}^{\nicefrac{1}{2}}}
				,
				\quad
				\MOM_n\geq\beta^n\MOM_0+\sum_{k=1}^n \bscl_n\beta^{n-k}\pr
				{g_k(\Theta_{k-1})}^2
				,
			\end{equation}
			and
			$\vass{g_0(\Theta_{-1})}\leq \Cst\prb{\vass{\Theta_0^{(i)}}+\fc}$.
			Then 
			\begin{multline}
				\label{eq:a_priori_to_prove_cor:momentum2:tilde} 
				\sup_{ n \in \N_0 }
				\vass{\Theta_n^{(i)}}\leq
				\fc
				+3\max\pRbbb{
					\vass{\Theta_0^{(i)}},
					\fc+\frac{\fc\alpha\Cst}{(1-\alpha)\cst},
					\fc
					+
					\PRbbb{\sup_{n\in\N}\gamma_n}
					\prbbb{\frac{\pr{1-\alpha}\vass{g_0(\Theta_{-1})}}{\eps+\PR{\MOM_0}^{\nicefrac{1}{2}}}
						\\+\frac{\max\{1,\Cst\}(2+\alpha)\beta^{\nicefrac{1}{2}}}{\PR{\textstyle\inf_{n\in\N}\bscl_n}^{\nicefrac{1}{2}}\cst\pr{\beta^{\nicefrac{1}{2}}-\alpha}}}
				}
				.
			\end{multline}
		\end{athm}

		\begin{aproof}
			Throughout this proof for every $n\in\N_0$ let $f_n\colon\R\to\R$ satisfy for all $\theta\in\R$ that
			\begin{equation}
				\label{eq:dim:red:for:incr}
				f_n(\theta)=g_n\prb{\Theta_{n-1}^{(1)},\Theta_{n-1}^{(2)},\dots,\Theta^{(i-1)}_{n-1},\theta,\Theta^{(i+1)}_{n-1},\dots,\Theta^{(\pars)}_{n-1}}
				.
			\end{equation}
			\startnewargseq
			\argument{\eqref{eq:setup:gen_grad:2:tilde};}{that $\fc\geq 0$}.
			\argument{
				\eqref{eq:setup:gen_grad:2:tilde};
				\eqref{eq:dim:red:for:incr}
			}{that for all $n\in\N$, $\theta\in\R$ it holds that
			\begin{equation}
				\label{eq:setup:gen_grad:sing:dim:ineq}
				\pr{\theta-\fc}
				\pr{\cst+(\Cst-\cst)\indicator{(-\infty,\fc]}(\theta)}
				\leq
				f_n(\theta)
				\leq
				\pr{\theta+\fc}
				\pr{\cst+(\Cst-\cst)\indicator{[-\fc,\infty)}(\theta)}
				.
			\end{equation}
		}
			\argument{\eqref{eq:dim:red:for:incr}}{for all $n\in\N$ that
			\begin{equation}
				\llabel{eq:setup:gen_grad:sing:dim:ineq2}
				f_n\prb{\Theta_{n-1}^{(i)}}
				=g_n(\Theta_{n-1})
				\qqandqq
				\vass{f_0\prb{\Theta_{-1}^{(i)}}}=\vass{g_0(\Theta_{-1})}\leq \Cst\prb{\vass{\Theta_0^{(i)}}+\fc}
				.
			\end{equation}
		}
			\argument{\lref{eq:setup:gen_grad:sing:dim:ineq2};
			\eqref{eq:recursion_cor:momentum2:tilde} ;}{
			that for all $n\in\N$ it holds that
			\begin{equation}
				\llabel{eq:recursion_cor:momentum2:tilde3} 
				\begin{split}
					\textstyle
					\MOM_n
					\geq
					\beta^n\MOM_0+\sum_{k=1}^n \bscl_n\beta^{n-k}\pr
					{g_k(\Theta_{k-1})}^2
					=
					\beta^n\MOM_0+\sum_{k=1}^n \bscl_n\beta^{n-k}\prb
					{f_k\prb{\Theta_{k-1}^{(i)}}}^2
					.
				\end{split}
			\end{equation}
		}
			\argument{
			\eqref{eq:recursion_cor:momentum2:tilde}; 
			\lref{eq:setup:gen_grad:sing:dim:ineq2} 
		}
			{that for all $n\in\N$ it holds that
			\begin{equation}
				\llabel{eq:multdim:ver}
				\begin{split}
					\Theta_n^{(i)}-\Theta_{ n - 1 }^{(i)}
					&= 
					-
					\frac{\gamma_n\PR{\sum_{k=0}^n (1-\alpha)\alpha^{n-k}g_k\pr
							{\Theta_{k-1}}}}{\eps+\PR{\MOM_n}^{\nicefrac{1}{2}}}
					= 
					-
					\frac{\gamma_n\PRb{ \sum_{k=0}^n (1-\alpha)\alpha^{n-k}f_k\prb{\Theta_{ k - 1 }^{(i)}}}}{\eps+\PR{\MOM_n}^{\nicefrac{1}{2}}}
					.
				\end{split}
			\end{equation}
		}
		\argument{\lref{eq:multdim:ver};
			\eqref{eq:recursion_cor:momentum2:tilde};
			\eqref{eq:setup:gen_grad:sing:dim:ineq};
			\lref{eq:setup:gen_grad:sing:dim:ineq2};
			\lref{eq:recursion_cor:momentum2:tilde3};
			\cref{cor:a_priori_bound_gen:momentum} (applied with 
			$\eps\curvearrowleft\eps$,
			$\cst\curvearrowleft\cst$,
			$\Cst\curvearrowleft\Cst$,
			$\alpha\curvearrowleft\alpha$,
			$\beta\curvearrowleft\beta$,
			$\fc\curvearrowleft\fc$,
			$(g_n)_{n\in\N_0}\curvearrowleft (f_n)_{n\in\N_0}$,
			$\bscl\curvearrowleft\bscl$,
			$\gamma\curvearrowleft\gamma$,
			$\MOM\curvearrowleft\MOM$,
			$\Theta\curvearrowleft\Theta^{(i)}$
			in the notation of \cref{cor:a_priori_bound_gen:momentum});}{that
						\begin{multline}
								\llabel{it:upper:bound:theta} 
				\sup_{ n \in \N_0 }
				\vass{\Theta_n^{(i)}}\leq
				\fc
				+3\max\pRbbb{
					\vass{\Theta_0^{(i)}},
					\frac{(\alpha\Cst+(1-\alpha)\cst) \fc}{(1-\alpha)\cst},
					\fc
					+
					\PRbbb{\sup_{n\in\N}\gamma_n}
					\prbbb{\frac{\pr{1-\alpha}\vass{f_0(\Theta_{-1}^{(i)})}}{\eps+\PR{\MOM_0}^{\nicefrac{1}{2}}}
						\\+\frac{\max\{1,\Cst\}(2+\alpha)\beta^{\nicefrac{1}{2}}}{\PR{\textstyle\inf_{n\in\N}\bscl_n}^{\nicefrac{1}{2}}\cst\pr{\beta^{\nicefrac{1}{2}}-\alpha}}}
				}
				.
			\end{multline}
		}
			\argument{
				\lref{it:upper:bound:theta};
				\lref{eq:setup:gen_grad:sing:dim:ineq2}}{\eqref{eq:a_priori_to_prove_cor:momentum2:tilde}}.
		\end{aproof}

		\subsection{A priori bounds for Adam for simple quadratic optimization problems}
		
		\begin{athm}{lemma}{momentum:representation}
			Let $\pars\in\N$, 
			let $\mom\colon\N_0\to\R^{\pars}$ be a function, and for every $n\in\N$ let $g_n\in\R^\pars$, $\beta_n\in\R$ satisfy
			\begin{equation}
				\label{setup:basic:geom:mom}
				\mom_n=\beta_n \mom_{n-1}+ g_n.
			\end{equation}
			Then it holds for all $n\in\N_0$ that
			$
			\mom_n
			=\PRb{\prod_{j=1}^{n}\beta_j}\mom_0
			+\sum_{k=1}^{n}\PRb{\prod_{j=k+1}^{n}\beta_j}g_k
			$.
		\end{athm}
		
		\begin{aproof}
			\argument{
				\eqref{setup:basic:geom:mom}}{that
				\begin{equation}
					\llabel{eq:base:momentum:geom}
					\begin{split}
						\mom_0
						\textstyle
						=\mom_0
						\qqandqq
						\mom_1
						=\textstyle\beta_1 \mom_0+ g_1
						.
					\end{split}
				\end{equation}
			}
			\argument{\eqref{setup:basic:geom:mom}}{that for all $n\in\N$ with $\mom_n=\PRb{\prod_{j=1}^n\beta_j}\mom_0+\sum_{k=1}^n\PRb{\prod_{j=k+1}^{n}\beta_j}g_k$ it holds that
				\begin{equation}
					\llabel{eq:base:momentum:geom:2}
					\begin{split}
						\mom_{n+1}
						=\textstyle\beta_{n+1} \mom_n+g_{n+1}
						&=\textstyle\beta_{n+1} \prb{\PRb{\prod_{j=1}^n\beta_j}\mom_0+\sum_{k=1}^n\PRb{\prod_{j=k+1}^n\beta_j}g_k}+g_{n+1}
						\\&=\textstyle\PRb{\prod_{j=1}^{n+1}\beta_j}\mom_0
						+\sum_{k=1}^{n}\PRb{\prod_{j=k+1}^{n+1}\beta_j}g_k+g_{n+1}
						\\&=\textstyle
						\PRb{\prod_{j=1}^{n+1}\beta_j}\mom_0
						+\sum_{k=1}^{n+1}\PRb{\prod_{j=k+1}^{n+1}\beta_j}g_k
						.
					\end{split}
				\end{equation}
			}
			\argument{\lref{eq:base:momentum:geom:2};
				\lref{eq:base:momentum:geom};
				induction}{that for all $n\in\N_0$ it holds that
				$
				\mom_n=\PRb{\prod_{j=1}^{n}\beta_j}\mom_0
				+\sum_{k=1}^{n}\PRb{\prod_{j=k+1}^{n}\beta_j}g_k
				$}.
			\finishproofthus
		\end{aproof}

		\begin{athm}{prop}{lem:ADAM:form}
			Let 
			$\lambda\colon\N\to[0,\infty)$,
			$\batch\colon\N\to\N$, and
			$\gamma\colon\N\to[0,\infty)$ satisfy
			\begin{equation}
				\llabel{eq:asymptotic:LR:nonstoch}
				\textstyle
				\inf_{n\in\N}\lambda_{n}>0
				\qqandqq
				\limsup_{n\to \infty}\PR{\gamma_n+\lambda_{n}} < \infty,
			\end{equation}
			let $\alpha\in[0,1)$, $\beta\in(\alpha^2,1)$, $\eps\in(0,\infty)$,
			let
			$\vartheta\colon\N\to\R$,
			$\mom\colon\N_0\to\R$,
			$\MOM\colon\N_0\to[0,\infty)$,
			$\MOMs\colon\N_0\to[0,\infty)$,
			and
			$ \Theta \colon\N_0\to \R$ 
			satisfy for all
			$n\in\N$ that
			\begin{equation}
				\llabel{eq:ADAM:1}
				\begin{split} 
					\mom_n&\textstyle=\alpha \mom_{n-1}+(1-\alpha)\lambda_{n}\prb{\Theta_{n-1}-\vartheta_n},
				\end{split}
			\end{equation}
			\begin{equation}
				\llabel{eq:ADAM:2}
				\begin{split}
					\MOM_{n}&\textstyle=\beta \MOM_{n-1}+(1-\beta)\PRb{\lambda_{n}\prb{\Theta_{n-1}-\vartheta_n}}^2,
				\end{split}
			\end{equation}
			\begin{equation}
				\llabel{eq:ADAM:3}
				\begin{split}
					\MOMs_0= \MOM_0,
					\qquad
					\MOMs_n\geq \MOM_n,
					\qqandqq
					\Theta_n
					= 
					\Theta_{ n - 1 }
					-
					\gamma_n\PRb{\eps+\PR{\MOMs_{n}}^{\nicefrac{1}{2}}}^{-1}\PRbbb{\frac{\mom_n}{1-\alpha^n}},
				\end{split}
			\end{equation}
			and let $\Cst,\fc\in\R$ satisfy
			\begin{equation}
				\llabel{eq:setup:bounds:nonstoch}
				\textstyle\sup_{n\in\N}\vass{\vartheta_{n}}\leq \fc
				\qqandqq
				\vass{\mom_0}
				\leq
				\Cst(1-\alpha)\pr{\vass{\Theta_0}+\fc}
				.
			\end{equation}
			Then
			\begin{multline}
				\llabel{eq:result:UB:pre:stoch:ADAM}
				\sup_{ n \in \N_0 }
				\vass{\Theta_n}
				\leq
				\fc
				+3\max\pRbbb{\vass{\Theta_0},
					\PRbbb{1+
						\frac{\alpha\PR{\sup_{n\in\N}\lambda_{n}} }{(1-\alpha)\PR{\inf_{n\in\N}\lambda_{n}}}}\fc,
					\fc+
					\PRbbb{\sup_{n\in\N}\gamma_n}
					\prbbb{\frac{\Cst(\vass{\Theta_0}+\fc)}{\eps}
						\\
						+
						\frac{(2+\alpha)\beta^{\nicefrac{1}{2}}\max\pR{1,\sup_{n\in\N}\lambda_{n}}}{(1-\alpha)\pr{1-\beta}^{\nicefrac{1}{2}}\pr{\beta^{\nicefrac{1}{2}}-\alpha}\PR{\inf_{n\in\N}\lambda_{n}}}
					}
				}.
			\end{multline}
			
		\end{athm}

		\begin{aproof}
			Throughout this proof let $g_n\colon\R\to\R$, $n\in\N_0$, satisfy for all $n\in\N$, $\theta\in\R$ that
			\begin{equation}
				\llabel{eq:setup:grad:fcts}
				g_0(\theta)=(1-\alpha)^{-1}\mom_0
				\qqandqq
				g_n(\theta)=\lambda_{n}\pr{\theta-\vartheta_n}
				.
			\end{equation}
			\startnewargseq
			\argument{
				\lref{eq:ADAM:1};
				\lref{eq:ADAM:3};
				\lref{eq:setup:grad:fcts};
				\cref{momentum:representation};}{that for all $n\in\N$ it holds that
				\begin{equation}
					\llabel{eq:Theta:ADAM:for:stable}
					\begin{split}
						\Theta_n
						&= 
						\Theta_{ n - 1 }
						-
						\frac{\gamma_n\prb{\prb{\frac{\alpha^n}{1-\alpha^n}}\mom_0+\prb{\frac{1-\alpha}{1-\alpha^n}}\sum_{k=1}^n \alpha^{n-k}\lambda_{k}\pr{\Theta_{k-1}-\vartheta_{k}}}}{\eps+\PR{\MOMs_n}^{\nicefrac{1}{2}}}
						\\&= 
						\Theta_{ n - 1 }
						-
						\frac{\prb{\frac{\gamma_n}{1-\alpha^n}}\prb{\alpha^n\mom_0+\PRb{\sum_{k=1}^n \alpha^{n-k}(1-\alpha)g_{k}\pr{\Theta_{k-1}}}}}{\eps+\PR{\MOMs_n}^{\nicefrac{1}{2}}}
						\\&= 
						\Theta_{ n - 1 }
						-
						\frac{\prb{\frac{\gamma_n}{1-\alpha^n}}\PRb{\sum_{k=0}^n \alpha^{n-k}(1-\alpha)g_{k}\pr{\Theta_{k-1}}}}{\eps+\PR{\MOMs_n}^{\nicefrac{1}{2}}}
						.
					\end{split}
				\end{equation}
			}
			\argument{\lref{eq:asymptotic:LR:nonstoch};\lref{eq:setup:bounds:nonstoch};}{that for all
				$\theta\in\R$,
				$n\in\N$ it holds that
				\begin{equation}
					\llabel{eq:gradient:boundaries}
					\begin{split}
						&	\textstyle
						\pr{\theta-\fc}
						\prb{\PR{\inf_{k\in\N}\lambda_{k}}+\prb{\PR{\sup_{k\in\N}\lambda_{k}}-\PR{\inf_{k\in\N}\lambda_{k}}}\indicator{(-\infty,\fc]}(\theta)}
						\\&\textstyle\leq
						\lambda_{n}\pr{\theta-\fc}
						\textstyle\leq
						\lambda_{n}\pr{\theta-\vartheta_{n}}
						\textstyle\leq
						\lambda_{n}\pr{\theta+\fc}
						\\&\textstyle\leq
						\pr{\theta+\fc}
						\prb{\PR{\inf_{k\in\N}\lambda_{k}}+\prb{\PR{\sup_{k\in\N}\lambda_{k}}-\PR{\inf_{k\in\N}\lambda_{k}}}\indicator{[-\fc,\infty)}(\theta)}
						.
					\end{split}
				\end{equation}
			}
			\argument{
				\lref{eq:gradient:boundaries};
				\lref{eq:setup:grad:fcts};
			}{that for all $n\in\N$, $\theta\in\R$ it holds that
				\begin{equation}
					\llabel{eq:gradient:bounds}
					\begin{split}
						\textstyle
						&\textstyle\pr{\theta-\fc}
						\prb{\PR{\inf_{k\in\N}\lambda_{k}}+\prb{\PR{\sup_{k\in\N}\lambda_{k}}-\PR{\inf_{k\in\N}\lambda_{k}}}\indicator{(-\infty,\fc]}(\theta)}
						\\&\textstyle\leq
						g_n(\theta)
						\leq
						\pr{\theta+\fc}
						\prb{\PR{\inf_{k\in\N}\lambda_{k}}+\prb{\PR{\sup_{k\in\N}\lambda_{k}}-\PR{\inf_{k\in\N}\lambda_{k}}}\indicator{[-\fc,\infty)}(\theta)}
						.
					\end{split}
				\end{equation}
			}
			\argument{
				\lref{eq:ADAM:2};
				\lref{eq:ADAM:3};
				\lref{eq:setup:grad:fcts};
				\cref{momentum:representation};
				the fact that $\beta\in(0,1)$;}{that for all $n\in\N$ it holds that
				\begin{equation}
					\llabel{eq:MOM:LB}
					\begin{split}
						\textstyle
						\MOMs_n\geq\MOM_n
						&=\textstyle\beta^n\MOM_{0}+\prb{1-\beta}\sum_{k=1}^n \beta^{n-k}\PR{\lambda_{k}\pr{\Theta_{k-1}-\vartheta_k}}^2
						\\&=\textstyle\beta^n\MOMs_{0}+\prb{1-\beta}\sum_{k=1}^n \beta^{n-k}\PR{\lambda_{k}\pr{\Theta_{k-1}-\vartheta_k}}^2
						.
					\end{split}
				\end{equation}
			}
			\argument{
				\lref{eq:setup:bounds:nonstoch};
			}{that
				\begin{equation}
					\llabel{eq:starting:cond}
					g_0\prb{\Theta_{\max\{-1,0\}}}
					=
					(1-\alpha)^{-1}\vass{\mom_{0}}
					\leq
					\Cst(\vass{\Theta_0}+\fc)
					.
				\end{equation}
			}
			\argument{
				\lref{eq:starting:cond};
				\lref{eq:setup:bounds:nonstoch};
				\lref{eq:Theta:ADAM:for:stable};
				\lref{eq:gradient:bounds};
				\lref{eq:MOM:LB};
				\cref{cor:a_priori_bound_gen:momentum} (applied with
				$\eps\curvearrowleft\eps$,
				$\cst \curvearrowleft \inf_{k\in\N}\lambda_{k}$,
				$\Cst \curvearrowleft \sup_{k\in\N}\lambda_{k}$,
				$\alpha \curvearrowleft \alpha$,
				$\beta \curvearrowleft \beta$,
				$\fc \curvearrowleft \fc$,
				$(g_n)_{n\in\N_0}\curvearrowleft (g_n)_{n\in\N_0}
				$,
				$\bscl\curvearrowleft \prb{\N\ni n\mapsto 1-\beta\in (0,\infty)}$, 
				$\gamma \curvearrowleft \prb{\N\ni n\mapsto \frac{\gamma_n}{1-\alpha^n}\in[0,\infty)}$,
				$\MOM \curvearrowleft 
				\MOMs
					$,
				$\Theta\curvearrowleft\pr{\Z\ni n\mapsto \Theta_{\max\{n,0\}}\in\R}$
				in the notation of \cref{cor:a_priori_bound_gen:momentum})}{that
				\begin{multline}
					\llabel{it:upper:bound:theta:ADAM}
					\sup_{ n \in \N_0 }
					\vass{\Theta_n}
					\leq
					\fc
					+3\max\pRbbb{\vass{\Theta_0},
						\PRbbb{1+
							\frac{\alpha\PR{\sup_{n\in\N}\lambda_{n}} }{(1-\alpha)\PR{\inf_{n\in\N}\lambda_{n}}}}\fc,
						\fc+
						\PRbbb{\sup_{n\in\N}\frac{\gamma_n}{1-\alpha^n}}
						\prbbb{\frac{(1-\alpha)\vass{g_0(\Theta_{\max\{-1,0\}})}}{\eps+\PR{\MOMs_{0}}^{\nicefrac{1}{2}}}
							\\+
							\frac{\max\pR{1,\sup_{n\in\N}\lambda_{n}}(2+\alpha)\beta^{\nicefrac{1}{2}}}{\pr{1-\beta}^{\nicefrac{1}{2}}\PR{\inf_{n\in\N}\lambda_{n}}\pr{\beta^{\nicefrac{1}{2}}-\alpha}}
						}
					}.
				\end{multline}
			}
			\argument{\lref{eq:starting:cond};}{
				\begin{equation}
					\begin{split}
						&\PRbbb{\sup_{n\in\N}\frac{\gamma_n}{1-\alpha^n}}
						\prbbb{\frac{(1-\alpha)\vass{g_0(\Theta_{\max\{-1,0\}})}}{\eps+\PR{\MOMs_{0}}^{\nicefrac{1}{2}}}
							+
							\frac{\max\pR{1,\sup_{n\in\N}\lambda_{n}}(2+\alpha)\beta^{\nicefrac{1}{2}}}{\pr{1-\beta}^{\nicefrac{1}{2}}\PR{\inf_{n\in\N}\lambda_{n}}\pr{\beta^{\nicefrac{1}{2}}-\alpha}}
						}
						\\&\leq\PRbbb{\sup_{n\in\N}\frac{\gamma_n}{1-\alpha}}
						\prbbb{\frac{(1-\alpha)\Cst(\vass{\Theta_0}+\fc)}{\eps}
							+
							\frac{\max\pR{1,\sup_{n\in\N}\lambda_{n}}(2+\alpha)\beta^{\nicefrac{1}{2}}}{\pr{1-\beta}^{\nicefrac{1}{2}}\PR{\inf_{n\in\N}\lambda_{n}}\pr{\beta^{\nicefrac{1}{2}}-\alpha}}
						}
						\\&=\PRbbb{\sup_{n\in\N}\gamma_n}
						\prbbb{\frac{\Cst(\vass{\Theta_0}+\fc)}{\eps}
							+
							\frac{(2+\alpha)\max\pR{1,\sup_{n\in\N}\lambda_{n}}}{(1-\alpha)\pr{1-\beta}^{\nicefrac{1}{2}}\pr{1-\alpha\beta^{-\nicefrac{1}{2}}}\PR{\inf_{n\in\N}\lambda_{n}}}
						}
						.
					\end{split}
				\end{equation}
			}
			\argument{
				\lref{it:upper:bound:theta:ADAM};
				\lref{eq:starting:cond};}{that
				\begin{multline}
					\llabel{it:upper:bound:ADAM}
					\sup_{ n \in \N_0 }
					\vass{\Theta_n}
					\leq
					\fc
					+3\max\pRbbb{\vass{\Theta_0},
						\PRbbb{1+
							\frac{\alpha\PR{\sup_{n\in\N}\lambda_{n}} }{(1-\alpha)\PR{\inf_{n\in\N}\lambda_{n}}}}\fc,
						\fc+
						\PRbbb{\sup_{n\in\N}\gamma_n}
						\prbbb{\frac{\Cst(\vass{\Theta_0}+\fc)}{\eps}
							\\
							+
							\frac{(2+\alpha)\max\pR{1,\sup_{n\in\N}\lambda_{n}}}{(1-\alpha)\pr{1-\beta}^{\nicefrac{1}{2}}\pr{1-\alpha\beta^{-\nicefrac{1}{2}}}\PR{\inf_{n\in\N}\lambda_{n}}}
						}
					}.
				\end{multline}
			}
			\argument{\lref{it:upper:bound:ADAM};}{\lref{eq:result:UB:pre:stoch:ADAM}}.
		\end{aproof}

		\begin{athm}{lemma}{cor:ADAM:form}
			Let $N\in\N$,
			let
			$\lambda\colon\N\to[0,\infty)$,
			$\batch\colon\N\to\N$, and
			$\gamma\colon\N\to[0,\infty)$ satisfy
			\begin{equation}
				\llabel{eq:asymptotic:LR:nonstoch}
				\textstyle
			\inf_{n\in\N}\lambda_{N+n}>0,
				\qquad
				\sup_{n\in\N}\gamma_n>0,
				\qqandqq
				\limsup_{n\to \infty}\PR{\gamma_n+\lambda_{n}} < \infty,
			\end{equation}
			let $\alpha\in[0,1)$, $\beta\in(\alpha^2,1)$, $\eps\in(0,\infty)$,
			let
			$\vartheta\colon\N\to\R$,
			$\mom\colon\N_0\to\R$,
			$\MOM\colon\N_0\to[0,\infty)$,
			and
			$ \Theta \colon\N_0\to \R$ 
			satisfy for all
			$n\in\N$ that
			\begin{equation}
				\llabel{eq:ADAM:1}
				\begin{split} 
					\mom_n&\textstyle=\alpha \mom_{n-1}+(1-\alpha)\lambda_{n}\prb{\Theta_{n-1}-\vartheta_n},
				\end{split}
			\end{equation}
			\begin{equation}
				\llabel{eq:ADAM:2}
				\begin{split}
					\MOM_{n}&\textstyle=\beta \MOM_{n-1}+(1-\beta)\PRb{\lambda_{n}\prb{\Theta_{n-1}-\vartheta_n}}^2,
				\end{split}
			\end{equation}
			\begin{equation}
				\llabel{eq:ADAM:3}
				\begin{split}
					\text{and}
					\qquad
					\Theta_n
					= 
					\Theta_{ n - 1 }
					-
					\gamma_n\PRbbb{\eps+\PRbbb{ \frac{\MOM_{n}}{1-\beta^n}}^{\nicefrac{1}{2}}}^{-1}\PRbbb{\frac{\mom_n}{1-\alpha^n}},
				\end{split}
			\end{equation}
			let $\Cst,\fc,D\in\R$ satisfy
			\begin{equation}
				\llabel{eq:setup:bounds:nonstoch}
				\textstyle\sup_{n\in\N}\vass{\vartheta_{n}}\leq \fc,
				\qquad
				\qquad
				\vass{\mom_0}
				\leq
				\Cst(1-\alpha)\pr{\vass{\Theta_0}+\fc}
				,
			\end{equation}
			\begin{equation}
				\qqandqq
				\qquad
				D\geq\prbbb{\frac{\PR{1+\sup_{n\in\N}\gamma_{n}}^{N+1}}{\sup_{n\in\N}\gamma_n}}
				(1-\alpha)^N
				\eps^{1-N}
				\PRbbb{1+\sup_{n\in\N}\lambda_n}^{N}
				\prbbb{\frac{\Cst+2}{\Cst}}.
			\end{equation}
			Then
			\begin{equation}
				\begin{split}
					\sup_{ n \in \N_0 }
					\vass{\Theta_{n}}
					\leq
					4D\prbbb{
						\frac{(3+\alpha)(1+\fc)\prb{1+\Cst(1+\eps^{-1})}^2\pr{1+\sup_{n\in\N}\lambda_{n}}}{\pr{1-\beta}^{\nicefrac{1}{2}}\pr{\beta^{\nicefrac{1}{2}}-\alpha}^2\PR{\inf_{n\in\N}\lambda_{N+n}}}}
					(\vass{\Theta_0}+1)
					.
				\end{split}
			\end{equation}
		\end{athm}
		
		\begin{aproof}
			Throughout this proof
			assume without loss of generality that $\max\{n\in\N\colon \lambda_n=0\}\in\N$ and 
			let $N\in\N$, $\fC\in\R$, $S,R\in(1,\infty)$ satisfy 
			\begin{equation}
				\llabel{eq:setup1}
				\fC=\frac{\pr{R^{N}-1}\pr{\Cst+2}}{(R-1)\Cst(1-\alpha)},
				\qquad
				\qquad
				R=\prbbb{1+\frac{\PR{\sup_{k\in\N}\gamma_k}}{\eps(1-\alpha)}}\textstyle\max\pRb{1,\sup_{k\in\N}\lambda_k}
				,
			\end{equation}
			\begin{equation}
				\qqandqq
				\qquad
				\llabel{eq:setup2}
				S=\frac{(2+\alpha)\beta^{\nicefrac{1}{2}}\max\pR{1,\sup_{n\in\N}\lambda_{N+n}}}{(1-\alpha)\pr{1-\beta}^{\nicefrac{1}{2}}\pr{\beta^{\nicefrac{1}{2}}-\alpha}\PR{\inf_{n\in\N}\lambda_{N+n}}}
				.
			\end{equation}
			\argument{
				\cref{momentum:representation};}{that for all $n\in\N_0$ it holds that
				\begin{equation}
					\llabel{eq:mom:bound}
					\begin{split}
						\vass{\mom_{n+1}}
						&\textstyle\leq\max\pRb{\vass{\mom_{n}},\lambda_{n+1}\pr{\vass{\Theta_{n}}+\fc}}
						\leq\textstyle\max\pRb{1,\sup_{k\in\N}\lambda_k}\prb{\max\pRb{\vass{\mom_{n}},\vass{\Theta_{n}}}+\fc}
					\end{split}
				\end{equation}
			}
			\argument{
				\lref{eq:mom:bound};
				\lref{eq:ADAM:3}}{that for all $n\in\N$ it holds that
				\begin{equation}
					\label{eq:UB:Thetas}
					\begin{split}
						\vass{\Theta_n}
						&\leq
						\vass[\bigg]{\Theta_{ n - 1 }
							-
							\gamma_n\PRbbb{\eps+\PRbbb{ \frac{\MOM_{n}}{1-\beta^n}}^{\nicefrac{1}{2}}}^{-1}\PRbbb{\frac{\mom_n}{1-\alpha^n}}}
						\\&\leq\vass{\Theta_{ n - 1 }}
						+
						\prbbb{\frac{\PR{\sup_{k\in\N}\gamma_k}}{\eps(1-\alpha)}}\vass{\mom_n}
						\\&\leq\vass{\Theta_{ n - 1 }}
						+
						\prbbb{\frac{\PR{\sup_{k\in\N}\gamma_k}}{\eps(1-\alpha)}}\textstyle\max\pRb{1,\sup_{k\in\N}\lambda_k}\prb{\max\pRb{\vass{\mom_{n-1}},\vass{\Theta_{n-1}}}+\fc}
						\\&\leq
						\prbbb{1+\frac{\PR{\sup_{k\in\N}\gamma_k}}{\eps(1-\alpha)}}\textstyle\max\pRb{1,\sup_{k\in\N}\lambda_k}\prb{\max\pRb{\vass{\mom_{n-1}},\vass{\Theta_{n-1}}}+\fc}
						\\&=
						R\prb{\max\pRb{\vass{\mom_{n-1}},\vass{\Theta_{n-1}}}+\fc}
						.
					\end{split}
				\end{equation}
			}
			\argument{\lref{eq:mom:bound};}{that for all $n\in\N$ it holds that
				\begin{equation}
					\llabel{eq:recursion:UB}
					\max\{\vass{\Theta_n},\vass{\mom_n}\}
					\leq 
					R\prb{\max\pRb{\vass{\mom_{n-1}},\vass{\Theta_{n-1}}}+\fc}
					.
			\end{equation}}
			\argument{
				\lref{eq:recursion:UB};
				induction
			}{
				that for all $n\in\N$ it holds that
				\begin{equation}
					\llabel{eq:late:ADAM}
					\begin{split}
						\textstyle
						\max\{\vass{\Theta_n},\vass{\mom_n}\}
						&\leq R\prb{\max\pRb{\vass{\mom_{n-1}},\vass{\Theta_{n-1}}}+\fc}
						\\&\leq R\prb{R\prb{\max\pRb{\vass{\mom_{n-2}},\vass{\Theta_{n-2}}}+\fc}+\fc}
						\\&= R^2\prb{\max\pRb{\vass{\mom_{n-2}},\vass{\Theta_{n-2}}}+(1+R)\fc}
						\\&\leq\dots
						\\&\leq\textstyle R^n\prb{\max\pRb{\vass{\mom_{0}},\vass{\Theta_{0}}}+\PR{\sum_{k=0}^{n-1}R^k}\fc}
						\\&\leq\textstyle R^n\prb{\max\pRb{\Cst(1-\alpha)\pr{\vass{\Theta_0}+\fc},\vass{\Theta_{0}}}+\PR{\sum_{k=0}^{n-1}R^k}\fc}
						\\&\leq\textstyle R^n\prb{\max\pRb{\Cst(1-\alpha),1}\pr{\vass{\Theta_0}+\fc}+\PR{\sum_{k=0}^{n-1}R^k}\fc}
						\\&\leq\textstyle \PR{\sum_{k=0}^{n}R^k}\prb{\max\pRb{\Cst(1-\alpha),1}\pr{\vass{\Theta_0}+\fc}+\fc}
						\\&\leq\textstyle \PR{\sum_{k=0}^{n}R^k}\max\pRb{\Cst(1-\alpha)+1,2}\pr{\vass{\Theta_0}+\fc}
						.
					\end{split}
				\end{equation}
			}
			\argument{\lref{eq:late:ADAM};}{
				\begin{equation}
					\llabel{eq:shifted:starting:cond}
					\begin{split}
						\frac{\vass{\mom_N}}{\Cst(1-\alpha)}
						&\leq
						\frac{\PR{\sum_{k=0}^{N}R^k}\max\pRb{\Cst(1-\alpha)+1,2}\pr{\vass{\Theta_0}+\fc}}{\Cst(1-\alpha)}
						=\fC\pr{\vass{\Theta_0}+\fc}
						\leq
						\vass{\Theta_N}+\fC\pr{\vass{\Theta_0}+\fc}
						.
					\end{split}
				\end{equation}
			}
			\argument{\lref{eq:late:ADAM};}{that for all $n\in\N_0\cap[0,N)$ it holds that
				\begin{equation}
					\llabel{eq:UB:first:finitely:many}
					\begin{split}
						\vass{\Theta_n}
						&\textstyle\leq \PR{\sum_{k=0}^{n}R^k}\max\pRb{\Cst(1-\alpha)+1,2}\pr{\vass{\Theta_0}+\fc}
						\\&\textstyle\leq \PR{\sum_{k=0}^{N-1}R^k}(\Cst+2)\pr{\vass{\Theta_0}+\fc}
						\\&=
						\frac{ \Cst(1-\alpha)\pr{R^{N}-1}(\Cst+2)\pr{\vass{\Theta_0}+\fc}}{ \Cst(1-\alpha)\pr{R-1}}
						\textstyle= \Cst(1-\alpha)\fC\pr{\vass{\Theta_0}+\fc}
						.
					\end{split}
			\end{equation}}
		\argument{\lref{eq:UB:first:finitely:many};}{that
		\begin{equation}
			\llabel{eq:UB:first:finitely:many:2}
			\begin{split}
			&\textstyle\max\pRb{\vass{\Theta_0},\vass{\Theta_1},\dots,\vass{\Theta_{N-1}}}
			\\&\leq
			\Cst(1-\alpha)\fC\pr{\vass{\Theta_0}+\fc}
			\\&\leq\textstyle 4\PR{\sup_{n\in\N}\max\{1,\gamma_{n}\}}\prb{\max\pRb{1,\tfrac{\Cst}{\eps}}+S}(1+\Cst)\fC(\vass{\Theta_0}+\max\{1,\fc\})
			.
			\end{split}
			\end{equation}
		}
			\argument{
				\lref{eq:UB:first:finitely:many};
				\lref{eq:shifted:starting:cond};
				\cref{lem:ADAM:form} (applied with
				$\lambda \curvearrowleft \pr{\N\ni n\mapsto \lambda_{N+n}\in[0,\infty)}$,
				$\batch \curvearrowleft \pr{\N\ni n\mapsto \batch_{N+n}\in\N}$,
				$\gamma \curvearrowleft \pr{\N\ni n\mapsto \gamma_{N+n}\in[0,\infty)}$,
				$\alpha \curvearrowleft \alpha$,
				$\beta \curvearrowleft \beta$,
				$\eps \curvearrowleft \eps$,
				$\vartheta \curvearrowleft \pr{\N\ni n\mapsto \vartheta_{N+n}\in\R}$,
				$\mom \curvearrowleft \pr{\N_0\ni n\mapsto \mom_{N+n}\in\R}$,
				$\MOM \curvearrowleft \pr{\N_0\ni n\mapsto \MOM_{N+n}\in[0,\infty)}$,
				$\Theta \curvearrowleft \pr{\N_0\ni n\mapsto \Theta_{N+n}\in\R}$,
				$\Cst \curvearrowleft \Cst$,
				$\fc \curvearrowleft \fC\pr{\vass{\Theta_0}+\fc}$,
				);}{that
				\begin{multline}
					\llabel{eq:result:UB:pre:stoch:ADAM}
					\sup_{ n \in \N_0 }
					\vass{\Theta_{N+n}}
					\leq
					\fC
					+3\max\pRbbb{\vass{\Theta_N},
						\PRbbb{1+
							\frac{\alpha\PR{\sup_{n\in\N}\lambda_{N+n}} }{(1-\alpha)\PR{\inf_{n\in\N}\lambda_{N+n}}}}\fC\pr{\vass{\Theta_0}+\fc},
						\fC\pr{\vass{\Theta_0}+\fc}
						\\	+
						\PRbbb{\sup_{n\in\N}\gamma_{N+n}}
						\prbbb{\frac{\Cst(\vass{\Theta_N}+\fC\pr{\vass{\Theta_0}+\fc})}{\eps}
							+
							\frac{(2+\alpha)\beta^{\nicefrac{1}{2}}\max\pR{1,\sup_{n\in\N}\lambda_{N+n}}}{(1-\alpha)\pr{1-\beta}^{\nicefrac{1}{2}}\pr{\beta^{\nicefrac{1}{2}}-\alpha}\PR{\inf_{n\in\N}\lambda_{N+n}}}
						}
					}.
				\end{multline}
			}
			\argument{
				\lref{eq:setup2};
				\lref{eq:result:UB:pre:stoch:ADAM};
				\lref{eq:UB:first:finitely:many};}{that
				\begin{equation}
					\llabel{eq:UB:for:last:infty:members}
					\begin{split}
						&\textstyle\sup_{ n \in \N_0 }
						\vass{\Theta_{N+n}}
						\\&\leq
						\textstyle \fC+3\max\{1,1+S,\PR{\sup_{n\in\N}\gamma_{n}}\prb{\tfrac{\Cst}{\eps}+S}\}(\vass{\Theta_N}+\fC(\vass{\Theta_0}+\max\{1,\fc\}))
						\\&\leq
						\textstyle \fC+3\PR{\sup_{n\in\N}\max\{1,\gamma_{n}\}}\prb{\max\pRb{1,\tfrac{\Cst}{\eps}}+S}(\vass{\Theta_N}+\fC(\vass{\Theta_0}+\max\{1,\fc\}))
						\\&\leq
						\textstyle \fC+3\PR{\sup_{n\in\N}\max\{1,\gamma_{n}\}}\prb{\max\pRb{1,\tfrac{\Cst}{\eps}}+S}(1+\Cst(1-\alpha))\fC(\vass{\Theta_0}+\max\{1,\fc\})
						\\&\leq
						\textstyle 4\PR{\sup_{n\in\N}\max\{1,\gamma_{n}\}}\prb{\max\pRb{1,\tfrac{\Cst}{\eps}}+S}(1+\Cst)\fC(\vass{\Theta_0}+\max\{1,\fc\})
						.
					\end{split}
				\end{equation}
			}
			\argument{\lref{eq:setup2};
			the fact that $\fC\geq 0$ and $\beta\leq 1$;}{that
			\begin{equation}
				\begin{split}
					\llabel{eq:constants}
					&\prb{\max\pRb{1,\tfrac{\Cst}{\eps}}+S}(1+\Cst)
					(\vass{\Theta_0}+\max\{1,\fc\})
					\\&=
					\prbbb{\max\pRb{1,\tfrac{\Cst}{\eps}}+\frac{(2+\alpha)\beta^{\nicefrac{1}{2}}\max\pR{1,\sup_{n\in\N}\lambda_{N+n}}}{\pr{1-\beta}^{\nicefrac{1}{2}}\pr{\beta^{\nicefrac{1}{2}}-\alpha}^2\PR{\inf_{n\in\N}\lambda_{N+n}}}}(1+\Cst)
					(\vass{\Theta_0}+\max\{1,\fc\})
					\\&\leq
					\prbbb{
						\prb{1+\Cst(1+\eps^{-1})}^2
						+
						\frac{(2+\alpha)(1+\Cst)\max\pR{1,\sup_{n\in\N}\lambda_{N+n}}}{\pr{1-\beta}^{\nicefrac{1}{2}}\pr{\beta^{\nicefrac{1}{2}}-\alpha}^2\PR{\inf_{n\in\N}\lambda_{N+n}}}}
						(1+\fc)(\vass{\Theta_0}+1)
					\\&\leq
					\prbbb{
						\frac{(3+\alpha)\prb{1+\Cst(1+\eps^{-1})}^2\pr{1+\sup_{n\in\N}\lambda_{n}}}{\pr{1-\beta}^{\nicefrac{1}{2}}\pr{\beta^{\nicefrac{1}{2}}-\alpha}^2\PR{\inf_{n\in\N}\lambda_{N+n}}}}
						(1+\fc)(\vass{\Theta_0}+1)
					\\&=
					\prbbb{
						\frac{(3+\alpha)(1+\fc)\prb{1+\Cst(1+\eps^{-1})}^2\pr{1+\sup_{n\in\N}\lambda_{n}}}{\pr{1-\beta}^{\nicefrac{1}{2}}\pr{\beta^{\nicefrac{1}{2}}-\alpha}^2\PR{\inf_{n\in\N}\lambda_{N+n}}}}
					(\vass{\Theta_0}+1)
					.
				\end{split}
			\end{equation}
		}
		\argument{
			\lref{eq:setup1};
			the fact that $R\geq 1$;}{
		\begin{equation}
			\llabel{eq:D:constants}
			\begin{split}
				&\PRbbb{\sup_{n\in\N}\max\{1,\gamma_{n}\}}\fC
				\\&=
				\PRbbb{\sup_{n\in\N}\max\{1,\gamma_{n}\}}
				\prbbb{\frac{\pr{R^{N}-1}\pr{\Cst+2}}{(R-1)\Cst(1-\alpha)}}
				\\&\leq
				\PRbbb{\sup_{n\in\N}\max\{1,\gamma_{n}\}}
				\prbbb{\frac{R^{N}\pr{\Cst+2}}{\pr{R-\max\pR{1,\sup_{n\in\N}\lambda_n}}\Cst(1-\alpha)}}
				\\&=
				\PRbbb{\sup_{n\in\N}\max\{1,\gamma_{n}\}}
				\prbbb{\frac{\PRb{\prb{1+\PR{\sup_{n\in\N}\gamma_n}}\pr{\eps(1-\alpha)}^{-1}\textstyle\max\pRb{1,\sup_{n\in\N}\lambda_n}}^{N}\pr{\Cst+2}}{\pr{\PR{\sup_{n\in\N}\gamma_n}\pr{\eps(1-\alpha)}^{-1}}\Cst(1-\alpha)}}
				\\&\leq
				\prbbb{\frac{\PRb{1+\sup_{n\in\N}\gamma_{n}}^{N+1}}{\sup_{n\in\N}\gamma_n}}
				(1-\alpha)^N
				\eps^{1-N}
				\PRbbb{1+\sup_{n\in\N}\lambda_n}^{N}
				\prbbb{\frac{\Cst+2}{\Cst}}
				=D
				.
			\end{split}
			\end{equation}
		}
			\argument{
				\lref{eq:D:constants};
				\lref{eq:UB:first:finitely:many:2};
				\lref{eq:UB:for:last:infty:members};
				\lref{eq:constants};
				}{that
				\begin{equation}
					\begin{split}
						\textstyle\sup_{ n \in \N_0 }
						\vass{\Theta_{n}}
						&\textstyle=\max\pRb{\vass{\Theta_0},\vass{\Theta_1},\dots,\vass{\Theta_{N-1}},\sup_{ n \in \N_0 }\vass{\Theta_{N+n}}}
						\\&\leq
						\textstyle 4\PR{\sup_{n\in\N}\max\{1,\gamma_{n}\}}\prb{\max\pRb{1,\tfrac{\Cst}{\eps}}+S}(1+\Cst)\fC(\vass{\Theta_0}+\max\{1,\fc\})
						\\&\leq
						\textstyle 4D\prb{\max\pRb{1,\tfrac{\Cst}{\eps}}+S}(1+\Cst)(\vass{\Theta_0}+\max\{1,\fc\})
						\\&\leq
						 4D\prbbb{
							\frac{(3+\alpha)(1+\fc)\prb{1+\Cst(1+\eps^{-1})}^2\pr{1+\sup_{n\in\N}\lambda_{n}}}{\pr{1-\beta}^{\nicefrac{1}{2}}\pr{\beta^{\nicefrac{1}{2}}-\alpha}^2\PR{\inf_{n\in\N}\lambda_{N+n}}}}
						(\vass{\Theta_0}+1)
						.
					\end{split}
				\end{equation}
			}
		\end{aproof}
		
		\begin{athm}{theorem}{lem:ADAM:bounded}
			Let $d\in\N$, 
			let
			$\lambda=(\lambda^{(1)},\dots,\lambda^{(d)})\colon\N\to[0,\infty)^d$,
			$\batch\colon\N\to\N$, and
			$\gamma\colon\N\to[0,\infty)$ satisfy for all $i\in\{1,2,\dots,d\}$ that
			\begin{equation}
				\llabel{eq:asymptotic:LR}
				\textstyle
				\liminf_{n\to\infty}\lambda_{n}^{(i)}>0
				\qqandqq
				\limsup_{n\to \infty}\prb{\gamma_n+\lambda_{n}^{(i)}} < \infty,
			\end{equation}
			let $(\Omega,\cF,\P)$ be a probability space,
			let $\fc\in[0,\infty)$,
			for every $n,i,j\in\N$
			let $X_{n,j}^{(i)}\colon\Omega\to[-\fc,\fc]$ be a random variable,
			let $\alpha\in[0,1)$, $\beta\in(\alpha^2,1)$, $\eps\in(0,\infty)$,
			let
			$\mom=\prb{\mom^{(1)},\dots,\mom^{(d)}}\colon\N_0\times\Omega\to\R^d$,
			$\MOM=\prb{\MOM^{(1)},\dots,\MOM^{(d)}}\colon\N_0\times\Omega\to[0,\infty)^d$, and
			$ \Theta =\prb{\Theta^{(1)},\dots,\Theta^{(d)}}\colon\N_0\times\Omega \to \R^d$ 
			be stochastic processes which satisfy for all
			$n\in\N$, $i\in\{1,2,\dots,d\}$ that
			\begin{equation}
				\llabel{eq:ADAM:1}
				\begin{split} 
					\mom_{n}^{(i)}&\textstyle=\alpha \mom_{n-1}^{(i)}+(1-\alpha)\PRbb{\frac{\lambda_{n}^{(i)}}{\batch_n}\sum_{j=1}^{\batch_n}\prb{\Theta_{n-1}^{(i)}-X_{n,j}^{(i)}}},
				\end{split}
			\end{equation}
			\begin{equation}
				\llabel{eq:ADAM:2}
				\begin{split}
					\MOM_{n}^{(i)}&\textstyle=\beta \MOM_{n-1}^{(i)}+(1-\beta)\PRbb{\frac{\lambda_{n}^{(i)}}{\batch_n}\sum_{j=1}^{\batch_n}\prb{\Theta_{n-1}^{(i)}-X_{n,j}^{(i)}}}^2,
				\end{split}
			\end{equation}
			\begin{equation}
				\llabel{eq:ADAM:3}
				\begin{split}
					\Theta_n^{(i)}
					= 
					\Theta_{ n - 1 }^{(i)}
					-
					\gamma_n\PRbbb{\eps+\PRbbb{ \frac{\MOM_{n}^{(i)}}{1-\beta^n}}^{\nicefrac{1}{2}}}^{-1}\PRbbb{\frac{\mom_n^{(i)}}{1-\alpha^n}},
				\end{split}
			\end{equation}
			and
			$				\vass{\mom_0^{(i)}}
			\leq
			\fc\vass{\Theta_0^{(i)}}+\fc$.
			Then there exists $\fC\in\R$ such that $\sup_{n\in\N_0}\norm{\Theta_n}\leq \fC\norm{\Theta_0}+\fC$.
		\end{athm}
		
		\begin{aproof}
			Throughout this proof assume without loss of generality that for all $i\in\{1,2,\dots,d\}$ it holds that
			$\fc\geq 1$,
			$\vass{\mom_0^{(i)}}
			\leq
			\fc(1-\alpha)\pr{\vass{\Theta_0^{(i)}}+\fc}$,
			$\sup_{n\in\N}\gamma_n>0$, and let $N\in\N$, $D\in\R$ satisfy $\min_{i\in\{1,2,\dots,d\}}\PRb{\inf_{n\in\N}\lambda_{N+n}^{(i)}}>0$ and
							\begin{equation}
					D=\prbbb{\frac{\PR{1+\sup_{n\in\N}\gamma_{n}}^{N+1}}{\sup_{n\in\N}\gamma_n}}
					(1-\alpha)^N
					\eps^{1-N}
					\PRbbb{1+\max_{i\in\{1,2,\dots,d\}}\prbbb{\sup_{n\in\N}\lambda_n^{(i)}}}^{N}
					\prbbb{\frac{\fc+2}{\fc}}.
				\end{equation}
			\argument{\lref{eq:ADAM:1};}{that for all 
				$n\in\N$, 
				$i\in\{1,2,\dots,d\}$,
				$\omega\in\Omega$ it holds that
				\begin{equation}
					\llabel{eq:momentum:path}
					\begin{split}
						\mom_n^{(i)}(\omega)
						&\textstyle=\alpha \mom_{n-1}^{(i)}(\omega)+(1-\alpha)\PRb{\frac{\lambda_n^{(i)}}{\batch_n}\sum_{j=1}^{\batch_n}\prb{\Theta_{n-1}^{(i)}(\omega)-X_{n,j}^{(i)}(\omega)}}
						\\&=\textstyle
						\alpha\mom_{n-1}^{(i)}(\omega)+(1-\alpha)\lambda_{n}^{(i)}\prb{\Theta_{n-1}^{(i)}(\omega)-\PRb{\frac{1}{\batch_n}\sum_{j=1}^{\batch_n}X_{n,j}^{(i)}(\omega)}}
						.
					\end{split}
				\end{equation}
			}
			\argument{\lref{eq:ADAM:2};}{that for all 
				$n\in\N$, 
				$i\in\{1,2,\dots,d\}$,
				$\omega\in\Omega$
				 it holds that
				\begin{equation}
					\llabel{eq:momentum:path2}
					\begin{split}
						\MOM_n^{(i)}(\omega)
						&\textstyle=\beta \MOM_{n-1}^{(i)}(\omega)+(1-\beta)\PRb{\frac{\lambda_{n}^{(i)}}{\batch_n}\sum_{j=1}^{\batch_n}\prb{\Theta_{n-1}^{(i)}(\omega)-X_{n,j}^{(i)}(\omega)}}^2
						\\&=\textstyle
						\beta\MOM_{n-1}^{(i)}(\omega)+(1-\beta)\PRb{\lambda_{n}^{(i)}\prb{\Theta_{n-1}^{(i)}(\omega)-\PRb{\frac{1}{\batch_n}\sum_{j=1}^{\batch_n}X_{n,j}^{(i)}(\omega)}}}^2
						.
					\end{split}
				\end{equation}
			}
			\argument{
				the the fact that for all $n,i,j\in\N$ it holds that $\vass{X_{n,j}^{(i)}}\leq\fc$;}{that for all 
				$i\in\{1,2,\dots,d\}$, 
				$\omega\in\Omega$ it holds that
				\begin{equation}
					\llabel{eq:grad:on:data:bounded}
					\begin{split}
						\textstyle
						\sup_{n\in\N}
						\vass{\frac{1}{\batch_n}\sum_{j=1}^{\batch_n}X_{n,j}^{(i)}(\omega)}
						\leq\sup_{n\in\N}
						\PRb{\frac{1}{\batch_n}\sum_{j=1}^{\batch_n}\vass{X_{n,j}^{(i)}(\omega)}}
						\leq\sup_{n\in\N}
						\PRb{\frac{1}{\batch_n}\sum_{j=1}^{\batch_n}\fc}
						=
						\fc
						.
					\end{split}
				\end{equation}
			}
			\argument{
				\lref{eq:grad:on:data:bounded};
				\lref{eq:momentum:path};
				\lref{eq:momentum:path2};
				\cref{cor:ADAM:form} (applied for every $\omega\in\Omega$ with 
				$\lambda\curvearrowleft \lambda^{(i)}$,
				$\batch\curvearrowleft \batch$,
				$\gamma\curvearrowleft \gamma$,
				$\alpha\curvearrowleft \alpha$,
				$\beta\curvearrowleft \beta$,
				$\eps\curvearrowleft \eps$,
				$\vartheta\curvearrowleft \pr{\N\ni n\mapsto \frac{1}{\batch_n}\sum_{j=1}^{\batch_n}X_{n,j}^{(i)}(\omega)\in\R}$,
				$\mom\curvearrowleft \pr{\N_0\ni n\mapsto \mom_n^{(i)}(\omega)\in\R}$,
				$\MOM\curvearrowleft \pr{\N_0\ni n\mapsto \MOM_n^{(i)}(\omega)\in[0,\infty)}$,
				$\Theta\curvearrowleft \pr{\N_0\ni n\mapsto \Theta_n^{(i)}(\omega)\in\R}$,
				$\Cst\curvearrowleft \fc$,
				$\fc\curvearrowleft \fc$
				for $i\in\{1,2,\dots,d\}$
				in the notation of \cref{cor:ADAM:form});
			}{that for all 
			$i\in\{1,2,\dots,d\}$,
			$\omega\in\Omega$ it holds that
				\begin{equation}
					\llabel{eq:result:UB:pre:stoch:ADAM}
					\sup_{ n \in \N_0 }
					\vass{\Theta_n^{(i)}(\omega)}
					\leq
						4D\prbbb{
						\frac{(3+\alpha)(1+\fc)\prb{1+\fc(1+\eps^{-1})}^2\pr{1+\sup_{n\in\N}\lambda_{n}^{(i)}}}{\pr{1-\beta}^{\nicefrac{1}{2}}\pr{\beta^{\nicefrac{1}{2}}-\alpha}^2\PR{\inf_{n\in\N}\lambda_{N+n}^{(i)}}}}
					\prb{\vass{\Theta_0^{(i)}(\omega)}+1}
				\end{equation}
			}
			\argument{
				\lref{eq:result:UB:pre:stoch:ADAM};
				\lref{eq:asymptotic:LR};
				the fact that $\alpha^2<\beta<1$ and $\min_{i\in\{1,2,\dots,d\}}\PRb{\inf_{n\in\N}\lambda_{N+n}^{(i)}}>0$;}{that there exists $c\in\R$ which satisfies that
				\begin{equation}
					\llabel{eq:constants}
					\begin{split}
						&\textstyle
						\sup_{ n \in \N_0 }
						\norm{\Theta_n}
						\\&\leq
						\textstyle
						\sup_{ n \in \N_0 }
						\PRb{\max_{i\in\{1,2,\dots,d\}}\sqrt{d}\vass{\Theta_n^{(i)}}}
						\\&=
						\textstyle
						\max_{i\in\{1,2,\dots,d\}}\sqrt{d}
						\PRb{\sup_{ n \in \N_0 }
						\vass{\Theta_n^{(i)}}}
						\\&\leq
						\max_{i\in\{1,2,\dots,d\}}
						\sqrt{d}
						\PRbbb{4D\prbbb{
							\frac{(3+\alpha)(1+\fc)\prb{1+\fc(1+\eps^{-1})}^2\prb{1+\sup_{n\in\N}\lambda_{n}^{(i)}}}{\pr{1-\beta}^{\nicefrac{1}{2}}\pr{\beta^{\nicefrac{1}{2}}-\alpha}^2\PRb{\inf_{n\in\N}\lambda_{N+n}^{(i)}}}}
						\prb{\vass{\Theta_0}^{(i)}+1}}
						\\&\leq
						4D\sqrt{d}
						\PRbbb{\max_{i\in\{1,2,\dots,d\}}
							\prbbb{
								\frac{(3+\alpha)(1+\fc)\prb{1+\fc(1+\eps^{-1})}^2\prb{1+\sup_{n\in\N}\lambda_{n}^{(i)}}}{\pr{1-\beta}^{\nicefrac{1}{2}}\pr{\beta^{\nicefrac{1}{2}}-\alpha}^2\PRb{\inf_{n\in\N}\lambda_{N+n}^{(i)}}}}
							}\prb{\norm{\Theta_0}+1}
						\\&=
						4D\sqrt{d}c
						\pr{\norm{\Theta_0}+1}
						.
					\end{split}
			\end{equation}}
			\argument{\lref{eq:constants};}{that there exists $\fC\in\R$ such that $\sup_{n\in\N_0}\norm{\Theta_n}\leq \fC\pr{\norm{\Theta_0}+1}$}.
		\end{aproof}

		\section{A priori bounds for the momentum optimizer}
		\label{section:apriori:momentum}
		
		In this section we specify in \cref{prop:mom-stable:diag} below the set of tupels of learning rates and eigenvalues of the Hessian for which the momentum method does not explode but stays bounded when applied to the simple class of quadratic \OPs\ in \cref{eq:min:prob} in \cref{subsection:1:1} above.
		This will allow us to explicitly depict the stability region (cf.\ \cref{def:stab:region} in \cref{subsection:1:1}) of the momentum optimizer, that is, to establish \cref{it:MOM:perm} in \cref{intro:thm:1.1} above (cf.~\cref{subsection:5.4} below).
		Our proof of \cref{prop:mom-stable:diag} employs the common concept of the spectral radius which we recall here in \cref{def:sr} below.
		The convergence part within the stability region in \cref{prop:mom-stable:diag} (cf.~\cref{prop:mom-stable}) is already well-known in the literature \cite{TrungMomLR}.

		\subsection{Asymptotic analysis for coupled systems}
		
				\begin{definition}[Spectral radius]
			\label{def:sr}
			We denote by $\sr\colon\pr{\cup_{n\in\N} \C^{ n \times n }}\to\R$ the function which satisfies for all $n\in\N$, $A \in \C^{ n \times n }$ that
			\begin{equation}
				\sr ( A ) = 
				\max\{ r \in \R \colon \pr{ \exists\, \mu\in \{ \lambda \in \C \colon \vass{\lambda} = r \}, v \in \C^n \backslash \{0\} \colon A v = \mu v } \}
				.
			\end{equation}
		\end{definition}
		
		\cfclear
		
		\begin{athm}{lemma}{lem:mat:seq:EV:small}
			Let $d\in\N$ and
			let 
			$\m\colon\N_0\to\C^d$, 
			$\Theta\colon\N_0\to\C^d$, and 
			$A\colon\N_0\to\C^{2d\times 2d}$ satisfy for all $n\in\N$ that
			\begin{equation}
				\llabel{eq:n_th:vec}
				\begin{pmatrix}
					\m_n \\ \Theta_n 
				\end{pmatrix}
				=
				A_n
				\begin{pmatrix}
					\m_{n-1} \\ \Theta_{n-1}
				\end{pmatrix},
				\quad
				\textstyle
				\limsup_{k\to\infty}\prb{\sup_{v\in\C^{d}\backslash\{0\}}\frac{\norm{A_kv-A_0v}}{\norm{v}}}=0,
				\qandq
				\sr(A_0)<1
			\end{equation}
			\cfload.
			Then $\limsup_{n\to\infty}\norm{\Theta_n}=0$.
		\end{athm}
		
		\setcounter{footnote}{1}
		\begin{aproof}
			Throughout this proof let $\mnorm{\cdot}\colon \C^{ d \times d}\to\R$ satisfy for all
			$B\in \C^{ d \times d }$ 
			that 
			\begin{equation}
				\llabel{eq:mnorm}
				\textstyle
				\mnorm{B} = 
				\sup_{v\in\C^{d}\backslash\{0\}}\prb{\frac{\norm{Bv}}{\norm{v}}}.
			\end{equation}
			\startnewargseq
			\argument{\kern-0.9ex, \eg, \cite[Lemma~2.4]{dereich2025sharphigherorderconvergence}
				(applied with 
				$d\curvearrowleft 2d$,
				$M\curvearrowleft A_0$,
				$(A_n)_{n\in\N}\curvearrowleft(A_n)_{n\in\N}$,
				$\delta\curvearrowleft\frac{1-\sr(A_0)}{2}$
				in the notation of \cite[Lemma~2.4]{dereich2025sharphigherorderconvergence});
				\lref{eq:n_th:vec};
				\lref{eq:mnorm}}{that there exist $c\in(0,\infty)$, $N\in\N$ which satisfy for all $m,n\in\N$ with $m>n\geq N$ that
				\begin{equation}
					\llabel{eq:Cauchy:norm}
					\mnorm{\sprod_{i=n+1}^m A_i}\leq c\prb{\sr(A_0)+\tfrac{1-\sr(A_0)}{2}}^{m-n}
					=
					c\prb{\tfrac{1+\sr(A_0)}{2}}^{m-n}
					.
			\end{equation}}
			\argument{
				\lref{eq:Cauchy:norm};
				\lref{eq:n_th:vec};
				\lref{eq:mnorm}}{that
				\begin{equation}
					\llabel{eq:Cauchy:norm2}
					\begin{split}
						\textstyle\limsup_{n \to \infty}\mnorm{\sprod_{i=1}^n A_i}
						&\textstyle\leq
						\limsup_{n \to \infty}\prb{\mnorm{\sprod_{i=1}^N A_i}\mnorm{\sprod_{i=N+1}^{\max\{N+1,n\}} A_i}}
						\\&\textstyle=
						\mnorm{\sprod_{i=1}^N A_i}
						\limsup_{n \to \infty}
						\mnorm{\sprod_{i=N+1}^{\max\{N+1,n\}} A_i}
						\\&\textstyle\leq
						\mnorm{\sprod_{i=1}^N A_i}
						\limsup_{n \to \infty}
						c\prb{\tfrac{1+\sr(A_0)}{2}}^{\max\{N+1,n\}-N}
						=0
						.
					\end{split}
			\end{equation}}
			\argument{
				\lref{eq:Cauchy:norm2};
				\lref{eq:mnorm};}{that
				\begin{equation}
					\begin{split}
						\limsup_{n \to \infty}\norm{\Theta_n}
						\leq
						\limsup_{n \to \infty}
						\norm[\bigg]{\begin{pmatrix}
								\m_n \\ \Theta_n 
						\end{pmatrix}}
						&=
						\limsup_{n \to \infty}
						\norm[\bigg]{A_n
							\begin{pmatrix}
								\m_{n-1} \\ \Theta_{n-1}
						\end{pmatrix}}
						\\&=
						\limsup_{n \to \infty}
						\norm[\bigg]{A_nA_{n-1}\cdots A_1
							\begin{pmatrix}
								\m_{0} \\ \Theta_{0}
						\end{pmatrix}}
						\\&\leq			
						\limsup_{n \to \infty}
						\mnorm{A_nA_{n-1}\cdots A_1}
						\norm[\bigg]{
							\begin{pmatrix}
								\m_{0} \\ \Theta_{0}
						\end{pmatrix}}
						\\&=
						\norm[\bigg]{
							\begin{pmatrix}
								\m_{0} \\ \Theta_{0}
						\end{pmatrix}}
						\limsup_{n \to \infty}\mnorm{\sprod_{i=1}^n A_i}
						=0
						.
					\end{split}
			\end{equation}}
		\end{aproof}

		\cfclear
		\begin{athm}{lemma}{lem:mat:seq:EV}
			Let 
			$v_1=(v_{1,1},v_{1,2})$, $v_2=(v_{2,1},v_{2,2})\in\C^2$, 
			$A\in\C^{2\times 2}$,
			$\mu_1,\mu_2\in\C$
			satisfy for all $i\in\{1,2\}$ that			
			\begin{equation}
				\llabel{eq:EVandEV:setup}
				\indicator{[0,1]}\pr{\sr\pr{A}}+\vass{v_{1,2}}>0,
				\qquad
				\sr\pr{A}=\max\{\vass{\mu_1},\vass{\mu_2}\},
				\qqandqq
				A v_i = \mu_i v_i
			\end{equation}
			and let 
			$\m\colon\N_0\to\C$ and
			$\Theta\colon\N_0\to\C$
			satisfy for all $n\in\N_0$ that
			\begin{equation}
				\llabel{eq:n_th:vec}
				\max\{1,\vass{\mu_1}\}>\vass{\mu_2}
				\qqandqq
				\begin{pmatrix}
					\m_n \\ \Theta_n 
				\end{pmatrix}
				=
				A^n
				\begin{pmatrix}
					\m_0\\ \Theta_{0}
				\end{pmatrix}
				=\textstyle
				A^n(v_1+v_2)
			\end{equation}
			\cfload.
			Then
			\begin{equation}
				\llabel{eq:result}
				\liminf_{n \to \infty} \vass{\Theta_n} = \limsup_{n \to \infty} \vass{\Theta_n } 
				=
				\begin{cases}
					0 & \colon \sr\pr{A}<1\\
					\vass{v_{1,2}} & \colon \sr\pr{A}=1\\
					\infty & \colon \sr\pr{A}>1.
				\end{cases}
			\end{equation}
		\end{athm}

		\begin{aproof}
			In our proof of \lref{eq:result} we distinguish between the cases $\sr\pr{A}<1$, $\sr\pr{A}=1$, and $\sr\pr{A}>1$.
			We first \prove \lref{eq:result} in the case 
			\begin{equation}
				\llabel{eq:small:sr}
				\sr\pr{A}<1.
			\end{equation}
			\argument{\lref{eq:small:sr};\cref{lem:mat:seq:EV:small} (applied with 
				$d\curvearrowleft 1$,
				$\m\curvearrowleft \m$,
				$\Theta\curvearrowleft\Theta$,
				$A\curvearrowleft \pr{\N_0\ni n\mapsto A\in\C^{2\times 2}}$) in the notation of \cref{lem:mat:seq:EV:small}}{that
				\begin{equation}
					\llabel{eq:small:sr:result}
					\limsup_{n \to \infty}\vass{\Theta_n}=0.
			\end{equation}}
			\argument{\lref{eq:small:sr:result}}{\lref{eq:result} in the case  $\sr\pr{A}<1$}.
			Next we \prove \lref{eq:result} in the case 
			\begin{equation}
				\llabel{eq:middle:sr}
				\sr\pr{A}=1.
			\end{equation}
			\startnewargseq
			\argument{
				\lref{eq:n_th:vec};
				\lref{eq:middle:sr};
			}{that
				\begin{equation}
					\llabel{eq:sr:to:evs}
					1=\sr\pr{A}=\max\{\vass{\mu_1},\vass{\mu_2}\}=\vass{\mu_1}>\vass{\mu_2}.
			\end{equation}}
			\argument{\lref{eq:sr:to:evs};\lref{eq:n_th:vec}}{that for all $z\in\{-1,1\}$ it holds that
				\begin{equation}
					\llabel{eq:middle:sr:asymp:val}
					\begin{split}
						z\PRbbb{\limsup_{n \to \infty}z\norm[\bigg]{\begin{pmatrix}
									\m_n \\ \Theta_n
								\end{pmatrix}
								-(\mu_1)^nv_1}
						}
						&=z\PRbbb{\limsup_{n \to \infty}z\norm{A^n(v_1+v_2)
								-(\mu_1)^nv_1}
						}
						\\&=z\PRbbb{\limsup_{n \to \infty}z\norm{(\mu_1)^nv_1+(\mu_2)^nv_2
								-(\mu_1)^nv_1}
						}
						\\&=z\PRbbb{\limsup_{n \to \infty}z\norm{(\mu_2)^nv_2}
						}
						= z\PRbbb{\limsup_{n \to \infty}z\vass{\mu_2}^n
						}\norm{v_2}=0
						.
					\end{split}
			\end{equation}}
			\argument{\lref{eq:sr:to:evs};\lref{eq:middle:sr:asymp:val}}{that for all $z\in\{-1,1\}$ it holds that
				\begin{equation}
					\llabel{eq:middle:sr:asymp:val:res}
					\textstyle
					z\PRb{ \limsup_{n \to \infty} z\vass{\Theta_n }}
					= z\PRb{ \limsup_{n \to \infty} z\vass{(\mu_1)^nv_{1,2}}}
					= z\PRb{ \limsup_{n \to \infty} z\vass{v_{1,2}}}
					=\vass{v_{1,2}}
					.
			\end{equation}}
			\argument{\lref{eq:middle:sr:asymp:val:res}}{\lref{eq:result} in the case $\sr\pr{A}=1$}.
			Next we \prove \lref{eq:result} in the case 
			\begin{equation}
				\llabel{eq:big:sr}
				\sr\pr{A}>1.
			\end{equation}
			\startnewargseq
			\argument{
				\lref{eq:EVandEV:setup};
				\lref{eq:n_th:vec};
				\lref{eq:big:sr}}{that
				\begin{equation}
					\llabel{eq:big:sr:asymp:val:res}
					\begin{split}
						\textstyle
						\liminf_{n \to \infty}\vass{
							\Theta_n
						}
						&\textstyle=
						\liminf_{n \to \infty}\vass{
							(\mu_1)^nv_{1,2}+(\mu_2)^nv_{2,2}
						}
						\\&\textstyle\geq
						\liminf_{n \to \infty}
						\PRb{
							\vass{(\mu_1)^nv_{1,2}
							}
							-
							\vass{
								(\mu_2)^nv_{2,2}
							}
						}
						\\&\textstyle=
						\liminf_{n \to \infty}
						\prb{\vass{\mu_1}^n\PRb{\vass{
									v_{1,2}
								}-\vass{\mu_2(\mu_1)^{-1}}^n\vass{v_{2,2}
								}
						}}
						\\&\textstyle=
						\prb{\lim_{n \to \infty}
							\vass{\mu_1}^n}
						\prb{\lim_{n \to \infty}\PRb{\vass{
									v_{1,2}
								}-\vass{\mu_2(\mu_1)^{-1}}^n\vass{	v_{1,2}
								}
						}}
						\\&\textstyle=
						\prb{\lim_{n \to \infty}
							\vass{\mu_1}^n}
						\vass{
							v_{1,2}
						}
						=\infty
						.
					\end{split}
			\end{equation}}
			\argument{\lref{eq:big:sr:asymp:val:res}}{\lref{eq:result} in the case $\sr\pr{A}>1$}.
		\end{aproof}
		
		 \subsection{One step analysis for the momentum optimizer}
		
		\cfclear
		\begin{athm}{lemma}{lem:EV}
			Let $\gamma, \cK \in (0, \infty)$,
			$\alpha \in (0, 1 )$, $A\in\R^{2\times 2}$ satisfy
			\begin{equation}
				\llabel{eq:matrix}
				A = \begin{pmatrix}
					\alpha & 2(1 - \alpha) \cK \\
					- \gamma \alpha & 1 - 2  ( 1 - \alpha ) \gamma\cK
				\end{pmatrix}
			\end{equation}
			and let $\mu_-,\mu_+\in\C$ satisfy
			\begin{equation}
				\llabel{eq:EV}
				\mu_\pm =	\frac{1 + \alpha - 2  ( 1 - \alpha ) \gamma\cK  }{2}
				\pm \sqrt{ \rbr*{ \frac{1 + \alpha - 2 ( 1 - \alpha )\gamma \cK }{2} } ^2 - \alpha }
				.
			\end{equation}
			Then
			\begin{enumerate}[label=(\roman*)]
				\item\llabel{it:evs} it holds that $\cu{ z \in \C \colon \br{ \exists \, \nu \in \C^2 \backslash \cu{0} \colon A \nu = z \nu } }
				=
				\{\mu_-,\mu_+\}$
				and
				\item\llabel{it:sr} it holds that
				\begin{equation}
					\llabel{eq:cond:for:result}
					\max\{1,\vass{\mu_-}\}>\vass{\mu_+}
					\qqandqq
					\sr(A)\in
					\begin{cases}
						[0,1) &\colon \gamma \cK<\frac{1 +\alpha}{1-\alpha}\\
						\{1\} &\colon \gamma \cK=\frac{1 +\alpha}{1-\alpha}\\
						(1,\infty) &\colon \gamma \cK>\frac{1 +\alpha}{1-\alpha}\\
					\end{cases}
				\end{equation}
			\end{enumerate}
			\cfout.
		\end{athm}

		\begin{aproof}
			\argument{
				\lref{eq:matrix};}{that
				\begin{equation}
					\llabel{eq:vanishing:det:0}
					\cu{ z \in \C \colon \br{ \exists \, \nu \in \C^2 \backslash \cu{0} \colon A \nu = z \nu } }
					=
					\pRbbb{z\in\C\colon\operatorname{det}
						\begin{pmatrix}
							\alpha-z& 2(1 - \alpha) \cK \\
							- \gamma \alpha & 1 - 2  ( 1 - \alpha ) \gamma\cK-z
						\end{pmatrix}
						=0}.
				\end{equation}
			}
			\argument{\lref{eq:vanishing:det:0}}{
				that 				\begin{equation}
					\llabel{eq:vanishing:det:1}
					\begin{split}
						&\cu{ z \in \C \colon \br{ \exists \, \nu \in \C^2 \backslash \cu{0} \colon A \nu = z \nu } }
						\\&=\pR{z\in\C\colon
							z^2
							-
							(1+\alpha-2(1-\alpha)\gamma\cK)z
							+
							\alpha
							-
							2\alpha(1-\alpha)\gamma\cK
							+
							2\alpha(1-\alpha)\gamma\cK
							=
							0}
						\\&=\pR{z\in\C\colon
							z^2
							-
							(1+\alpha-2(1-\alpha)\gamma\cK)z
							+
							\alpha
							=
							0}
						.
					\end{split}
			\end{equation}}
			\argument{
				\lref{eq:vanishing:det:1};
				\lref{eq:EV};}{that
				\begin{equation}
					\llabel{eq:EV:of:A}
					\cu{ z \in \C \colon \br{ \exists \, \nu \in \C^2 \backslash \cu{0} \colon A \nu = z \nu } }
					=
					\{\mu_-,\mu_+\}
					\qqandqq
					\sr(A)=\max\{\vass{\mu_-},\vass{\mu_+}\}
					.
				\end{equation}
			}
			\argument{\lref{eq:EV:of:A}}{\lref{it:evs}}.
			\argument{the fact that $\alpha<1$}{that
				\begin{equation}
					\llabel{eq:location:distinguisher}
					\begin{split}
						\frac{1+\sqrt{\alpha}}{2(1-\sqrt{\alpha})}
						=
						\frac{\pr{1+\sqrt{\alpha}}^2}{2(1-\alpha)}
						&=
						\frac{\pr{1+\sqrt{\alpha}}^2-2(1+\alpha)+2(1+\alpha)}{2(1-\alpha)}
						\\&=
						\frac{-(1-\sqrt{\alpha})^2+2(1+\alpha)}{2(1-\alpha)}
						\leq
						\frac{1+\alpha}{1-\alpha}
						.
					\end{split}
				\end{equation}
			}
			In our proof of \lref{it:sr} we distinguish between the cases
			$\gamma \cK\leq\frac{1 - \sqrt{\alpha}}{2\pr{1 + \sqrt{\alpha}}}$,
			$\frac{1 - \sqrt{\alpha}}{2\pr{1 + \sqrt{\alpha}}}<\gamma \cK< \frac{1 + \sqrt{\alpha}}{2\pr{1 - \sqrt{\alpha}}}$,
			$\frac{1 + \sqrt{\alpha}}{2\pr{1 - \sqrt{\alpha}}}\leq\gamma \cK< \frac{1 + \alpha}{1-\alpha}$,
			$\gamma \cK=\frac{1 + \alpha}{1-\alpha}$, and
			$\gamma \cK>\frac{1 + \alpha}{1-\alpha}$.
			We first \prove[ps]\ \lref{it:sr} in the case 
			\begin{equation}
				\llabel{eq:1st:case}
				\gamma \cK\leq\frac{1 - \sqrt{\alpha}}{2\pr{1 + \sqrt{\alpha}}}
				.
			\end{equation}
			\argument{\lref{eq:1st:case};}{that
				\begin{equation}
					\llabel{eq:1st:case:2}
					\frac{1 + \alpha - 2 ( 1 - \alpha )\gamma \cK }{2} 
					\geq
					\frac{1 + \alpha - ( 1 - \sqrt{\alpha} )^2 }{2}
					=
					\frac{2\sqrt{\alpha}  }{2}
					=
					\sqrt{\alpha}
					.
			\end{equation}}
			\argument{\lref{eq:1st:case:2};}{that
				\begin{equation}
					\llabel{eq:1st:case:3}
					\rbr*{ \frac{1 + \alpha - 2 ( 1 - \alpha )\gamma \cK }{2} } ^2 - \alpha \ge  0
					.
			\end{equation}}
			\argument{
				\lref{eq:EV};
				\lref{eq:1st:case:3};
			}{that
				\begin{equation}
					\llabel{eq:1st:case:4}
					\{\mu_-,\mu_+\}
					\subseteq
					\R
					\qqandqq
					\mu_-\leq\mu_+.
				\end{equation}
			}
			\argument{
				\lref{eq:1st:case:4};
				\lref{eq:EV};}{that
				\begin{equation}
					\llabel{eq:1st:case:5}
					\begin{split}
						\mu_+ 
						&=
						\frac{1 + \alpha - 2  ( 1 - \alpha ) \gamma\cK  }{2}
						+ \sqrt{ \rbr*{ \frac{1 + \alpha - 2( 1 - \alpha )\gamma\cK }{2} } ^2 - \alpha }
						\\&=
						\frac{1 + \alpha - 2  ( 1 - \alpha ) \gamma\cK  }{2}
						+ \sqrt{ \rbr*{ \frac{1 - \alpha + 2( 1 - \alpha )\gamma\cK }{2} +\alpha-2\gamma( 1 - \alpha )\cK } ^2 - \alpha }
						\\&= 	\frac{1 + \alpha - 2  ( 1 - \alpha ) \gamma\cK  }{2}
						+ \sqrt{ \rbr*{ \frac{1 - \alpha + 2  ( 1 - \alpha ) \gamma\cK }{2} } ^2+\alpha - 2  ( 1 - \alpha ) \gamma\cK -\alpha} \\
						& < \frac{1 + \alpha - 2  ( 1 - \alpha ) \gamma\cK  }{2} + \frac{\vass{1 - \alpha + 2  ( 1 - \alpha ) \gamma\cK}}{2}  = 1
						.
					\end{split} 
				\end{equation}
			}
			\argument{\lref{eq:1st:case};}{that
				\begin{equation}
					\llabel{eq:1st:case:6}
					\begin{split}
						&\rbr*{ \frac{1 + \alpha - 2 ( 1 - \alpha )\gamma \cK }{2} } ^2 - \alpha  - \rbr*{\frac{3 + \alpha - 2 ( 1 - \alpha )\gamma \cK }{2}}^2
						\\&=
						\frac{(1 + \alpha)^2 - (3+\alpha)^2}{4}
						+\frac{4(3+\alpha-(1+\alpha)) ( 1 - \alpha )\gamma \cK}{4}-\alpha 
						\\&=
						-2-\alpha
						+
						2 ( 1 - \alpha )\gamma \cK-\alpha 
						\\&\leq
						-2-2\alpha
						+
						(1-\sqrt{\alpha})^2
						=-1-2\sqrt{\alpha}-\alpha<0
						.
					\end{split}
			\end{equation}}
			\argument{
				\lref{eq:EV};
				\lref{eq:1st:case:6};
				\lref{eq:1st:case:3};
			}{that
				\begin{equation}
					\llabel{eq:1st:case:7}
					\begin{split}
						\mu_-
						&=	\frac{1 + \alpha - 2  ( 1 - \alpha ) \gamma\cK  }{2}
						- \sqrt{ \rbr*{ \frac{1 + \alpha - 2( 1 - \alpha )\gamma\cK }{2} } ^2 - \alpha }
						\\&>
						\frac{1 + \alpha - 2  ( 1 - \alpha ) \gamma\cK  }{2}
						- 
						\frac{3 + \alpha - 2  ( 1 - \alpha ) \gamma\cK  }{2}
						=-1
						.
					\end{split}
				\end{equation}
			}
			\argument{
				\lref{eq:EV:of:A};
				\lref{eq:1st:case:4};
				\lref{eq:1st:case:7};
				\lref{eq:1st:case:5};}{
				\begin{equation}
					\llabel{eq:1st:case:8}
					\sr(A)=\max\{\vass{\mu_-},\vass{\mu_+}\}<1.
			\end{equation}}
			\argument{\lref{eq:1st:case:8};\lref{eq:location:distinguisher}}[verbs=ps]{\lref{it:sr} in the case $\gamma \cK\leq\frac{1 - \sqrt{\alpha}}{2\pr{1 + \sqrt{\alpha}}}$}.
			In the next step we \prove[ps]\ \lref{it:sr} in the case 
			\begin{equation}
				\llabel{eq:2nd:case}
				\frac{1 - \sqrt{\alpha}}{2\pr{1 + \sqrt{\alpha}}}<\gamma \cK< \frac{1 + \sqrt{\alpha}}{2\pr{1 - \sqrt{\alpha}}}
				.
			\end{equation}
			\startnewargseq
			\argument{\lref{eq:2nd:case};}{that 
				\begin{equation}
					\llabel{eq:2nd:case:1}
					-\sqrt{\alpha}=\frac{1 + \alpha -(1+\sqrt{\alpha})^2}{2}
					<\frac{1 + \alpha - 2 ( 1 - \alpha )\gamma \cK }{2}
					< \frac{1 + \alpha -(1-\sqrt{\alpha})^2}{2}
					=\sqrt{\alpha}
					.
			\end{equation}}
			\argument{\lref{eq:2nd:case:1};}{that 
				\begin{equation}
					\llabel{eq:2nd:case:2}
					\prbbb{\frac{1 + \alpha - 2 ( 1 - \alpha )\gamma \cK }{2}}^2 
					<
					\alpha
					.
			\end{equation}}
			\argument{
				\lref{eq:EV};
				\lref{eq:EV:of:A};
				\lref{eq:2nd:case:2};
			}{that
				\begin{equation}
					\llabel{eq:2nd:case:3}
					\{\mu_-,\mu_+\}
					\subseteq
					\C\backslash\R
					\qqandqq
					\sr(A)=\max\{\vass{\mu_-},\vass{\mu_+}\}=\sqrt{\alpha}<1.
				\end{equation}
			}
			\argument{\lref{eq:2nd:case:3};\lref{eq:location:distinguisher};}[verbs=ps]{\lref{it:sr} in the case $\frac{1 - \sqrt{\alpha}}{2\pr{1 + \sqrt{\alpha}}}<\gamma \cK< \frac{1 + \sqrt{\alpha}}{2\pr{1 - \sqrt{\alpha}}}$}.
			In the next step we \prove[ps]\ \lref{it:sr} in the case 
			\begin{equation}
				\llabel{eq:3rd:case}
				\frac{1 + \sqrt{\alpha}}{2\pr{1 - \sqrt{\alpha}}}\leq\gamma \cK< \frac{1 + \alpha}{1-\alpha}
				.
			\end{equation}
			\startnewargseq
			\argument{\lref{eq:3rd:case};}{that
				\begin{equation}
					\llabel{eq:3rd:case:2}
					\frac{1 + \alpha - 2 ( 1 - \alpha )\gamma \cK }{2} 
					\leq
					\frac{1 + \alpha - ( 1 + \sqrt{\alpha} )^2 }{2}
					=
					-\frac{2\sqrt{\alpha}  }{2}
					=
					-\sqrt{\alpha}
					.
			\end{equation}}
			\argument{\lref{eq:3rd:case:2};}{that
				\begin{equation}
					\llabel{eq:3rd:case:3}
					\rbr*{ \frac{1 + \alpha - 2 ( 1 - \alpha )\gamma \cK }{2} } ^2 - \alpha \ge  0
					.
			\end{equation}}
			\argument{
				\lref{eq:EV};
				\lref{eq:3rd:case:3};
			}{that
				\begin{equation}
					\llabel{eq:3rd:case:4}
					\{\mu_-,\mu_+\}
					\subseteq
					\R
					\qqandqq
					\mu_-\leq\mu_+.
				\end{equation}
			}
			\argument{
				\lref{eq:3rd:case:4};
				\lref{eq:EV};}{that
				\begin{equation}
					\llabel{eq:3rd:case:5}
					\begin{split}
						\mu_+ 
						&=
						\frac{1 + \alpha - 2  ( 1 - \alpha ) \gamma\cK  }{2}
						+ \sqrt{ \rbr*{ \frac{1 + \alpha - 2( 1 - \alpha )\gamma\cK }{2} } ^2 - \alpha }
						\\&=
						\frac{1 + \alpha - 2  ( 1 - \alpha ) \gamma\cK  }{2}
						+ \sqrt{ \rbr*{ \frac{1 - \alpha + 2( 1 - \alpha )\gamma\cK }{2} +\alpha-2\gamma( 1 - \alpha )\cK } ^2 - \alpha }
						\\&= 	\frac{1 + \alpha - 2  ( 1 - \alpha ) \gamma\cK  }{2}
						+ \sqrt{ \rbr*{ \frac{1 - \alpha + 2  ( 1 - \alpha ) \gamma\cK }{2} } ^2+\alpha - 2  ( 1 - \alpha ) \gamma\cK -\alpha} \\
						& < \frac{1 + \alpha - 2  ( 1 - \alpha ) \gamma\cK  }{2} + \frac{\vass{1 - \alpha + 2  ( 1 - \alpha ) \gamma\cK}}{2}  = 1
						.
					\end{split} 
				\end{equation}
			}
			\argument{\lref{eq:3rd:case};}{that
				\begin{equation}
					\llabel{eq:3rd:case:6}
					\begin{split}
						&\rbr*{ \frac{1 + \alpha - 2 ( 1 - \alpha )\gamma \cK }{2} } ^2 - \alpha  - \rbr*{\frac{3 + \alpha - 2 ( 1 - \alpha )\gamma \cK }{2}}^2
						\\&=
						\frac{(1 + \alpha)^2 - (3+\alpha)^2}{4}
						+\frac{4(3+\alpha-(1+\alpha)) ( 1 - \alpha )\gamma \cK}{4}-\alpha 
						\\&=
						-2-\alpha
						+
						2 ( 1 - \alpha )\gamma \cK-\alpha 
						\\&<
						-2-2\alpha
						+2(1+\alpha)
						=0
						.
					\end{split}
			\end{equation}}
			\argument{
				\lref{eq:EV};
				\lref{eq:3rd:case:6};
				\lref{eq:3rd:case:3};
			}{that
				\begin{equation}
					\llabel{eq:3rd:case:7}
					\begin{split}
						\mu_-
						&=	\frac{1 + \alpha - 2  ( 1 - \alpha ) \gamma\cK  }{2}
						- \sqrt{ \rbr*{ \frac{1 + \alpha - 2( 1 - \alpha )\gamma\cK }{2} } ^2 - \alpha }
						\\&>
						\frac{1 + \alpha - 2  ( 1 - \alpha ) \gamma\cK  }{2}
						- 
						\frac{3 + \alpha - 2  ( 1 - \alpha ) \gamma\cK  }{2}
						=-1
						.
					\end{split}
				\end{equation}
			}
			\argument{
				\lref{eq:EV:of:A};
				\lref{eq:3rd:case:4};
				\lref{eq:3rd:case:7};
				\lref{eq:3rd:case:5};}{that
				\begin{equation}
					\llabel{eq:3rd:case:8}
					\sr(A)=\max\{\vass{\mu_-},\vass{\mu_+}\}<1.
			\end{equation}}
			\argument{\lref{eq:3rd:case:8};}[verbs=ps]{\lref{it:sr} in the case of $\frac{1 + \sqrt{\alpha}}{2\pr{1 - \sqrt{\alpha}}}\leq\gamma \cK< \frac{1 + \alpha}{1-\alpha}$}.
			In the next step we \prove[ps]\ \lref{it:sr} in the case  
			\begin{equation}
				\llabel{eq:4th:case}
				\gamma \cK=\frac{1 + \alpha}{1-\alpha}
				.
			\end{equation}
			\argument{
				\lref{eq:4th:case};
				\lref{eq:EV};}{that
				\begin{equation}
					\begin{split}
						\llabel{eq:4th:case:1}
						\mu_\pm 
						&=	\frac{1 + \alpha - 2  ( 1 - \alpha ) \gamma\cK  }{2}
						\pm \sqrt{ \rbr*{ \frac{1 + \alpha - 2( 1 - \alpha )\gamma\cK }{2} } ^2 - \alpha }
						\\&=	\frac{-1 - \alpha }{2}
						\pm \sqrt{ \rbr*{ \frac{1 + \alpha}{2} } ^2 - \alpha }
						\\&=	\frac{-1 - \alpha }{2}
						\pm \sqrt{ \rbr*{ \frac{1 - \alpha}{2} } ^2 }
						\\&=	\frac{-1 - \alpha }{2}
						\pm   \frac{1 - \alpha}{2}
						.
					\end{split}
			\end{equation}}
			\argument{\lref{eq:4th:case:1};\lref{eq:EV:of:A}}{that
				\begin{equation}
					\llabel{eq:4th:case:2}
					\mu_-=-1,
					\qquad
					\mu_+=-\alpha,
					\qqandqq
					\sr(A)=\max\{1,\alpha\}=1
					.
			\end{equation}}
			\argument{\lref{eq:4th:case:2};}[verbs=ps]{\lref{it:sr} in the case  $\gamma \cK=\frac{1 + \alpha}{1-\alpha}$}.
			In the next step we \prove[ps]\ \lref{eq:cond:for:result} in the case  
			\begin{equation}
				\llabel{eq:5th:case}
				\gamma \cK>\frac{1 + \alpha}{1-\alpha}
				.
			\end{equation}
			\startnewargseq
			\argument{\lref{eq:5th:case};}{that
				\begin{equation}
					\llabel{eq:5th:case:2}
					\frac{1 + \alpha - 2 ( 1 - \alpha )\gamma \cK }{2} 
					<
					\frac{1 + \alpha - 2(1+\alpha) }{2}
					=
					-\frac{1+\alpha }{2}
					=
					-\sqrt{\alpha}
					-\frac{(1-\sqrt{\alpha})^2}{2}
					<
					-\sqrt{\alpha}
					.
			\end{equation}}
			\argument{\lref{eq:5th:case:2};}{that
				\begin{equation}
					\llabel{eq:5th:case:3}
					\rbr*{ \frac{1 + \alpha - 2 ( 1 - \alpha )\gamma \cK }{2} } ^2 - \alpha >  0
					.
			\end{equation}}
			\argument{
				\lref{eq:EV};
				\lref{eq:EV:of:A};
				\lref{eq:5th:case:3};
			}{that
				\begin{equation}
					\llabel{eq:5th:case:4}
					\{\mu_-,\mu_+\}
					\subseteq
					\R
					\qqandqq
					\mu_-<\mu_+.
				\end{equation}
			}
			\argument{
				\lref{eq:5th:case:4};
				\lref{eq:EV};
				\lref{eq:5th:case};}{that
				\begin{equation}
					\llabel{eq:5th:case:5}
					\begin{split}
						\mu_+ 
						&=
						\frac{1 + \alpha - 2  ( 1 - \alpha ) \gamma\cK  }{2}
						+ \sqrt{ \rbr*{ \frac{1 + \alpha - 2( 1 - \alpha )\gamma\cK }{2} } ^2 - \alpha }
						\\&=
						\frac{1 + \alpha - 2  ( 1 - \alpha ) \gamma\cK  }{2}
						+ \sqrt{ \rbr*{ \frac{1 - \alpha + 2( 1 - \alpha )\gamma\cK }{2} +\alpha-2\gamma( 1 - \alpha )\cK } ^2 - \alpha }
						\\&= 	\frac{1 + \alpha - 2  ( 1 - \alpha ) \gamma\cK  }{2}
						+ \sqrt{ \rbr*{ \frac{1 - \alpha + 2  ( 1 - \alpha ) \gamma\cK }{2} } ^2+\alpha - 2  ( 1 - \alpha ) \gamma\cK -\alpha} \\
						& < \frac{1 + \alpha - 2  ( 1 - \alpha ) \gamma\cK  }{2} + \frac{\vass{1 - \alpha + 2  ( 1 - \alpha ) \gamma\cK}}{2}  = 1
						.
					\end{split} 
				\end{equation}
			}
			\argument{\lref{eq:5th:case};}{that
				\begin{equation}
					\llabel{eq:5th:case:6}
					\begin{split}
						&\rbr*{ \frac{1 + \alpha - 2 ( 1 - \alpha )\gamma \cK }{2} } ^2 - \alpha  - \rbr*{\frac{3 + \alpha - 2 ( 1 - \alpha )\gamma \cK }{2}}^2
						\\&=
						\frac{(1 + \alpha)^2 - (3+\alpha)^2}{4}
						+\frac{4(3+\alpha-(1+\alpha)) ( 1 - \alpha )\gamma \cK}{4}-\alpha 
						\\&=
						-2-\alpha
						+
						2 ( 1 - \alpha )\gamma \cK-\alpha 
						>
						-2-2\alpha
						+
						2(1+\alpha)
						=0
						.
					\end{split}
			\end{equation}}
			\argument{
				\lref{eq:EV};
				\lref{eq:5th:case:3};
				\lref{eq:5th:case:6};
			}{that
				\begin{equation}
					\llabel{eq:5th:case:7}
					\begin{split}
						\mu_-
						&=	\frac{1 + \alpha - 2  ( 1 - \alpha ) \gamma\cK  }{2}
						- \sqrt{ \rbr*{ \frac{1 + \alpha - 2( 1 - \alpha )\gamma\cK }{2} } ^2 - \alpha }
						\\&<
						\frac{1 + \alpha - 2  ( 1 - \alpha ) \gamma\cK  }{2}
						- \sqrt{ \rbr*{ \frac{3 + \alpha - 2( 1 - \alpha )\gamma\cK }{2} } ^2  }
						\\&=
						\frac{1 + \alpha - 2  ( 1 - \alpha ) \gamma\cK  }{2}
						-
						\frac{\vass{3 + \alpha - 2  ( 1 - \alpha ) \gamma\cK  }}{2}
						\\&=
						\frac{1 + \alpha - 2  ( 1 - \alpha ) \gamma\cK  }{2}
						-
						\vass[\bigg]{1+
							\frac{1 + \alpha - 2  ( 1 - \alpha ) \gamma\cK  }{2}}
						\leq-1
						.
					\end{split}
				\end{equation}
			}
			\argument{
				\lref{eq:EV:of:A};
				\lref{eq:5th:case:7};
				\lref{eq:5th:case:6};}{that
				\begin{equation}
					\llabel{eq:5th:case:8}
					\begin{split}
						\vass{\mu_-}>1,
						\qquad
						\vass{\mu_-}>\vass{\mu_+},
						\qqandqq
						\sr(A)=\max\{\vass{\mu_-},\vass{\mu_+}\}>1
						.
					\end{split}
			\end{equation}}
			\argument{\lref{eq:5th:case:8};}[verbs=ps]{\lref{it:sr} in the case $\gamma \cK>\frac{1 + \alpha}{1-\alpha}$}.
		\end{aproof}

		\subsection{Asymptotic analysis for the momentum optimizer with constant learning rates}

		\begin{athm}{prop}{prop:mom-stable}
			Let $\fd \in \N$,
			$\gamma, \cK \in (0, \infty)$,
			$\alpha \in (0, 1 )$,
			$\lp \in \R^\fd$,
			let
			$\fl \colon \R^\fd \to \R$ satisfy for all $\theta \in \R^\fd$ that
			\begin{equation}
				\fl ( \theta ) = \cK\norm{\theta - \lp} ^2,
			\end{equation}
			and let $\Theta \colon \N_0 \to \R^\fd$
			and
			$ \m \colon \N_0 \to \R^\fd$ satisfy for all $n \in \N$ that
			\begin{equation}
				\llabel{eq:start:and:setup}
				\m_0=0,
				\quad
				\m_n = \alpha \m_{n-1} + (1-\alpha) ( \nabla \fl ) ( \Theta_{n-1} ) ,
				\quad
				\Theta_0 \ne \lp,
				\qandq  \Theta_n = \Theta_{n-1} - \gamma \m_n.
			\end{equation}
			Then
			\begin{equation}
				\llabel{eq:result}
				\liminf_{n \to \infty} \norm{\Theta_n - \lp} = \limsup_{n \to \infty} \norm{\Theta_n - \lp } 
				\in
				\begin{cases}
					\{0\} & \colon \gamma\cK < \tfrac{  1 + \alpha  }{  1 - \alpha} \\
					( 0 , \infty ) & \colon \gamma\cK =  \tfrac{ 1 + \alpha  }{ 1 - \alpha  }\\
					\{\infty\} & \colon \gamma\cK >  \tfrac{1 + \alpha }{   1 - \alpha }.
				\end{cases}
			\end{equation}
		\end{athm}

		\begin{aproof}
			Throughout this proof assume without loss of generality that $\fd = 1$ and $\lp = 0$
			and let $A\in\C^{2\times 2}$, $\mu_-,\mu_+\in\C$ satisfy
			\begin{equation}
				\llabel{eq:matrix}
				A = \begin{pmatrix}
					\alpha & 2(1 - \alpha) \cK \\
					- \gamma \alpha & 1 - 2  ( 1 - \alpha ) \gamma\cK
				\end{pmatrix}
			\end{equation}
			\begin{equation}
				\llabel{eq:EV}
				\qqandqq
				\mu_\pm =	\frac{1 + \alpha - 2  ( 1 - \alpha ) \gamma\cK  }{2}
				\pm \sqrt{ \rbr*{ \frac{1 + \alpha - 2 ( 1 - \alpha )\gamma \cK }{2} } ^2 - \alpha }
				.
			\end{equation}
			\argument{\lref{eq:matrix};\lref{eq:EV}; \cref{lem:EV} (applied with
				$\gamma\curvearrowleft\gamma$,
				$\cK\curvearrowleft\cK$,
				$\alpha\curvearrowleft\alpha$,
				$A\curvearrowleft A$,
				$\mu-\curvearrowleft\mu_-$,
				$\mu_+\curvearrowleft\mu_+$ in the notation of \cref{lem:EV})}{that
				\begin{equation}
					\llabel{eq:EV:of:A}
					\cu{ z \in \C \colon \br{ \exists \, \nu \in \C^2 \backslash \cu{0} \colon A \nu = z \nu } }
					=
					\{\mu_-,\mu_+\},
				\end{equation}
				\begin{equation}
					\llabel{eq:cond:for:result:2}
					\max\{1,\vass{\mu_-}\}>\vass{\mu_+},
					\qqandqq
					\sr(A)\in
					\begin{cases}
						[0,1) &\colon \gamma \cK<\frac{1 +\alpha}{1-\alpha}\\
						\{1\} &\colon \gamma \cK=\frac{1 +\alpha}{1-\alpha}\\
						(1,\infty) &\colon \gamma \cK>\frac{1 +\alpha}{1-\alpha}
					\end{cases}
			\end{equation}}
			\argument{the fact that $\forall \, \theta \in \R \colon \nabla \fl ( \theta ) = 2\cK \theta$}{
				\llabel{eq:recursion}
				for all $n \in \N$ that
				$\m_n = \alpha \m_{n-1} + 2(1 - \alpha ) \cK \Theta_{n-1}$
				and
				$\Theta_n = \Theta_{n-1} - \gamma \alpha \m_{n-1} - 2  ( 1 - \alpha )\gamma\cK \Theta_{n-1}$.
			} 
			\argument{
				\lref{eq:matrix};
				\lref{eq:recursion};
			}{that for all $n \in \N$ it holds that
				\begin{equation}
					\llabel{eq:recursion:with:matrix}
					\begin{pmatrix}
						\m_n \\ \Theta_n 
					\end{pmatrix}
					= A \begin{pmatrix}
						\m_{n-1} \\
						\Theta_{n-1} 
					\end{pmatrix}.
				\end{equation}
			}
			\argument{\lref{eq:recursion:with:matrix};}{for all $n \in \N$ that
				\begin{equation}
					\llabel{eq:rec}
					\begin{pmatrix}
						\m_n \\ \Theta_n 
					\end{pmatrix}
					= A ^n \begin{pmatrix}
						\m_{0} \\
						\Theta_{0} 
					\end{pmatrix}.
				\end{equation}
			}
			\argument{
				\lref{eq:matrix};}{that for all $\lambda\in\C$ it holds that
				\begin{equation}
					\begin{split}
						\llabel{eq:no:EV}
						\norm[\bigg]{A\begin{pmatrix}
								0 \\ 1
							\end{pmatrix}
							-
							\lambda\begin{pmatrix}
								0 \\ 1
						\end{pmatrix}} 
						&=
						\norm[\bigg]{\begin{pmatrix}
								2(1-\alpha)\cK \\ 1-2(1-\alpha)\gamma\cK-\lambda
							\end{pmatrix}
						}
						\geq
						2(1-\alpha)\cK>0
						\qqandqq
					\end{split}
				\end{equation}
				\begin{equation}					
					\norm[\bigg]{A\begin{pmatrix}
							1 \\ 0
						\end{pmatrix}
						-
						\lambda\begin{pmatrix}
							1 \\ 0
					\end{pmatrix}} 
					=
					\norm[\bigg]{\begin{pmatrix}
							\alpha-\lambda \\ -\gamma\alpha
						\end{pmatrix}
					}
					\geq
					\gamma\alpha>0
					.
			\end{equation}}
			\argument{
				\lref{eq:start:and:setup};
				\lref{eq:no:EV};
				\lref{eq:EV:of:A}}{that there exist $v=(v_1,v_2)$, $w=(w_1,w_2)\in\{(x,y)\in\C^2\colon xy\neq 0\}$ which satisfy that
				\begin{equation}
					\begin{split}
						\llabel{eq:start:EV:repr}
						\begin{pmatrix}
							0 \\ \Theta_0
						\end{pmatrix}
						=
						\begin{pmatrix}
							\m_0 \\ \Theta_0
						\end{pmatrix}
						=v+w,
						\qquad
						Av=\mu_+v,
						\qqandqq
						Aw=\mu_-w
						.
					\end{split}
				\end{equation}
			}
			\argument{
				\lref{eq:rec};
				\lref{eq:start:EV:repr};}{that for all $n\in\N$ it holds that
				\begin{equation}
					\llabel{eq:n_th:vec}
					\begin{pmatrix}
						\m_n \\ \Theta_n 
					\end{pmatrix}
					=
					A^n\begin{pmatrix}
						\m_0 \\ \Theta_0
					\end{pmatrix}
					=
					A^n
					\begin{pmatrix}
						0 \\ \Theta_0
					\end{pmatrix}
					=
					A^n\pr{v+w}
					.
			\end{equation}}
			\argument{
				\lref{eq:n_th:vec};
				\lref{eq:cond:for:result:2};
				\lref{eq:start:EV:repr};
				\cref{lem:mat:seq:EV} (applied with 
				$\mu_1\curvearrowleft\mu_-$,
				$\mu_2\curvearrowleft\mu_+$,
				$v_1\curvearrowleft w$,
				$v_2\curvearrowleft v$,
				$A\curvearrowleft A$,
				$\m\curvearrowleft\m$,
				$\Theta\curvearrowleft\Theta$,
				in the notation of \cref{lem:mat:seq:EV})}{
				\begin{equation}
					\llabel{eq:result:appl}
					\liminf_{n \to \infty} \vass{\Theta_n} = \limsup_{n \to \infty} \vass{\Theta_n } 
					=
					\begin{cases}
						0 & \colon \sr\pr{A}<1\\
						\vass{w_{2}} & \colon \sr\pr{A}=1\\
						\infty & \colon \sr\pr{A}>1.
					\end{cases}
				\end{equation}
			}
			\argument{
				\lref{eq:result:appl};
				\lref{eq:start:EV:repr};}{\lref{eq:result}}.
		\end{aproof}

			\begin{athm}{cor}{prop:mom-stable:diag}
		Let $\fd \in \N$,
		$\gamma, \lambda_1,\lambda_2,\dots,\lambda_d \in [0, \infty)$,
		$\alpha \in (0, 1 )$,
		$\lp=(\lp_1,\dots,\lp_d) \in \R^\fd$,
		let
		$\fl \colon \R^\fd \to \R$ satisfy for all $\theta=(\theta_1,\dots,\theta_d) \in \R^\fd$ that
		\begin{equation}
			\llabel{eq:loss:diag:mom}
			\fl ( \theta ) = \textstyle\sum_{i=1}^d\lambda_i\pr{\theta_i - \lp_i}^2,
		\end{equation}
		and let $\Theta \colon \N_0 \to \R^\fd$
		and
		$ \m \colon \N_0 \to \R^\fd$ satisfy for all $n \in \N$ that
		\begin{equation}
			\llabel{eq:start:and:setup:diag}
			\m_0=0,
			\quad
			\m_n = \alpha \m_{n-1} + (1-\alpha) ( \nabla \fl ) ( \Theta_{n-1} ) ,
			\quad
			\Theta_0 \ne \lp,
			\qandq  \Theta_n = \Theta_{n-1} - \gamma \m_n.
		\end{equation}
		Then it holds that $\sup_{n \in\N} \norm{\Theta_n }<\infty$ if and only if 
		$
					 \gamma\max\{\lambda_1,\lambda_2,\dots,\lambda_d\} \leq \tfrac{  1 + \alpha  }{  1 - \alpha}
		$.
	\end{athm}

	\begin{aproof}
		Throughout this proof for every $i\in\{1,2,\dots,d\}$ let $p_i\colon\R^d\to\R$ satisfy for all $x=(x_1,\dots,x_d)\in\R^d$ that 
		\begin{equation}
			\llabel{eq:def:proj}
		p_i(x)=x_i.
		\end{equation}
		\argument{
			\lref{eq:start:and:setup:diag};
			\lref{eq:def:proj};
			\cref{momentum:representation}}{that for all $i\in\{1,2,\dots,d\}$, $n\in\N$ it holds that
			\begin{equation}
				\llabel{eq:rep:mom:expl}
				\begin{split}
				p_i(\m_n)
				\textstyle=p_i\prb{(1-\alpha)\sum_{k=1}^n\alpha^{n-k} ( \nabla \fl ) ( \Theta_{k-1} )}
				&\textstyle=(1-\alpha)\sum_{k=1}^n\alpha^{n-k} p_i\prb{( \nabla \fl ) ( \Theta_{k-1} )}
				\\&\textstyle=(1-\alpha)\sum_{k=1}^n\alpha^{n-k} 2\lambda_i p_i ( \Theta_{k-1} -\lp)
				.
				\end{split}
		\end{equation}}
		\argument{
			\lref{eq:loss:diag:mom};
			\lref{eq:start:and:setup:diag};
			\lref{eq:def:proj}}{that for all $i\in\{1,2,\dots,d\}$, $n\in\N$ with $\gamma\lambda_i=0$ it holds that
			\begin{equation}
				\llabel{eq:vanishing:case}
				\begin{split}
				p_i(\Theta_n)
				=p_i(\Theta_{n-1}-\gamma\m_n)
				&=p_i(\Theta_{n-1})-\gamma p_i(\m_n)
				\\&\textstyle=p_i(\Theta_{n-1})-2\gamma\lambda_i (1-\alpha)\sum_{k=1}^n\alpha^{n-k}  p_i ( \Theta_{k-1} -\lp)
				\\&=p_i(\Theta_{n-1}).
				\end{split}
				\end{equation}}
			\argument{\cref{prop:mom-stable} (applied with 
				$d\curvearrowleft 1$,
				$\gamma\curvearrowleft \gamma$,
				$\cK\curvearrowleft \lambda_i$,
				$\alpha\curvearrowleft \alpha$,
				$\lp\curvearrowleft \lp_i$,
				$\fl\curvearrowleft \pr{\R\ni\theta\mapsto \lambda_i(\theta-\lp_i)^2\in\R}$,
				$\Theta\curvearrowleft \pr{\N_0\ni n\mapsto p_i(\Theta_n)\in\R}$,
				$\m\curvearrowleft \pr{\N_0\ni n\mapsto p_i(\m_n)\in\R}$ for $i\in\{1,2,\dots,d\}$
				 in the notation of \cref{prop:mom-stable})}{that for all $i\in\{1,2,\dots,d\}$ with $\gamma\lambda_i>0$ it holds that
				 \begin{equation}
				 	\llabel{eq:result:prop:mom}
				 	\liminf_{n \to \infty} \vass{p_i(\Theta_n - \lp)} 
				 	= \limsup_{n \to \infty} \vass{p_i(\Theta_n - \lp) } 
				 	\in
				 	\begin{cases}
				 		\{0\} & \colon \gamma\lambda_i < \tfrac{  1 + \alpha  }{  1 - \alpha} \\
				 		( 0 , \infty ) & \colon \gamma\lambda_i =  \tfrac{ 1 + \alpha  }{ 1 - \alpha  }\\
				 		\{\infty\} & \colon \gamma\lambda_i >  \tfrac{1 + \alpha }{   1 - \alpha }.
				 	\end{cases}
				 	\end{equation}}
			 	\argument{
			 		\lref{eq:vanishing:case};
			 		\lref{eq:result:prop:mom};}{that 
			 		\begin{equation}
			 			\llabel{eq:result}
			 			\textstyle
			 			\sup_{n\in\N}\norm{\Theta_n}^2
			 			=
			 			\textstyle
			 			\sup_{n\in\N}\prb{\sum_{i=1}^d \vass{p_i(\Theta_n)}^2}
			 			\in
			 			\begin{cases}
			 				[0,\infty)
			 				& \colon \gamma\max\{\lambda_1,\lambda_2,\dots,\lambda_d\} \leq \tfrac{  1 + \alpha  }{  1 - \alpha} \\
			 				\{\infty\} 
			 				& \colon \gamma\max\{\lambda_1,\lambda_2,\dots,\lambda_d\} >  \tfrac{1 + \alpha }{   1 - \alpha }.
			 				\end{cases}
			 			\end{equation}}
		 	\argument{\lref{eq:result}}{that $\sup_{n \in\N} \norm{\Theta_n }<\infty$ if and only if 
		 		$
		 		\gamma\max\{\lambda_1,\lambda_2,\dots,\lambda_d\} \leq \tfrac{  1 + \alpha  }{  1 - \alpha}
		 		$}.	\end{aproof}

		\subsection{Asymptotic analysis for the momentum optimizer with convergent learning rates}
		
		\begin{athm}{prop}{prop:mom_conv_stable}
			Let 
			$\fd \in \N$,
			$\lp \in \R^\fd$,
			$\gamma, \cK \in (0, \infty)$,
			$\alpha \in (0, 1 )$,
			let
			$\fl \colon \R^\fd \to \R$ satisfy for all $\theta \in \R^\fd$ that
			\begin{equation}
				\llabel{eq:setup:cvg:LR}
				\fl ( \theta ) = \cK \norm{\theta - \lp} ^2
				\qqandqq
				\gamma\cK < \tfrac{ 1 + \alpha  }{  1 - \alpha },
			\end{equation}
			let $a\colon\N\to(0,\infty)$ and $\lambda\colon\N\to(0,\infty)$
			satisfy $\limsup_{n \to \infty }\pr{\vass{ a_n - \alpha}+\vass{\lambda_n - \gamma}}=0$,
			and let $\Theta\colon \N_0 \to \R^\fd$ and $\m \colon \N_0 \to \R^\fd$ satisfy for all $n \in \N$ that
			\begin{equation}
				\llabel{eq:recursion:cvg:LR}
				\m_0=0,
				\qquad
				\m_n = a_n \m_{n-1} + ( 1 - a_n ) ( \nabla \fl ) ( \Theta_{n-1} ),
			\end{equation}
			\begin{equation}
				\Theta_0 \ne \lp,
				\qqandqq  
				\Theta_n = \Theta_{n-1} - \lambda_n \m_n.
			\end{equation}
			Then 
			$\limsup_{n \to \infty} \norm{\Theta_n - \lp }  = 0$.
		\end{athm}
		
		\begin{cproof}{prop:mom_conv_stable}
			Throughout this proof assume without loss of generality that $\fd = 1$ and $\lp = 0$
			and let $A\colon\N_0\to\C^{2\times 2}$ satisfy for all $n\in\N$ that
			\begin{equation}
				\llabel{eq:matrix:Nest}
				A_n = \begin{pmatrix}
					a _n & 2(1 - a_n ) \cK \\
					- \lambda _n a _n & 1 -  2( 1 - a _n )\lambda _n \cK
				\end{pmatrix}
				\qqandqq
				A_0 = \begin{pmatrix}
					\alpha & 2(1 - \alpha ) \cK \\
					- \gamma  \alpha  & 1 -  2( 1 - \alpha  )\gamma  \cK
				\end{pmatrix}
				.
			\end{equation}
			\argument{\lref{eq:matrix:Nest};\cref{lem:EV};}{
				\begin{equation}
					\llabel{eq:small:sr}
					\sr(A_0)<1
			\end{equation}}
			\argument{the fact that $\forall \, \theta \in \R \colon \nabla \fl ( \theta ) = 2\cK \theta$}{
				\llabel{eq:recursion}
				for all $n \in \N$ that
				$\m_n = \alpha_n \m_{n-1} + 2(1 - \alpha_n ) \cK \Theta_{n-1}$
				and
				$\Theta_n = \Theta_{n-1} - \gamma_n \alpha_n \m_{n-1} - 2  ( 1 - \alpha_n )\gamma_n\cK \Theta_{n-1}$.
			} 
			\argument{
				\lref{eq:recursion};
				\lref{eq:recursion:cvg:LR};
			}{that for all $n \in \N$ it holds that
				\begin{equation}
					\llabel{eq:recursion:with:matrix}
					\begin{pmatrix}
						\m_n \\ \Theta_n 
					\end{pmatrix}
					= A_n \begin{pmatrix}
						\m_{n-1} \\
						\Theta_{n-1} 
					\end{pmatrix}.
				\end{equation}
			}
			\argument{the assumption that $\lim_{n \to \infty } a_n = \alpha$ and $\lim_{n \to \infty } \lambda_n = \gamma$}{that 
				\begin{equation}
					\llabel{eq:conv:matr}
					\begin{split}
					&\limsup_{n\to\infty}\prbb{\sup_{v\in\C^{2}\backslash\{0\}}\frac{\norm{A_kv-A_0v}}{\norm{v}}}
					\\&=\limsup_{n\to\infty}\prbbb{\sup_{v\in\{w\in\C^{2}\colon\norm{w}=1\}}\norm{(A_n-A_0)v}}
					\\&=\limsup_{n\to\infty}\prbbb{\sup_{v\in\{w\in\C^{2}\colon\norm{w}=1\}}\norm[\bigg]{
							\begin{pmatrix}
								a _n-\alpha & 2(1 - a_n ) \cK-2(1 - \alpha) \cK \\
								- \lambda _n a _n + \gamma \alpha&  -  2( 1 - a _n )\lambda _n \cK+  2( 1 - \alpha )\gamma \cK)
							\end{pmatrix}v}}
					\\&=0
					\end{split}
			\end{equation}}
			\argument{
				\lref{eq:small:sr};
				\lref{eq:recursion:with:matrix};
				\lref{eq:conv:matr};
				\cref{lem:mat:seq:EV:small};}{that
				$
				\limsup_{n \to \infty}\vass{\Theta_n}=0$}.
		\end{cproof}

		\begin{athm}{cor}{cor:mom-stable}[Bias-adjusted momentum]
				Let $\fd \in \N$,
				$\gamma, \cK \in (0, \infty)$,
				$\alpha \in (0, 1 )$,
				$\lp \in \R^\fd$, let
				$\fl \colon \R^\fd \to \R$ satisfy for all $\theta \in \R^\fd$ that
				\begin{equation}
					\fl ( \theta ) = \cK \norm{\theta - \lp} ^2,
				\end{equation}
				and let $\Theta\colon \N_0 \to \R^\fd$ and $\m \colon \N_0 \to \R^\fd$ satisfy for all $n \in \N$ that
				\begin{equation}
					\m_0=0,
					\qquad
					\m_n = \alpha \m_{n-1} + (1-\alpha) ( \nabla \fl ) ( \Theta_{n-1} ) ,
				\end{equation}
				\begin{equation}
					\Theta_0 \ne \lp,
					\qqandqq  \Theta_n = \Theta_{n-1} - \gamma (1 - \alpha^n ) ^{-1} \m_n.
				\end{equation}
				Then
				$\limsup_{n \to \infty} \norm{\Theta_n - \lp }  = 0$.
		\end{athm}
		
		\begin{aproof}
			\argument{the fact that $\lim_{n\to\infty}((1-\alpha^n)^{-1}\gamma)=\gamma$;\cref{prop:mom_conv_stable} (applied with 
				$d\curvearrowleft d$,
				$\gamma\curvearrowleft \gamma$,
				$\cK\curvearrowleft \cK$,
				$\alpha\curvearrowleft\alpha$,
				$\fl\curvearrowleft\fl$,
				$a\curvearrowleft\pr{\N\ni n\mapsto \alpha\in(0,\infty)}$,
				$\lambda\curvearrowleft\pr{\N\ni n\mapsto (1-\alpha^n)^{-1}\gamma\in(0,\infty)}$,
				$\Theta \curvearrowleft \Theta$,
				$\m\curvearrowleft \m$
				in the notation of \cref{prop:mom_conv_stable});}{that
				\begin{equation}
					\limsup_{n \to \infty} \norm{\Theta_n - \lp }  = 0.
			\end{equation}}
		\end{aproof}
		
		\section{A priori bounds for the Nesterov optimizer}
		\label{section:apriori:Nest}

		In this section we specify in \cref{prop:Nest-stable:diag} below the set of tupels of learning rates and eigenvalues of the Hessian for which the Nesterov method does not explode but stays bounded when applied to the simple class of quadratic \OPs\ in \cref{eq:min:prob} in \cref{subsection:1:1} above.
		This will allow us to explicitly specify the stability region (cf.\ \cref{def:stab:region} in \cref{subsection:1:1}) of the Nesterov optimizer, that is, to establish \cref{it:Nest:perm} in \cref{intro:thm:1.1} above (cf.~\cref{subsection:5.5}).
		The convergence parts within the stability region in \cref{prop:Nest-stable:diag} (cf.\ \cref{prop:nesterov-stable}) is already well-known in the literature \cite{TrungMomLR}.
		
		\subsection{One step analysis for the Nesterov optimizer}
		
		\begin{athm}{lemma}{lem:EV:Nest}
			Let $\gamma, \cK \in (0, \infty)$,
			$\alpha \in (0, 1 )$, $A\in\C^{2\times 2}$ satisfy
			\begin{equation}
				\llabel{eq:matrix:2}
				A = \begin{pmatrix}
					(1+\alpha)(1-2\gamma\cK) & -\alpha \\
					1 - 2 \gamma\cK & 0
				\end{pmatrix}
			\end{equation}
			and let $\mu_-,\mu_+\in\C$ satisfy
			\begin{equation}
				\llabel{eq:EV}
				\mu_\pm =\frac{( 1 + \alpha )( 1 - 2\gamma \cK ) }{2}
				\pm \sqrt{ \rbr*{ \frac{ ( 1 + \alpha )( 1 - 2\gamma \cK ) }{2} } ^2 - \alpha ( 1 - 2\gamma \cK ) }
				.
			\end{equation}
			Then
			\begin{enumerate}[label=(\roman*)]
				\item\llabel{it:evs} it holds that $\cu{ z \in \C \colon \br{ \exists \, \nu \in \C^2 \backslash \cu{0} \colon A \nu = z \nu } }
				=
				\{\mu_-,\mu_+\}$
				and
				\item\llabel{it:sr} it holds that
				\begin{equation}
					\llabel{eq:cond:for:result:3}
					\max\{1,\vass{\mu_-}\}>\vass{\mu_+}
					\qqandqq
					\sr(A)\in
					\begin{cases}
						[0,1) &\colon \gamma \cK<\frac{1 +\alpha}{1+2\alpha}\\
						\{1\} &\colon \gamma \cK=\frac{1 +\alpha}{1+2\alpha}\\
						(1,\infty) &\colon \gamma \cK>\frac{1 +\alpha}{1+2\alpha}.\\
					\end{cases}
				\end{equation}
			\end{enumerate}
		\end{athm}

		\begin{aproof}
			\argument{
				\lref{eq:matrix:2};}{that
				\begin{equation}
					\llabel{eq:vanishing:det:0}
					\cu{ z \in \C \colon \br{ \exists \, \nu \in \C^2 \backslash \cu{0} \colon A \nu = z \nu } }
					=
					\pRbbb{z\in\C\colon\operatorname{det}
						\begin{pmatrix}
							(1+\alpha)(1-2\gamma\cK)-z & -\alpha \\
							1 - 2 \gamma\cK & -z
						\end{pmatrix}
						=0}.
				\end{equation}
			}
			\argument{\lref{eq:vanishing:det:0}}{
				that 				\begin{equation}
					\llabel{eq:vanishing:det:1}
					\begin{split}
						&\cu{ z \in \C \colon \br{ \exists \, \nu \in \C^2 \backslash \cu{0} \colon A \nu = z \nu } }
						\\&=\pR{z\in\C\colon
							z^2
							-
							(1+\alpha)(1-2\gamma\cK)z
							+
							\alpha(1-2\gamma\cK)
							=
							0}
						.
					\end{split}
			\end{equation}}
			\argument{
				\lref{eq:vanishing:det:1};
				\lref{eq:EV};}{that
				\begin{equation}
					\llabel{eq:EV:of:A:Nest:mat}
					\cu{ z \in \C \colon \br{ \exists \, \nu \in \C^2 \backslash \cu{0} \colon A \nu = z \nu } }
					=
					\{\mu_-,\mu_+\}
					\qqandqq
					\sr(A)=\max\{\vass{\mu_-},\vass{\mu_+}\}
					.
				\end{equation}
			}
			\argument{\lref{eq:EV:of:A}}{\lref{it:evs}}.
			\argument{the fact that $0<\alpha<1$}{that
				\begin{equation}
					\llabel{eq:location:distinguisher}
					\begin{split}
						\frac{\pr{1-\alpha}^2}{2\pr{1+\alpha}^2}
						=
						\frac{(1-\alpha)^2}{2+4\alpha+2\alpha^2}
						\leq
						\frac{(1-\alpha)^2}{1+2\alpha}
						<
						\frac{1+\alpha}{1+2\alpha}
						<
						\frac{3+\alpha}{1+\alpha}
						.
					\end{split}
				\end{equation}
			}
			In our proof of \lref{it:sr} we distinguish between the cases
			$\gamma \cK\leq\frac{\pr{1-\alpha}^2}{2\pr{1+\alpha}^2}$,
			$\frac{\pr{1-\alpha}^2}{2\pr{1+\alpha}^2}<\gamma \cK< \frac{1}{2}$,
			$\frac{1}{2}\leq\gamma \cK< \frac{1+\alpha}{1+2\alpha}$,
			$\gamma\cK=\frac{1+\alpha}{1+2\alpha}$, and
			$\gamma \cK>\frac{1+\alpha}{1+2\alpha}$.
			We first \prove[ps]\ \lref{it:sr} in the case 
			\begin{equation}
				\llabel{eq:1st:case}
				\gamma \cK\leq\frac{\pr{1-\alpha}^2}{2\pr{1+\alpha}^2}
				.
			\end{equation}
			\argument{\lref{eq:1st:case};}{that
				\begin{equation}
					\llabel{eq:1st:case:2}
					\begin{split}
						\frac{ ( 1 + \alpha )^2( 1 - 2\gamma \cK ) }{4}
						\geq
						\frac{ ( 1 + \alpha )^2- (1-\alpha)^2 }{4}
						&=
						\frac{4\alpha}{4}
						=
						\alpha
						.
					\end{split}
			\end{equation}}
			\argument{\lref{eq:1st:case:2};\lref{eq:1st:case};the fact that $\alpha>0$}{that
				\begin{equation}
					\llabel{eq:1st:case:3}
					\rbr*{ \frac{ ( 1 + \alpha )( 1 - 2\gamma \cK ) }{2} } ^2  - \alpha ( 1 - 2\gamma \cK ) \ge  0
					.
			\end{equation}}
			\argument{
				\lref{eq:EV};
				\lref{eq:1st:case:3};
			}{that
				\begin{equation}
					\llabel{eq:1st:case:4}
					\{\mu_-,\mu_+\}
					\subseteq
					\R
					\qqandqq
					\mu_-\leq\mu_+.
				\end{equation}
			}
			\argument{
				\lref{eq:1st:case};
				\lref{eq:1st:case:4};
				\lref{eq:EV};}{that
				\begin{equation}
					\llabel{eq:1st:case:5}
					\begin{split}
						\mu_+ 
						&=
						\frac{( 1 + \alpha )( 1 - 2\gamma \cK ) }{2}
						+ \sqrt{ \rbr*{ \frac{ ( 1 + \alpha )( 1 - 2\gamma \cK ) }{2} } ^2 - \alpha ( 1 - 2\gamma \cK ) }
						\\&<
						\frac{( 1 + \alpha )( 1 - 2\gamma \cK ) }{2}
						+ \sqrt{ \rbr*{ \frac{ ( 1 + \alpha )( 1 - 2\gamma \cK ) }{2} } ^2 - (1+\alpha) ( 1 - 2\gamma \cK )+1 }
						\\&= \frac{( 1 + \alpha )( 1 - 2\gamma \cK ) }{2} + \vass*{\frac{( 1 + \alpha )( 1 - 2\gamma \cK ) }{2} -1}
						\\&= \frac{( 1 + \alpha )( 1 - 2\gamma \cK ) }{2} + \rbr*{ 1 - \frac{( 1 + \alpha )( 1 - 2\gamma \cK ) }{2} } = 1
						.
					\end{split} 
				\end{equation}
			}
			\argument{
				\lref{eq:EV};
				\lref{eq:1st:case};
				\lref{eq:1st:case:4};
				\lref{eq:1st:case:3};
			}{that
				\begin{equation}
					\llabel{eq:1st:case:7}
					\begin{split}
						\mu_-
						&=	
						\frac{( 1 + \alpha )( 1 - 2\gamma \cK ) }{2}
						- \sqrt{ \rbr*{ \frac{ ( 1 + \alpha )( 1 - 2\gamma \cK ) }{2} } ^2 - \alpha ( 1 - 2\gamma \cK ) }
						\\&>
						\frac{( 1 + \alpha )( 1 - 2\gamma \cK ) }{2}
						- \sqrt{ \rbr*{ \frac{ ( 1 + \alpha )( 1 - 2\gamma \cK ) }{2} } ^2 }
						\\&=
						\frac{( 1 + \alpha )( 1 - 2\gamma \cK ) }{2}
						- 
						\frac{( 1 + \alpha )( 1 - 2\gamma \cK ) }{2}
						=0
						.
					\end{split}
				\end{equation}
			}
			\argument{
				\lref{eq:vanishing:det:1};
				\lref{eq:1st:case:4};
				\lref{eq:1st:case:7};
				\lref{eq:1st:case:5};}{that
				\begin{equation}
					\llabel{eq:1st:case:8}
					\sr(A)=\max\{\vass{\mu_-},\vass{\mu_+}\}<1.
			\end{equation}}
			\argument{\lref{eq:1st:case:8}}[verbs=ps]{\lref{it:sr} in the case  $\gamma \cK\leq\frac{\pr{1-\alpha}^2}{2\pr{1+\alpha}^2}$}.
			In the next step we \prove[ps]\ \lref{it:sr} in the case 
			\begin{equation}
				\llabel{eq:2nd:case}
				\frac{\pr{1-\alpha}^2}{2\pr{1+\alpha}^2}<\gamma \cK< \frac{1}{2}
				.
			\end{equation}
			\startnewargseq
			\argument{\lref{eq:2nd:case};}{that 
				\begin{equation}
					\llabel{eq:2nd:case:1}
					\alpha
					=\frac{ ( 1 + \alpha )^2- (1-\alpha)^2 }{4}
					=\frac{ ( 1 + \alpha )^2\prb{1- \frac{\pr{1-\alpha}^2}{\pr{1+\alpha}^2}} }{4}
					>\frac{ ( 1 + \alpha )^2( 1 - 2\gamma \cK ) }{4}
					>\frac{ ( 1 + \alpha )^2( 1 - 1 ) }{4}=0
					.
			\end{equation}}
			\argument{\lref{eq:2nd:case:1};\lref{eq:2nd:case}}{that 
				\begin{equation}
					\llabel{eq:2nd:case:2}
					\rbr*{ \frac{ ( 1 + \alpha )( 1 - 2\gamma \cK ) }{2} } ^2 - \alpha ( 1 - 2\gamma \cK )
					=
					( 1 - 2\gamma \cK )
					\prbb{\frac{ ( 1 + \alpha )^2( 1 - 2\gamma \cK ) }{4}-\alpha}
					<0
					.
			\end{equation}}
			\argument{
				\lref{eq:EV};
				\lref{eq:EV:of:A:Nest:mat};
				\lref{eq:2nd:case:2};
				\lref{eq:2nd:case};
			}{that
				\begin{equation}
					\llabel{eq:2nd:case:3}
					\{\mu_-,\mu_+\}
					\subseteq
					\C\backslash\R
					\qqandqq
					\sr(A)=\max\{\vass{\mu_-},\vass{\mu_+}\}=\alpha ( 1 - 2\gamma \cK )<1.
				\end{equation}
			}
			\argument{\lref{eq:2nd:case:3};\lref{eq:location:distinguisher};}[verbs=ps]{\lref{it:sr} in the case $\frac{\pr{1-\alpha}^2}{2\pr{1+\alpha}^2}<\gamma \cK< \frac{1}{2}$}.
			In the next step we \prove[ps]\ \lref{it:sr} in the case
			\begin{equation}
				\llabel{eq:3rd:case}
				\frac{1}{2}\leq\gamma \cK< \frac{1+\alpha}{1+2\alpha}
				.
			\end{equation}
			\startnewargseq
			\argument{\lref{eq:3rd:case};}{that
				\begin{equation}
					\llabel{eq:3rd:case:1}
					\rbr*{ \frac{ ( 1 + \alpha )( 1 - 2\gamma \cK ) }{2} } ^2  - \alpha ( 1 - 2\gamma \cK ) 
					=
					\prbb{\frac{ ( 1 + \alpha )^2( 2\gamma \cK -1) }{4}  + \alpha }( 2\gamma \cK -1)\ge  0
					.
			\end{equation}}
			\argument{
				\lref{eq:EV};
				\lref{eq:3rd:case:1};
			}{that
				\begin{equation}
					\llabel{eq:3rd:case:4}
					\{\mu_-,\mu_+\}
					\subseteq
					\R
					\qqandqq
					\mu_-\leq\mu_+.
				\end{equation}
			}
			\argument{
				\lref{eq:3rd:case:4};
				\lref{eq:EV};
			}{that
				\begin{equation}
					\llabel{eq:3rd:case:5}
					\begin{split}
						\mu_{\pm}
						&=
						\frac{( 1 + \alpha )( 1 - 2\gamma \cK ) }{2}
						\pm \sqrt{ \rbr*{ \frac{ ( 1 + \alpha )( 1 - 2\gamma \cK ) }{2} } ^2 - \alpha ( 1 - 2\gamma \cK ) }
						\\&=
						\frac{( 1 + \alpha )( 1 - 2\gamma \cK ) }{2}
						\pm \sqrt{ \rbr*{ \frac{ ( 1 + \alpha )( 2\gamma \cK -1) }{2} } ^2 + \alpha (  2\gamma \cK -1) }
						\\&=
						\frac{( 1 + \alpha )( 1 - 2\gamma \cK ) }{2}
						\pm\sqrt{ \rbr*{ \frac{ ( 1 + \alpha )( 2\gamma \cK -1) }{2} } ^2 + (1+\alpha) (  2\gamma \cK -1)-(2\gamma\cK-1) }
						\\&=
						\frac{( 1 + \alpha )( 1 - 2\gamma \cK ) }{2}
						\pm \sqrt{ \rbr*{ \frac{ ( 1 + \alpha )( 2\gamma \cK -1) }{2} +1} ^2 -2\gamma\cK }
						.
					\end{split} 
				\end{equation}
			}
			\argument{
				\lref{eq:3rd:case};
				\lref{eq:3rd:case:5};}{that
				\begin{equation}
					\llabel{eq:3rd:case:5.5}
					\begin{split}
						\mu_+ 
						&=
						\frac{( 1 + \alpha )( 1 - 2\gamma \cK ) }{2}
						+ \sqrt{ \rbr*{ \frac{ ( 1 + \alpha )( 2\gamma \cK -1) }{2} +1} ^2 -2\gamma\cK }
						\\&<
						\frac{( 1 + \alpha )( 1 - 2\gamma \cK ) }{2}
						+ \sqrt{ \rbr*{ \frac{ ( 1 + \alpha )( 2\gamma \cK -1) }{2} +1} ^2 }
						\\&=
						\frac{( 1 + \alpha )( 1 - 2\gamma \cK ) }{2}
						+  \frac{ ( 1 + \alpha )( 2\gamma \cK -1) }{2} +1 =1
						.
					\end{split} 
				\end{equation}
			}
			\argument{\lref{eq:3rd:case};the fact that $\gamma\cK<\frac{1+\alpha}{1+2\alpha}$}{that
				\begin{equation}
					\llabel{eq:3rd:case:6}
					\begin{split}
						&\rbr*{ \frac{ ( 1 + \alpha )( 2\gamma \cK -1) }{2} +1} ^2 -2\gamma\cK 
						\\&=\frac{ ( 1 + \alpha )^2( 2\gamma \cK -1)^2 }{4}+ ( 1 + \alpha )( 2\gamma \cK -1) +1 -2\gamma\cK 
						\\&=\frac{ ( 1 + \alpha )^2( 2\gamma \cK -1)^2 }{4}
						- ( 1 + \alpha )( 2\gamma \cK -1)
						+ (1+2\alpha) ( 2\gamma \cK -1)  
						\\&<
						\frac{ ( 1 + \alpha )^2( 2\gamma \cK -1)^2 }{4}- ( 1 + \alpha )( 2\gamma \cK -1) +1 
						=\prbb{ \frac{ ( 1 + \alpha )( 2\gamma \cK-1 ) }{2}-1}^2
						.
					\end{split}
			\end{equation}}
			\argument{
				\lref{eq:3rd:case:6};
				\lref{eq:3rd:case:5};
			}{that
				\begin{equation}
					\llabel{eq:3rd:case:7}
					\begin{split}
						\mu_-
						&=	\frac{( 1 + \alpha )( 1 - 2\gamma \cK ) }{2}
						- \sqrt{ \rbr*{ \frac{ ( 1 + \alpha )( 2\gamma \cK -1) }{2} +1} ^2 -2\gamma\cK }
						\\&>
						\frac{( 1 + \alpha )( 1 - 2\gamma \cK ) }{2}
						- \sqrt{ \rbr*{ \frac{ ( 1 + \alpha )( 2\gamma \cK -1) }{2} -1} ^2 }
						\\&=
						\frac{( 1 + \alpha )( 1 - 2\gamma \cK ) }{2}
						- \vass*{\frac{ ( 1 + \alpha )( 2\gamma \cK-1 ) }{2}-1}
						\\&=
						\frac{( 1 + \alpha )( 1 - 2\gamma \cK ) }{2}
						+ \frac{ ( 1 + \alpha )( 2\gamma \cK-1 ) }{2}-1
						=-1
						.
					\end{split}
				\end{equation}
			}
			\argument{
				\lref{eq:EV:of:A:Nest:mat};
				\lref{eq:3rd:case:5.5};
				\lref{eq:3rd:case:7};
			}{
				\begin{equation}
					\llabel{eq:3rd:case:8}
					\sr(A)=\max\{\vass{\mu_-},\vass{\mu_+}\}<1.
			\end{equation}}
			\argument{\lref{eq:3rd:case:8};}[verbs=ps]{\lref{it:sr} in the case $\frac{1}{2}\leq\gamma \cK< \frac{1+\alpha}{1+2\alpha}$}.
			In the next step we \prove[ps]\ \lref{it:sr} in the case
			\begin{equation}
				\llabel{eq:4th:case}
				\gamma \cK=\frac{1+\alpha}{1+2\alpha}
				.
			\end{equation}
			\argument{
				\lref{eq:4th:case};
				\lref{eq:EV};}{that
				\begin{equation}
					\begin{split}
						\llabel{eq:4th:case:1}
						\mu_\pm 
						&=	\frac{( 1 + \alpha )( 1 - 2\gamma \cK ) }{2}
						\pm \sqrt{ \rbr*{ \frac{ ( 1 + \alpha )( 1 - 2\gamma \cK ) }{2} } ^2 - \alpha ( 1 - 2\gamma \cK ) }
						\\&=\frac{ -1 - \alpha  }{2( 1 +2\alpha )}
						\pm \sqrt{ \rbr*{\frac{ 1 + \alpha  }{2( 1 +2\alpha )}} ^2 + \frac{\alpha}{1+2\alpha} }
						\\&=\frac{ -1 - \alpha  }{2( 1 +2\alpha )}
						\pm \sqrt{ \frac{ (1 + \alpha)^2  }{4(1+2\alpha)^2} + \frac{4\alpha(1+2\alpha)}{4(1+2\alpha)^2} }
						\\&=\frac{ -1 - \alpha  }{2( 1 +2\alpha )}
						\pm \sqrt{ \frac{ 1+2\alpha+\alpha^2+4\alpha+8\alpha^2 }{4(1+2\alpha)^2}  }
						\\&=\frac{ -1 - \alpha  }{2( 1 +2\alpha )}
						\pm \sqrt{ \frac{ 1+6\alpha+9\alpha^2 }{4(1+2\alpha)^2}  }
						\\&=\frac{ -1 - \alpha  }{2( 1 +2\alpha )}
						\pm \sqrt{ \frac{ (1+3\alpha)^2 }{4(1+2\alpha)^2}  }
						=\frac{ -1 - \alpha  }{2( 1 +2\alpha )}
						\pm  \frac{ 1+3\alpha }{2(1+2\alpha)}  
						.
					\end{split}
			\end{equation}}
			\argument{\lref{eq:4th:case:1};\lref{eq:EV:of:A:Nest:mat}}{that
				\begin{equation}
					\llabel{eq:4th:case:2}
					\mu_-=-1,
					\qquad
					\mu_+=\tfrac{\alpha}{1+2\alpha}<1,
					\qqandqq
					\sr(A)=\max\pRb{1,\tfrac{\alpha}{1+2\alpha}}=1
					.
			\end{equation}}
			\argument{\lref{eq:4th:case:2};}[verbs=ps]{\lref{it:sr} in the case $\gamma \cK=\frac{1+\alpha}{1+2\alpha}$}.
			In the next step we \prove[ps]\ \lref{eq:cond:for:result:3} in the case 
			\begin{equation}
				\llabel{eq:5th:case}
				\gamma \cK>\frac{1+\alpha}{1+2\alpha}
				.
			\end{equation}
			\startnewargseq
			\argument{\lref{eq:5th:case}; and the fact that $\frac{1+\alpha}{1+2\alpha}>\frac{1}{2}$;}{that
				\begin{equation}
					\llabel{eq:5th:case:2}
					\rbr*{ \frac{ ( 1 + \alpha )( 1 - 2\gamma \cK ) }{2} } ^2  - \alpha ( 1 - 2\gamma \cK ) 
					=
					\prbb{\frac{ ( 1 + \alpha )^2( 2\gamma \cK -1) }{4}  + \alpha }( 2\gamma \cK -1)>  0
					.
			\end{equation}}
			\argument{
				\lref{eq:EV};
				\lref{eq:EV:of:A:Nest:mat};
				\lref{eq:5th:case:2};
			}{that
				\begin{equation}
					\llabel{eq:5th:case:4}
					\{\mu_-,\mu_+\}
					\subseteq
					\R
					\qqandqq
					\mu_-<\mu_+.
				\end{equation}
			}
			\argument{
				\lref{eq:5th:case:4};
				\lref{eq:EV};
			}{that
				\begin{equation}
					\llabel{eq:5th:case:5}
					\begin{split}
						\mu_{\pm}
						&=
						\frac{( 1 + \alpha )( 1 - 2\gamma \cK ) }{2}
						\pm \sqrt{ \rbr*{ \frac{ ( 1 + \alpha )( 1 - 2\gamma \cK ) }{2} } ^2 - \alpha ( 1 - 2\gamma \cK ) }
						\\&=
						\frac{( 1 + \alpha )( 1 - 2\gamma \cK ) }{2}
						\pm \sqrt{ \rbr*{ \frac{ ( 1 + \alpha )( 2\gamma \cK -1) }{2} } ^2 + \alpha (  2\gamma \cK -1) }
						\\&=
						\frac{( 1 + \alpha )( 1 - 2\gamma \cK ) }{2}
						\pm\sqrt{ \rbr*{ \frac{ ( 1 + \alpha )( 2\gamma \cK -1) }{2} } ^2 + (1+\alpha) (  2\gamma \cK -1)-(2\gamma\cK-1) }
						\\&=
						\frac{( 1 + \alpha )( 1 - 2\gamma \cK ) }{2}
						\pm \sqrt{ \rbr*{ \frac{ ( 1 + \alpha )( 2\gamma \cK -1) }{2} +1} ^2 -2\gamma\cK }
						.
					\end{split} 
				\end{equation}
			}
			\argument{
				\lref{eq:5th:case};
				\lref{eq:5th:case:5};}{that
				\begin{equation}
					\llabel{eq:5th:case:5.5}
					\begin{split}
						\mu_+ 
						&=
						\frac{( 1 + \alpha )( 1 - 2\gamma \cK ) }{2}
						+ \sqrt{ \rbr*{ \frac{ ( 1 + \alpha )( 2\gamma \cK -1) }{2} +1} ^2 -2\gamma\cK }
						\\&<
						\frac{( 1 + \alpha )( 1 - 2\gamma \cK ) }{2}
						+ \sqrt{ \rbr*{ \frac{ ( 1 + \alpha )( 2\gamma \cK -1) }{2} +1} ^2 }
						\\&=
						\frac{( 1 + \alpha )( 1 - 2\gamma \cK ) }{2}
						+  \frac{ ( 1 + \alpha )( 2\gamma \cK -1) }{2} +1 =1
						.
					\end{split} 
				\end{equation}
			}
			\argument{\lref{eq:5th:case};the fact that $\gamma\cK>\frac{1+\alpha}{1+2\alpha}$}{that
				\begin{equation}
					\llabel{eq:5th:case:6}
					\begin{split}
						&\rbr*{ \frac{ ( 1 + \alpha )( 2\gamma \cK -1) }{2} +1} ^2 -2\gamma\cK 
						\\&=\frac{ ( 1 + \alpha )^2( 2\gamma \cK -1)^2 }{4}+ ( 1 + \alpha )( 2\gamma \cK -1) +1 -2\gamma\cK 
						\\&=\frac{ ( 1 + \alpha )^2( 2\gamma \cK -1)^2 }{4}
						- ( 1 + \alpha )( 2\gamma \cK -1)
						+ (1+2\alpha) ( 2\gamma \cK -1)  
						\\&>
						\frac{ ( 1 + \alpha )^2( 2\gamma \cK -1)^2 }{4}- ( 1 + \alpha )( 2\gamma \cK -1) +1 
						=\prbb{ \frac{ ( 1 + \alpha )( 2\gamma \cK-1 ) }{2}-1}^2
						.
					\end{split}
			\end{equation}}
			\argument{
				\lref{eq:5th:case:6};
				\lref{eq:5th:case:5};
			}{that
				\begin{equation}
					\llabel{eq:5th:case:7}
					\begin{split}
						\mu_-
						&=	\frac{( 1 + \alpha )( 1 - 2\gamma \cK ) }{2}
						- \sqrt{ \rbr*{ \frac{ ( 1 + \alpha )( 2\gamma \cK -1) }{2} +1} ^2 -2\gamma\cK }
						\\&<
						\frac{( 1 + \alpha )( 1 - 2\gamma \cK ) }{2}
						- \sqrt{ \rbr*{ \frac{ ( 1 + \alpha )( 2\gamma \cK -1) }{2} -1} ^2 }
						\\&=
						\frac{( 1 + \alpha )( 1 - 2\gamma \cK ) }{2}
						- \vass*{\frac{ ( 1 + \alpha )( 2\gamma \cK-1 ) }{2}-1}
						\\&=
						\frac{( 1 + \alpha )( 1 - 2\gamma \cK ) }{2}
						+ \frac{ ( 1 + \alpha )( 2\gamma \cK-1 ) }{2}-1
						=-1
						.
					\end{split}
				\end{equation}
			}
			\argument{
				\lref{eq:EV:of:A};
				\lref{eq:5th:case:7};
				\lref{eq:5th:case:6};}{that
				\begin{equation}
					\llabel{eq:5th:case:8}
					\begin{split}
						\vass{\mu_-}>1,
						\qquad
						\vass{\mu_-}>\vass{\mu_+},
						\qqandqq
						\sr(A)=\max\{\vass{\mu_-},\vass{\mu_+}\}>1
						.
					\end{split}
			\end{equation}}
			\argument{\lref{eq:5th:case:8};}[verbs=ps]{\lref{it:sr} in the case $\gamma \cK>\frac{1+\alpha}{1+2\alpha}$}.
		\end{aproof}
		
		\subsection{Asymptotic analysis for the Nesterov optimizer}
		
		\begin{prop}
			\label{prop:nesterov-stable}
			Let 
			$\fd \in \N$,
			$\lp \in \R^\fd$,
			$\gamma, \cK \in (0, \infty)$,
			$\alpha \in (0, 1 )$,
			let
			$\fl \colon \R^\fd \to \R$ satisfy for all $\theta \in \R^\fd$ that
			\begin{equation}
				\fl ( \theta ) = \cK\norm{\theta - \lp} ^2,
			\end{equation}
			and let $\Theta\colon \N_0 \to \R^\fd$ and $\m \colon \N_0 \to \R^\fd$ satisfy for all $n \in \N$ that
			\begin{equation}
				\llabel{eq:nest:setup}
				\m_n = (1 + \alpha ) \Theta_n - \alpha \Theta_{n-1},
				\qquad
				\Theta_n = \m_{n-1} - \gamma ( \nabla \fl ) ( \m_{n-1} ),
				\qqandqq
				\Theta_0 \ne \lp.
			\end{equation}
			Then
			\begin{equation}
				\llabel{eq:prop:Nest:res}
				\liminf_{n \to \infty} \norm{\Theta_n - \lp} = \limsup_{n \to \infty} \norm{\Theta_n - \lp } 
				\in
				\begin{cases}
					\{0\} & \colon \gamma\cK < \frac{1 +\alpha}{1+2\alpha} \\
					( 0 , \infty ) & \colon \gamma\cK = \frac{1 +\alpha}{1+2\alpha} \\
					\{\infty\} & \colon \gamma\cK >  \frac{1 +\alpha}{1+2\alpha}  .
				\end{cases}
			\end{equation}
		\end{prop}

		\begin{cproof}{prop:nesterov-stable}
			Throughout this proof assume without loss of generality that $\fd = 1$ and $\lp = 0$
			and let $A\in\C^{2\times 2}$, $\mu_-,\mu_+\in\C$ satisfy
			\begin{equation}
				\llabel{eq:matrix:3}
				A = \begin{pmatrix}
					(1+\alpha)(1-2\gamma\cK) & -\alpha \\
					1 - 2 \gamma\cK & 0
				\end{pmatrix}
			\end{equation}
			\begin{equation}
				\llabel{eq:EV:Nest}
				\qqandqq
				\mu_\pm =\frac{( 1 + \alpha )( 1 - 2\gamma \cK ) }{2}
				\pm \sqrt{ \rbr*{ \frac{ ( 1 + \alpha )( 1 - 2\gamma \cK ) }{2} } ^2 - \alpha ( 1 - 2\gamma \cK ) }
				.
			\end{equation}
			\argument{\lref{eq:matrix:3};\lref{eq:EV:Nest}; \cref{lem:EV:Nest} (applied with
				$\gamma\curvearrowleft\gamma$,
				$\cK\curvearrowleft\cK$,
				$\alpha\curvearrowleft\alpha$,
				$A\curvearrowleft A$,
				$\mu-\curvearrowleft\mu_-$,
				$\mu_+\curvearrowleft\mu_+$ in the notation of \cref{lem:EV:Nest})}{that
				\begin{equation}
					\llabel{eq:EV:of:A}
					\cu{ z \in \C \colon \br{ \exists \, \nu \in \C^2 \backslash \cu{0} \colon A \nu = z \nu } }
					=
					\{\mu_-,\mu_+\},
				\end{equation}
				\begin{equation}
					\llabel{eq:cond:for:result:4}
					\max\{1,\vass{\mu_-}\}>\vass{\mu_+},
					\qqandqq
					\sr(A)\in
					\begin{cases}
						[0,1) &\colon \gamma \cK<\frac{1 +\alpha}{1+2\alpha}\\
						\{1\} &\colon \gamma \cK=\frac{1 +\alpha}{1+2\alpha}\\
						(1,\infty) &\colon \gamma \cK>\frac{1 +\alpha}{1+2\alpha}.\\
					\end{cases}
			\end{equation}}
			\argument{\lref{eq:nest:setup};the fact that $\forall \, \theta \in \R \colon (\nabla \fl) ( \theta ) = 2\cK \theta$}{for all $n\in\N$ that
				\begin{equation}
					\llabel{eq:theta:recursion}
					\begin{split}
						\Theta_n = \m_{n-1} - \gamma ( \nabla \fl ) ( \m_{n-1} )
						= \m_{n-1} - 2\gamma\cK  \m_{n-1} 
						= (1-2\gamma\cK)\m_{n-1}
						.
					\end{split}
			\end{equation}}
			\argument{\lref{eq:theta:recursion};\lref{eq:nest:setup}}{for all $n\in\N$ that
				\begin{equation}
					\llabel{eq:mom:recursion}
					\begin{split}
						\m_n 
						= (1 + \alpha ) \Theta_n - \alpha \Theta_{n-1}
						&= (1+\alpha)(1-2\gamma\cK)\m_{n-1}- \alpha \Theta_{n-1}
						.
					\end{split}
			\end{equation}}
			\argument{
				\lref{eq:mom:recursion};
				\lref{eq:nest:setup};
				\lref{eq:matrix:3};
				\lref{eq:theta:recursion};
			}{that for all $n \in \N$ it holds that
				\begin{equation}
					\llabel{eq:recursion:with:matrix}
					\begin{pmatrix}
						\m_n \\ \Theta_n 
					\end{pmatrix}
					= A \begin{pmatrix}
						\m_{n-1} \\
						\Theta_{n-1} 
					\end{pmatrix}.
				\end{equation}
			}
			\argument{\lref{eq:recursion:with:matrix};}{for all $n \in \N$ that
				\begin{equation}
					\llabel{eq:rec}
					\begin{pmatrix}
						\m_n \\ \Theta_n 
					\end{pmatrix}
					= A ^n \begin{pmatrix}
						\m_{0} \\
						\Theta_{0} 
					\end{pmatrix}.
				\end{equation}
			}
			\argument{
				\lref{eq:matrix:3};}{that for all $\lambda,x\in\C$ it holds that
				\begin{equation}
					\llabel{eq:no:EV}					
					\norm[\bigg]{A\begin{pmatrix}
							x \\ 0
						\end{pmatrix}
						-
						\lambda\begin{pmatrix}
							x \\ 0
					\end{pmatrix}} 
					=
					\norm[\bigg]{\begin{pmatrix}
							(1+\alpha)(1-2\gamma\cK)x-\lambda x \\ (1-2\gamma\cK)x
						\end{pmatrix}
					}
					\geq
					\vass{1-2\gamma\cK}\vass{x}
					.
			\end{equation}}
			\argument{
				\lref{eq:cond:for:result:4};
				the fact that $\tfrac{2(1+\alpha)}{1+2\alpha}>1$
			}{that for all $x\in\C\backslash\{0\}$ it holds that
				\begin{equation}
					\llabel{eq:non:van:sec:coord}
					\begin{split}
						\indicator{[0,1]}\pr{\sr\pr{A}}+
						\indicator{(1,\infty)}\pr{\sr\pr{A}}\vass{1-2\gamma\cK}\vass{x}
						&\geq
						\indicator{[0,1]}\pr{\sr\pr{A}}+
						\indicator{(1,\infty)}\pr{\sr\pr{A}}\vass{1-\tfrac{2(1+\alpha)}{1+2\alpha}}\vass{x}
						\\&=
						\indicator{[0,1]}\pr{\sr\pr{A}}+
						\indicator{(1,\infty)}\pr{\sr\pr{A}}\prb{\tfrac{\alpha}{1+2\alpha}}\vass{x}
						>0
						.
					\end{split}
			\end{equation}}
			\argument{
				\lref{eq:EV:of:A}}{that there exist $v=(v_1,v_2)$, $w=(w_1,w_2)\in\C^2\backslash\{0\}$ which satisfy that
				\begin{equation}
					\begin{split}
						\llabel{eq:start:EV:repr}
						\begin{pmatrix}
							\m_0 \\ \Theta_0
						\end{pmatrix}
						=v+w,
						\qquad
						Av=\mu_+v,
						\qqandqq
						Aw=\mu_-w
						.
					\end{split}
				\end{equation}
			}
			\argument{
				\lref{eq:rec};
				\lref{eq:start:EV:repr};}{that for all $n\in\N$ it holds that
				\begin{equation}
					\llabel{eq:n_th:vec}
					\begin{pmatrix}
						\m_n \\ \Theta_n 
					\end{pmatrix}
					=
					A^n\begin{pmatrix}
						\m_0 \\ \Theta_0
					\end{pmatrix}
					=
					A^n\pr{v+w}
					.
			\end{equation}}
			\argument{
				\lref{eq:start:EV:repr};
				\lref{eq:non:van:sec:coord};
				\lref{eq:no:EV};}{that
				\begin{equation}
					\llabel{eq:non:van:sec:co}
					\indicator{[0,1]}\pr{\sr\pr{A}}+\vass{w_2}>0.
			\end{equation}}
			\argument{
				\lref{eq:n_th:vec};
				\lref{eq:non:van:sec:co};
				\lref{eq:cond:for:result:4};
				\lref{eq:start:EV:repr};
				\cref{lem:mat:seq:EV} (applied with 
				$\mu_1\curvearrowleft\mu_-$,
				$\mu_2\curvearrowleft\mu_+$,
				$v_1\curvearrowleft w$,
				$v_2\curvearrowleft v$,
				$A\curvearrowleft A$,
				$\m\curvearrowleft\m$,
				$\Theta\curvearrowleft\Theta$,
				in the notation of \cref{lem:mat:seq:EV})}{
				\begin{equation}
					\llabel{eq:result:appl}
					\liminf_{n \to \infty} \vass{\Theta_n} = \limsup_{n \to \infty} \vass{\Theta_n } 
					=
					\begin{cases}
						0 & \colon \sr\pr{A}<1\\
						\vass{w_{2}} & \colon \sr\pr{A}=1\\
						\infty & \colon \sr\pr{A}>1.
					\end{cases}
				\end{equation}
			}
			\argument{
				\lref{eq:result:appl};
				\lref{eq:non:van:sec:co};}{\lref{eq:prop:Nest:res}}.
		\end{cproof}
		
		\begin{athm}{cor}{prop:Nest-stable:diag}
			Let $\fd \in \N$,
			$\gamma, \lambda_1,\lambda_2,\dots,\lambda_d \in [0, \infty)$,
			$\alpha \in (0, 1 )$,
			$\lp=(\lp_1,\dots,\lp_d) \in \R^\fd$,
			let
			$\fl \colon \R^\fd \to \R$ satisfy for all $\theta=(\theta_1,\dots,\theta_d) \in \R^\fd$ that
			\begin{equation}
				\llabel{eq:loss:diag:mom}
				\fl ( \theta ) = \textstyle\sum_{i=1}^d\lambda_i\pr{\theta_i - \lp_i}^2,
			\end{equation}
			and let $\Theta \colon \N_0 \to \R^\fd$
			and
			$ \m \colon \N_0 \to \R^\fd$ satisfy for all $n \in \N$ that
			\begin{equation}
				\llabel{eq:start:and:setup:diag}
				\m_n = (1 + \alpha ) \Theta_n - \alpha \Theta_{n-1},
				\qquad
				\Theta_n = \m_{n-1} - \gamma ( \nabla \fl ) ( \m_{n-1} ),
				\qqandqq
				\Theta_0 \ne \lp.
			\end{equation}
			Then it holds that $\sup_{n \in\N} \norm{\Theta_n }<\infty$ if and only if 
			$
			\gamma\max\{\lambda_1,\lambda_2,\dots,\lambda_d\} \leq \tfrac{  1 + \alpha  }{  1 +2 \alpha}
			$.
		\end{athm}
		
		\begin{aproof}
			Throughout this proof for every $i\in\{1,2,\dots,d\}$ let $p_i\colon\R^d\to\R$ satisfy for all $x=(x_1,\dots,x_d)\in\R^d$ that 
			\begin{equation}
				\llabel{eq:def:proj}
				p_i(x)=x_i.
			\end{equation}
			\argument{
				\lref{eq:loss:diag:mom};
				\lref{eq:start:and:setup:diag};
				\lref{eq:def:proj}}{that for all $i\in\{1,2,\dots,d\}$, $n\in\N$ with $\gamma\lambda_i=0$ it holds that
				\begin{equation}
					\llabel{eq:vanishing:case}
					\begin{split}
						p_i(\Theta_n)
						=p_i(\Theta_{n-1}-\gamma( \nabla \fl ) ( \m_{n-1} ))
						&=p_i(\Theta_{n-1})-\gamma p_i(( \nabla \fl )(\m_{n-1}))
						\\&=p_i(\Theta_{n-1})-2\gamma \lambda_i p_i(\m_{n-1}-\lp)
						=p_i(\Theta_{n-1}).					\end{split}
			\end{equation}}
			\argument{\cref{prop:nesterov-stable} (applied with 
				$d\curvearrowleft 1$,
				$\gamma\curvearrowleft \gamma$,
				$\cK\curvearrowleft \lambda_i$,
				$\alpha\curvearrowleft \alpha$,
				$\lp\curvearrowleft \lp_i$,
				$\fl\curvearrowleft \pr{\R\ni\theta\mapsto \lambda_i(\theta-\lp_i)^2\in\R}$,
				$\Theta\curvearrowleft \pr{\N_0\ni n\mapsto p_i(\Theta_n)\in\R}$,
				$\m\curvearrowleft \pr{\N_0\ni n\mapsto p_i(\m_n)\in\R}$ for $i\in\{1,2,\dots,d\}$
				in the notation of \cref{prop:nesterov-stable})}{that for all $i\in\{1,2,\dots,d\}$ with $\gamma\lambda_i>0$ it holds that
				\begin{equation}
					\llabel{eq:result:prop:mom}
					\liminf_{n \to \infty} \vass{p_i(\Theta_n - \lp)} 
					= \limsup_{n \to \infty} \vass{p_i(\Theta_n - \lp) } 
					\in
					\begin{cases}
						\{0\} & \colon \gamma\lambda_i < \tfrac{  1 + \alpha  }{  1 +2 \alpha} \\
						( 0 , \infty ) & \colon \gamma\lambda_i =  \tfrac{ 1 + \alpha  }{ 1 +2 \alpha  }\\
						\{\infty\} & \colon \gamma\lambda_i >  \tfrac{1 + \alpha }{   1 +2 \alpha }.
					\end{cases}
			\end{equation}}
			\argument{
				\lref{eq:vanishing:case};
				\lref{eq:result:prop:mom};}{that 
				\begin{equation}
					\llabel{eq:result}
					\textstyle
					\sup_{n\in\N}\norm{\Theta_n}^2
					=
					\textstyle
					\sup_{n\in\N}\prb{\sum_{i=1}^d \vass{p_i(\Theta_n)}^2}
					\in
					\begin{cases}
						[0,\infty)
						& \colon \gamma\max\{\lambda_1,\lambda_2,\dots,\lambda_d\} \leq \tfrac{  1 + \alpha  }{  1 +2 \alpha} \\
						\{\infty\} 
						& \colon \gamma\max\{\lambda_1,\lambda_2,\dots,\lambda_d\} >  \tfrac{1 + \alpha }{   1 +2 \alpha }.
					\end{cases}
			\end{equation}}
			\argument{\lref{eq:result}}{that $\sup_{n \in\N} \norm{\Theta_n }<\infty$ if and only if 
				$
				\gamma\max\{\lambda_1,\lambda_2,\dots,\lambda_d\} \leq \tfrac{  1 + \alpha  }{  1 +2 \alpha}
				$}.	\end{aproof}

		\section{Asymptotical stability for gradient based methods}
		\label{section:5}
		
		In this section we combine the findings from 
		\cref{section:Aprioribounds,section:apriori:momentum,section:apriori:Nest} above to explicitly specify in
	 \cref{final:thm:1.1} below the stability region (see \cref{subsection:1:1} above) for 
	 the Nesterov optimizer (cf.\ \cref{prop:Nest-stable:diag}),
	 the \GD\ optimizer (cf., \eg, \cite[Theorem~6.1.12]{jentzen2023mathematical}),
	 the momentum optimizer (cf.\ \cref{prop:mom-stable:diag}),
	 the \Adam\ optimizer (cf.\ \cref{lem:ADAM:bounded}), and the
	 \RMSprop\ optimizer (cf.\ \cref{lem:ADAM:bounded}).
	 \Cref{intro:thm:1.1} in \cref{subsection:1:1} is an immediate consequence of \cref{final:thm:1.1}.
	 
	 In \cref{def:asymptotically:stable} below we introduce the notion of an asymptotically stable optimizer and in \cref{def:uniformly:stable,def:strongly:asymptotically:stable,def:super:strongly:asymptotically:stable} below we present related stability notions for optimization methods.
	 In the elementary results in \cref{stab:def:relation} and \cref{lem:basic:stable:properties} below we present basic properties and relations between these stability concepts.

		\subsection{Introduction of asymptotic stability}
		\label{subsection:5.1}
		
		\begin{definition}[Asymptotically stable]
			\label{def:asymptotically:stable}
			Let $d \in \N$,
			let $\Phi_n\colon(\R^d)^{n}\to\R^d$, $n \in \N$, be functions,
			and let $\mathcal{A}\subseteq [0,\infty)^{d+1}$ be a set.
			Then we say that $(\Phi_n)_{n\in \N}$ is $\mathcal{A}$-asymptotically stable if and only if it holds
		 	for every 
			$\gamma\in\R$,
			$\vartheta,\lambda=(\lambda_1,\dots,\lambda_d)\in\R^d$
			and every 
			$\Theta\colon\N_0  \to \R^d$ with
			$(\gamma,\lambda_1,\dots,\lambda_d)\in \mathcal{A}$
			and
			\begin{equation}
				\label{eq:asymptotically:stable}
				\forall\, n \in \N\colon \Theta_{n}=\Theta_{n-1}
				-
				\gamma\Phi_n\prb{
				\operatorname{diag}(\lambda)(\Theta_0-\vartheta),
				\operatorname{diag}(\lambda)(\Theta_1-\vartheta),
				\dots,
				\operatorname{diag}(\lambda)(\Theta_{n-1}-\vartheta)}
			\end{equation}
		that
	$\limsup_{n\to \infty} \norm{\Theta_n}< \infty$.
		\end{definition}
		
		\cfclear
		\begin{athm}{lemma}{stab:def:relation}
			Let $d \in \N$,
			let $\Phi_n\colon(\R^d)^{n}\to\R^d$, $n \in \N$, be functions,
			let $B=\{\mathcal{B}\subseteq [0,\infty)^{d+1}\colon \allowbreak\pr{\Phi_n}_{n\in\N}\text{ is $\mathcal{B}$-asymptotically stable\cfadd{def:asymptotically:stable}}\}$
			and let $\mathcal{A}\subseteq [0,\infty)^{d+1}$ satisfy
			\begin{equation}
				\llabel{eq:stab:def:rel}
				\mathcal{A}=\cup_{\mathcal{B}\in B}\mathcal{B}
			\end{equation}
				\cfload.
			Then the stability region of $\pr{\Phi_n}_{n\in\N}$ is $\mathcal{A}$\cfadd{def:stab:region} \cfout.
			\end{athm}
			
			\begin{aproof}
				Throughout this proof let $\mathcal{C}\subseteq[0,\infty)^{d+1}$ satisfy that the stability region of $\pr{\Phi_n}_{n\in\N}$ is $\mathcal{C}$\cfadd{def:stab:region} \cfload.
				\argument{\eqref{eq:asymptotically:stable}}{that for
				all 
				$\mathcal{B}\subseteq[0,\infty)^{d+1}$, 
				$\gamma\in[0,\infty)$, 
				$\lambda=(\lambda_1,\dots,\lambda_d)\in[0,\infty)^d$,
				$\vartheta\in\R^d$ with
				$(\gamma,\lambda_1,\dots,\lambda_d)\in\mathcal{B}$ and
				$\pr{\Phi_n}_{n\in\N}$ \stable{\mathcal{B}}and all 
				$\Theta\colon\N_0  \to \R^d$ with
				\begin{equation}
					\label{eq:stab:reg:1}
					\forall\, n \in \N\colon \Theta_{n}=\Theta_{n-1}
					-
					\gamma\Phi_n\prb{
						\operatorname{diag}(\lambda)(\Theta_0-\vartheta),
						\operatorname{diag}(\lambda)(\Theta_1-\vartheta),
						\dots,
						\operatorname{diag}(\lambda)(\Theta_{n-1}-\vartheta)}
				\end{equation}
				it holds that
				\begin{equation}
					\llabel{eq:stab:reg:2}
					\limsup_{n\to \infty} \norm{\Theta_n}
					<\infty.
				\end{equation}
			}
			\argument{
				\lref{eq:stab:reg:2};
				\lref{eq:stab:def:rel}}{$\mathcal{A}\subseteq \mathcal{C}$}.
			\argument{
				\eqref{eq:def:stab:reg:1};
				\eqref{eq:def:stab:reg:2}}{that for
				all 
				$\gamma\in[0,\infty)$, 
				$\lambda=(\lambda_1,\dots,\lambda_d)\in[0,\infty)^d$,
				$\vartheta\in\R^d$ with
				$(\gamma,\lambda_1,\dots,\lambda_d)\in\mathcal{C}$ and all 
				$\Theta\colon\N_0  \to \R^d$ with
				\begin{equation}
					\label{eq:stab:reg:1.2}
					\forall\, n \in \N\colon \Theta_{n}=\Theta_{n-1}
					-
					\gamma\Phi_n\prb{
						\operatorname{diag}(\lambda)(\Theta_0-\vartheta),
						\operatorname{diag}(\lambda)(\Theta_1-\vartheta),
						\dots,
						\operatorname{diag}(\lambda)(\Theta_{n-1}-\vartheta)}
				\end{equation}
				it holds that
				\begin{equation}
					\llabel{eq:stab:reg:2.2}
					\limsup_{n\to \infty} \norm{\Theta_n}
					<\infty.
					\end{equation}}
					\argument{
						\lref{eq:stab:reg:2.2};
						\lref{eq:stab:def:rel}}{$\mathcal{C}\subseteq \mathcal{A}$}.
			\end{aproof}

	\begin{definition}[Uniformly stable]
			\label{def:uniformly:stable}
			Let $d \in \N$ and
			let $\Phi_n\colon(\R^d)^{n}\to\R^d$, $n \in \N$, be functions.
			Then we say that $(\Phi_n)_{n\in \N}$ is uniformly stable if and only if it holds that
			$(\Phi_n)_{n\in \N}$ \stable{[0,\infty)^{d+1}}\cfload.
	\end{definition}

			\begin{definition}[Strongly asymptotically stable]
			\label{def:strongly:asymptotically:stable}
			Let $d \in \N$,
			let $\Phi_n\colon(\R^d)^{n}\to\R^d$, $n \in \N$, be functions,
			and let $\mathcal{A}\subseteq [0,\infty)^{d+1}$ be a set.
			Then we say that $(\Phi_n)_{n\in \N}$ is strongly $\mathcal{A}$-asymptotically stable if and only if it holds
			for every
			$\vartheta,\lambda=(\lambda_1,\dots,\lambda_d)\in\R^d$,
			every
			$\gamma
			\colon\N\to\R$, 
			and every
			$\Theta\colon\N_0 \to \R^d$ with
			$\cup_{n\in\N}\{(\gamma_n,\lambda_1,\dots,\lambda_d)\}\subseteq \mathcal{A}$
			and
			\begin{equation}
				\label{eq:strongly:asymptotically:stable}
				\forall\, n \in \N\colon \Theta_{n}=\Theta_{n-1}
				-
				\gamma_n\Phi_n(
				\operatorname{diag}(\lambda)(\Theta_0-\vartheta),
				\operatorname{diag}(\lambda)(\Theta_1-\vartheta),
				\dots,
				\operatorname{diag}(\lambda)(\Theta_{n-1}-\vartheta))
			\end{equation}
			that
			$\limsup_{n\to \infty} \norm{\Theta_n}< \infty$.
		\end{definition}
		
			\begin{definition}[Super strongly asymptotically stable]
			\label{def:super:strongly:asymptotically:stable}
			Let $d \in \N$,
			let $\Phi_n\colon(\R^d)^{n}\to\R^d$, $n \in \N$, be functions,
			and let $\mathcal{A}\subseteq [0,\infty)^{d+1}$ be a set.
			Then we say that $(\Phi_n)_{n\in \N}$ is super strongly $\mathcal{A}$-asymptotically stable if and only if it holds
			for every
			$\vartheta\in\R^d$,
			every
			$\gamma\colon\N\to\R$,
			every
			$\lambda=(\lambda_n)_{n\in\N}=((\lambda_{n,1},\dots,\lambda_{n,d}))_{n\in\N}\colon\N\to\R^d$, and every
			$\Theta\colon\N_0  \to \R^d$ with
			$\cup_{n\in\N}\{(\gamma_n,\lambda_{n,1},\dots,\lambda_{n,d})\}\allowbreak\subseteq \mathcal{A}$
			and
			\begin{equation}
				\label{eq:super:strongly:asymptotically:stable}
				\forall\, n \in \N\colon \Theta_{n}=\Theta_{n-1}
				-
				\gamma_n\Phi_n(
				\operatorname{diag}(\lambda_1)(\Theta_0-\vartheta),
				\operatorname{diag}(\lambda_2)(\Theta_1-\vartheta),
				\dots,
				\operatorname{diag}(\lambda_n)(\Theta_{n-1}-\vartheta))
			\end{equation}
			that
			$\limsup_{n\to \infty} \norm{\Theta_n}< \infty$.
		\end{definition}

		\begin{athm}{lemma}{lem:basic:stable:properties}
			Let $d\in\N$,
			for every $n\in\N$
			let $\mathcal{A}_n\subseteq [0,\infty)^{d+1}$ be a set, and
			let $\Phi_n\colon(\R^d)^{n}\to\R^d$, $n \in \N$, satisfy 
			\begin{equation}
				\label{eq:asympt}
				\mathcal{A}_1=\pRb{\mathcal{B}\subseteq [0,\infty)^{d+1}\colon (\Phi_n)_{n \in \N} \text{ \stable{\mathcal{B}}}\!\!},
				\vspace{-0.15cm}
			\end{equation} 
			\begin{equation}
				\label{eq:sasympt}
				\mathcal{A}_2=\pRb{\mathcal{B}\subseteq [0,\infty)^{d+1}\colon (\Phi_n)_{n \in \N} \text{ \sstable{\mathcal{B}}}\!\!},
			\end{equation} 
			\begin{equation}
				\label{eq:ssasympt}
				\qandq
				\mathcal{A}_3=\pRb{\mathcal{B}\subseteq [0,\infty)^{d+1}\colon (\Phi_n)_{n \in \N} \text{ \ssstable{\mathcal{B}}}\!\!}
			\end{equation} 
			\cfload. Then
			\begin{enumerate}[label=(\roman*)]
				\item\label{it:1} it holds that $\mathcal{A}_3\subseteq \mathcal{A}_2\subseteq\mathcal{A}_1$,
				\item\label{it:2} it holds that $\emptyset\in\mathcal{A}_3$,
				\item\label{it:closed:subsets} it holds for all $i\in\{1,2,3\}$, $\mathcal{B}\in\mathcal{A}_i$, $\mathcal{C}\subseteq\mathcal{B}$ that $\mathcal{C}\in\mathcal{A}_i$, 
				\item\label{it:4} it holds for all 
				$\mathcal{B},\mathcal{C}\in\mathcal{A}_1$ that $\pr{\mathcal{B}\cup\mathcal{C}}\in\mathcal{A}_1$, and
				\item\label{it:5} it holds for all $i\in\{1,2\}$, $\mathcal{B}\in\mathcal{A}_i$, $\lambda\in\mathcal{B}$ that $\{\lambda\}\in\mathcal{A}_{i+1}$.
			\end{enumerate}
			\end{athm}
			
			\begin{aproof}
				\argument{
					\cref{eq:asymptotically:stable};
					\cref{eq:strongly:asymptotically:stable};
					\cref{eq:asympt};
					\cref{eq:sasympt}}{that 
					\begin{equation}
						\llabel{eq:it1:first:incl}
						\mathcal{A}_2\subseteq\mathcal{A}_1.
					\end{equation}}
				\argument{
					\cref{eq:strongly:asymptotically:stable};
					\cref{eq:super:strongly:asymptotically:stable};
					\cref{eq:sasympt};
					\cref{eq:ssasympt}}{that 
					\begin{equation}
						\llabel{eq:it1:second:incl}
						\mathcal{A}_3\subseteq\mathcal{A}_2.
				\end{equation}}
			\argument{
				\lref{eq:it1:first:incl};
				\lref{eq:it1:second:incl}}{\cref{it:1}}.
				\startnewargseq
				\argument{
					\cref{eq:super:strongly:asymptotically:stable};
					\cref{eq:ssasympt};}{\cref{it:2}}.
				\startnewargseq
				\argument{
					\cref{eq:asymptotically:stable};
					\cref{eq:asympt};}{that for all $\mathcal{B}\in\mathcal{A}_1$, $\mathcal{C}\subseteq\mathcal{B}$ it holds that
					\begin{equation}
						\llabel{eq:closure:s}
						\mathcal{C}\in\mathcal{A}_1.
					\end{equation}
				}
				\argument{
					\cref{eq:strongly:asymptotically:stable};
					\cref{eq:sasympt};}{that for all $\mathcal{B}\in\mathcal{A}_2$, $\mathcal{C}\subseteq\mathcal{B}$ it holds that
					\begin{equation}
						\llabel{eq:closure:ss}
						\mathcal{C}\in\mathcal{A}_2.
					\end{equation}
				}
				\argument{
					\cref{eq:super:strongly:asymptotically:stable};
					\cref{eq:ssasympt};}{that for all $\mathcal{B}\in\mathcal{A}_3$, $\mathcal{C}\subseteq\mathcal{B}$ it holds that
					\begin{equation}
						\llabel{eq:closure:sss}
						\mathcal{C}\in\mathcal{A}_3.
					\end{equation}
					}
					\argument{
						\lref{eq:closure:s};
						\lref{eq:closure:ss};
						\lref{eq:closure:sss};}{\cref{it:closed:subsets}}.
				\startnewargseq
				\argument{
					\cref{eq:asymptotically:stable};
					\cref{eq:asympt};}{that for all $\mathcal{B},\mathcal{C}\in\mathcal{A}_1$ it holds that
					\begin{equation}
						\llabel{eq:closure:union}
						\pr{\mathcal{B}\cup\mathcal{C}}\in\mathcal{A}_1.
					\end{equation}
				}
				\argument{\lref{eq:closure:union};}{\cref{it:4}}.
				\startnewargseq
				\argument{\cref{it:1};\cref{it:closed:subsets};}{\cref{it:5}}.
			\end{aproof}
		
		\subsection{Asymptotic stability of the Adam and the RMSprop optimizer}
		\label{subsection:5.2}
		
		In the following notion, \cref{def:ADAM} below, we recall the introduction of the $\alpha$-$\beta$-$\eps$-\Adam\ optimizer \cite{KingmaBa2024_Adam} using the general functions $\Phi_n\colon(\R^d)^{n}\to\R^d$, $n \in \N$, from \cref{def:stab:region} above (cf., \eg, \cite[Definitions~4.1]{dereich2025sharphigherorderconvergence}).

	\begin{definition}[Adam optimizer]
		\label{def:ADAM}
		Let $d\in\N$, $\alpha,\beta\in[0,1)$, $\eps\in(0,\infty)$ and
		let $\Phi_n=\pr{\Phi_n^{(1)},\dots\allowbreak,\Phi_n^{(d)}}\colon\allowbreak(\R^d)^{n}\to\R^d$, $n \in \N$, be functions.
		Then we say that $(\Phi_n)_{n \in \N}$ is the $\alpha$-$\beta$-$\eps$-\Adam\ optimizer on $\R^d$
		(we say that $(\Phi_n)_{n \in \N}$ is the $\alpha$-$\beta$-$\eps$-\Adam\ optimizer)
		 if and only if it holds
		for all 
		$n\in\N$,
		$i\in\{1,2,\dots,d\}$,
		$g_1=(\grad{1}{1},\dots,\grad{1}{d})$,
		$g_2=(\grad{2}{1},\dots,\grad{2}{d})$,
		$\dots$,
		$g_{n}=(\grad{n}{1},\dots,\grad{n}{d})\in\R^d$
		that
		\begin{equation}
			\label{def:eq:ADAM}
			\Phi_n^{(i)}(g_1,g_2,\dots,g_n)
			=
			\frac{\prb{\frac{1-\alpha}{1-\alpha^n}}\sum_{k=1}^n \alpha^{n-k}\grad{k}{i}}{\eps+\PRb{\prb{\frac{1-\beta}{1-\beta^n}}\sum_{k=1}^n\beta^{n-k}\vass{\grad{k}{i}}^2}^{\nicefrac{1}{2}}}
			.
		\end{equation}
	\end{definition}

			\cfclear
	\begin{athm}{lemma}{lem:ADAM:classic:representation}
		Let $d\in\N$, $\alpha,\beta\in[0,1)$, $\eps\in(0,\infty)$ and
		let $\Phi_n=\pr{\Phi_n^{(1)},\dots,\Phi_n^{(d)}}\colon(\R^d)^{n}\to\R^d$, $n \in \N$, be functions.
		Then $(\Phi_n)_{n \in \N}$ is the $\alpha$-$\beta$-$\eps$-\Adam\ optimizer on $\R^d$\cfadd{def:ADAM} if and only if it holds
		for 
		every
		$g=(g_n)_{n\in\N}=\pr{\pr{\grad{n}{1},\dots,\grad{n}{d}}}_{n\in\N}\colon\N\to\R^d$
		that there exist
		$\mom=(\mom_n)_{n\in\N_0}=\pr{\pr{\mom_n^{(1)},\dots,\mom_n^{(d)}}}_{n\in\N_0}\colon\N_0\to\R^d$
		and 
		$\MOM=(\MOM_n)_{n\in\N_0}=\pr{\pr{\MOM_n^{(1)},\dots,\MOM_n^{(d)}}}_{n\in\N_0}\colon\N_0\to\R^d$ 
		such that for all $n\in\N$, $i\in\{1,2,\dots,d\}$ it holds that
		\begin{equation}
			\label{eq:ADAM:1}
			\begin{split}
				\mom_0=0, 
				\quad\qquad 
				\mom_n&\textstyle=\alpha \mom_{n-1}+(1-\alpha)g_n,
			\end{split}
		\end{equation}
		\begin{equation}
			\label{eq:ADAM:2}
			\begin{split}
				\MOM^{(i)}_0=0,
				\quad\qquad
				\MOM_{n}^{(i)}&\textstyle=\beta \MOM_{n-1}^{(i)}+(1-\beta)\vass
				{\grad{n}{i}}^2,
			\end{split}
		\end{equation}
		\begin{equation}
			\label{eq:ADAM:3}
			\begin{split}
				\text{and}
				\qquad
				\Phi_n^{(i)}(g_1,g_2,\dots,g_n)
				=\PRbbb{\eps+\PRbbb{ \frac{\MOM_{n}^{(i)}}{1-\beta^n}}^{\nicefrac{1}{2}}}^{-1}\PRbbb{\frac{\mom_n^{(i)}}{1-\alpha^n}}
			\end{split}
		\end{equation}
		\cfload.
	\end{athm}
	\cfclear

	\begin{aproof}
		\argument{\cref{momentum:representation} (applied with
			$\pars\curvearrowleft d$,
			$\mom\curvearrowleft \pr{\N_0\ni n\mapsto \mom_n\in\R^d}$,
			$(g_n)_{n\in\N}\curvearrowleft \pr{\N\ni n\mapsto
				\prb{(1-\delta)\prb{g_n^{(1)}}^j,
					(1-\delta)\prb
					{g_n^{(2)}}^j,
					\dots,
					(1-\delta)\prb
					{g_n^{(d)}}^j
				}\in\R^d}$,
			$(\beta_n)_{n\in\N}\curvearrowleft \pr{\N\ni n\mapsto\delta\in\R}$
			for $j\in\{1,2\}$, $\delta\in\{\alpha,\beta\}$
			in the notation of \cref{momentum:representation});
		}{that for every $j\in\{1,2\}$, $\delta\in\{\alpha,\beta\}$,
		every
			$(g_n)_{n\in\N}=\pr{\pr{\grad{n}{1},\dots,\grad{n}{d}}}_{n\in\N}\colon\N\to\R^d$,
			and	every
			$(\mom_n)_{n\in\N_0}=\pr{\pr{\mom_n^{(1)},\dots,\mom_n^{(d)}}}_{n\in\N_0}\colon\N_0\to\R^d$ which satisfy for all 
			$n\in\N$,  
			$i\in\{1,2,\dots,d\}$ that
			\begin{equation}
				\llabel{eq:ADAM:12}
				\begin{split}
					\mom_{0}=0
					\qqandqq
					\mom_{n}^{(i)}&\textstyle=\delta \mom_{n}^{(i)}+(1-\delta)\prb
					{g_n^{(i)}}^j,
				\end{split}
			\end{equation}
			it holds for all $n\in\N$, $i\in\{1,2,\dots,d\}$ that
			\begin{equation}
				\begin{split}
					\textstyle
					\mom_n=(1-\delta)\sum_{k=1}^n\delta^{n-k}\prb
					{g_k^{(i)}}^j
					.
				\end{split}
			\end{equation}
		}
		\argument{\lref{eq:ADAM:12};
			\eqref{def:eq:ADAM};}{
			that $(\Phi_n)_{n \in \N}$ is the $\alpha$-$\beta$-$\eps$-\Adam\ optimizer on $\R^d$ if and only if 
			it holds
			for 
			every
			$(g_n)_{n\in\N}=\pr{\pr{\grad{n}{1},\dots,\grad{n}{d}}}_{n\in\N}\colon\N\to\R^d$
			that there exist
			$(\mom_n)_{n\in\N_0}=\pr{\pr{\mom_n^{(1)},\dots,\mom_n^{(d)}}}_{n\in\N_0}\colon\N_0\to\R^d$
			and 
			$(\MOM_n)_{n\in\N_0}=\pr{\pr{\MOM_n^{(1)},\dots,\MOM_n^{(d)}}}_{n\in\N_0}\colon\N_0\to\R^d$ 
			 such that for all $n\in\N$, $i\in\{1,2,\dots,d\}$ it holds that
			\begin{equation}
				\begin{split}
					\mom_0=0, 
					\qquad\qquad 
					\mom_n&\textstyle=\alpha \mom_{n-1}+(1-\alpha)g_n,
				\end{split}
			\end{equation}
			\begin{equation}
				\begin{split}
					\MOM^{(i)}_0=0,
					\qquad\qquad
					\MOM_{n}^{(i)}&\textstyle=\beta \MOM_{n-1}^{(i)}+(1-\beta)\vass
					{\grad{n}{i}}^2,
				\end{split}
			\end{equation}
			\begin{equation}
				\begin{split}
					\text{and}
					\qquad\qquad
					\Phi_n^{(i)}(g_1,g_2,\dots,g_n)
					&=\PRbbb{\eps+\PRbbb{ \frac{\MOM_{n}^{(i)}}{1-\beta^n}}^{\nicefrac{1}{2}}}^{-1}\PRbbb{\frac{\mom_n^{(i)}}{1-\alpha^n}}
					\\&=
					\frac{\prb{\frac{1-\alpha}{1-\alpha^n}}\sum_{k=1}^n \alpha^{n-k}\grad{k}{i}}{\eps+\PRb{\prb{\frac{1-\beta}{1-\beta^n}}\sum_{k=1}^n\beta^{n-k}\vass{\grad{k}{i}}^2}^{\nicefrac{1}{2}}}
					.
				\end{split}
			\end{equation}		
		}
	\end{aproof}

				\begin{athm}{cor}{cor:stoch:ADAM:bounded}
		Let 
		$d\in\N$, $\alpha\in[0,1)$, $\beta\in(\alpha^2,1)$, $\eps\in(0,\infty)$,
		let $\Phi_n=\pr{\Phi_n^{(1)},\dots,\Phi_n^{(d)}}\colon(\R^d)^{n}\to\R^d$, $n \in \N$,\cfadd{def:ADAM} be the $\alpha$-$\beta$-$\eps$-\Adam\ optimizer, 
		let
		$\gamma\colon\N\to[0,\infty)$ be bounded,
		let
		$\batch\colon\N\to\N$ be a function,
		let $\lambda=(\lambda_1,\dots,\lambda_d)\in[0,\infty)^d$, $\fc\in[0,\infty)$,
		let $(\Omega,\cF,\P)$ be a probability space,
		for every $n,j\in\N$
		let $X_{n,j}=\pr{X_{n,j}^{(1)},\dots,X_{n,j}^{(d)}}\colon\Omega\to[-\fc,\fc]^d$ be a random variable,
		and
		let
		$ \Grad=\pr{\Grad^{(1)},\dots,\Grad^{(d)}}\colon\N\times\Omega \to \R^d$ 				
		and
		$ \Theta \colon\N_0\times\Omega \to \R^d$ 
		satisfy for all
		$n\in\N$ that
		\begin{equation}
			\llabel{eq:ADAM:3}
			\begin{split}
				\textstyle
				\Grad_n=
				\frac{1}{\batch_n}\sum_{j=1}^{\batch_n}\diag(\lambda)\pr{\Theta_{n-1}-X_{n,j}}
				\qqandqq
				\Theta_n
				= 
				\Theta_{ n - 1 }
				-
				\gamma_n\Phi_n\prb{\Grad_1,\Grad_2,\dots,\Grad_n}
			\end{split}
		\end{equation}
		\cfload.
		Then there exists $c\in\R$ such that $\sup_{n\in\N_0}\norm{\Theta_n}\leq c\norm{\Theta_0}+ c$.
	\end{athm}
	
	\begin{aproof}
		Throughout this proof assume without loss of generality that $\lambda\in(0,\infty)^d$.
		\argument{the assumption that $\pr{\Phi_n}_{n \in \N}$\cfadd{def:ADAM} is the $\alpha$-$\beta$-$\eps$-\Adam\ optimizer;
			\cref{lem:ADAM:classic:representation} (applied with
			$d\curvearrowleft d$,
			$\alpha\curvearrowleft \alpha$,
			$\beta\curvearrowleft \beta$,
			$\eps\curvearrowleft \eps$,
			$(\Phi_n)_{n\in\N}\curvearrowleft(\Phi_n)_{n\in\N}$,
			$(g_n)_{n\in\N}\curvearrowleft(\Grad_n)_{n\in\N}$
			in the notation of \cref{lem:ADAM:classic:representation});
		}{that there exist 	$\mom=\prb{\mom^{(1)},\dots,\mom^{(d)}}\colon\N_0\times\Omega\to\R^d$ and
			$\MOM=\prb{\MOM^{(1)},\dots,\MOM^{(d)}}\colon\N_0\times\Omega\to[0,\infty)^d$ which satisfy for all $n\in\N$, $i\in\{1,2,\dots,d\}$ that
			\begin{equation}
				\llabel{it:ADAM:1}
				\begin{split}
					\mom_0=0, 
					\quad\qquad 
					\mom_n&\textstyle=\alpha \mom_{n-1}+(1-\alpha)\Grad_n,
				\end{split}
			\end{equation}
			\begin{equation}
				\llabel{it:ADAM:2}
				\begin{split}
					\MOM^{(i)}_0=0,
					\quad\qquad
					\MOM_{n}^{(i)}&\textstyle=\beta \MOM_{n-1}^{(i)}+(1-\beta)\vass
					{\Grad_{n}^{(i)}}^2,
				\end{split}
			\end{equation}
			\begin{equation}
				\llabel{it:ADAM:3}
				\begin{split}
					\text{and}
					\qquad
					\Phi_n^{(i)}(\Grad_{1},\Grad_{2},\dots,\Grad_{n})
					=\PRbbb{\eps+\PRbbb{ \frac{\MOM_{n}^{(i)}}{1-\beta^n}}^{\nicefrac{1}{2}}}^{-1}\PRbbb{\frac{\mom_n^{(i)}}{1-\alpha^n}}
					.
				\end{split}
		\end{equation}}
		\argument{\lref{it:ADAM:3};\lref{eq:ADAM:3}}{that for all $n\in\N$, $i\in\{1,2,\dots,d\}$ it holds that
			\begin{equation}
				\llabel{eq:ADAM:1.i}
				\begin{split} 
					\mom_{n}^{(i)}&\textstyle=\alpha \mom_{n-1}^{(i)}+(1-\alpha)\PRbb{\frac{\lambda_{i}}{\batch_n}\sum_{j=1}^{\batch_n}\prb{\Theta_{n-1}^{(i)}-X_{n,j}^{(i)}}},
				\end{split}
			\end{equation}
			\begin{equation}
				\llabel{eq:ADAM:2.i}
				\begin{split}
					\MOM_{n}^{(i)}&\textstyle=\beta \MOM_{n-1}^{(i)}+(1-\beta)\PRbb{\frac{\lambda_{i}}{\batch_n}\sum_{j=1}^{\batch_n}\prb{\Theta_{n-1}^{(i)}-X_{n,j}^{(i)}}}^2,
					\qqandqq
				\end{split}
			\end{equation}
			\begin{equation}
				\llabel{eq:ADAM:3.i}
				\begin{split}
					\Theta_n^{(i)}
										= 
										\Theta_{ n - 1 }^{(i)}
										-
										\gamma_n\Phi_n^{(i)}(\Grad_{1},\Grad_{2},\dots,\Grad_{n})
					= 
					\Theta_{ n - 1 }^{(i)}
					-
					\gamma_n\PRbbb{\eps+\PRbbb{ \frac{\MOM_{n}^{(i)}}{1-\beta^n}}^{\nicefrac{1}{2}}}^{-1}\PRbbb{\frac{\mom_n^{(i)}}{1-\alpha^n}}
					.
				\end{split}
		\end{equation}}
		\argument{\lref{it:ADAM:1};\lref{eq:ADAM:3.i};\cref{lem:ADAM:bounded} (applied with 
			$d\curvearrowleft d$,
			$\lambda\curvearrowleft \pr{\N\ni n\mapsto\lambda\in[0,\infty)^d}$,
			$\batch\curvearrowleft \batch$,
			$\gamma\curvearrowleft \gamma$,
			$\fc\curvearrowleft \fc$,
			$\prb{X_{n,j}^{(i)}}_{(n,i,j)\in\N^3}\curvearrowleft \prb{X_{n,j}^{(\min\{i,d\})}}_{(n,i,j)\in\N^3}$,
			$\alpha\curvearrowleft \alpha$,
			$\beta\curvearrowleft \beta$,
			$\eps\curvearrowleft \eps$,
			$\mom\curvearrowleft \mom$,
			$\MOM\curvearrowleft \MOM$,
			$\Theta\curvearrowleft \Theta$
			in the notation of \cref{lem:ADAM:bounded})}{that there exists $c\in\R$ such that $\sup_{n\in\N_0}\norm{\Theta_n}\leq c\norm{\Theta_0}+ c$}.
	\end{aproof}

	In the following notion, \cref{def:RMSprop} below, we recall the introduction of the $\beta$-$\eps$-\RMSprop\ optimizer (cf.\ \cite{Hinton24_RMSprop} and, \eg,  \cite[Definitions~6.6.5 and 7.7.3]{jentzen2023mathematical}) using the general functions $\Phi_n\colon(\R^d)^{n}\to\R^d$, $n \in \N$, from \cref{def:stab:region} above (cf., \eg, \cite[Definitions~2.2]{dereich2025sharphigherorderconvergence}).

		\begin{definition}[RMSprop optimizer]
		\label{def:RMSprop}
		Let $d\in\N$, $\beta\in[0,1)$, $\eps\in(0,\infty)$ and
		let $\Phi_n=\pr{\Phi_n^{(1)},\dots,\Phi_n^{(d)}}\colon\allowbreak(\R^d)^{n}\to\R^d$, $n \in \N$, be functions.
		Then we say that $(\Phi_n)_{n \in \N}$ is the $\beta$-$\eps$-\RMSprop\ optimizer on $\R^d$ 
		(we say that $(\Phi_n)_{n \in \N}$ is the $\beta$-$\eps$-\RMSprop\ optimizer)
		if and only if it holds
		for all 
		$n\in\N$,
		$i\in\{1,2,\dots,d\}$,
		$g_1=(\grad{1}{1},\dots,\grad{1}{d})$,
		$g_2=(\grad{2}{1},\dots,\grad{2}{d})$,
		$\dots$,
		$g_{n}=(\grad{n}{1},\dots,\grad{n}{d})\in\R^d$
		that
		\begin{equation}
			\label{def:eq:RMSprop}
			\Phi_n^{(i)}(g_1,g_2,\dots,g_n)
			=
			\frac{\grad{n}{i}}{\eps+\PRb{\prb{\frac{1-\beta}{1-\beta^n}}\sum_{k=1}^n\beta^{n-k}\vass{\grad{k}{i}}^2}^{\nicefrac{1}{2}}}
			.
		\end{equation}
	\end{definition}
	 In the following statement, \cref{RMS:is:ADAM} below, we  briefly recall the elementary fact that for every $\beta\in[0,1)$, $\eps\in(0,\infty)$ we have that the $\beta$-$\eps$-\RMSprop\ optimizer (see \cref{def:RMSprop} above) coincides with the $0$-$\beta$-$\eps$-\Adam\ method (see \cref{def:ADAM} above).

	\cfclear
		\begin{athm}{lemma}{RMS:is:ADAM}
		Let $d\in\N$, $\beta\in[0,1)$, $\eps\in(0,\infty)$ and
		let $\Phi_n\colon\allowbreak(\R^d)^{n}\to\R^d$, $n \in \N$, be the
		$0$-$\beta$-$\eps$-\Adam\ optimizer on $\R^d$\cfadd{def:ADAM} \cfload. 
		Then $(\Phi_n)_{n \in \N}$ is the
		$\beta$-$\eps$-\RMSprop\ optimizer on $\R^d$\cfadd{def:RMSprop} \cfout. 
	\end{athm}
		
	\begin{aproof}
		\argument{
			\eqref{def:eq:ADAM};
			the assumption that 
			$\pr{\Phi_n}_{n\in\N}$ is the $0$-$\beta$-$\eps$-\Adam\ optimizer;
			}{
		for all 
		$n\in\N$,
		$i\in\{1,2,\dots,d\}$,
		$g_1=(\grad{1}{1},\dots,\grad{1}{d})$,
		$g_2=(\grad{2}{1},\dots,\grad{2}{d})$,
		$\dots$,
		$g_{n}=(\grad{n}{1},\dots,\grad{n}{d})\in\R^d$
		that
		\begin{equation}
			\llabel{eq:adam:RMSprop:is:adam}
			\begin{split}
			\Phi_n^{(i)}(g_1,g_2,\dots,g_n)
			&=
			\frac{\grad{n}{i}}{\eps+\PRb{\prb{\frac{1-\beta}{1-\beta^n}}\sum_{k=1}^n\beta^{n-k}\vass{\grad{k}{i}}^2}^{\nicefrac{1}{2}}}
			.
			\end{split}
		\end{equation}
	}
	\argument{\eqref{def:eq:RMSprop};\lref{eq:adam:RMSprop:is:adam}}{that $(\Phi_n)_{n \in \N}$ is the
		$\beta$-$\eps$-\RMSprop\ optimizer on $\R^d$ \cfadd{def:RMSprop}}
	\end{aproof}

		\cfclear
		\begin{athm}{lemma}{lem:ADAM:stability}
			Let $d\in\N$, $\alpha\in[0,1)$, $\beta\in(\alpha^2,1)$, $\eps\in(0,\infty)$, $\vartheta=(\vartheta_1,\dots,\vartheta_d)\in\R^d$,
			 let $\Phi_n\colon(\R^d)^{n}\to\R^d$, $n \in \N$, be the $\alpha$-$\beta$-$\eps$-\Adam\ optimizer on $\R^d$\cfadd{def:ADAM},
			 and let 
			 $\lambda=(\lambda_n)_{n\in\N}=(\pr{\lambda_{n,1},\dots,\lambda_{n,d}})_{n\in\N}\colon\N\to[0,\infty)^d$,
			 $\gamma\colon\N\to[0,\infty)$, and
			 $\Theta=(\Theta_n)_{n\in\N_0}=\pr{\pr{\Theta_n^{(1)},\dots,\Theta_n^{(d)}}}_{n\in\N_0}\colon\allowbreak\N_0\to\R^d$ satisfy for all $n\in\N$ that
			 \begin{align}
			 	\llabel{eq:Theta:ADAM:for:stability}
			 	\Theta_n
			 	&=\Theta_{n-1}-\gamma_n\Phi_n(
			 	\operatorname{diag}(\lambda_1)(\Theta_0-\vartheta),
			 	\operatorname{diag}(\lambda_2)(\Theta_1-\vartheta),
			 	\dots,
			 	\operatorname{diag}(\lambda_n)(\Theta_{n-1}-\vartheta))
			 \end{align}
			 \cfload.
			Then 
			\begin{enumerate}[label=(\roman*)]
				\item \label{it:fixed:zero}it holds for all $m\in\N_0$, $i\in\{1,2,\dots,d\}$ with $\sum_{n\in\N}\lambda_{n,i}=0$ that $\Theta_m^{(i)}=\Theta_0^{(i)}$,
				\item \label{it:fixed:nonzero}it holds for all $i\in\{1,2,\dots,d\}$, $R\in[1,\infty)$ with $\cup_{n\in\N}\{(\gamma_n,\lambda_{n,i})\}\subseteq[0,R]\times[R^{-1},R]$ that $\sup_{n\in\N}\vass{\Theta_n^{(i)}}<\infty$,
				\item \label{it:gen:statement}it holds for all $R\in[1,\infty)$ with $\cup_{n\in\N}\{(\gamma_n,\lambda_n)\}\subseteq[0,R]\times[R^{-1},R]^d$ that $\sup_{n\in\N} \norm{\Theta_n}< \infty$,
				\item \label{it:fixed:for:stable}it holds for all $v\in[0,\infty)\times[0,\infty)^{d}$ with $\cup_{n\in\N}\{(\gamma_n,\lambda_{n})\}=\{v\}$ that $\sup_{n\in\N} \norm{\Theta_n}< \infty$,
				\item \label{it:for:sstable}it holds for all $R\in[0,\infty)$ with $\cup_{n\in\N}\{(\gamma_n,\lambda_{n})\}\subseteq[0,R]\times\{\lambda_1\}\subseteq[0,R]^{d+1}$ that $\sup_{n\in\N} \norm{\Theta_n}< \infty$, and
				\item \label{it:for:ssstable}it holds for all $R\in[1,\infty)$ with $\cup_{n\in\N}\{(\gamma_n,\lambda_{n})\}\subseteq[0,R]\times[R^{-1},R]^{d}$ that $\sup_{n\in\N} \norm{\Theta_n}< \infty$.
			\end{enumerate}
		\end{athm}
		
		\begin{aproof}
			\argument{\cref{it:fixed:zero,it:fixed:nonzero};}[plural]{\cref{it:fixed:for:stable,it:for:sstable}}.
			\startnewargseq\argument{\cref{it:fixed:nonzero};}{\cref{it:gen:statement}}.
			\startnewargseq\argument{\cref{it:gen:statement}}{\cref{it:for:ssstable}}.
			\startnewargseq\argument{
				\lref{eq:Theta:ADAM:for:stability};
				\cref{lem:ADAM:classic:representation};
				the assumption that $(\Phi_n)_{n \in \N}$ is the $\alpha$-$\beta$-$\eps$-\Adam\ on $\R^d$\cfadd{def:ADAM}}
				{there exist		
					$\mom_n=(\mom_n^{(1)},\dots,\mom_n^{(d)})\in\R^d$, $n\in\N_0$, and 
					$\MOM_n=(\MOM_n^{(1)},\dots,\MOM_n^{(d)})\in\R^d$, $n\in\N_0$, such that for all $n\in\N$, $i\in\{1,2,\dots,d\}$ it holds that
									\begin{equation}
						\llabel{eq:ADAM:1}
						\begin{split}
							\mom_0=0, 
							\quad\qquad 
							\mom_n&\textstyle=\alpha \mom_{n-1}+(1-\alpha)\operatorname{diag}(\lambda_n)(\Theta_{n-1}-\vartheta),
						\end{split}
					\end{equation}
					\begin{equation}
						\llabel{eq:ADAM:2}
						\begin{split}
							\MOM^{(i)}_0=0,
							\quad\qquad
							\MOM_{n}^{(i)}&\textstyle=\beta \MOM_{n-1}^{(i)}+(1-\beta)\PRb
							{\lambda_{n,i}\prb{\Theta_{n-1}^{(i)}-\vartheta_i}}^2,
						\end{split}
					\end{equation}
					\begin{equation}
						\llabel{eq:ADAM:3}
						\begin{split}
							\text{and}
							\qquad
							\Theta_n^{(i)}
							=
							\Theta_{n-1}^{(i)}
							-\gamma_n
							\PRbbb{\eps+\PRbbb{ \frac{\MOM_{n}^{(i)}}{1-\beta^n}}^{\nicefrac{1}{2}}}^{-1}\PRbbb{\frac{\mom_n^{(i)}}{1-\alpha^n}}
							.
						\end{split}
					\end{equation}
				}
			\argument{
				\lref{eq:ADAM:3};
				}{\cref{it:fixed:zero}}.
			\Nobs that
				\lref{eq:ADAM:1},
				\lref{eq:ADAM:2},
				\lref{eq:ADAM:3},
				the fact that for all $i\in\{1,2,\dots,d\}$, $R\in[1,\infty)$ with $\forall\,n\in\N\colon(\gamma_n\in[0,R])\wedge\pr{\lambda_{n,i}\in[R^{-1},R]}$ it holds that
			\begin{equation}
				\llabel{eq:LR:and:Lambda:bounds}
				\textstyle
				\inf_{n\in\N}\lambda_{n,i}\geq R^{-1} > 0
				\qqandqq
				\limsup_{n\to \infty}\PR{\gamma_n+\lambda_{n,i}}\leq R +R < \infty,
				\end{equation}
				and
				\cref{lem:ADAM:bounded}
				 (applied for every $i\in\{1,2,\dots,d\}$ with
				 $d\curvearrowleft 1$,
				$\lambda \curvearrowleft \pr{\N\ni n\mapsto \lambda_{n,i}\in\R}$,
				$\batch\curvearrowleft \pr{\N\ni n\mapsto 1\in\N}$,
				$\gamma \curvearrowleft \gamma$,
				$\alpha \curvearrowleft \alpha$,
				$\beta \curvearrowleft \beta$,
				$\eps \curvearrowleft \eps$,
				$(X_{n,j}^{(i)})_{(n,i,j)\in\N^3} \curvearrowleft \pr{\N\ni n\mapsto \vartheta_i\in\R}$,
				$\mom \curvearrowleft \mom$,
				$\MOM \curvearrowleft \MOM$,
				$\Theta \curvearrowleft \Theta$
				in the notation of \cref{lem:ADAM:bounded})
				\prove
				that for all $i\in\{1,2,\dots,d\}$, $R\in[1,\infty)$ with $\forall\,n\in\N\colon(\gamma_n\in[0,R])\wedge\pr{\lambda_{n,i}\in[R^{-1},R]}$ there exists $\fC\in\R$ such that
				\begin{equation}
					\llabel{eq:ADAM:ith:bound}
					\begin{split}
						\textstyle\sup_{n\in\N_0}\vass{\Theta_n^{(i)}}\leq \fC\prb{\vass{\Theta_0^{(i)}}+1}<\infty.
					\end{split}
				\end{equation}
		This \proves \cref{it:fixed:nonzero}.
		\end{aproof}

		\cfclear
		\begin{athm}{cor}{cor:ADAM:stable}
			Let $d\in\N$, $\alpha\in[0,1)$, $\beta\in(\alpha^2,1)$, $\eps\in(0,\infty)$ and let $\Phi_n
			\colon(\R^d)^{n}\to\R^d$, $n \in \N$, be the $\alpha$-$\beta$-$\eps$-\Adam\ optimizer on $\R^d$\cfadd{def:ADAM} \cfload.
			Then $(\Phi_n)_{n \in \N}$ \ustable\cfout.
		\end{athm}
		\cfclear
		
		\begin{aproof}
			Throughout this proof let $\gamma\in[0,\infty)$, $\vartheta\in\R^d$, $\lambda\in[0,\infty)^d$ and let $\Theta\colon\N_0  \to \R^d$ satisfy for all $n\in\N$ that
			\begin{equation}
				\llabel{eq:ADAM:recursion}
				\Theta_{n}=\Theta_{n-1}
				-
				\gamma\Phi_n(
				\operatorname{diag}(\lambda)(\Theta_0-\vartheta),
				\operatorname{diag}(\lambda)(\Theta_1-\vartheta),
				\dots,
				\operatorname{diag}(\lambda)(\Theta_{n-1}-\vartheta))
				.
			\end{equation}
			\argument{
				\lref{eq:ADAM:recursion};
				\cref{it:fixed:for:stable} in \cref{lem:ADAM:stability} (applied with 
				$d\curvearrowleft d$,
				$\alpha\curvearrowleft\alpha$,
				$\beta \curvearrowleft \beta$,
				$(\Phi_n)_{n\in\N}\curvearrowleft(\Phi_n)_{n\in\N}$,
				$\lambda \curvearrowleft \pr{\N\ni n \mapsto \lambda\in[0,\infty)^d}$,
				$\gamma \curvearrowleft \pr{\N\ni n\mapsto \gamma\in[0,\infty)}$,
				$\Theta \curvearrowleft \Theta$,
				$v\curvearrowleft (\gamma,\lambda)$
				in the notation of \cref{lem:ADAM:stability})}{\llabel{L1}that $\limsup_{n\to \infty}\norm{\Theta_n}<\infty$}.
				\argument{\lref{L1}}{that $(\Phi_n)_{n \in \N}$ \ustable}
		\end{aproof}

			\cfclear
		\begin{athm}{cor}{cor:ADAM:stable:clear}
			Let $d\in\N$, $\alpha\in[0,1)$, $\beta\in(\alpha^2,1)$, $\eps\in(0,\infty)$,
			let $\mathcal{A}\subseteq[0,\infty)^{d+1}$ be a set, and
			 let $\Phi_n\colon(\R^d)^{n}\to\R^d$, $n \in \N$, be the $\alpha$-$\beta$-$\eps$-\Adam\ optimizer on $\R^d$\cfadd{def:ADAM} \cfload.
			Then $(\Phi_n)_{n \in \N}$ \stable{\mathcal{A}}\cfout.
		\end{athm}
		\cfclear
		
		\begin{aproof}
			\argument{\cref{cor:ADAM:stable}}{\llabel{U:stable}that $(\Phi_n)_{n \in \N}$ \ustable}
			\argument{\lref{U:stable};}{\llabel{spec:stable}that $(\Phi_n)_{n\in \N}$ \stable{[0,\infty)^{d+1}}}
			\argument{\lref{spec:stable};\cref{it:closed:subsets} in \cref{lem:basic:stable:properties};the assumption that $\mathcal{A}\subseteq[0,\infty)^{d+1}$}{that $(\Phi_n)_{n \in \N}$ \stable{\mathcal{A}}}\!\!.
		\end{aproof}

		\begin{athm}{cor}{cor:ADAM:sstable}
			Let $d\in\N$, $\alpha\in[0,1)$, $\beta\in(\alpha^2,1)$, $\eps\in(0,\infty)$,
			let $\mathcal{A}\subseteq[0,\infty)^{d+1}$ be bounded, and let $\Phi_n
			\colon(\R^d)^{n}\to\R^d$, $n \in \N$, be the $\alpha$-$\beta$-$\eps$-\Adam\ optimizer on $\R^d$\cfadd{def:ADAM} \cfload.
			Then $(\Phi_n)_{n \in \N}$ \sstable{\mathcal{A}}\cfout.
		\end{athm}
		\cfclear
		
				\begin{aproof}
			Throughout this proof let $\vartheta\in\R^d$, $R\in(0,\infty)$, $\lambda\in[0,R]^d$ satisfy $\mathcal{A}\subseteq [0,R]^{d+1}$ and let $\gamma\colon\N\to[0,R]$ and $\Theta\colon\N_0  \to \R^d$ satisfy for all $n\in\N$ that
			\begin{equation}
				\llabel{eq:ADAM:recursion}
				\Theta_{n}=\Theta_{n-1}
				-
				\gamma_n\Phi_n(
				\operatorname{diag}(\lambda)(\Theta_0-\vartheta),
				\operatorname{diag}(\lambda)(\Theta_1-\vartheta),
				\dots,
				\operatorname{diag}(\lambda)(\Theta_{n-1}-\vartheta))
				.
			\end{equation}
			\argument{
				\lref{eq:ADAM:recursion};
				the fact that $(\Phi_n)_{n \in \N}$ is the $\alpha$-$\beta$-$\eps$-\Adam\ optimizer on $\R^d$;
				\cref{it:for:sstable} in \cref{lem:ADAM:stability} (applied with 
				$d\curvearrowleft d$,
				$\alpha\curvearrowleft\alpha$,
				$\beta \curvearrowleft \beta$,
				$(\Phi_n)_{n\in\N}\curvearrowleft(\Phi_n)_{n\in\N}$,
				$\lambda \curvearrowleft \pr{\N\ni n \mapsto \lambda\in[0,\infty)^d}$,
				$\gamma \curvearrowleft \gamma$,
				$\Theta \curvearrowleft \Theta$
				in the notation of \cref{lem:ADAM:stability})}[verbs=s]{\llabel{L1}that $\limsup_{n\to \infty}\norm{\Theta_n}<\infty$}.
			\argument{\lref{L1}}{\llabel{L2}that $(\Phi_n)_{n \in \N}$ \sstable{[0,R]^{d+1}}}
			\argument{
				\lref{L2};
				\cref{it:closed:subsets} in \cref{lem:basic:stable:properties};
				the fact that $\mathcal{A}\subseteq[0,R]^{d+1}$;}{that $(\Phi_n)_{n \in \N}$ \ssstable{\mathcal{A}}}
		\end{aproof}

		\begin{athm}{cor}{cor:ADAM:ssstable}
			Let $d\in\N$, $\alpha\in[0,1)$, $\beta\in(\alpha^2,1)$, $\eps\in(0,\infty)$,
			let $\mathcal{A}\subseteq[0,\infty)\times(0,\infty)^{d}$ be compact,
			and let $\Phi_n=\pr{\Phi_n^{(1)},\dots,\Phi_n^{(d)}}\colon(\R^d)^{n}\to\R^d$, $n \in \N$, be the $\alpha$-$\beta$-$\eps$-\Adam\ optimizer on $\R^d$\cfadd{def:ADAM} \cfload.
			Then $(\Phi_n)_{n \in \N}$ \ssstable{\mathcal{A}}\cfout.
		\end{athm}
		
						\begin{aproof}
			Throughout this proof let $\vartheta\in\R^d$, $R\in(0,\infty)$, $\lambda\in[0,R]^d$ satisfy $\mathcal{A}\subseteq\pr{[0,R]\times[R^{-1},R]^{d}}$ and let $\gamma\colon\N\to[0,R]$ and $\Theta\colon\N_0  \to \R^d$ satisfy for all $n\in\N$ that
			\begin{equation}
				\llabel{eq:ADAM:recursion}
				\Theta_{n}=\Theta_{n-1}
				-
				\gamma_n\Phi_n(
				\operatorname{diag}(\lambda_1)(\Theta_0-\vartheta),
				\operatorname{diag}(\lambda_2)(\Theta_1-\vartheta),
				\dots,
				\operatorname{diag}(\lambda_n)(\Theta_{n-1}-\vartheta))
				.
			\end{equation}
			\argument{
				\lref{eq:ADAM:recursion};
				\cref{it:for:ssstable} in \cref{lem:ADAM:stability} (applied with 
				$d\curvearrowleft d$,
				$\alpha\curvearrowleft\alpha$,
				$\beta \curvearrowleft \beta$,
				$(\Phi_n)_{n\in\N}\curvearrowleft(\Phi_n)_{n\in\N}$,
				$\lambda \curvearrowleft \lambda$,
				$\gamma \curvearrowleft \gamma$,
				$\Theta \curvearrowleft \Theta$
				in the notation of \cref{lem:ADAM:stability})}{\llabel{L1}that $\limsup_{n\to \infty}\norm{\Theta_n}<\infty$}.
			\argument{\lref{L1}}{\llabel{L2}that $(\Phi_n)_{n \in \N}$ \ssstable{\pr{[0,R]\times[R^{-1},R]^{d}}}\cfload}.
			\argument{
				\lref{L2};
				\cref{lem:basic:stable:properties};
				the fact that $\mathcal{A}\subseteq\pr{[0,R]\times[R^{-1},R]^{d}}$;}{that $(\Phi_n)_{n \in \N}$ \ssstable{\mathcal{A}}}\!\!.
		\end{aproof}

		\subsection{Asymptotic stability of the standard GD optimizer}
		\label{subsection:5.3}
		
			In the following notion, \cref{def:SGD} below, we recall the introduction of the \GD\ optimizer using the general functions $\Phi_n\colon(\R^d)^{n}\to\R^d$, $n \in \N$, from \cref{def:stab:region} above (cf., \eg, \cite[Definitions~2.1]{dereich2025sharphigherorderconvergence}).
		
		\begin{definition}[GD optimizer]
			\label{def:SGD}
			Let $d\in\N$ and
			let $\Phi_n\colon\allowbreak(\R^d)^{n}\to\R^d$, $n \in \N$, be functions.
			Then we say that $(\Phi_n)_{n \in \N}$ is the \GD\ optimizer on $\R^d$
			(we say that $(\Phi_n)_{n \in \N}$ is the \GD\ optimizer) if and only if it holds
			for all 
			$n\in\N$,
			$g_1,g_2,\dots,g_{n}\in\R^d$
			that
			\begin{equation}
				\label{def:eq:SGD}
				\Phi_n(g_1,g_2,\dots,g_n)
				=
				g_n
				.
			\end{equation}
		\end{definition}
		
			\begin{athm}{cor}{cor:SGD:stable}
			Let $d\in\N$,
			let $\mathcal{A}\subseteq[0,\infty)^{d+1}$ be a set, and let $\Phi_n\colon(\R^d)^{n}\to\R^d$, $n \in \N$, be the \GD\ optimizer on $\R^d$\cfadd{def:SGD} \cfload.
			Then $(\Phi_n)_{n \in \N}$ \stable{\mathcal{A}}\!\! if and only if it holds for all $(\lambda_0,\lambda_1,\dots,\lambda_{d})\in\mathcal{A}$ that \begin{equation}
				\llabel{eq:GD:asymp-stable}
				\textstyle\max_{i\in\{1,2,\dots,d\}}\pr{\lambda_0\lambda_{i}}\leq2
			\end{equation}\cfout.
		\end{athm}
		\cfclear
		
		\begin{aproof}
			Throughout this proof 
			let $\mathcal{B}\subseteq[0,\infty)^{d+1}$ satisfy $\mathcal{B}=\{(\lambda_0,\lambda_1,\dots,\lambda_{d})\in[0,\infty)^{d+1}\colon \max_{i\in\{1,2,\dots,d\}}(\lambda_0\lambda_i)\leq 2\}$,
			let $\vartheta=(\vartheta_1,\dots,\vartheta_d)\in\R^d$,
			for every 
			$i\in\{1,2,\dots,d\}$, $\lambda=(\lambda_1,\dots,\lambda_{d})\in[0,\infty)^{d}$ let $\mathscr{L}_i^{\lambda}\colon\R\to\R$ satisfy for all $\theta\in\R$ that
			\begin{equation}
				\llabel{eq:loss:quadr:stretched}
				\mathscr{L}_i^{\lambda}(\theta)=\tfrac{\lambda_i}{2}(\theta-\vartheta_i)^2,
			\end{equation}
			and let $\Theta
			=\prb{\Theta_{n}^{\gamma,\lambda}}_{(\gamma,\lambda,n)\in[0,\infty)\times[0,\infty)^{d}\times\N_0}
			=\prb{\Theta_{n,1}^{\gamma,\lambda},\dots,\Theta_{n,d}^{\gamma,\lambda}}_{(\gamma,\lambda,n)\in[0,\infty)\times[0,\infty)^{d}\times\N_0}\colon[0,\infty)\times[0,\infty)^{d}\times\N_0\to\R^d$ satisfy for all
			$\gamma\in[0,\infty)$,
			$\lambda=(\lambda_1,\dots,\lambda_{d})\in[0,\infty)^{d}$, 
			$n\in\N$, 
			$i\in\{1,2,\dots,d\}$ that
			\begin{equation}
				\llabel{eq:all:dim:SGD}
				\begin{split}
					\Theta_{n}^{\gamma,\lambda}
					&=
					\Theta_{n-1}^{\gamma,\lambda}
					-
					\gamma\Phi_n
					\prb{
						\operatorname{diag}(\lambda)\prb{\Theta_{0}^{\gamma,\lambda}-\vartheta},
						\operatorname{diag}(\lambda)\prb{\Theta_{1}^{\gamma,\lambda}-\vartheta},
						\dots,
						\operatorname{diag}(\lambda)\prb{\Theta_{n-1}^{\gamma,\lambda}-\vartheta}}
					.
				\end{split}
			\end{equation}
			\startnewargseq
			\argument{
				\lref{eq:loss:quadr:stretched};
				\lref{eq:all:dim:SGD};
				the assumption that $(\Phi_n)_{n \in \N}$ is the \GD\ optimizer on $\R^d$\cfadd{def:SGD};
			}{
				that for all 
				$\gamma\in[0,\infty)$,
				$\lambda=(\lambda_1,\dots,\lambda_{d})\in[0,\infty)^{d}$, 
				$n\in\N$, 
				$i\in\{1,2,\dots,d\}$ it holds that
				\begin{equation}
					\llabel{eq:repr:like:book}
					\Theta_{n,i}^{\gamma,\lambda}
					=
					\Theta_{n-1,i}^{\gamma,\lambda}
					-
					\gamma\lambda_i\prb{\Theta_{n-1,i}^{\gamma,\lambda}-\vartheta_i}
					=
					\Theta_{n-1,i}^{\gamma,\lambda}
					-
					\gamma\prb{\nabla\mathscr{L}_i^{\lambda}}\prb{\Theta_{n-1,i}^{\gamma,\lambda}}
			\end{equation}}
			\!\argument{\lref{eq:repr:like:book};
				\cite[Theorem~6.1.12]{jentzen2023mathematical}
				(applied with 
				$\pars\curvearrowleft 1$,
				$\alpha\curvearrowleft\lambda_i$,
				$\gamma \curvearrowleft \gamma$,
				$\vartheta\curvearrowleft\vartheta_i$,
				$\xi\curvearrowleft \Theta_{0,i}^{\gamma,\lambda}$,
				$\mathscr{L}\curvearrowleft\mathscr{L}_i^{\lambda}$,
				$\Theta\curvearrowleft\prb{\N_0\ni n\mapsto\Theta_{n,i}^{\gamma,\lambda}\in\R}$ 
				for 
				$\gamma\in[0,\infty)$,
				$\lambda=(\lambda_1,\dots,\lambda_d)\in[0,\infty)^d$,
				$i\in\{1,2,\dots,d\}$
				in the notation of \cite[Theorem~6.1.12]{jentzen2023mathematical})}{that for all 
				$\gamma\in[0,\infty)$,
				$\lambda=(\lambda_1,\dots,\lambda_{d})\in[0,\infty)^{d}$,  
				$i\in\{1,2,\dots,d\}$
				with $\Theta_{0,i}^{\gamma,\lambda}\neq\vartheta_i$ and $\lambda_i>0$ it holds that
				\begin{equation}
					\llabel{eq:BUCH:SGD:result}
					\liminf_{n\to\infty}\vass{\Theta_{n,i}^{\gamma,\lambda}-\vartheta_i}
					=\limsup_{n \to \infty}\vass{\Theta_{n,i}^{\gamma,\lambda}-\vartheta_i}
					=
					\begin{cases}
						0&\colon\gamma\lambda_i\in(0,2)\\
						\vass{\Theta_{0,i}^{\gamma,\lambda}-\vartheta_i}&\colon\gamma\lambda_i\in\{0,2\}\\
						\infty&\colon\gamma\lambda_i\in(2,\infty). 
					\end{cases}
			\end{equation}}
			\argument{
				\lref{eq:repr:like:book};
				\lref{eq:BUCH:SGD:result};}{that for all 
				$\gamma\in[0,\infty)$,
				$\lambda=(\lambda_1,\dots,\lambda_{d})\in[0,\infty)^{d}$,  
				$i\in\{1,2,\dots,d\}$ 
				with $\gamma\lambda_i\leq 2$ it holds that
				\begin{equation}
					\llabel{eq:BUCH:SGD:result:2}
					\limsup_{n \to \infty}\vass{\Theta_{n,i}^{\gamma,\lambda}-\vartheta_i}<\infty.
			\end{equation}}
			\argument{\lref{eq:BUCH:SGD:result:2};}{that for all 
				$\gamma\in[0,\infty)$,
				$\lambda=(\lambda_1,\dots,\lambda_{d})\in[0,\infty)^{d}$
				with $\max_{i\in\{1,2,\dots,d\}}(\gamma\lambda_i)\leq 2$
				it holds that
				\begin{equation}
					\llabel{eq:stable:case}
					\begin{split}
						\textstyle
						\limsup_{n\to\infty}\norm{\Theta_{n}^{\gamma,\lambda}-\vartheta}^2
						&\textstyle=
						\limsup_{n\to\infty}\PRb{\sum_{i=1}^d\vass{\Theta_{n,i}^{\gamma,\lambda}-\vartheta_i}^2}
						\\&\textstyle\leq
						\sum_{i=1}^d \PRb{\limsup_{n\to\infty} \vass{\Theta_{n,i}^{\gamma,\lambda}-\vartheta_i}^2}
						\\&\textstyle=
						\sum_{i=1}^d \PRb{\limsup_{n\to\infty} \vass{\Theta_{n,i}^{\gamma,\lambda}-\vartheta_i}}^2
						<\infty.
					\end{split}
			\end{equation}}
			\argument{
				\lref{eq:all:dim:SGD};
				\lref{eq:stable:case};
			}{that 
			$(\Phi_n)_{n \in \N}$ \stable{\mathcal{B}}\!}
			\argument{\lref{eq:BUCH:SGD:result};}{that for all
				$\gamma\in[0,\infty)$,
				$\lambda=(\lambda_1,\dots,\lambda_{d})\in[0,\infty)^{d}$ with
				$\Theta_{0,i}^{\gamma,\lambda}\neq\vartheta_i$ and $\max_{i\in\{1,2,\dots,d\}}\pr{\gamma\lambda_i}>2$ it holds that
				\begin{equation}
					\llabel{eq:unstable:case}
					\begin{split}
						\textstyle
						\limsup_{n\to\infty}\norm{\Theta_{n}^{\gamma,\lambda}-\vartheta}
						&\textstyle\geq
						\limsup_{n\to\infty} \PRb{\max_{i\in\{1,2,\dots,d\}}\vass{\Theta_{n,i}^{\gamma,\lambda}-\vartheta_i}}
						\\&\textstyle=
						\max_{i\in\{1,2,\dots,d\}}
						\PRb{\limsup_{n\to\infty} \vass{\Theta_{n,i}^{\gamma,\lambda}-\vartheta_i}}
						=
						\infty.
					\end{split}
				\end{equation}
			}
			\argument{
				\lref{eq:all:dim:SGD};
				\lref{eq:unstable:case};
				\cref{lem:basic:stable:properties}}{that $(\Phi_n)_{n \in \N}$ \stable{\mathcal{A}}\!\! if and only if $\mathcal{A}\subseteq\mathcal{B}$}.
		\end{aproof}
		
		\subsection{Asymptotic stability of the momentum optimizer}
		\label{subsection:5.4}
		
		In the following notion, \cref{def:MOM} below, we recall the introduction of the $\alpha$-momentum optimizer \cite{Poljak_momentum_SGD} using the general functions $\Phi_n\colon(\R^d)^{n}\to\R^d$, $n \in \N$, from \cref{def:stab:region} above (cf., \eg, \cite[Definitions~3.1]{dereich2025sharphigherorderconvergence}).
		
			\begin{definition}[Momentum optimizer]
			\label{def:MOM}
			Let $d\in\N$, $\alpha\in[0,1]$ and
			let $\Phi_n\colon\allowbreak(\R^d)^{n}\to\R^d$, $n \in \N$, be functions.
			Then we say that $(\Phi_n)_{n \in \N}$ is the $\alpha$-momentum optimizer on $\R^d$ 
			(we say that $(\Phi_n)_{n \in \N}$ is the $\alpha$-momentum optimizer)
			if and only if it holds
			for all 
			$n\in\N$,
			$g_1,g_2,\dots,g_{n}\in\R^d$
			that
			\begin{equation}
				\label{def:eq:MOM}
				\Phi_n(g_1,g_2,\dots,g_n)
				=
				\pr{1-\alpha}\sum_{k=1}^n \alpha^{n-k}g_k
				.
			\end{equation}
		\end{definition}

		In the following statement, \cref{cor:SGD:is:MOM} below, we  briefly recall the elementary fact that the \GD\ optimizer (see \cref{def:SGD} above) coincides with the $0$-momentum method (see \cref{def:MOM} above).
		
		\begin{athm}{lemma}{cor:SGD:is:MOM}
			Let $d\in\N$ and
			let $\Phi_n\colon(\R^d)^{n}\to\R^d$, $n \in \N$, be the $0$-momentum optimizer on $\R^d$\cfadd{def:MOM} \cfload.
			Then $(\Phi_n)_{n \in \N}$ is the \GD\ optimizer on $\R^d$\cfadd{def:SGD} \cfout.
			\end{athm}
			
			\begin{aproof}
				\argument{
					\eqref{def:eq:MOM};
					the assumption that 
					$\pr{\Phi_n}_{n\in\N}$ is the $0$-momentum optimizer}{for all 
					$n\in\N$,
					$g_1,g_2,\dots,g_{n}\in\R^d$
					that
					\begin{equation}
						\llabel{eq:1.sgd:mom}
						\Phi_n(g_1,g_2,\dots,g_n)
						=
						g_n
						.
					\end{equation}
			}
			\argument{\lref{eq:1.sgd:mom};\eqref{def:eq:SGD};}{that $(\Phi_n)_{n \in \N}$ is the \GD\ optimizer on $\R^d$ \cfadd{def:SGD}}
			\end{aproof}
		
		\cfclear
		\begin{lemma}
			\label{lem:MOM:classic:representation}
			Let $d\in\N$, $\alpha\in[0,1]$ and
			let $\Phi_n\colon(\R^d)^{n}\to\R^d$, $n \in \N$, be functions.
			Then $(\Phi_n)_{n \in \N}$ is the $\alpha$-momentum optimizer on $\R^d$\cfadd{def:MOM} if and only if it holds that
			for all
			$g=(g_n)_{n\in\N}\colon\N\to\R^d$
			there exists		
			$\mom=(\mom_n)_{n\in\N_0}\colon\N_0\to\R^d$
			such that for all $n\in\N$ it holds that
			\begin{equation}
				\llabel{eq:MOM:0}
				\begin{split}
					\mom_0=0, 
					\qquad 
					\mom_n&\textstyle=\alpha \mom_{n-1}+(1-\alpha)g_n,
					\qqandqq
					\Phi_n(g_1,g_2,\dots,g_n)
					=\mom_n
				\end{split}
			\end{equation}
			\cfload.
		\end{lemma}
		\cfclear

		\begin{aproof}
			\argument{
				\cref{momentum:representation} (applied with
				$\pars\curvearrowleft d$,
				$\mom\curvearrowleft \mom$,
				$(g_n)_{n\in\N}\curvearrowleft(1-\alpha)g_n$,
				$(\beta_n)_{n\in\N}\curvearrowleft\pr{\N\ni n\mapsto \alpha\in\R}$	
				in the notation of \cref{momentum:representation});}{that 
				for all	$g_n\in\R^d$, $n\in\N$,
				and		
				$\mom_n\in\R^d$, $n\in\N_0$,
				which satisfy for all $n\in\N$ that
				\begin{equation}
					\llabel{eq:MOM:1}
					\begin{split}
						\mom_0=0
						\qqandqq
						\mom_n&\textstyle=\alpha \mom_{n-1}+(1-\alpha)g_n
					\end{split}
				\end{equation}
				it holds for all
				 $n\in\N$ that
			\begin{equation}
				\llabel{eq:alt:repr:mom}
				\mom_n
				\textstyle
				=\alpha^n\mom_0+\sum_{k=1}^n\alpha^{n-k}(1-\alpha)g_k
				=(1-\alpha)\sum_{k=1}^n\alpha^{n-k}g_k
				.
				\end{equation}}
				\argument{
					\lref{eq:alt:repr:mom};
					\eqref{def:eq:MOM};
					\lref{eq:MOM:0};}{that $(\Phi_n)_{n \in \N}$ is the $\alpha$-momentum optimizer on $\R^d$\cfadd{def:MOM} if and only if it holds that
					for all 
					$g_n\in\R^d$, $n\in\N$,
					there exist		
					$\mom_n\in\R^d$, $n\in\N_0$,
					such that for all $n\in\N$ it holds that
					\begin{equation}
						\llabel{eq:MOM:1.2}
						\begin{split}
							\mom_0=0, 
							\qquad 
							\mom_n&\textstyle=\alpha \mom_{n-1}+(1-\alpha)g_n,
							\qqandqq
							\Phi_n(g_1,g_2,\dots,g_n)
							=\mom_n
						\end{split}
						\end{equation}}
		\end{aproof}

		\begin{athm}{cor}{cor:MOM:stable}
			Let $d\in\N$, $\alpha\in[0,1)$,
			let $\mathcal{A}\subseteq[0,\infty)^{d+1}$ be a set, and let $\Phi_n\colon(\R^d)^{n}\to\R^d$, $n \in \N$, be the $\alpha$-momentum optimizer on $\R^d$\cfadd{def:MOM} \cfload.
			Then $(\Phi_n)_{n \in \N}$ \stable{\mathcal{A}}\!\! if and only if it holds for all $(\lambda_0,\lambda_1,\dots,\lambda_{d})\in\mathcal{A}$ that \begin{equation}
				\textstyle\max_{i\in\{1,2,\dots,d\}}\pr{\lambda_0\lambda_{i}}\leq2\PRb{\frac{1+\alpha}{1-\alpha}}
			\end{equation}\cfout.
		\end{athm}
		\cfclear

		\begin{aproof}
				Throughout this proof 
				assume without loss of generality that $\alpha>0$ (cf.\ \cref{cor:SGD:stable} and \cref{cor:SGD:is:MOM})
			let $\vartheta=(\vartheta_1,\dots,\vartheta_d)\in\R^d$,
			for every 
			 $\lambda=(\lambda_1,\dots,\lambda_{d})\in[0,\infty)^{d}$ let $\mathscr{L}^{\lambda}\colon\R^d\to\R$ satisfy for all $\theta=(\theta_1,\dots,\theta_d)\in\R^d$ that
			\begin{equation}
				\llabel{eq:loss:quadr:stretched}
				\textstyle
				\mathscr{L}^{\lambda}(\theta)=\sum_{i=1}^d\tfrac{\lambda_i}{2}(\theta_i-\vartheta_i)^2,
			\end{equation}
			for every 
			$\gamma\in[0,\infty)$,
			$\lambda\in[0,\infty)^{d}$
			let
			$\Theta^{\gamma,\lambda}=\prb{\Theta_{n}^{\gamma,\lambda}}_{n\in\N_0}\colon\N_0\to\R^d$ satisfy for all
			$n\in\N$ that
			\begin{equation}
				\llabel{eq:all:dim:mom}
				\begin{split}
					\Theta_{n}^{\gamma,\lambda}
					&=
					\Theta_{n-1}^{\gamma,\lambda}
					-
					\gamma\Phi_n
					\prb{
						\operatorname{diag}(\lambda)\prb{\Theta_{0}^{\gamma,\lambda}-\vartheta},
						\operatorname{diag}(\lambda)\prb{\Theta_{1}^{\gamma,\lambda}-\vartheta},
						\dots,
						\operatorname{diag}(\lambda)\prb{\Theta_{n-1}^{\gamma,\lambda}-\vartheta}}
					,
				\end{split}
			\end{equation}
			and for every 
			$\gamma\in[0,\infty)$,
			$\lambda\in[0,\infty)^{d}$
			let $\m_{n}^{\gamma,\lambda}\in\R^d$, $n\in\N_0$, satisfy for all $n\in\N$ that
			\begin{equation}
				\llabel{eq:mom:def:for:mom}
				\m_{0}^{\gamma,\lambda}=0
				\qqandqq
				\m_{n}^{\gamma,\lambda}=\alpha\m_{n-1}^{\gamma,\lambda}+(1-\alpha)\prb{\nabla\mathscr{L}^{\lambda}}\prb{\Theta_{n-1}^{\gamma,\lambda}}.
			\end{equation}
			\startnewargseq
				\argument{
					\lref{eq:mom:def:for:mom};
					\lref{eq:loss:quadr:stretched};
					\lref{eq:all:dim:mom};
					the assumption that $(\Phi_n)_{n \in \N}$ is the $\alpha$-momentum optimizer on $\R^d$\cfadd{def:MOM};
					the fact that for all $\theta\in\R^d$, $\lambda\in[0,\infty)^d$ it holds that $\prb{\nabla\mathscr{L}^{\lambda}}(\theta)=\operatorname{diag}(\lambda)\prb{\theta-\vartheta}$}{that for all 
				$\gamma\in[0,\infty)$,
				$\lambda=(\lambda_1,\dots,\lambda_{d})\in[0,\infty)^{d}$, 
				$n\in\N$ it holds that
				\begin{equation}
					\llabel{eq:all:dim:MOM:2}
					\begin{split}
						\Theta_{n}^{\gamma,\lambda}
						&=\textstyle
						\Theta_{n-1}^{\gamma,\lambda}
						-
						\gamma (1-\alpha)\sum_{k=1}^n\alpha^{n-k}\operatorname{diag}(\lambda)\prb{\Theta_{k-1}^{\gamma,\lambda}-\vartheta}
						\\&=\textstyle
						\Theta_{n-1}^{\gamma,\lambda}
						-
						\gamma (1-\alpha)\sum_{k=1}^n\alpha^{n-k}\prb{\nabla\mathscr{L}^{\lambda}}\prb{\Theta_{k-1}^{\gamma,\lambda}}
						=\textstyle
						\Theta_{n-1}^{\gamma,\lambda}
						-
						\gamma \m_{n}^{\gamma,\lambda}
						.
					\end{split}
			\end{equation}}
			\argument{
				\lref{eq:all:dim:MOM:2};
				\cref{prop:mom-stable:diag}
				(applied with 
				$d\curvearrowleft d$,
				$\gamma \curvearrowleft \gamma$,
				$(\lambda_1,\lambda_2,\dots,\lambda_d)\curvearrowleft \pr{\frac{\lambda_1}{2},\frac{\lambda_2}{2},\dots\allowbreak,\frac{\lambda_d}{2}}$,
				$\alpha\curvearrowleft\alpha$,
				$\vartheta\curvearrowleft\vartheta$,
				$\fl\curvearrowleft\mathscr{L}^{\lambda}$,
				$\Theta\curvearrowleft\prb{\N_0\ni n\mapsto\Theta_{n}^{\gamma,\lambda}\in\R^d}$,
				$\m\curvearrowleft\prb{\N_0\ni n\mapsto\m_{n}^{\gamma,\lambda}\in\R^d}$ 
				for 
				$\gamma\in[0,\infty)$,
				$\lambda=(\lambda_1,\dots,\lambda_d)\in[0,\infty)^d$
				in the notation of \cref{prop:mom-stable:diag})}{that for all 
				$\gamma\in[0,\infty)$,
				$\lambda=(\lambda_1,\dots,\lambda_{d})\in[0,\infty)^{d}$
				with $\Theta_{0}^{\gamma,\lambda}\neq\vartheta$ it holds that
				\begin{equation}
					\llabel{eq:BUCH:SGD:result}
					\sup_{n \in\N}\norm{\Theta_{n}^{\gamma,\lambda}}
					\in
					\begin{cases}
						[0,\infty)&\colon \gamma\max\{\lambda_1,\lambda_2,\dots,\lambda_d\}\leq2\PR{\frac{1+\alpha}{1-\alpha}}\\
						\{\infty\}&\colon\gamma\max\{\lambda_1,\lambda_2,\dots,\lambda_d\}>2\PR{\frac{1+\alpha}{1-\alpha}}. 
					\end{cases}
			\end{equation}}
		\argument{
			\lref{eq:all:dim:mom};
			\lref{eq:all:dim:MOM:2};
			\lref{eq:BUCH:SGD:result};}[verbs=p]{that $(\Phi_n)_{n \in \N}$ \stable{\mathcal{A}}\!\! if and only if it holds for all
			$(\lambda_0,\lambda_1,\allowbreak\dots,\lambda_{d})\in\mathcal{A}$ that $\max_{i\in\{1,2,\dots,d\}}(\lambda_0\lambda_i)\leq 2\PR{\frac{1+\alpha}{1-\alpha}}$ }
		\end{aproof}

\begin{athm}{cor}{cor:stoch:MOM:bounded}
	Let 
	$d\in\N$, $\alpha\in[0,1)$, $\lambda=(\lambda_1,\dots,\lambda_d)\in[0,\infty)^d$,
	let $\Phi_n\colon(\R^d)^{n}\to\R^d$, $n \in \N$,\cfadd{def:MOM} be the $\alpha$-momentum optimizer, 
	let
	$\gamma\colon\N\to[0,\infty)$ satisfy that $\sup_{n\in\N}\gamma_n\leq \frac{1-\alpha}{\pr{1+2\alpha}\max\{1,\lambda_1,\lambda_2,\dots,\lambda_n\}}$,
	let
	$\batch\colon\N\to\N$ be a function,
	let $\fc\in[0,\infty)$,
	let $(\Omega,\cF,\P)$ be a probability space,
	for every $n,j\in\N$
	let $X_{n,j}\colon\Omega\to[-\fc,\fc]^d$ be a random variable,
	and
	let
	$ \Grad=\pr{\Grad^{(1)},\dots,\Grad^{(d)}}\colon\N\times\Omega \to \R^d$ 				
	and
	$ \Theta \colon\N_0\times\Omega \to \R^d$ 
	satisfy for all
	$n\in\N$ that
	\begin{equation}
		\llabel{eq:MOM:3}
		\begin{split}
			\textstyle
			\Grad_n=
			\frac{1}{\batch_n}\sum_{j=1}^{\batch_n}\diag(\lambda)\pr{\Theta_{n-1}-X_{n,j}}
			\qqandqq
			\Theta_n
			= 
			\Theta_{ n - 1 }
			-
			\gamma_n\Phi_n\prb{\Grad_1,\Grad_2,\dots,\Grad_n}
		\end{split}
	\end{equation}
	\cfload.
	Then there exists $c\in\R$ such that $\sup_{n\in\N_0}\norm{\Theta_n}\leq c\norm{\Theta_0}+ c$.
\end{athm}

\begin{aproof}
	Throughout this proof assume without loss of generality that $d=1$ and $\lambda\neq0$,
	let $\cst,\Cst\in(0,\infty)$ satisfy 
	$\cst=\frac{(1+2\alpha)\lambda}{1-\alpha}$ and $\Cst=\max\{1,\cst\}$.
	\argument{\lref{eq:MOM:3}}{that for all $n\in\N$ it holds that
		\begin{equation}
			\llabel{eq:gradient:boundaries}
			\begin{split}
				\textstyle
				\lambda\pr{\Theta_{n-1}-\fc}
				\textstyle=
				\frac{1}{\batch_n}\PRb{\sum_{j=1}^{\batch_n}\lambda\pr{\Theta_{n-1}-\fc}}
				&\textstyle\leq
				\frac{1}{\batch_n}\PRb{\sum_{j=1}^{\batch_n}\lambda\pr{\Theta_{n-1}-X_{n,j}}}
				=
				\Grad_k
				\\&\textstyle\leq
				\frac{1}{\batch_n}\PRb{\sum_{j=1}^{\batch_n}\lambda\pr{\Theta_{n-1}+\fc}}
				=
				\lambda\pr{\Theta_{n-1}+\fc}
				.
			\end{split}
	\end{equation}}
	\argument{
		\eqref{def:eq:MOM};
		the assumtion that $(\Phi_n)_{n\in\N}$ is the $\alpha$-momentum optimizer}{that for all $n\in\N$ it holds that
		\begin{equation}
			\llabel{eq:repr:change}
			\begin{split}
				\Theta_n
				=\Theta_{n-1}-\gamma_n\Phi_n\prb{\Grad_1,\Grad_2,\dots,\Grad_n}
				&=\textstyle\Theta_{n-1}-\gamma_n \PRb{\sum_{k=1}^n \pr{1-\alpha}\alpha^{n-k}\Grad_k}
				.
			\end{split}
	\end{equation}}
	\argument{\cref{prop:a_priori_bound_one_dim:Adam:1} (applied with
		$\alpha\curvearrowleft\alpha$,
		$\fc\curvearrowleft\fc$,
		$\nu\curvearrowleft\lambda$,
		$\mu\curvearrowleft\lambda$,
		$N\curvearrowleft 1$,
		$M\curvearrowleft M$
		for $M\in\N$
		in the notation of \cref{prop:a_priori_bound_one_dim:Adam:1})}{that for all $M\in\N$ it holds that
		\begin{equation}
			\llabel{eq:finite:steps:bound}
			\max_{n\in\N\cap[1,M]}\vass{\Theta_n}
			\leq
			4\fc+\frac{3\fc\alpha\lambda}{(1-\alpha)\lambda}
			+3\vass{\Theta_{0}}
			\leq
			\prbbb{4\fc+\frac{3\fc\alpha}{1-\alpha}
				+3}\pr{\vass{\Theta_{0}}+1}
	\end{equation}}
	\argument{\lref{eq:finite:steps:bound}}{that there exist $c\in\R$ that
		\begin{equation}
			\textstyle
			\sup_{n\in\N}\vass{\Theta_n}
			\leq c\pr{\vass{\Theta_{0}}+1}
			=c\vass{\Theta_{0}}+c
			.
	\end{equation}}
\end{aproof}

		\subsection{Asymptotic stability of the Nesterov optimizer}
		\label{subsection:5.5}
		
		In the literature there are several slightly modified variants how to describe the Nesterov optimizer that can be easily transferred to each other (cf.~\cite{Nesterov1983} and, \eg, \cite[Sections~6.4 and 7.5]{jentzen2023mathematical} and \cite{Sutskever2013}). In the next notion, \cref{def:Nest} below, we recall one of these variants of the Nesterov optimizer (cf., \eg, \cite[Definition~6.4.22]{jentzen2023mathematical}) and in \cref{lem:Buch_Nest:agrees} and \cref{lem:Buch_Nest:agrees:simp} below we briefly recall how this variant of is related to other variants of the Nesterov method in the literature.

	\begin{definition}[Nesterov optimizer]
		\label{def:Nest}
		Let $d\in\N$, $\alpha\in[0,1]$ and
		let 
		$\Phi_n\colon\allowbreak
		(\R^d)^{n}\to\R^d$, $n \in \N$, be functions.
		Then we say that $(\Phi_n)_{n \in \N}$ is the $\alpha$-Nesterov optimizer on $\R^d$ 
		(we say that $(\Phi_n)_{n \in \N}$ is the $\alpha$-Nesterov optimizer)
		if and only if it holds
		for all 
		$n\in\N$,
		$g_1,g_2,\dots,g_n\in\R^d$
		that
		\begin{equation}
			\label{def:eq:Nest}
			\Phi_{n}(g_1,g_2,\dots,g_n)
			=
			g_n+\sum_{k=1}^n \alpha^{n+1-k}g_k
			.
		\end{equation}
	\end{definition}

	\begin{athm}{lemma}{lem:Buch_Nest:agrees}
		Let $d\in\N$, $\alpha\in(0,1)$, $\gamma,\Gamma\in(0,\infty)$ satisfy $\Gamma=\gamma(1-\alpha)$,
		let
		$\mathcal{G}\colon\R^d\to\R^d$,
		$\Theta\colon\N_0\to\R^d$, and $\m\colon\N_0\to\R^d$ satisfy for all $n\in\N$ that
		\begin{equation}
			\llabel{eq:Nest:book}
			\m_0=0,
			\quad
			\m_n = \alpha\m_{n-1} + (1-\alpha) \mathcal{G}\pr{\Theta_{n-1}-\gamma\alpha\m_{n-1}},
			\qandq
			\Theta_n = \Theta_{n-1} - \gamma \m_n,
		\end{equation}
		and let 
			$\Psi\colon\N_0\to\R^d$ and $M\colon\N_0\to\R^d$ satisfy for all $n\in\N$ that
		\begin{equation}
			\llabel{eq:Nest:here}
			\Psi_0=M_0=\Theta_0,
			\quad
			M_n = (1 + a) \Psi_n - \alpha\Psi_{n-1},
			\qqandqq
			\Psi_n = M_{n-1} - \Gamma\mathcal{G}(M_{n-1}).
		\end{equation}
		Then $\Theta=\Psi$.
	\end{athm}
	
	\begin{aproof}
				\argument{\lref{eq:Nest:book};\lref{eq:Nest:here}}{that
			\begin{equation}
				\llabel{eq:n=1}
				\Psi_1
				=M_0-\Gamma\mathcal{G}(M_0)
				=\Theta_0-\gamma(\alpha\m_0+(1-\alpha)\mathcal{G}(\Theta_0))
				=\Theta_0-\gamma\m_1=\Theta_1.
			\end{equation}
		}
		\argument{\lref{eq:Nest:book};}{that for all $n\in\N$ it holds that
		\begin{equation}
			\llabel{eq:help:recursion}
			(1+\alpha)\Theta_{n}-\alpha\Theta_{n-1}
			=\Theta_n-\alpha(\Theta_{n-1}-\Theta_n)
			=\Theta_n-\alpha\gamma\m_n.
			\end{equation}}
		\argument{
			\lref{eq:help:recursion};
			\lref{eq:Nest:book};
			\lref{eq:Nest:here};}{that for all $n\in\N\cap(1,\infty)$ with $\forall\, m\in\N\cap(1,n)\colon\Theta_m=\Psi_m$ it holds that
		\begin{equation}
			\llabel{eq:other:ns}
			\begin{split}
				\Psi_n
				&=M_{n-1} - \Gamma \mathcal{G} ( M_{n-1} )
				\\&=(1 + \alpha) \Psi_{n-1} - \alpha\Psi_{n-2} - \Gamma \mathcal{G} ( (1 + \alpha) \Psi_{n-1} - \alpha\Psi_{n-2} )
				\\&=(1 + \alpha) \Theta_{n-1} - \alpha\Theta_{n-2} - \gamma_n\pr{1-\alpha} \mathcal{G} ( (1 + \alpha) \Theta_{n-1} - \alpha\Theta_{n-2} )
				\\&=\Theta_{n-1} -\alpha\gamma\m_{n-1} - \gamma\pr{1-\alpha} \mathcal{G} (\Theta_{n-1} -\gamma\alpha\m_{n-1} )
				\\&=\Theta_{n-1} -\alpha\gamma\m_{n-1} - \gamma\pr{1-\alpha} \mathcal{G} (\Theta_{n-1} -\gamma\alpha\m_{n-1} )
				\\&=\Theta_{n-1} -\gamma\pr{\alpha\m_{n-1} + \pr{1-\alpha} \mathcal{G} (\Theta_{n-1} -\gamma\alpha\m_{n-1} )}
				=\Theta_{n-1} -\gamma\m_n=\Theta_n
				.
			\end{split}
		\end{equation}
	}
	\argument{
		\lref{eq:n=1};
		\lref{eq:other:ns};}{that for all $n\in\N$ it holds that $\Theta_n=\Psi_n$}.
	\end{aproof}

		\newcommand{\mm}{\mathfrak{m}}
	\begin{athm}{lemma}{lem:Buch_Nest:agrees:simp}
		Let 
		$d\in\N$, 
		$\alpha\in(0,1)$, 
		$\gamma,\Gamma\in(0,\infty)$, 
		$\mathcal{G}\in C(\R^d,\R^d)$ satisfy $\Gamma=\gamma(1-\alpha)$, 
		let $\Theta\colon\N_0\to\R^d$ and $\m\colon\N_0\to\R^d$ satisfy for all $n\in\N$ that
		\begin{equation}
			\llabel{eq:Nest:book:V1}
			\m_0=0,
			\quad
			\m_n=\alpha\m_{n-1}+(1-\alpha)\mathcal{G}(\Theta_{n-1}-\gamma\alpha\m_{n-1}),
			\qandq
			\Theta_n=\Theta_{n-1}-\gamma\alpha\m_n,
		\end{equation}
		let $\Psi\colon\N_0\to\R^d$ and $\mm\colon\N_0\to\R^d$ satisfy for all $n\in\N$ that
		\begin{equation}
			\llabel{eq:Nest:book}
			\mm_0=0,
			\quad
			\mm_n = \alpha\mm_{n-1} + \mathcal{G}\pr{\Psi_{n-1}},
			\qandq
			\Psi_n = \Psi_{n-1} 
			-\Gamma\alpha\mm_n
			- \Gamma\mathcal{G}(\Psi_{n-1}),
		\end{equation}
		and let 
		$\varTheta\colon\N_0\to\R^d$ and $M\colon\N_0\to\R^d$ satisfy for all $n\in\N$ that
		\begin{equation}
			\llabel{eq:Nest:here}
			M_0=\varTheta_0=\Psi_0=\Theta_0,
			\quad
			M_n = (1 + \alpha) \varTheta_n - \alpha\varTheta_{n-1},
			\qandq
			\varTheta_n = M_{n-1} - \gamma\mathcal{G}(M_{n-1}).
		\end{equation}
		Then 
		\begin{enumerate}[label=(\roman*)]
			\item\llabel{it:relation} it holds for all $n\in\N_0$ that $\Theta_n=\Psi_n+\gamma\alpha\m_n=\varTheta_n$ and $(1-\alpha)\mm_n=\m_n$,
			\item\llabel{it:asymp} it holds that $\sup_{n\in\N}\norm{\Theta_n}<\infty$ if and only if $\sup_{n\in\N}\norm{\Psi_n}<\infty$, and
			\item\llabel{it:repr} it holds for all $n\in\N$ that $\Psi_n=\Psi_{n-1}-\Gamma\PRb{\mathcal{G}(\Psi_{n-1})+\sum_{k=1}^n\alpha^{n+1-k}\mathcal{G}(\Psi_{k-1})}$.
		\end{enumerate}
	\end{athm}
	
	\begin{aproof}
		\argument{
			\lref{eq:Nest:book:V1};
			\lref{eq:Nest:book};\!\!, \eg,
			\cite[Lemma~6.4.21]{jentzen2023mathematical} (applied with
			$\m\curvearrowleft\m$,
			$\Theta\curvearrowleft\Theta$,
			$\gamma\curvearrowleft\pr{\N\ni n\mapsto \gamma\in[0,\infty)}$,
			$\Psi\curvearrowleft\Psi$,
			$\beta\curvearrowleft\pr{\N\ni n\mapsto \alpha\in[0,\infty)}$,
			$\delta\curvearrowleft\pr{\N\ni n\mapsto \Gamma\in[0,\infty)}$,
			in the notation of \cite[Lemma~7.5.2]{jentzen2023mathematical})}{that for all $n\in\N$ it holds that
		\begin{equation}
			\llabel{eq:booklike:setup}
					\Psi_n=\Theta_n-\gamma\alpha\m_n
					\qqandqq
					(1-\alpha)\mm_n=\m_n.
			\end{equation}}
		\argument{
			\lref{eq:Nest:book:V1};
			\lref{eq:Nest:here};
			\cref{lem:Buch_Nest:agrees} (applied with 
			$d\curvearrowleft d$,
			$\alpha\curvearrowleft \alpha$,
			$\gamma\curvearrowleft \gamma$,
			$\Gamma\curvearrowleft d$,
			$\mathcal{G}\curvearrowleft \mathcal{G}$,
			$\Theta\curvearrowleft \Theta$,
			$\m\curvearrowleft \m$,
			$\Psi\curvearrowleft \varTheta$,
			$M\curvearrowleft M$
			in the notation of \cref{lem:Buch_Nest:agrees})}{that for all $n\in\N_0$ it holds that
		\begin{equation}
			\llabel{eq:prev:L}
			\Theta_n=\varTheta_n.
		\end{equation}
	}
	\argument{\lref{eq:prev:L};\lref{eq:booklike:setup}}{\lref{it:relation}}.
	\argument{\lref{eq:Nest:book:V1}}{that 
	\begin{equation}
		\llabel{eq:Theta:mom:bounded}
		\textstyle
		\sup_{n\in\N}\norm{\gamma\alpha\m_n}
		=
		\sup_{n\in\N}\norm{\Theta_{n}-\Theta_{n-1}}
		\leq
		2 \PRb{\sup_{n\in\N}\norm{\Theta_n}} 
		.
		\end{equation}}
	\argument{
		\lref{eq:booklike:setup};
		\lref{eq:Theta:mom:bounded};
		}{that
		\begin{equation}
			\llabel{eq:concl:Psi:bounded}
			\begin{split}
			\textstyle
			\sup_{n\in\N}\norm{\Psi_n}
			=
			\sup_{n\in\N}\norm{\Theta_n-\gamma\alpha\m_n}
			&\textstyle\leq \PRb{\sup_{n\in\N}\norm{\Theta_n}}+\PRb{\sup_{n\in\N}\norm{\gamma\alpha\m_n}}
			\\&\textstyle\leq
			3 \PRb{\sup_{n\in\N}\norm{\Theta_n} }
			.
			\end{split}
			\end{equation}}
		\argument{\lref{eq:Nest:book};}{that 
			\begin{equation}
				\llabel{eq:Psi:mom:bounded}
				\begin{split}
				\textstyle
				\sup_{n\in\N}\norm{\gamma\alpha\m_n}
				=
				\sup_{n\in\N}\norm{\gamma\alpha(1-\alpha)\mm_n}
				&=\textstyle
				\sup_{n\in\N}\norm{\Gamma\alpha\mm_n}
				\\&=\textstyle
				\sup_{n\in\N}\norm{\Psi_{n-1}-\Psi_n-\Gamma\mathcal{G}(\Psi_{n-1})}
				\\&\leq\textstyle
				\sup_{n\in\N}\PRb{\norm{\Psi_{n-1}}+\norm{\Psi_n}+\Gamma\norm{\mathcal{G}(\Psi_{n-1})}}
				\\&\leq\textstyle
				2\PRb{\sup_{n\in\N}\norm{\Psi_{n}}}+\Gamma\PRb{\sup_{n\in\N}\norm{\mathcal{G}(\Psi_{n})}}
				.
				\end{split}
		\end{equation}}
		\argument{
			\lref{eq:booklike:setup};
			\lref{eq:Psi:mom:bounded};
		}{that
			\begin{equation}
				\llabel{eq:concl:Theta:bounded}
				\begin{split}
				\textstyle
				\sup_{n\in\N}\norm{\Theta_n}
				=
				\sup_{n\in\N}\norm{\Psi_n+\gamma\alpha\m_n}
				&\textstyle\leq \PRb{\sup_{n\in\N}\norm{\Psi_n}}+\PRb{\sup_{n\in\N}\norm{\gamma\alpha\m_n}}
				\\&\textstyle\leq
				3\PRb{\sup_{n\in\N}\norm{\Psi_n}}
				+
				\Gamma\PRb{\sup_{n\in\N}\norm{\mathcal{G}(\Psi_{n})}}
				.
				\end{split}
		\end{equation}}
	\argument{\lref{eq:concl:Theta:bounded};\lref{eq:concl:Psi:bounded};}{that
	\begin{equation}
		\llabel{eq:it:2:sol}
		\sup_{n\in\N}\norm{\Psi_n}
		\leq
		3 \PRbbb{\sup_{n\in\N}\norm{\Theta_n}}
		\qqandqq
		\sup_{n\in\N}\norm{\Theta_n}
		\leq
		3\PRbbb{\sup_{n\in\N}\norm{\Psi_n}}
		+
		\Gamma\PRbbb{\sup_{n\in\N}\norm{\mathcal{G}(\Psi_{n})}}
		.
	\end{equation}}
		\argument{\lref{eq:it:2:sol};the fact that $\mathcal{G}\in C(\R^d,\R^d)$;}{\lref{it:asymp}}.
		\argument{\lref{eq:Nest:book};induction}{that for all $n\in\N$ it holds that
		\begin{equation}
			\llabel{eq:explicit:mom}
			\begin{split}
			\mm_n
			&=\alpha\mm_{n-1}+\mathcal{G}(\Psi_{n-1})
			\\&=\alpha^2\mm_{n-2}+\alpha\mathcal{G}(\Psi_{n-2})+\mathcal{G}(\Psi_{n-1})
			\\&=\dots
			\\&\textstyle=\alpha^n\mm_{0}+\sum_{k=1}^n\alpha^{n-k}\mathcal{G}(\Psi_{k-1})
			=\sum_{k=1}^n\alpha^{n-k}\mathcal{G}(\Psi_{k-1}).
			\end{split}
			\end{equation}}
		\argument{\lref{eq:explicit:mom};\lref{eq:Nest:book}}{that for all $n\in\N$ it holds that
		\begin{equation}
			\llabel{eq:repr:Nest}
			\begin{split}
			\Psi_n
			=\Psi_{n-1}-\Gamma\alpha\mm_n-\Gamma\mathcal{G}(\Psi_{n-1})
			&\textstyle=\Psi_{n-1}-\Gamma\alpha\PRb{\sum_{k=1}^n\alpha^{n-k}\mathcal{G}(\Psi_{k-1})}-\Gamma\mathcal{G}(\Psi_{n-1})
			\\&\textstyle=\Psi_{n-1}-\Gamma\PRb{\mathcal{G}(\Psi_{n-1})+\sum_{k=1}^n\alpha^{n+1-k}\mathcal{G}(\Psi_{k-1})}	
			.			
			\end{split}
			\end{equation}}
		\argument{\lref{eq:repr:Nest}}{\lref{it:repr}}.
	\end{aproof}

			In \cref{def:Nest} above we recall for every $\alpha\in[0,1]$ the concept of the $\alpha$-Nesterov optimizer. In the following statement, \cref{cor:0-Nest:is:SGD} below, we  briefly recall the elementary fact that the $0$-Nesterov optimizer coincides with the standard \GD\ method.

			\begin{athm}{lemma}{cor:0-Nest:is:SGD}
				Let $d\in\N$ and
				let $\Phi_n\colon(\R^d)^{n}\to\R^d$, $n \in \N$, be the $0$-Nesterov optimizer on $\R^d$\cfadd{def:Nest} \cfload.
				Then $(\Phi_n)_{n \in \N}$ is the \GD\ optimizer on $\R^d$\cfadd{def:SGD} \cfout.
	\end{athm}
	
		\begin{aproof}
		\argument{\eqref{def:eq:SGD};\eqref{def:eq:Nest}}{for all 
			$n\in\N$,
			$g_1,g_2,\dots,g_{n}\in\R^d$
			that
			\begin{equation}				
				\Psi_n(g_1,g_2,\dots,g_n)
				=
				g_n
				=
				\Phi_n(g_1,g_2,\dots,g_n)
				.
			\end{equation}
		}
	\end{aproof}
	
		\begin{athm}{cor}{cor:Nesterov:stable}
		Let $d\in\N$, $\alpha\in[0,1)$,
		let $\mathcal{A}\subseteq[0,\infty)^{d+1}$ be a set, and let $\Phi_n\colon(\R^d)^{n}\to\R^d$, $n \in \N$, be the $\alpha$-Nesterov optimizer\cfadd{def:Nest} \cfload.
		Then $(\Phi_n)_{n \in \N}$ \stable{\mathcal{A}}\!\! if and only if it holds for all $(\lambda_0,\lambda_1,\dots,\lambda_{d})\in\mathcal{A}$ that \begin{equation}
			\textstyle\max_{i\in\{1,2,\dots,d\}}\pr{\lambda_0\lambda_{i}}\leq2\PRb{\frac{1-\alpha^2 }{1+2\alpha}}
		\end{equation}\cfout.
	\end{athm}
	\cfclear

		\begin{aproof}
		Throughout this proof 
		assume without loss of generality that $\alpha>0$ (cf.\ \cref{cor:SGD:stable} and \cref{cor:0-Nest:is:SGD})
		let $\vartheta=(\vartheta_1,\dots,\vartheta_d)\in\R^d$,
		for every 
		 $\lambda=(\lambda_1,\dots,\lambda_{d})\in[0,\infty)^{d}$ let $\mathscr{L}^{\lambda}\colon\R^d\to\R$ satisfy for all $\theta=(\theta_1,\dots,\theta_d)\in\R^d$ that
		\begin{equation}
			\llabel{eq:loss:quadr:stretched}
			\textstyle
			\mathscr{L}^{\lambda}(\theta)=\sum_{i=1}^d\tfrac{\lambda_i}{2}(\theta-\vartheta_i)^2,
		\end{equation}
		for every 
		$\gamma\in[0,\infty)$,
		$\lambda\in[0,\infty)^{d}$
		let
		$\Psi_{n}^{\gamma,\lambda}\in\R^d$ satisfy for all
		$n\in\N$ that
		\begin{equation}
			\llabel{eq:all:dim:Nest}
			\begin{split}
				\Psi_{n}^{\gamma,\lambda}
				&=
				\Psi_{n-1}^{\gamma,\lambda}
				-
				\gamma\Phi_n
				\prb{
					\operatorname{diag}(\lambda)\prb{\Psi_{0}^{\gamma,\lambda}-\vartheta},
					\operatorname{diag}(\lambda)\prb{\Psi_{1}^{\gamma,\lambda}-\vartheta},
					\dots,
					\operatorname{diag}(\lambda)\prb{\Psi_{n-1}^{\gamma,\lambda}-\vartheta}}
				,
			\end{split}
		\end{equation}
		for every 
		$\gamma\in[0,\infty)$,
		$\lambda\in[0,\infty)^{d}$
		let $\mm_{n}^{\gamma,\lambda}\in\R^d$, $n\in\N_0$, satisfy for all $n\in\N$ that
		\begin{equation}
			\llabel{eq:Nest:def:for:Nest}
			\mm_{0}^{\gamma,\lambda}=0
			\qqandqq
			\mm_{n}^{\gamma,\lambda}=\alpha\mm_{n-1}^{\gamma,\lambda}+\prb{\nabla\mathscr{L}^{\lambda}}\prb{\Psi_{n-1}^{\gamma,\lambda}},
		\end{equation}
		and for every 
		$\gamma\in[0,\infty)$,
		$\lambda\in[0,\infty)^{d}$
		let $\Theta_{n}^{\gamma,\lambda}\in\R^d$, $n\in\N_0$, and $M_{n}^{\gamma,\lambda}\in\R^d$, $n\in\N_0$, satisfy for all $n\in\N$ that
		\begin{equation}
			\llabel{eq:Nest:def:for:Nest:2}
			M_{0}^{\gamma,\lambda}=\Theta_{0}^{\gamma,\lambda}=\Psi_{0}^{\gamma,\lambda},
			\qquad\qquad
			M_{n}^{\gamma,\lambda}=(1+\alpha)	\Theta_{n}^{\gamma,\lambda}-\Theta_{n-1}^{\gamma,\lambda},
		\end{equation}
		\begin{equation}
			\llabel{eq:Nest:def:for:Nest:3}
			\text{and}\qquad\qquad
			\Theta_{n}^{\gamma,\lambda}=M_{n-1}^{\gamma,\lambda}-\gamma\prb{\nabla\mathscr{L}^{\lambda}}\prb{M_{n-1}^{\gamma,\lambda}}
			.
		\end{equation}
		\startnewargseq
		\argument{
			\lref{eq:Nest:def:for:Nest:2};
			\lref{eq:Nest:def:for:Nest:3};
			\cref{prop:Nest-stable:diag} (applied with
			$d\curvearrowleft d$,
			$(\lambda_1,\lambda_2,\dots,\lambda_d)\curvearrowleft \pr{\frac{\lambda_1}{2},\frac{\lambda_2}{2},\dots,\frac{\lambda_d}{2}}$,
			$\alpha\curvearrowleft\alpha$,
			$\gamma \curvearrowleft \gamma$,
			$\vartheta\curvearrowleft\vartheta$,
			$\fl\curvearrowleft\mathscr{L}^{\lambda}$,
			$\Theta\curvearrowleft\prb{\N_0\ni n\mapsto\Theta_{n}^{\gamma,\lambda}\in\R^d}$,
			$\m\curvearrowleft\prb{\N_0\ni n\mapsto\m_{n}^{\gamma,\lambda}\in\R^d}$ 
			for 
			$\gamma\in[0,\infty)$,
			$\lambda=(\lambda_1,\dots,\lambda_d)\in[0,\infty)^d$
			in the notation of \cref{prop:Nest-stable:diag})}{that for all
			$\gamma\in[0,\infty)$,
			$\lambda=(\lambda_1,\dots,\lambda_{d})\in[0,\infty)^{d}$ with  $\Theta_{0}^{\gamma,\lambda}\neq\vartheta$ it holds that
				\begin{equation}
					\llabel{eq:prop:Nest:res}
										\textstyle
					\sup_{n\in\N}\norm{\Theta_{n}^{\gamma,\lambda}}
					\in
					\begin{cases}
						[0,\infty)
						& \colon \gamma\max\{\lambda_1,\lambda_2,\dots,\lambda_d\} \leq 2\PRb{\tfrac{1 + \alpha }{   1 +2 \alpha }} \\
						\{\infty\} 
						& \colon \gamma\max\{\lambda_1,\lambda_2,\dots,\lambda_d\} >  2\PRb{\tfrac{1 + \alpha }{   1 +2 \alpha }}.
					\end{cases}
				\end{equation}
			}
					\argument{
				\lref{eq:Nest:def:for:Nest:2};
				\lref{eq:Nest:def:for:Nest:3};
				\lref{eq:prop:Nest:res};
				}{that for all 
				$\gamma\in[0,\infty)$,
				$\lambda=(\lambda_1,\dots,\lambda_{d})\in[0,\infty)^{d}$
				with $\Theta_{0}^{\gamma,\lambda}\neq\vartheta$ it holds that
				\begin{equation}
					\llabel{eq:BUCH:SGD:result:2}
					\textstyle
					\sup_{n\in\N}\norm{\Theta_{n}^{\gamma(1-\alpha)^{-1},\lambda}}
					\in
					\begin{cases}
						[0,\infty)
						& \colon \gamma\max\{\lambda_1,\lambda_2,\dots,\lambda_d\} \leq 2\PRb{\tfrac{1 - \alpha^2 }{   1 +2 \alpha }} \\
						\{\infty\} 
						& \colon \gamma\max\{\lambda_1,\lambda_2,\dots,\lambda_d\} >  2\PRb{\tfrac{1 - \alpha^2 }{   1 +2 \alpha }}.
					\end{cases}
			\end{equation}}
		\argument{\lref{eq:Nest:def:for:Nest};\cref{momentum:representation}}{for all
			$\gamma\in[0,\infty)$,
			$\lambda\in[0,\infty)^{d}$, 
			$n\in\N$ that
			\begin{equation}
				\llabel{eq:Nest:explicit}
				\textstyle
				\mm_{n}^{\gamma,\lambda}=\sum_{k=1}^n\alpha^{n-k}\prb{\nabla\mathscr{L}^{\lambda}}\prb{\Psi_{k-1}^{\gamma,\lambda}}.
		\end{equation}}
		\argument{\lref{eq:Nest:explicit};\lref{eq:loss:quadr:stretched};\lref{eq:all:dim:Nest};the assumption that $(\Phi_n)_{n \in \N}$ is the $\alpha$-Nesterov optimizer on $\R^d$\cfadd{def:Nest};}{that for all 
			$\gamma\in[0,\infty)$,
			$\lambda\in[0,\infty)^{d}$, 
			$n\in\N$ it holds that
			\begin{equation}
				\llabel{eq:all:dim:Nest:2}
				\begin{split}
					\Psi_{n}^{\gamma,\lambda}
					&=\textstyle
					\Psi_{n-1}^{\gamma,\lambda}
					-
					\gamma \PRb{\operatorname{diag}(\lambda)\prb{\Psi_{n-1}^{\gamma,\lambda}-\vartheta}
					+\sum_{k=1}^n \alpha^{n+1-k}\operatorname{diag}(\lambda)\prb{\Psi_{k-1}^{\gamma,\lambda}-\vartheta}}
					\\&=\textstyle
					\Psi_{n-1}^{\gamma,\lambda}
					-
					\gamma \PRb{\prb{\nabla\mathscr{L}^{\lambda}}\prb{\Psi_{n-1}^{\gamma,\lambda}}
						+\sum_{k=1}^n \alpha^{n+1-k}\prb{\nabla\mathscr{L}^{\lambda}}\prb{\Psi_{k-1}^{\gamma,\lambda}}}
					\\&=\textstyle
						\Psi_{n-1}^{\gamma,\lambda}
						-
						\gamma\alpha\PRb{\sum_{k=1}^n \alpha^{n-k}\prb{\nabla\mathscr{L}^{\lambda}}\prb{\Psi_{k-1}^{\gamma,\lambda}}}
						-\gamma\prb{\nabla\mathscr{L}^{\lambda}}\prb{\Psi_{n-1}^{\gamma,\lambda}}
					\\&=\textstyle
					\Psi_{n-1}^{\gamma,\lambda}
					-
					\gamma\alpha\mm_{n}^{\gamma,\lambda}
					-\gamma\prb{\nabla\mathscr{L}^{\lambda}}\prb{\Psi_{n-1}^{\gamma,\lambda}}
					.
				\end{split}
		\end{equation}}
	\argument{
		\lref{eq:all:dim:Nest:2};
		\lref{eq:Nest:def:for:Nest};
		\lref{eq:Nest:def:for:Nest:2};
		\lref{eq:Nest:def:for:Nest:3};
		\cref{lem:Buch_Nest:agrees:simp} (applied with 
		$\Psi\curvearrowleft\prb{\N_0\ni n\mapsto \Psi_{n}^{\gamma(1-\alpha),\lambda}\in\R^d}$,
		$\mm\curvearrowleft\prb{\N_0\ni n\mapsto \mm_{n}^{\gamma,\lambda}\in\R^d}$,
		$\Gamma\curvearrowleft \gamma(1-\alpha)$,
		$\varTheta\curvearrowleft \prb{\N_0\ni n\mapsto \Theta_{n}^{\gamma,\lambda}\in\R^d}$,
		$M\curvearrowleft \prb{\N_0\ni n\mapsto M_{n}^{\gamma,\lambda}\in\R^d}$,
		$\gamma\curvearrowleft\gamma$
		for
		$\gamma\in(0,\infty)$,
		$\lambda=(\lambda_1,\dots,\lambda_{d})\in[0,\infty)^{d}$
		in the notation of \cref{lem:Buch_Nest:agrees:simp})}{that for all
		$\gamma\in(0,\infty)$,
		$\lambda=(\lambda_1,\dots,\lambda_{d})\in[0,\infty)^{d}$ it holds that
		\begin{equation}
			\llabel{eq:transf:of:Nest}
			\textstyle
			\PRb{\limsup_{n\to\infty}\norm{\Theta_{n}^{\gamma,\lambda}}<\infty}
			\Leftrightarrow
			\PRb{\limsup_{n\to\infty}\norm{\Psi_{n}^{\gamma(1-\alpha),\lambda}}<\infty}
			.
		\end{equation}
		}
			\argument{\lref{eq:transf:of:Nest};the assumption that $\alpha<1$}{
			that for all
			$\gamma\in(0,\infty)$,
			$\lambda=(\lambda_1,\dots,\lambda_{d})\in[0,\infty)^{d}$
			 it holds that
			\begin{equation}
				\llabel{eq:transf:of:Nest.1}
				\textstyle
				\PRb{\limsup_{n\to\infty}\norm{\Theta_{n}^{\gamma(1-\alpha)^{-1},\lambda}}<\infty}
				\Leftrightarrow
				\PRb{\limsup_{n\to\infty}\norm{\Psi_{n}^{\gamma,\lambda}}<\infty}
				.
			\end{equation}
		}
		\argument{
			\lref{eq:BUCH:SGD:result:2};
			\lref{eq:transf:of:Nest.1};}{
			that for all
			$\gamma\in(0,\infty)$,
			$\lambda=(\lambda_1,\dots,\lambda_{d})\in[0,\infty)^{d}$, 
			with $\Theta_{0}^{\gamma,\lambda}\neq\vartheta$ it holds that
			\begin{equation}
				\llabel{eq:BUCH:Nest:result:2}
				\textstyle
				\sup_{n\in\N}\norm{\Psi_{n}^{\gamma,\lambda}}
				\in
				\begin{cases}
					[0,\infty)
					& \colon \gamma\max\{\lambda_1,\lambda_2,\dots,\lambda_d\} \leq 2\PRb{\tfrac{1 - \alpha^2 }{   1 +2 \alpha }} \\
					\{\infty\} 
					& \colon \gamma\max\{\lambda_1,\lambda_2,\dots,\lambda_d\} >  2\PRb{\tfrac{1 - \alpha^2 }{   1 +2 \alpha }}.
				\end{cases}
				\end{equation}}
				\argument{
		\lref{eq:Nest:def:for:Nest:2};
		\lref{eq:all:dim:Nest:2};
	}{that for all 
		$\gamma\in[0,\infty)$,
		$\lambda=(\lambda_1,\dots,\lambda_{d})\in[0,\infty)^{d}$,  
		$n\in\N$ with $\gamma(\Theta_{0}^{\gamma,\lambda}-\vartheta)=0$ it holds that
		\begin{equation}
			\llabel{eq:edge:cases}
			\Psi_{n}^{\gamma,\lambda}=\Psi_{0}^{\gamma,\lambda}.
	\end{equation}}
				\argument{
	\lref{eq:edge:cases};
	\lref{eq:BUCH:Nest:result:2};
}{that for all 
	$\gamma\in[0,\infty)$,
	$\lambda=(\lambda_1,\dots,\lambda_{d})\in[0,\infty)^{d}$ with $\Theta_{0}^{\gamma,\lambda}\neq\vartheta$ it holds that
	\begin{equation}
		\llabel{eq:all:cases}
		\textstyle
		\sup_{n\in\N}\norm{\Psi_{n}^{\gamma,\lambda}}
		\in
		\begin{cases}
			[0,\infty)
			& \colon \gamma\max\{\lambda_1,\lambda_2,\dots,\lambda_d\} \leq 2\PRb{\tfrac{1 - \alpha^2 }{   1 +2 \alpha }} \\
			\{\infty\} 
			& \colon \gamma\max\{\lambda_1,\lambda_2,\dots,\lambda_d\} >  2\PRb{\tfrac{1 - \alpha^2 }{   1 +2 \alpha }}.
		\end{cases}
\end{equation}}
		\argument{
			\lref{eq:all:dim:Nest};
			\lref{eq:all:cases};}[verbs=e]{that $(\Phi_n)_{n \in \N}$ \stable{\mathcal{A}}\!\! if and only if it holds for all
			$(\lambda_0,\lambda_1,\allowbreak\dots,\lambda_{d})\in\mathcal{A}$ that $\max_{i\in\{1,2,\dots,d\}}(\lambda_0\lambda_i)\leq 2\PRb{\frac{1-\alpha^2 }{1+2\alpha}}$ }
	\end{aproof}
		
\subsection{Asymptotic stability properties for deep learning optimizers}
\label{subsection:5.6}

\begin{figure}[H]
			\centering
			\begin{tikzpicture}[scale=.7,samples=2000]
				\tznode(-1.5,4){Largest eigenvalue $\max\{\lambda_1,\dots,\lambda_d\}$}[rotate=90]
				\tznode(10,-1.2){Learning rate $\gamma$}
				\tzticks{5,10,...,20} 
				{2,4,...,8} 
				\def\sgdgraph{2/\x} 
				\def\momentumgraph{38/\x} 
				\def\momentumgraphs{6/\x} 
				\fill [blue!30, domain=4.75:20, variable=\x] 
				plot ({\x}, {38/\x})
				-- (20,0)
				-- (0,0)
				-- (0,8)
				-- cycle;
				\fill [violet!50, domain=0.75:20, variable=\x]
				plot ({\x}, {6/\x})
				-- (20,0)
				-- (0,0)
				-- (0,8)
				-- cycle;
				\fill [red!30, domain=0.25:20, variable=\x]
				plot ({\x}, {2/\x})
				-- (20,0)
				-- (0,0)
				-- (0,8)
				-- cycle;
				\tzfn[red,thick]"sgdgraph"{\sgdgraph}[0.25:20]{}
				\tzfn[blue,thick]"momentumgraph"{\momentumgraph}[4.75:20]{}
				\tzfn[violet,thick]"momentumgraphs"{\momentumgraphs}[0.75:20]{}
				\tzpath+[pattern=north east lines,opacity=.1]
				(0,0)(0,8)(20,0)(0,-8);
				\tznode(11,6){$\alpha$-$\beta$-$\eps$-Adam}
				\tznode(1.7,0.3){GD/$0$-mom.}[red]
				\tznode(6,3){$0.9$-mom.}[blue]
				\tznode(3.5,1.1){$0.5$-mom.}[violet]
				\tzaxes(-0.5,-0.5)(20.5,8.5)
			\end{tikzpicture}
			\caption{In this figure we graphically represent for every 
				$\alpha\in[0,1)$, $\beta\in(\alpha^2,1)$, $\eps\in(0,\infty)$ the stability region of
				the \GD\ optimizer (the $0$-momentum optimizer), 
				the $0.5$-momentum optimizer,  
				the $0.9$-momentum optimizer, and
				the $\alpha$-$\beta$-$\eps$-\Adam\ optimizer (standard axis).}
			\label{fig:comp:without:Nets}
		\end{figure}
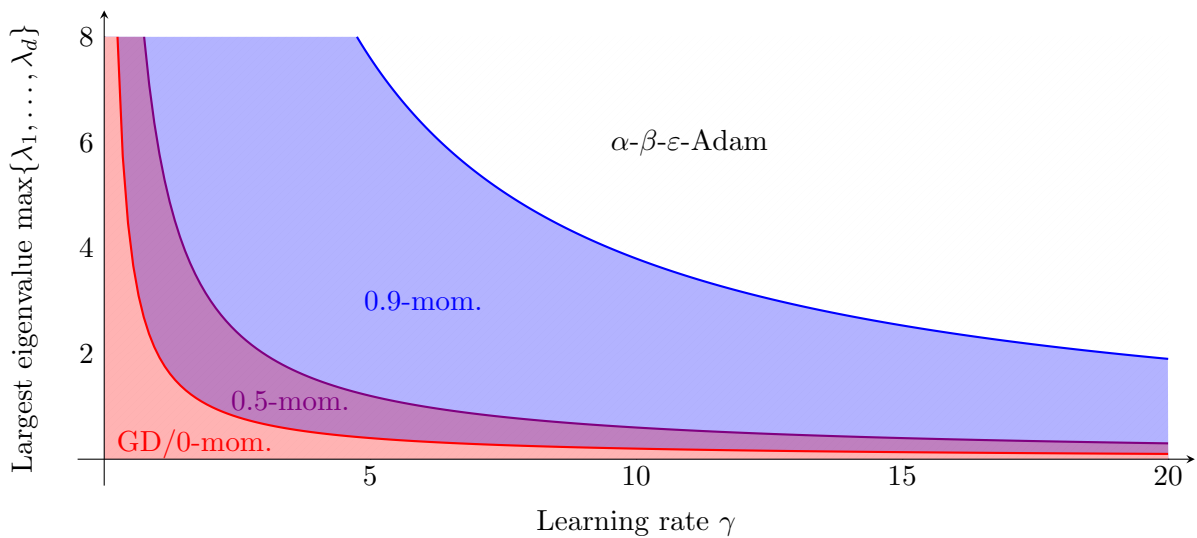

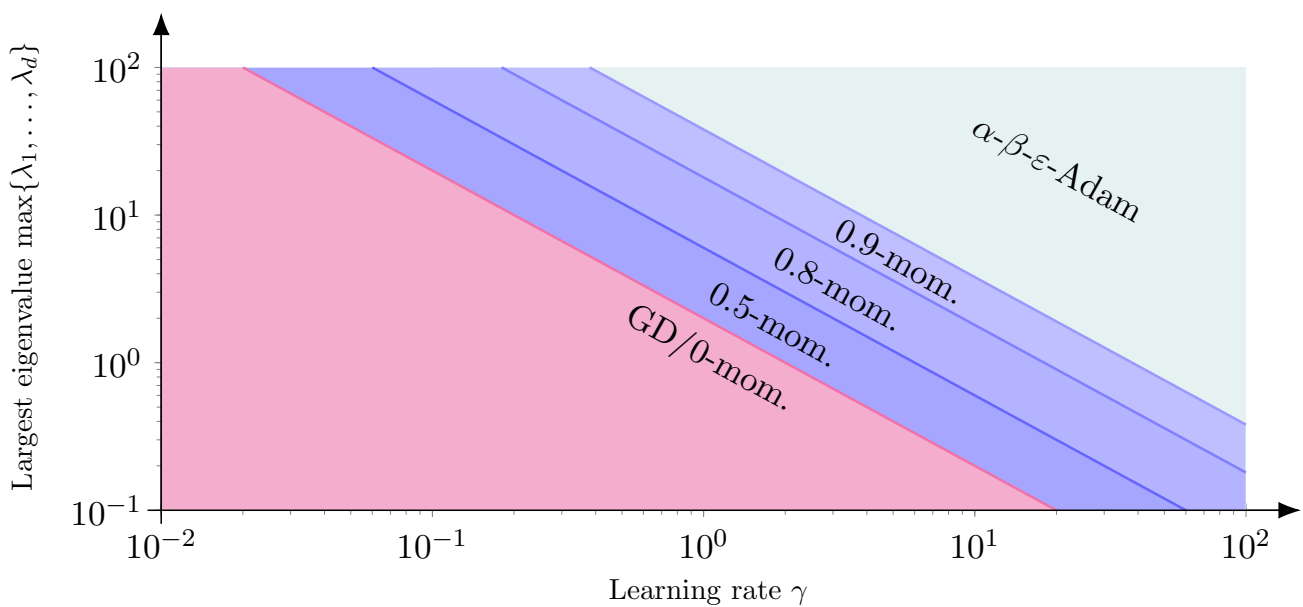
\begin{figure}[H]
			\centering
			\begin{tikzpicture}[scale=1.2]
				\begin{axis}[
					axis lines = middle,
					ymode=log,
					xmode=log,
					ymax=210,xmax=150,
					ymin=0.1,xmin=0.01,
					width=14cm,
					height=7cm,
					axis x line* = bottom,
					axis y line* = left,
					]
					\fill [teal!10, domain=0.1:100, variable=\x]
					(-4.6,4.6)
					-- (4.6,4.6)
					-- (4.6,-2.3)
					-- (-4.6,-2.3)
					-- cycle;
					\fill [blue!25, domain=0.1:100, variable=\x]
					(-2.3,4.6)
					-- (-0.95,4.6)
					-- (4.6,-0.95)
					-- (4.6,-2.3)
					-- (-2.3,-2.3)
					-- cycle;
					\fill [blue!30, domain=0.1:100, variable=\x]
					(-2.8,4.6)
					-- (-1.7,4.6)
					-- (4.6,-1.7)
					-- (4.6,-2.3)
					-- (-2.3,-2.3)
					-- cycle;
					\fill [blue!35, domain=0.1:100, variable=\x]
					(-4.6,4.6)
					-- (-2.8,4.6)
					-- (4.1,-2.3)
					-- (-4.6,-2.3)
					-- cycle;
					\fill [magenta!40, domain=0.1:100, variable=\x]
					(-4.6,4.6)
					-- (-3.9,4.6)
					-- (3,-2.3)
					-- (-4.6,-2.3)
					-- cycle;
					\addplot+[domain=0.38:100,mark=none,smooth,samples=200,blue!40,thick] (\x, {38/\x}); 
					\addplot+[domain=0.18:100,mark=none,smooth,samples=200,blue!50,thick] (\x, {18/\x}); 
					\addplot+[domain=0.06:100,mark=none,smooth,samples=200,blue!60,thick] (\x, {6/\x}); 
					\addplot+[domain=0.02:100,mark=none,smooth,samples=200,magenta!70,thick] (\x, {2/\x}); 
					\tznode(3,3){$\alpha$-$\beta$-$\eps$-Adam}[rotate=331]
					\tznode(1.67,1.58){$0.9$-mom.}[rotate=331]
					\tznode(1.15,1.15){$0.8$-mom.}[rotate=331]
					\tznode(0.6,0.6){$0.5$-mom.}[rotate=331]
					\tznode(0.05,0.05){GD/$0$-mom.}[rotate=331]
				\end{axis}
				\tznode(6,-0.9){Learning rate $\gamma$}
				\tznode(-1.5,2.7){Largest eigenvalue $\max\{\lambda_1,\dots,\lambda_d\}$}[rotate=90]
				\draw [line width=0.2mm, -{Latex[length=3mm]}]
				(-0.15,0) -- (12.5,0);
				\draw [line width=0.2mm, -{Latex[length=3mm]}]
				(0,-0.15) -- (0,5.5);
			\end{tikzpicture}
			\caption{In this figure we graphically represent for every 
				$\alpha\in[0,1)$, $\beta\in(\alpha^2,1)$, $\eps\in(0,\infty)$ the stability region of
				the \GD\ optimizer (the $0$-momentum optimizer), 
				the $0.5$-momentum optimizer, 
				the $0.8$-momentum optimizer, 
				the $0.9$-momentum optimizer, and
				the $\alpha$-$\beta$-$\eps$-\Adam\ optimizer (logarithmically scaled axis)
				.}
			\label{fig:comp:without:Nets:log}
		\end{figure}

		In \cref{cor:ADAM:stable},
		\cref{cor:SGD:stable}, and 
		\cref{cor:MOM:stable} above we specify the stability region of the 
		\Adam\ optimizer, the momentum optimizer, and the \GD\ optimizer. 
		In \cref{fig:comp:without:Nets} we graphically represent for every
		$\alpha\in[0,1)$, $\beta\in(\alpha^2,1)$, $\eps\in(0,\infty)$ the stability regions of
		the \GD\ optimizer (the $0$-momentum optimizer), 
		the $0.5$-momentum optimizer, 
		the $0.9$-momentum optimizer, and
		the $\alpha$-$\beta$-$\eps$-\Adam\ optimizer.
		In \cref{fig:comp:without:Nets:log} we graphically represent for every
		$\alpha\in[0,1)$, $\beta\in(\alpha^2,1)$, $\eps\in(0,\infty)$ the stability regions of
		the \GD\ optimizer (the $0$-momentum optimizer), 
		the $0.5$-momentum optimizer, 
		the $0.8$-momentum optimizer, 
		the $0.9$-momentum optimizer, and
		the $\alpha$-$\beta$-$\eps$-\Adam\ optimizer.

	\begin{athm}{theorem}{final:thm:1.1}
			Let $d\in\N$, $\alpha\in[0,1)$, $\beta\in(\alpha^2,1)$, $\eps\in(0,\infty)$. Then
			\begin{enumerate}[label=(\roman*)]
				\item \llabel{it:NEST}
				it holds for every $\alpha$-Nesterov optimizer $\Phi_n\colon(\R^d)^{n}\to\R^d$, $n \in \N$,\cfadd{def:Nest} 
				that the stability region\cfadd{def:stab:region} of $(\Phi_n)_{n\in \N}$ is
				\begin{equation}
					\pRb{\lambda=(\lambda_0,\lambda_1,\dots,\lambda_{d})\in[0,\infty)^{d+1}\colon \textstyle\max_{i\in\{1,2,\dots,d\}}\pr{\lambda_0\lambda_{i}}\leq2\PRb{\frac{1-\alpha^2}{1+2\alpha}}}\color{black},
				\end{equation}
				\vspace{-0.7cm}
				\item \llabel{it:SGD}
				it holds for every \GD\ optimizer $\Phi_n\colon(\R^d)^{n}\to\R^d$, $n \in \N$,\cfadd{def:SGD} 
				that the stability region\cfadd{def:stab:region} of $(\Phi_n)_{n\in \N}$ is
				\begin{equation}
					\pRb{\lambda=(\lambda_0,\lambda_1,\dots,\lambda_{d})\in[0,\infty)^{d+1}\colon \textstyle\max_{i\in\{1,2,\dots,d\}}\pr{\lambda_0\lambda_{i}}\leq2},
				\end{equation}
				\vspace{-0.7cm}
				\item \llabel{it:MOM}
				it holds for every $\alpha$-momentum optimizer $\Phi_n\colon(\R^d)^{n}\to\R^d$, $n \in \N$,\cfadd{def:MOM} 
				that the stability region\cfadd{def:stab:region} of $(\Phi_n)_{n\in \N}$ is
				\begin{equation}
					\pRb{\lambda=(\lambda_0,\lambda_1,\dots,\lambda_{d})\in[0,\infty)^{d+1}\colon \textstyle\max_{i\in\{1,2,\dots,d\}}\pr{\lambda_0\lambda_{i}}\leq2\PRb{\frac{1+\alpha}{1-\alpha}}},
				\end{equation}
				\vspace{-0.7cm}
				\item \llabel{it:RMS}
				it holds for every $\beta$-$\eps$-\RMSprop\ optimizer $\Phi_n\colon(\R^d)^{n}\to\R^d$, $n \in \N$,\cfadd{def:RMSprop}
				that the stability region\cfadd{def:stab:region} of $(\Phi_n)_{n\in \N}$ is
				$[0,\infty)^{d+1}$,
				\item \llabel{it:ADAM}
				it holds for every $\alpha$-$\beta$-$\eps$-\Adam\ optimizer  $\Phi_n\colon(\R^d)^{n}\to\R^d$, $n \in \N$,\cfadd{def:ADAM} 
				that the stability region\cfadd{def:stab:region} of $(\Phi_n)_{n\in \N}$ is
				$[0,\infty)^{d+1}$,
				\item \llabel{it:ADAM:s}
				it holds
				for 
				 every $\alpha$-$\beta$-$\eps$-\Adam\ optimizer  $\Phi_n\colon(\R^d)^{n}\to\R^d$, $n \in \N$,\cfadd{def:ADAM} and every bounded $\mathcal{A}\subseteq[0,\infty)^{d+1}$
				that $(\Phi_n)_{n \in \N}$ \sstable{\mathcal{A}}\!\!, and
				\item \llabel{it:ADAM:ss}
				it holds for
				every $\alpha$-$\beta$-$\eps$-\Adam\ optimizer  $\Phi_n\colon(\R^d)^{n}\to\R^d$, $n \in \N$,\cfadd{def:ADAM} and every compact $\mathcal{A}\subseteq[0,\infty)\times(0,\infty)^{d}$ 
				that $(\Phi_n)_{n \in \N}$ \ssstable{\mathcal{A}}
			\end{enumerate}
			\cfload.
	\end{athm}
	
	\begin{aproof}
	\argument{\cref{stab:def:relation};\cref{cor:Nesterov:stable}}{\lref{it:NEST}}.
	\argument{\cref{stab:def:relation};\cref{cor:SGD:stable}}{\lref{it:SGD}}.
	\argument{\cref{stab:def:relation};\cref{cor:MOM:stable}}{\lref{it:MOM}}.
	\argument{\cref{stab:def:relation};\cref{cor:ADAM:stable:clear}}{\lref{it:ADAM}}.
	\argument{\cref{RMS:is:ADAM};\lref{it:ADAM}}{\lref{it:RMS}}.
	\argument{\cref{cor:ADAM:sstable}}{\lref{it:ADAM:s}}.
	\argument{\cref{cor:ADAM:ssstable}}{\lref{it:ADAM:ss}}.
	\end{aproof}

		In \cref{fig:5}, \cref{fig:6}, and \cref{fig:7} below and \cref{fig:1} in \cref{section:1} above we graphically illustrate the fact that \Adam\ and \RMSprop\ are uniformly stable (cf.~\cref{def:uniformly:stable}) according to \cref{final:thm:1.1} above and
		the fact that momentum and \Adam\ attain higher order convergence rates according to \cite[Theorem~1.2]{dereich2025sharphigherorderconvergence}.

	\begin{figure}[H]
		\centering
		\begin{tikzpicture}[shorten >=1pt,-latex,draw=black!100, node distance=\layersep,auto]
			\tznode(-2.5,0.7){Convergence speed}[rotate=90]
			\tznode(3,-2){Stability}		
			\node (SGD) at (1,0) {GD};
			\node (momentum-SGD) at (1,1) {momentum};
			\node (Adam) at (4,0) {RMSprop};
			\node (RMSprop) at (4,1) {Adam};
			\draw [line width=0.2mm, -{Latex[length=3mm]}]
			(-2.5,-1.5) -- (8,-1.5);
			\draw [line width=0.2mm, -{Latex[length=3mm]}]
			(-2,-2) -- (-2,3);
		\end{tikzpicture}
		\caption{Graphical illustration of stability and convergence speed properties of the momentum optimizer,
			the \Adam\ optimizer,
			the \GD\ optimizer, and
			the \RMSprop\ optimizer
		}
		\label{fig:5}
	\end{figure}		
	
	\begin{figure}[H]
		\centering
		\begin{tikzpicture}[shorten >=1pt,-latex,draw=black!100, node distance=\layersep,auto]
			\tznode(-2.5,0.7){Convergence speed}[rotate=90]
			\tznode(3,-2){Stability}
			\tzplot[fill,blue!20,smooth cycle,thick,opacity=.5]{1} 
			(4,-0.5)(2.7,0.5)(4,1.5)(5.2,0.5);
			\tzplot[draw,blue,smooth cycle,thick]{1} 
			(4,-0.5)(2.7,0.5)(4,1.5)(5.2,0.5);
			
			\node [blue,text width=12em, text centered,font=\small](uniformly stable) at (6,-0.5) {uniformly stable};
			
			\node (SGD) at (1,0) {GD};
			\node (momentum-SGD) at (1,1) {momentum};
			\node (Adam) at (4,0) {RMSprop};
			\node (RMSprop) at (4,1) {Adam};
			\draw [line width=0.2mm, -{Latex[length=3mm]}]
			(-2.5,-1.5) -- (8,-1.5);
			\draw [line width=0.2mm, -{Latex[length=3mm]}]
			(-2,-2) -- (-2,3);
		\end{tikzpicture}
		\caption{Graphical illustration of stability and convergence speed properties of the momentum optimizer,
			the \Adam\ optimizer,
			the \GD\ optimizer, and
			the \RMSprop\ optimizer
		}
		\label{fig:6}
	\end{figure}

	\begin{figure}[H]
		\centering
		\begin{tikzpicture}[shorten >=1pt,-latex,draw=black!100, node distance=\layersep,auto]
			\tznode(-2.5,0.7){Convergence speed}[rotate=90]
			\tznode(3,-2){Stability}
			\tzplot[fill,red!20,smooth cycle,thick,opacity=.5]{1} 
			(2.3,0.5)(-1,1.1)(3,1.5)(5,1);
			\tzplot[draw,red,smooth cycle,thick]{1} 
			(2.3,0.5)(-1,1.1)(3,1.5)(5,1);

			\node [red,text width=12em, text centered,font=\small](optimal convergence rate) at (2,1.7) {optimal convergence rate};
			
			\node (SGD) at (1,0) {GD};
			\node (momentum-SGD) at (1,1) {momentum};
			\node (Adam) at (4,0) {RMSprop};
			\node (RMSprop) at (4,1) {Adam};
			\draw [line width=0.2mm, -{Latex[length=3mm]}]
			(-2.5,-1.5) -- (8,-1.5);
			\draw [line width=0.2mm, -{Latex[length=3mm]}]
			(-2,-2) -- (-2,3);
		\end{tikzpicture}
		\caption{Graphical illustration of stability and convergence speed properties of the momentum optimizer,
			the \Adam\ optimizer,
			the \GD\ optimizer, and
			the \RMSprop\ optimizer
		}
		\label{fig:7}
	\end{figure}

	\subsection*{Acknowledgements}
	This work has been partially funded by the European Union (ERC, MONTECARLO, 101045811). The views and the opinions expressed in this work are however those of the authors only and do not necessarily reflect those of the European Union or the European Research Council (ERC). Neither the European Union nor the granting authority can be held responsible for them. We also gratefully acknowledge the Cluster of Excellence EXC 2044-390685587, Mathematics Münster: Dynamics-Geometry-Structure funded by the Deutsche Forschungsgemeinschaft (DFG, German Research Foundation).
	Most of the specific formulations in the proofs of this work have been created using \cite{KuckuckTools}.

	\bibliographystyle{acm}
	\bibliography{bibfile}

	\end{document}